%% file: main2.tex
\documentclass[12pt,twoside,singlespace]{mitthesis}
\usepackage{lgrind,latexsym,amssymb,amsmath,amsthm,amsfonts,graphicx,float,bm,grffile,hyperref}
\usepackage{algorithmicx}
\usepackage{algpseudocode}
\newcommand{\pd}{\partial}

\newcommand{\R}{\mathbb{R}}

\newcommand{\<}{\langle}
\renewcommand{\>}{\rangle}

\newcommand{\eps}{\varepsilon}

\newcommand{\tr}{\mbox{tr}}

\newcommand{\beq}{\begin{equation}}
\newcommand{\eeq}{\end{equation}}
\renewcommand{\tilde}{\widetilde}
\renewcommand{\hat}{\widehat}
\newcommand{\bit}{\begin{itemize}}
\newcommand{\eit}{\end{itemize}}
\newcommand{\ben}{\begin{enumerate}}
\newcommand{\een}{\end{enumerate}}
\newcommand{\x}{\mathbf{x}}
\newcommand{\y}{\mathbf{y}}
\newcommand{\bc}{\mathbf{c}}
\newcommand{\Rm}{R_\text{max}}

\newtheorem{definition}{Definition}

\newtheorem{theorem}{Theorem}
\newtheorem{lemma}{Lemma}
\newtheorem{algorithm}{Algorithm}

\newtheorem{remark}[subsection]{Remark}

\begin{document}

\include{cover2}
\pagestyle{plain}

\include{contents}

\include{intro2}

\include{back2}
\include{probing2}
\include{plr2}

\include{conclusion2}

\appendix
\include{steps2}
\include{biblio}
\end{document}

%% file: cover2.tex
\title{Compressed Absorbing Boundary Conditions for the Helmholtz Equation}

\author{Rosalie B\'elanger-Rioux}
\prevdegrees{B.Sc., McGill University (2009)}
\department{Department of Mathematics}

\degree{Doctor of Philosophy}
\degreemonth{June}
\degreeyear{2014}
\thesisdate{May 8, 2014}

\supervisor{Laurent Demanet}{Assistant Professor}

\chairman{Michel Goemans}{Chairman, Department Committee on Graduate Theses}

\maketitle

\cleardoublepage
\setcounter{savepage}{\thepage}
\begin{abstractpage}
\input{abstract}
\end{abstractpage}

\cleardoublepage

\section*{Acknowledgments}

I would like to thank my advisor, Laurent Demanet, for his advice and support throughout my time as a graduate student.
\newline
\newline
Also, I would like to acknowledge professors I took classes from, including my thesis committee members Steven Johnson and Jacob White, and qualifying examiners Ruben Rosales and again Steven Jonhson. I appreciate your expertise, availability and help.
\newline
\newline
I would like to thank collaborators and researchers with whom I discussed various ideas over the years, it was a pleasure to learn and discover with you.
\newline
\newline
Thanks also to my colleagues here at MIT in the mathematics PhD program, including office mates. %
It was a pleasure to get through this together.
\newline
\newline
Thanks to my friends, here and abroad, who supported me. From my decision to come to MIT, all the way to graduation, it was good to have you.
\newline
\newline
Finally, I thank my loved ones for being there for me always, and especially when I needed them.

%% file: abstract.tex
%
%
%
Absorbing layers are sometimes required to be impractically thick in order to offer an accurate approximation of an absorbing boundary condition for the Helmholtz equation in a heterogeneous medium. It is always possible to reduce an absorbing layer to an operator at the boundary by layer-stripping elimination of the exterior unknowns, but the linear algebra involved is costly. We propose to bypass the elimination procedure, and directly fit the surface-to-surface operator in compressed form from a few exterior Helmholtz solves with random Dirichlet data. We obtain a concise description of the absorbing boundary condition, with a complexity that grows slowly (often, logarithmically) in the frequency parameter. We then obtain a fast (nearly linear in the dimension of the matrix) algorithm for the application of the absorbing boundary condition using partitioned low rank matrices. The result, modulo a precomputation, is a fast and memory-efficient compression scheme of an absorbing boundary condition for the Helmholtz equation.

%% file: contents.tex
\tableofcontents
\newpage
\listoffigures
\newpage
\listoftables

%% file: intro2.tex
\chapter{Introduction}\label{ch:intro}

This work investigates arbitrarily accurate realizations of absorbing (a.k.a. open, radiating) boundary conditions (ABC), including absorbing layers, for the 2D acoustic high-frequency Helmholtz equation in certain kinds of heterogeneous media. Instead of considering a specific modification of the partial differential equation, such as a perfectly matched layer, we study the broader question of compressibility of the nonlocal kernel that appears in the exact boundary integral form $Du=\pd_\nu u$ of the ABC. The operator $D$ is called the \emph{Dirichlet-to-Neumann} (DtN) map. This boundary integral viewpoint invites to rethink ABCs as a two-step numerical scheme, where 
\begin{enumerate}
\item a precomputation sets up an expansion of the kernel of the boundary integral equation, then 
\item a fast algorithm is used for each application of this integral kernel at the open boundaries in a Helmholtz solver.
\end{enumerate}
This two-step approach may pay off in scenarios when the precomputation is amortized over a large number of solves of the original equation with different data. This framework is, interestingly, half-way between a purely analytical or physical method and a purely numerical one. It uses both the theoretical grounding of analytic knowledge and the intuition from understanding the physics of the problem in order to obtain a useful solution.

The numerical realization of ABC typically involves absorbing layers that become impractical for difficult $c(\mathbf{x})$, or for high accuracy. We instead propose to realize the ABC by directly compressing the integral kernel of $D$, so that the computational cost of its setup and application would become competitive when (\ref{eq:HE}) is to be solved multiple times. Hence this paper is not concerned with the design of a new ABC, but rather with the reformulation of existing ABCs that otherwise require a lot of computational work per solve. In many situations of practical interest we show that it is possible to ``learn" the integral form of $D$, as a precomputation, from a small number of solves of the exterior problem with the expensive ABC. By ``small number", we mean a quantity essentially independent of the number of discretization points $N$ along one dimension -- in practice as small as 1 or as large as 50. We call this strategy matrix probing. To show matrix probing is efficient, we prove a result on approximating $D$ in a special case, with inverse powers multiplied by a complex exponential. This leads us to the successful design of a basis for a variety of heterogeneous media.

Once we obtain a matrix realization $\tilde{D}$ of the ABC from matrix probing, we can use it in a Helmholtz solver. However, a solver would use matrix-vector multiplications to apply the dense matrix $\tilde{D}$. Hence the second step of our numerical scheme: we compress $\tilde{D}$ using partitioned low-rank (PLR) matrices to acquire a fast matrix-vector product. This second step can only come after the first, since it is the first step that gives us access to the entries in $D$ and allows us to use the compression algorithms of interest to us. We know we can use hierarchical or partitioned-low-rank matrices to compress the DtN map because we prove the numerical low-rank property of off-diagonal blocks of the DtN map, and those algorithms exploit favorably low-rank blocks. Since PLR matrices are more flexible than hierarchical matrices, we use them to compress $\tilde{D}$ into $\overline{D}$ and apply it to vectors in complexity ranging from $O(N\log N)$ (more typical) to $O(N^{3/2})$ (worst case). %
The precomputation necessary to set up the PLR approximation is of similar complexity. This can be compared to the complexity of a dense matrix-vector product, which is $O(N^2)$.

In this introduction, we first motivate in section \ref{sec:motivate} the study of the Helmholtz equation in an unbounded domain by presenting important applications. We then give more details on the steps of our numerical scheme in section \ref{sec:struct}, which will make explicit the structure of this thesis.%

\section{Applications of the Helmholtz equation in an unbounded domain}\label{sec:motivate}

We consider the scalar Helmholtz equation in $\mathbb{R}^2$,
\begin{equation}\label{eq:HE0}
\Delta u(\x)+\frac{\omega^2}{c^2(\x)} u(\x) = f(\x), \qquad \x = (x_1, x_2),
\end{equation}
with compactly supported $f$, the \emph{source}. We seek the solution $u$. The function $c$ in \eqref{eq:HE0} is called the \emph{medium}, $\omega$ the frequency. %

We consider the unique solution $u$ to \eqref{eq:HE0} determined by the Sommerfeld radiation condition (SRC), which demands that the solution be outgoing. %
We call the problem of finding a solution $u$ to \eqref{eq:HE0} with the SRC the \emph{free-space problem}. There are many applications that require a solution to the free-space Helmholtz problem, or the related free-space Maxwell's equations problem. To solve the free-space problem numerically, we reformulate the problem to a bounded domain $\Omega$ in which we shall find the solution $u$. One must then impose an \emph{absorbing boundary condition} (ABC) on the boundary $\pd \Omega$. ABCs are designed to absorb waves impinging on $\pd \Omega$ so the waves do not reflect back in $\Omega$ and pollute the solution there.
ABCs in heterogeneous media are often too costly, but the two-step scheme of this thesis addresses this problem. %
Let us now discuss applications of this two-step scheme.%

The main application is wave-based imaging, an inverse problem for the wave equation
\begin{equation}\label{eq:WE}
 \Delta U(\x,t) - \frac{1}{c^2(\x)}\frac{\pd^2 U(\x,t)}{\pd t^2} = F(\x,t), \qquad
\end{equation}
and related equations. The Helmholtz equation is equivalent to the wave equation because we can decompose the solution $U$ and the source $F$ into time-harmonic components by a Fourier transform in time. %
 Another type of wave equation, the Maxwell's equations, can also be reformulated as Helmholtz equations on the various components of the electric and magnetic fields.%

An inverse problem is as follows: instead of trying to find the solution $u$ of \eqref{eq:HE0} given $\omega$, $f$ and $c$, we do the opposite. In other words, we are given the solution $u$ at a set of receiver locations for some $\omega$ and various sources $f$. We try to determine the medium $c$ from that information. We are usually interested in (or can only afford) $c(\x)$ for $\x$ in some part of the whole space, say some bounded domain $\Omega$, and absorbing boundary conditions are then necessary. To solve the inverse problem in $\Omega$, we need to use a lot of sources, say a thousand. The details of why many sources are useful, and how to solve the inverse problem, are not relevant to this thesis. What is relevant is that, in the course of solving the inverse problem for say the Helmholtz equation, it is needed to solve the free-space Helmholtz problem for all these sources, with ABCs, in heterogeneous media. We list here some imaging applications where the Helmholtz equation, or other types of wave equations, are used to solve inverse problems with ABCs in heterogeneous media, that is, where the numerical scheme in this thesis might prove useful.

\begin{itemize}
 \item In seismic imaging, we acquire knowledge of the rock formations under the earth's surface. That is, we want to know the medium $c$ in which the waves propagate. The domain $\Omega$ might be a region surrounding a seismic fault \cite{seismic} where one wants to assess earthquake hazards, or a place where one might like to find hydrocarbons, other minerals, or even geothermal energy \cite{geothermal}. ABCs are needed to simulate the much larger domain in which $\Omega$ is embedded, that is, the Earth. Sources (which might be acoustic or electromagnetic) and receivers might be placed on the surface of the Earth or inside a well, or might be towed by a boat or placed at the bottom of the sea. An earthquake can also be used as a source.
 \item Ultrasonic testing \cite{ultrasonic} is a form of non-destructive testing where very short ultrasonic pulses are sent inside an object. The object might be, say, a pipe that is being tested for cracks or damage from rust, or a weld being tested for defects. The received reflections or refractions from those ultrasonic pulses are used for diagnosis on the object of interest. $\Omega$ might be the object itself, or a part of it. Ultrasonic imaging \cite{3dmedimag} is also used in medicine to visualize fetuses, muscles, tendons and organs. The domain $\Omega$ is then the body part of interest.
 \item Synthetic-aperture radar imaging is used to visualize a scene by sending electromagnetic pulses from an antenna aboard a plane or satellite \cite{radar}. It is also used to detect the presence of an object \cite{borden} far away such as a planet, or through clouds or foliage.
\end{itemize}

An entirely different application of the Helmholtz equation resides in photonics \cite{JohnsonPhot}. The field of photonics studies the optical properties of materials. In particular, one tries to construct photonic crystals, a periodic medium with desired properties depending on the specific application. It is of course useful to first test photonic crystals numerically to observe their properties before actually building them. However, since crystals are assumed to be infinite, absorbing boundary conditions need to be used to reformulate the problem on a bounded domain. This is where our two-step numerical scheme can be used.

\section{A two-step numerical scheme for compressing ABCs}\label{sec:struct}

The next chapter, chapter \ref{ch:back}, will present theoretical facts about the Helmholtz equation and related concepts, which will be useful at various points in later chapters for developing the proposed framework. Then, chapter \ref{ch:probing} will present the first step of the scheme: matrix probing and how it is used for providing a rapidly converging expansion of the DtN map. Next, chapter \ref{ch:plr} will contain the second part of the proposed two-step scheme: the compression of the matrix probing expansion, for fast evaluation, using partitioned low-rank matrices. After the precomputations of matrix probing and PLR matrices, we finally obtain a fast and accurate compressed absorbing boundary condition. Chapter \ref{ch:conclusion} concludes this thesis with a review of the main new ideas presented, identification of further research directions and open problems, and also an overview of how the presented framework could be used in other contexts. A summary of the steps involved in the presented two-step numerical scheme is available in appendix \ref{ch:steps}.

\subsection{Background material}

We begin by introducing again the problem we wish to solve, the \emph{free-space problem} for the Helmholtz equation. This problem is defined on an unbounded domain such as $\mathbb{R}^2$, but can be solved by reformulating to a bounded domain we will call $\Omega$. One way to do this is to impose what we call the \emph{exterior Dirichlet-to-Neumann map} (DtN map) on $\pd \Omega$. The DtN map $D$ relates Dirichlet data $u$ on $\pd \Omega$ to the normal derivative $\pd_\nu u$ of $u$, where $\nu$ is the unit outward vector to $\pd \Omega$: $Du=\pd_\nu u$. This allows us to solve the Helmholtz equation on $\Omega$ to obtain a solution $u$ which coincides, inside $\Omega$, with the solution to the free-space problem.

We shall see that any \emph{absorbing boundary condition} (ABC) or Absorbing Layer (AL) can be used to obtain the DtN map. As we mentioned before, an ABC is a special condition on the boundary $\pd \Omega$ of $\Omega$, which should minimize reflections of waves reaching $\pd \Omega$. %
Different ABCs work differently however, and so we will review existing techniques for constructing ABCs for the Helmholtz equation, and see again how they are computationally expensive in a variable medium. This suggests finding a different, more efficient, way of obtaining the exterior DtN map from an ABC. To attack this problem, we consider the \emph{half-space} DtN map, known analytically in a constant medium. This half-space DtN map is actually quite similar to the exterior DtN map, at least when $\Omega$ is a rectangular domain. Indeed, the interface along which the DtN map is defined for the half-space problem is a straight infinite line $x_2=0$. For the exterior problem on a rectangle, the DtN map is defined on each of the straight edges of $\Omega$. Hence the restriction of the exterior DtN map to one such edge of $\Omega$, which say happens to be on $x_2=0$, behaves similarly to the restriction of the half-space DtN map to that same edge. The main difference between those two maps is created by scattering from corners of $\pd \Omega$. Both in chapter \ref{ch:probing} and \ref{ch:plr}, we will prove facts about the half-space DtN map which will inform our numerical scheme for the exterior DtN map. 

We end the chapter of background material with how the most straightforward way of obtaining the exterior DtN map from an ABC, \emph{layer-stripping}, is prohibitively slow, especially in variable media. This will also explain how, even if we have an efficient way of obtaining the DtN map from an ABC, applying it at every solve of the Helmholtz equation will be slow. In fact, once we have the map $D$ in $Du=\pd_\nu u$, we need to apply this map to vectors inside a Helmholtz solver. But $D$ has dimension approximately $N$, if $N$ is the number of points along each direction in $\Omega$, so matrix-vector products with $D$ have complexity $N^2$. This is why we have developed the two-step procedure presented in this thesis: first an expansion of the DtN map using matrix probing, then a fast algorithm using partitioned low-rank matrices for the application of the DtN map.

\subsection{Matrix probing}

Chapter \ref{ch:probing} is concerned with the first step of this procedure, namely, setting up an expansion of the exterior DtN map kernel, in a precomputation. This will pave the way for compression in step two, presented in chapter \ref{ch:plr}.

Matrix probing is used to find an expansion of a matrix $M$. For this, we assume an expansion of the type 
\[ M \approx \tilde{M} = \sum_{j=i}^{p} c_jB_j, \] 
where the \emph{basis matrices} $\left\{B_j\right\}$ are known, and we wish to find the coefficients $\left\{c_j\right\}$. We do not have access to $M$ itself, but only to products of $M$ with vectors. In particular, we can multiply $M$ with a random vector $z$ to obtain
\[ w=Mz \approx \sum_{j=1}^p c_j B_j z = \Psi_z \, \bc.\]
We can thus obtain the vector of coefficients $\bc$ by applying the pseudoinverse of $\Psi_z$ to $Mz=w$.

For matrix probing to be an efficient expansion scheme, we need to carefully choose the basis matrices $\left\{B_j\right\}$. Here, we use knowledge of the half-space DtN map to inform our choices: a result on approximating the half-space DtN map with a particular set of functions, inverse powers multiplied by a complex exponential. Hence our basis matrices are typically a discretization of the kernels
\begin{equation}%
B_j(x,y)= \frac{e^{ik|x-y|}}{(h+|x-y|)^{j/2}},
\end{equation}
where $h$ is our discretization parameter, $h=1/N$.  The need for a careful design of the basis matrices can however be a limitation of matrix probing. In this thesis, we have used also insights from geometrical optics to derive basis matrices which provide good convergence in a variety of cases. Nonetheless, in a periodic medium such as a photonic crystal, where the wavelength is as large as the features of the medium, the geometrical optics approximation breaks down. Instead, we use insights from the solution in a periodic medium, which we know behaves like a Bloch wave, to design basis matrices. However, the periodic factor of a Bloch wave does not lead to very efficient basis matrices since it is easily corrupted by numerical error. Another limitation is that a medium which has discontinuities created discontinuities as well in the DtN map, forcing again a more careful design of basis matrices.

Once we do have basis matrices, we can probe the DtN map $D$, block by block. Indeed, we choose $\Omega=[0,1]^2$, and so $\pd \Omega$ is the four edges of the square\footnote{The framework does carry over to polygonal domains easily, but we do not cover this here.}, numbered counter-clockwise. Hence, $D$ has a 4 by 4 structure corresponding to restrictions of the DtN map to each pair of edges. To obtain the product of the $(i_M,j_M)$ block $M$ of $D$ with a random vector, we need to solve what we call the \emph{exterior problem}: we put a random Dirichlet boundary condition on edge $j_M$ of $\pd \Omega$, solve the Helmholtz equation outside of $\Omega$ using an ABC, and take the normal derivative of the solution on edge $i_M$ of $\pd \Omega$.

In situations of practical interest, we obtain the DtN map as a precomputation using matrix probing with leading complexity of about 1 to 50 solves of the exterior problem with the expensive ABC -- a number of solves essentially independent of the number of discretization points $N$. A solve of the exterior problem is essentially equivalent to a solve of the original problem with that same expensive ABC.

We then present a careful study of using matrix probing to expand the DtN map in various media, use the matrix probing expansion to solve the Helmholtz equation, and document the complexity of the method.

In the next chapter, we will present a fast algorithm for applying the DtN map's expansion found by matrix probing.

\subsection{Partitioned low-rank matrices}

In chapter \ref{ch:plr}, we produce a fast algorithm for applying the DtN map $D$ to vectors, since this is an operation that a Helmholtz solver needs. Again, we needed matrix probing first, to obtain an explicit approximation of $D$, to be compressed in order to obtain a fast matrix-vector product. Indeed, we do not have direct access to the entries of $D$ at first, but rather we need to solve a costly problem, the exterior problem, every time we need a multiplication of $D$ with a vector. Now that we have an explicit representation of $D$ from matrix probing, with cost $O(N^2)$ for matrix-vector products, we can compress that representation to obtain a faster product.

We have mentioned before how the half-space DtN map is similar to the exterior DtN map. In chapter \ref{ch:plr}, we will prove that the half-space DtN map kernel $K$ in constant medium is numerically low rank. This means that, away from the singularity of $K$, the function $K(|x-y|)$ can be written as a short sum of functions of $x$ multiplied by functions of $y$: 
\[ K(|x-y|)=\sum_{j=1}^J \Psi_j(x)\chi_j(y) + E(x,y)\]
with error $E$ small. The number of terms $J$ depends logarithmically on the error tolerance and on the frequency $k$. This behavior carries over to some extent to the exterior DtN map. Thus we use a compression algorithm which can exploit the low-rank properties of blocks of $D$ that are not on the diagonal, that is, that are away from the singularity.

A well-known such compression framework is called \emph{hierarchical matrices}. Briefly, hierarchical matrices adaptively divide or compress diagonal blocks. We start by dividing the matrix in 4 blocks of half the original matrix's dimension. The two off-diagonal blocks are compressed: we express them by their singular value decomposition (SVD), truncated after $R$ terms. The two diagonal blocks are divided in four again, and we recurse: off-diagonal blocks are compressed, diagonal blocks are divided, etc. We do not divide a block if compressing it will result in less error than our error tolerance $\eps$ -- this is adaptivity. The parameters $R$ and $\eps$ are chosen judiciously to provide a fast matrix-vector multiplication with small error.

However, the hierarchical matrix framework can only apply to matrices with a singularity along the diagonal. This is not useful to us, since for example a block of $D$ corresponding to two consecutive edges of $\pd \Omega$ will have the singularity in a corner. We thus decide to use partitioned low-rank (PLR) matrices for compressing and applying $D$. PLR matrices are more flexible than hierarchical matrices: when we divide a block, or the original matrix, in 4 sub-blocks, any of those four sub-blocks can be divided again. Any block that is not divided any further is called a \emph{leaf}, and $\mathcal{B}$ is the set of all leaves. If a singularity is in a corner, then the PLR compression algorithm will automatically divide blocks close to that corner, but will compress farther blocks since they have lower numerical rank. Note that we use the randomized SVD \cite{randomSVD} to speed us the compression, so that its complexity is of order $O(NR^2|\mathcal{B}|)$, where $|\mathcal{B}|$ is often on the order of $\log N$ but will be $\sqrt{N}$ in the worst case. Similarly, the complexity of a matrix-vector product is usually $O(NR|\mathcal{B}|)$, which for $N\approx 1000$ provides a speed-up over a dense matrix-vector product of a factor of 30 to 100. We also show that the worst case complexity of a matrix-vector product in the PLR framework is $O(N^{3/2})$. This should be compared to the complexity of a dense matrix-vector product, which is $O(N^2)$.

We may then use PLR matrices to compress the DtN map, and use this compressed map in a Helmholtz solver. We verify the complexity of the method, and present results on the solution of the Helmholtz equation with a probed and compressed DtN map.

\subsection{Summary of steps}

In appendix \ref{ch:steps}, we present a summary of the various operations involved in implementing the presented two-step algorithm.

%% file: back2.tex
\chapter{Background material}\label{ch:back}

In this chapter, we review the necessary background material for this thesis.
We begin by introducing in section \ref{sec:HEfree} the problem we wish to solve, that is, the \emph{free-space problem} for the Helmholtz equation. This problem is defined on an unbounded domain. We will explain how \emph{absorbing boundary conditions} (ABCs) and the \emph{exterior Dirichlet-to-Neumann map} (DtN map) are related, and can be used to reformulate the free-space problem from an unbounded domain to a bounded one. This reformulation allows us to solve the free-space problem.

We then briefly review existing techniques for constructing ABCs for the Helmholtz equation in section \ref{sec:abc}. From this, we will understand why, when solving the Helmholtz equation in a heterogeneous medium, ABCs are computationally expensive. This suggests an approach where we compress an ABC, or the related exterior DtN map. In other words, we find a way to apply the ABC which allows for a faster solve of the free-space problem.

To attack this problem, we first look at the \emph{half-space} DtN map, which is known analytically in uniform medium. In some respect, this half-space DtN map is quite similar to the exterior DtN map. To understand this similarity, we introduce in section \ref{sec:HEhalf} the concepts of the \emph{exterior} and \emph{half-space} problems for the Helmholtz equation, and the related half-space DtN map. In the next chapters, we will prove facts about the half-space DtN map which will inform our compression scheme for the exterior DtN map. %

Before we explain our work, we end this chapter of background material with section \ref{sec:strip}. We first describe how we may eliminate unknowns in the Helmholtz discretized system to obtain a Riccati equation governing the half-space DtN map. We then describe how the most straightforward way of compressing the exterior DtN map, by eliminating the exterior unknowns in the Helmholtz discretized system, which we call \emph{layer-stripping}, is prohibitively slow. This will explain also how, even if we have a more efficient way of obtaining the DtN map, applying it at every solve of the Helmholtz equation might still be slow. This is why we have developed the two-step procedure presented in this thesis: first an expansion of the DtN map, then a fast algorithm for the application of the DtN map.

\section{The Helmholtz equation: free-space problem}\label{sec:HEfree}

We consider the scalar Helmholtz equation in $\mathbb{R}^2$,
\begin{equation}\label{eq:HE}
\Delta u(\x)+\frac{\omega^2}{c^2(\x)} u(\x) = f(\x), \qquad \x = (x_1, x_2),
\end{equation}
with compactly supported $f$. We shall call this $f$ the \emph{right-hand side} or the \emph{source}. Here, the solution we seek to this equation is $u$. The function $c$ in \eqref{eq:HE} is called the \emph{medium of propagation}, or simply \emph{medium}. When $c$ is a constant, we say the medium is \emph{uniform}. When $c$ varies, we say the medium is \emph{heterogeneous}. We call $\omega$ the frequency, and note (this shall be explained later) that a high frequency makes the problem harder to solve numerically.

Throughout this work we consider the unique solution $u$ to \eqref{eq:HE} determined by the Sommerfeld radiation condition (SRC) at infinity: when $c(\x)$ extends to a uniform $c$ outside of a bounded set\footnote{If the medium $c(\x)$ does not extend to a uniform value outside of a bounded set, it is possible one could use the limiting absorption principle to define uniqueness of the solution to the Helmholtz equation. Start from the wave equation solution $u(\x,t)$. We would like to take the Fourier transform in time of $u$ to obtain the Helmholtz solution. We may take the Laplace transform of $u$ instead to obtain $\hat{u}(\x,s)$, and ask that the path of integration in the $s$ variable approach the imaginary axis from the decaying side. This could then be used to define the Helmholtz solution.}, the SRC is \cite{McLean}
\begin{equation}\label{eq:src}
\lim_{r \rightarrow \infty} r^{1/2} \left( \frac{\pd u}{\pd r} - ik u \right) = 0, \qquad k = \frac{\omega}{c},
\end{equation}
where $r$ is the radial coordinate. We call the problem of finding a solution $u$ to \eqref{eq:HE} and \eqref{eq:src} the \emph{free-space problem}.

\subsection{Solution in a uniform medium using the Green's function}\label{sec:Gfsol}

When the medium $c$ is uniform, there exists an analytical solution to the free-space problem, using the \emph{Green's function}.
\begin{definition}\label{def:Green}
The free-space \emph{Green's function} for the Helmholtz equation \eqref{eq:HE} is the unique function $G(\x,\y)$ which solves the free-space problem with right-hand side a delta function, $f(\x)=\delta(|\x-\y|)$, for every fixed $\y$.
\end{definition}
It is well-known (p. 19 of \cite{bookChew}) that the Green's function $G$ of the uniform free-space problem \eqref{eq:HE}, \eqref{eq:src} is the following:
\begin{equation}\label{eq:Gfree}
G(\x,\y)=\frac{i}{4} H_0^{(1)}(k|\x-\y|)
\end{equation}
where $H_0^{(1)}$ is the Hankel function of zeroth order of the first kind. Then, one can compute the solution $u(\x)$ for any $\x \in \mathbf{R}^2$, with now any given right-hand side $f$ supported on some bounded domain, using the following formula (see p. 62 of \cite{FollandIntroPDE}):
\begin{equation}\label{eq:Gsol}
u(\x)=\int_{\mathbf{R}^2} G(|\x-\y|)f(\y) \ dy.
\end{equation}
To rapidly and accurately evaluate the previous expression at many $\x$'s is not easy: $G(\x,\y)$ has a singularity at $\x=\y$, and so care must be taken in the numerical evaluation of the integral in \eqref{eq:Gsol}. Since the issue of evaluating the solution to the free-space problem when the Green's function is known is not directly relevant to this thesis, we refer the reader to \cite{bem}.

\subsection{Solution in a heterogeneous medium using Absorbing Boundary Conditions}

For a number of practical applications of the Helmholtz equation, as was mentioned in the introductory chapter, one cannot find an analytical solution to this problem because the medium is not uniform, that is, $c$ is not a constant. In particular, the Green's function is not known, and one expects that calculating this Green's function numerically is just as hard as calculating the solution with numerical methods.

To obtain a numerical solution to an equation on an unbounded domain then, one must first reformulate this equation to obtain a possibly modified equation on a bounded domain. Hence, we pick a domain $\Omega$ in $\mathbb{R}^2$, with the compact support of $f$ contained in $\Omega$, such that $\Omega$ contains the area of interest, that is, where we care to obtain a solution. We now seek to reformulate the SRC on the boundary $\pd \Omega$, so that the resulting solution inside $\Omega$ matches that of the free-space problem. This leads us to define an \emph{Absorbing Boundary Condition}:

\begin{definition}
An \emph{Absorbing Boundary Condition} for the Helmholtz equation \eqref{eq:HE} is a condition on $\pd \Omega$, the boundary of a closed, bounded domain $\Omega \in \mathbb{R}^2$, which uniquely defines a solution to the Helmholtz equation restricted to $\Omega$, such that this unique solution matches the solution to the free-space problem \eqref{eq:HE}, \eqref{eq:src}.
\end{definition}

Clearly, if we can reformulate the SRC on the boundary $\pd \Omega$, so that the resulting solution inside $\Omega$ matches that of the free-space problem, we will obtain an ABC.

ABCs are extremely important to the numerical solution of the free-space problem because, as alluded to previously, they allow us to restrict the computational domain to a bounded domain $\Omega$, where a solution can be computed in finite time with finite memory. We will discuss ABCs in more detail in the next section, section \ref{sec:abc}.

We will explain in the next chapter the particular, quite rudimentary, solver we used in our numerical experiments. We note here that a better solver should be used for treating larger problems or obtaining better accuracy, for two reasons. First of all, as we explain in more detail in section \ref{sec:strip}, the cost of solving the Helmholtz problem with a standard solver in two dimensions is $O(N^4)$, which is prohibitive. Secondly, as we discuss in section \ref{sec:compABC}, a higher frequency means we need more points per wavelength in our discretization -- this is known as the \emph{pollution effect}. To treat larger problems, there exist better solvers such as the sweeping preconditioner of Engquist and Ying \cite{Hsweep,Msweep}, the shifted Laplacian preconditioner of Erlangga \cite{erlangga,er}, the domain decomposition method of Stolk \cite{stolk}, or the direct solver with spectral collocation of Martinsson, Gillman and Barnett \cite{dirfirst,dirstab}. The problem of appropriate numerical solvers for the Helmholtz equation in high frequency is very much a subject of ongoing research and not the purpose of this thesis, hence we do not discuss this further.

\subsection{The Dirichlet-to-Neumann map as an ABC}\label{sec:dtn}

We now seek to reformulate the SRC \eqref{eq:src} on $\pd \Omega$. There are many ways to do that numerically, as we shall see in section \ref{sec:abc}, but we wish to highlight this particular, analytical way because it introduces a fundamental concept, the \emph{Dirichlet-to-Neumann} map.

Let $G(\x,\y)$ be the Green's function for the free-space problem. Define the single and double layer potentials, respectively, on some closed contour $\Gamma$ by the following, for $\psi, \ \phi$ on $\Gamma$ (see details in \cite{McLean}, \cite{CK}):
\[
 S \psi (\x)=\int_{\Gamma} G(\x,\y) \ \psi(\y) \ dS_y, \qquad T \phi (\x) =\int_{\Gamma} \frac{\pd G}{\pd \nu_{\y}} (\x,\y) \ \phi(\y) \ dS_{\y},
\]
where $\nu$ is the outward pointing normal to the curve $\Gamma$, and $\x$ is not on $\Gamma$.
Now let $u^+$ satisfy the Helmholtz equation \eqref{eq:HE} in the exterior domain $\R \setminus \overline{\Omega}$, along with the SRC \eqref{eq:src}. Then Green's third identity is satisfied in the exterior domain: using $\Gamma= \pd \Omega$, we get
\begin{equation}\label{eq:GRF}
T u^+ - S \frac{\pd u}{\pd \nu}^+ =  u^+, \qquad \x \in \R^2 \setminus \overline{\Omega}.
\end{equation}
Finally, using the jump condition of the the double layer $T$, we obtain Green's identity on the boundary $\pd \Omega$:
\[
(T - \frac{1}{2} I ) \, u^+ - S \frac{\pd u}{\pd \nu}^+ =  0, \qquad \x \in \pd \Omega.
\]
When the single-layer potential $S$ is invertible\footnote{This is the case when there is no interior resonance at frequency $\omega$, which could be circumvented by the use of combined field integral equations as in \cite{CK}. The existence and regularity of $D$ ultimately do not depend on the invertibility of $S$.}, we can let $D = S^{-1} (T - \frac{1}{2} I )$, and equivalently write (dropping the $+$ in the notation)
\begin{equation}\label{eq:dtn-abc}
\frac{\pd u}{\pd \nu} = D u, \qquad \x \in \pd \Omega.
\end{equation}
The operator $D$ is called the \emph{exterior Dirichlet-to-Neumann map} (or DtN map), because it maps the Dirichlet data $u$ to the Neumann data $\pd u/\pd \nu$ with $\nu$ pointing outward. The DtN map is independent of the right-hand side $f$ of \eqref{eq:HE} as long as $f$ is supported in $\Omega$. The notion that \eqref{eq:dtn-abc} can serve as an exact ABC was made clear in a uniform medium, e.g., in \cite{engmaj} and in \cite{kelgiv}. Equation \eqref{eq:dtn-abc} continues to hold even when $c(\mathbf{x})$ is heterogenous in the vicinity of $\partial \Omega$, provided the correct (often unknown) Green's function is used. The medium is indeed heterogeneous near $\partial \Omega$ in many situations of practical interest, such as in geophysics.

The DtN map $D$ is symmetric. The proof of the symmetry of D was shown in a slightly different setting here \cite{symm} and can be adapted to our situation. Much more is known about DtN maps, such as the many boundedness and coercivity theorems between adequate fractional Sobolev spaces (mostly in free space, with various smoothness assumptions on the boundary). We did not attempt to leverage these properties of $D$ in the scheme presented here.

We only compress the exterior DtN map in this work, and often refer to it as the DtN map for simplicity, unless there could be confusion with another concept, for example with the \emph{half-space DtN map}. We shall talk more about the half-space DtN map soon, but first, we review in the upcoming section existing methods for discrete absorbing boundary conditions.%

\section{Discrete absorbing boundary conditions}\label{sec:abc}

There are many ways to realize an absorbing boundary condition for the Helmholtz equation, and we briefly describe the main ones in this section. We start with ABCs that are surface-to-surface, and move on to ABCs which involve surrounding the computational domain $\Omega$ by an absorbing layer. The later approach is often more desirable because the parameters of the layer can usually be adjusted to obtain a desired accuracy. We then discuss the complexity of ABCs in heterogeneous media.%

\subsection{Surface-to-surface ABCs}

An early seminal work in absorbing boundary conditions was from Engquist and Majda, who in \cite{engmaj} consider the half-space problem ($x_1 \geq 0$) for the wave equation, with uniform medium, and write down the form of a general wave packet traveling to the left (towards the negative $x_1$ direction).
From this, they calculate the boundary condition which exactly annihilates those wave packets, and obtain \footnote{We are omitting details here for brevity, including a discussion of pseudo-differential operators and their symbols.}: %
\begin{equation} \label{sym}
 d/dx - i\sqrt{\omega^2-\xi^2}
\end{equation}
where $(\omega,\xi)$ are the dual variables to $(y,t)$ in Fourier space. They can then approximate the square root in various ways in order to obtain usable, i.e. local in both space and time, boundary conditions, recalling that $i\omega$ corresponds to $\pd/\pd y$ and $i\xi$ corresponds to $\pd/\pd t$.

Hagstrom and Warburton in \cite{haglap} also consider the half-space problem, take the transverse Fourier-Laplace transforms of the solution and use a given Dirichlet data on the boundary of the half-space to obtain what they call a complete wave representation of the solution, valid away from the boundary. They then use this representation to obtain approximate local boundary condition sequences to be used as ABC's. Again, their method was developed for the uniform case.

Keller and Givoli, in \cite{kelgiv}, use a different technique: they assume a circular or spherical $\Omega$, and a uniform medium outside of this $\Omega$. They can then use the Green's function of the exterior problem, which is known for a circle, in order to know the solution anywhere outside $\Omega$, given the boundary data $u$ on $\pd \Omega$. They can then differentiate this solution in the radial (which is the normal) coordinate, and evaluate it on $\pd \Omega$ to obtain the exterior DtN map. They can now use this DtN map as an ABC in a Helmholtz solver. This technique can be seen as \emph{eliminating the exterior unknowns}: as we do not care to know the solution outside of $\Omega$, we can use the information we have on the exterior solution to reduce the system to one only on the inside of $\Omega$. This meaning of \emph{eliminating the exterior unknowns} shall become more obvious when we apply this to the discretized Helmholtz equation in section \ref{sec:strip}.

Somersalo et al. (\cite{somer}) also use the DtN map, this time to solve an interior problem related to the Helmholtz equation. They use a differential equation of Riccati type to produce the DtN map. When introduce the elimination of unknowns in section \ref{sec:strip}, we shall demonstrate the connection we have found between the \emph{half-space} problem and a Riccati equation for the DtN map.

The aforementioned ABCs are not meant to be a representative sample of all the existing techniques. Unfortunately, these techniques either do not apply to heterogeneous medium, or do not perform very well in that situation. In contrast, various types of absorbing layers can be used in heterogeneous media, with caution, and are more flexible.

\subsection{Absorbing layer ABCs}\label{sec:layers}

Another approach to ABCs is to surround the domain of interest by an \emph{absorbing layer}. %
While a layer should preferably be as thin as possible, to reduce computational complexity, its design involves at least two different factors: 1) waves that enter the layer must be significantly damped before they re-enter the computational domain, and 2) reflections created when waves cross the domain-layer interface must be minimized. The Perfectly Matched Layer of B\'erenger (called PML, see \cite{berenger}) is a convincing solution to this problem in a uniform acoustic medium. Its performance often carries through in a general heterogeneous acoustic medium $c(\mathbf{x})$, though its derivation strictly speaking does not.

PML consists of analytically continuing the solution to the complex plane for points inside the layer. This means we apply a coordinate change of the type $x \rightarrow x + i \int^x \sigma(\tilde{x}) \ d\tilde{x}$, say for a layer in the positive $x$ direction, with $\sigma=0$ in the interior $\Omega$, and $\sigma$ positive and increasing inside the layer. Hence the equations are unchanged in the interior, so the solution there is the desired solution to the Helmholtz equation, but waves inside the layer will be more and more damped, the deeper they go into the layer. This is because the solution is a superposition of complex exponentials, which become decaying exponentials under the change of variables. Then we simply put zero Dirichlet boundary conditions at the end of the layer. Whatever waves reflect there will be tiny when they come back out of the layer into the interior. For a heterogeneous medium, we may still define a layer-based scheme from a transformation of the spatial derivatives which mimics the one done for the PML in a uniform medium, by replacing the Laplacian operator $\Delta$ by some $\Delta_{layer}$ inside the PML, but this layer will not be perfectly matched anymore and is called a \emph{pseudo-PML} (pPML). In this case, reflections from the interface between $\Omega$  and the layer are usually not small. It has been shown in \cite{adiabatic} that, in some cases of interest to the optics community with nonuniform media, pPML for Maxwell's equations can still work, but the layer needs to be made very thick in order to minimize reflections at the interface. In this case, the Helmholtz equation has to be solved in a very large computational domain, where most of the work will consist in solving for the pPML. In fact, the layer might even cause the solution to grow exponentially inside it, instead of forcing it to decay (\cite{diazjoly}, \cite{back}), because the group and phase velocities have an opposite sign. %

An ABC scheme which is more stable by construction is the one of Appel\"o and Colonius \cite{appcol}. They use a smooth coordinate transform to reduce an unbounded domain to a bounded one with a slowing-down layer, and damp the spurious waves thus created by artificial viscosity (high-order undivided differences). The stability of this scheme follows from its construction, so that it can be used in problems for which the pPML is unstable. However, this method is not ideal because it requires discretizing higher and higher order space derivatives in order to obtain better and better results.

\subsection{Complexity of ABCs in heterogeneous media}\label{sec:compABC}

Unfortunately, discrete absorbing layers such as the pPML may need to be quite wide in practice, or may be otherwise computationally costly (because for example of high-order artificial viscosity in \cite{appcol}). Call $L$ this width (in meters). Although this is not a limitation of the framework presented in this paper, we discretize the Helmholtz operator in the most elementary way using the standard five-point difference stencil. Put $h = 1/N$ for the grid spacing, where $N$ is the number of points per dimension for the interior problem, inside the unit square $\Omega = [0,1]^2$. While $\Omega$ contains $N^2$ points, the total number of unknowns is $O\left((N+2w)^2\right)$ in the presence of the layer, where $w=L/h$ is its width in number of grid points. In a uniform medium, the PML width $L$ needed is a fraction of the wavelength, i.e. $L \sim \lambda=\frac{2\pi}{\omega} \sim \frac{1}{N}$, so that we need a constant number of points independently of $N$: $w=L/h=LN \sim 1$. However, in nonuniform media, the heterogeneity of $c(\mathbf{x})$ can limit the accuracy of the layer. If we consider an otherwise uniform medium with an embedded scatterer outside of $\Omega$, then the pPML will have to be large enough to enclose this scatterer, no matter $N$. For more general, heterogeneous media such as the ones considered in this paper, we often observe that convergence as a function of $L$ or $w$ is delayed compared to a uniform medium. That means that we have $L \sim L_0$ so that $w \sim NL_0$ or $w = O(N)$, as we assume in the sequel.

The authors of \cite{appcol} have not investigated the computational complexity of their layer for a given accuracy, but it is clear that a higher accuracy will require higher-order derivatives for the artificial viscosity, and those are quite costly. Fortunately, the framework to be developed over the next chapters also applies to the compression of such a layer, just as it does to any other ABC.

In the case of a second-order discretization, the rate at which one must increase $N$ in order to preserve a constant accuracy in the solution, as $\omega$ grows, is about $N =O(\omega^{1.5})$. This unfortunate phenomenon, called the \emph{pollution effect}, is well-known: it begs to \emph{increase} the resolution, or number of points per wavelength, of the scheme as $\omega$ grows  \cite{nvsom,BabPollut}. As we saw, the width of the pPML may be as wide as a constant value $L_0$ independent of $N$, hence its width generally needs to scale as $O(\omega^{1.5})$ grid points.

Next, we introduce the exterior and half-space problems for the Helmholtz equation. We explain how those are related, and how knowledge from the solution of one will help us with the solution to the other.

\section{The Helmholtz equation: exterior and half-space problems}\label{sec:HEhalf}

The previous two sections addressed the fact that we wish to obtain the exterior DtN map in order to approximate the free-space solution in $\Omega$, and how to do that using ABCs and the DtN map. However, ABCs can be computationally intensive. To obtain the exterior DtN map numerically in a feasible way, we will need solve in chapter \ref{ch:probing} the exterior problem, and so we define it here. The \emph{half-space problem} for the Helmholtz equation is also interesting to us because we can write down an analytical formula for the DtN map, and use that to gain knowledge that might prove to be more general and apply to the exterior DtN map. Hence we begin by explaining the exterior problem. Then we introduce the half-space problem and its DtN map, and why this might give us insights into the exterior DtN map. We then state important results which will be used in the next chapters.

\subsection{The exterior problem}\label{sec:extprob}

The exterior problem consists of solving the free-space problem, but outside of some domain $\Omega$, given a Dirichlet boundary condition $g$ on $\pd \Omega$ and the SRC \eqref{eq:src}. That is, the following has to hold:
\begin{equation}\label{eq:HEext}
\Delta u(\x)+\frac{\omega^2}{c^2(\x)} u(\x) = 0, \ \x = (x_1, x_2) \in \Omega^c,
\end{equation}
where $\Omega^c=\mathbb{R}^2 \setminus \Omega$, with the boundary data $u(\x)=g(\x)$ for $\x \in \pd \Omega$ and a given $g$. Again, we require the SRC to hold. We call this the \emph{exterior problem}, since we solve for the solution $u$ outside of the domain $\Omega$. We are interested in this problem because, if we can find the solution $u$, then we can take its derivative, normal to $\pd \Omega$, and obtain the exterior Dirichlet-to-Neumann map. This is how we will calculate the DtN map numerically in the next chapter. Then, of course, the DtN map can be used to solve the free-space problem reformulated on $\Omega$, which is our goal. This sounds like a circular way of solving the free-space problem. This is because, as we shall see in more detail in the next chapter, solving the exterior problem a few times will give us the exterior DtN map, which will speed up all free-space solves.

In the next chapter, we will use $\Omega=\left[0,1\right]^2$, a square of side 1. For practical computations, a rectangular domain made to fit tightly around the scatterer of interest is often used, especially if we have a thin and long scatterer. As we will see in section \ref{sec:abc}, numerous ABCs have been designed for the rectangular domain, so choosing a rectangular domain here is not a limitation. Then, the numerical DtN map is a matrix which, when multiplied by a vector of Dirichlet values $u$ on $\pd \Omega$, outputs Neumann values $\pd_\nu u$ on $\pd \Omega$. In particular, we can consider the submatrix of the DtN map corresponding to Dirichlet values on one particular side of $\Omega$, and Neumann values on that same side. As we shall see next, this submatrix should be quite similar to the half-space DtN map.

\subsection{The half-space problem}\label{sec:half}

We consider again the scalar Helmholtz equation, but this time we take $\Omega$ to be the top half-space $\mathbb{R}_+^2=\left\{ (x_1,x_2) : x_2 \geq 0 \right\}$, so that now the boundary $\pd \Omega$ is the $x_1$-axis, that is, when $x_2=0$. And we consider the exterior problem for this $\Omega$: 
\begin{equation}\label{eq:hsHE}
\Delta u(\x)+\frac{\omega^2}{c^2(\x)} u(\x) = 0, \qquad \x = (x_1, x_2), \qquad x_2 \leq 0
\end{equation}
with given boundary condition $g$:
\begin{equation}\label{eq:bchalf}
u(x_1,0)=g(x_1,0),
\end{equation}
requiring some decay on $g(x_1,0)$ as $x_1 \rightarrow \pm \infty$ and the SRC \eqref{eq:src} to hold in the bottom half-space $\mathbb{R}_-^2=\left\{ (x_1,x_2) : x_2 \leq 0 \right\}$.  We shall explain how to find this analytical DtN map for the half-space problem in uniform medium, but first, we discuss the relevance of the half-space DtN map for the exterior DtN map.

Let us use again $\Omega=\left[0,1\right]^2$ as we introduced earlier, and let us call $S_1$ the bottom side of $\pd \Omega$: $S_1=\left\{(x_1,x_2): 0 \leq x_1 \leq 1, x_2=0 \right\}$. For the exterior problem, we prescribe boundary values $g$ on $\pd \Omega$, solve the exterior problem and obtain values of $u_\text{ext}$ everywhere outside of $\Omega$. From those, we know the exterior DtN map $D_\text{ext}:\pd \Omega \rightarrow \pd \Omega$. Consider now the values of $u_\text{ext}$ we have just found, along the $x_1$ axis: we can use those to define the boundary condition \eqref{eq:bchalf} of the half-space problem. The solution $u_\text{half}$ we find for this half-space problem on the bottom half-space $\mathbb{R}_-^2$ coincides with $u_\text{ext}$ on that half-space. Similarly, the exterior DtN map $D_\text{ext}$ restricted to the bottom side of $\Omega$, $D_\text{ext}: S_1 \rightarrow S_1$ coincides with the half-space DtN map $D_\text{half}$, restricted to this same side $D_\text{half}: S_1 \rightarrow S_1$.

\begin{figure}[ht]
\begin{minipage}[t]{0.48\linewidth}
\includegraphics[scale=.45]{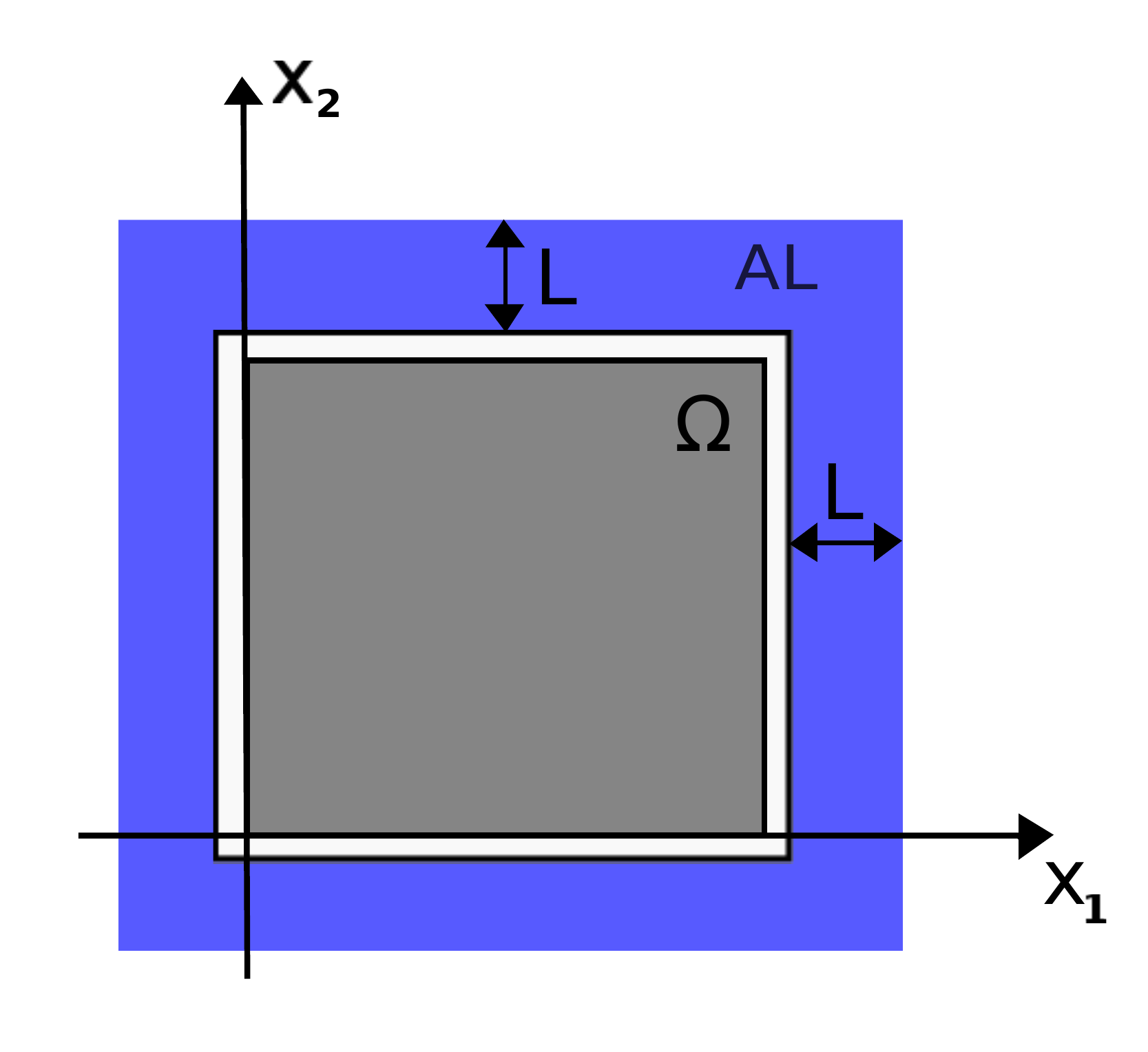}
\caption{The exterior problem: $\Omega$ is in grey, there is a thin white strip around $\Omega$, then an absorbing layer in blue.}
\label{fig:ext}
\end{minipage}
\hspace{0.1cm}
\begin{minipage}[t]{0.48\linewidth}
\includegraphics[scale=.45]{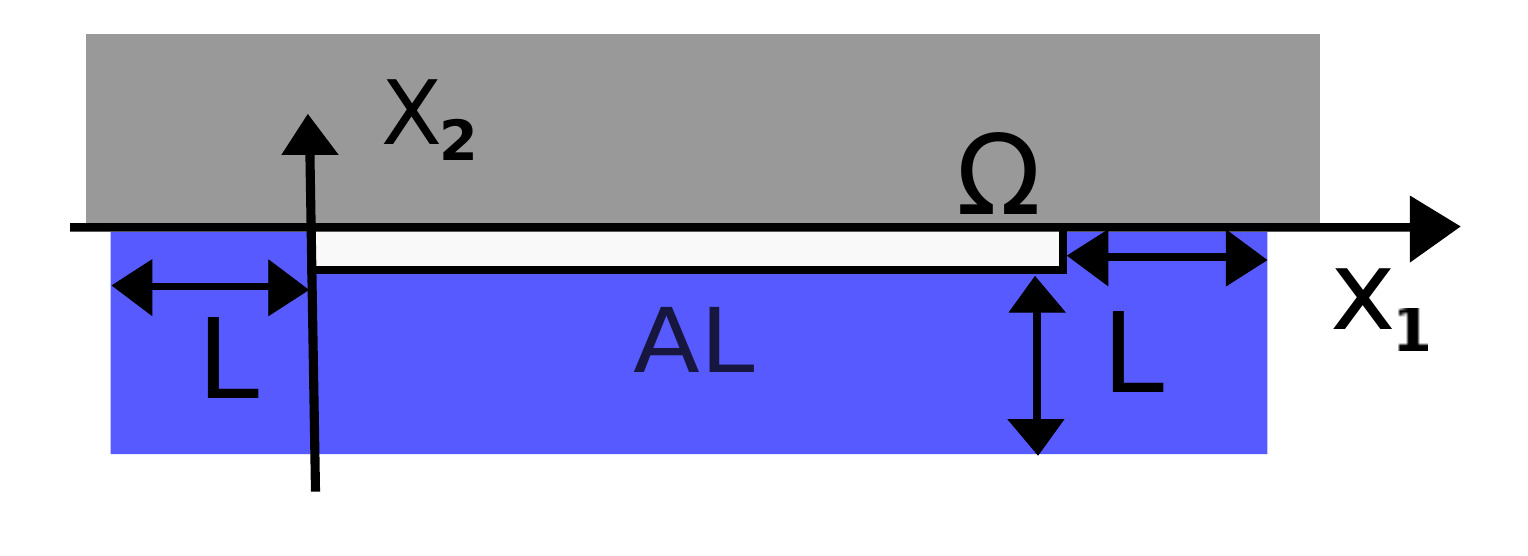}
\caption{The half-space problem: there is thin white strip below the bottom side of $\Omega$, and an absorbing layer (in blue) surrounds that strip on three sides.}
\label{fig:half}
\end{minipage}
\end{figure}

This relationship between the half-space and exterior DtN maps remains when solving those problems numerically. To solve the exterior problem numerically, we proceed very similarly to how we would for the free-space problem. As we saw in the previous section on ABCs, we may place an ABC on $\pd \Omega$ for the free-space problem. For the exterior problem, we place this ABC just a little outside of $\pd \Omega$, as in Figure \ref{fig:ext}. We then enforce the boundary condition $u(\x)=g(\x)$ for $\x$ on $\pd \Omega$, and solve the Helmholtz equation outside of $\Omega$, inside the domain delimited by the ABC. We thus obtain the solution $u$ on a thin strip just outside of $\Omega$, and we can compute from that $u$ the DtN map $\pd_\nu u$ on $\pd \Omega$. This is why we need to put the ABC just a little outside of $\Omega$, and not exactly on $\pd \Omega$.

To solve the half-space problem numerically, we will need again an ABC, in order to reformulate the problem on a smaller, bounded domain. We can prescribe values of $u$ only along the bottom edge $S_1$ as in Figure \ref{fig:half}, leave a thin strip just below that edge, surround that thin strip by an ABC on all three other sides, and solve for the solution $u$. We can then, from the solution $u$ in this thin strip, calculate an approximate half-space DtN map on $S_1$. And we see how this approximate half-space DtN map restricted to $S_1$ will be similar (but not exactly the same) to the exterior DtN map restricted to that same edge $S_1$, given the same boundary data on $S_1$ and 0 boundary data on the other edges of $\pd \Omega$. The main difference between the two maps is that in the exterior case, the two corners of $S_1$ will cause scattering, some of which will affect the solution in the top half-space $\mathbb{R}_+^2$ as well.%

It is because of this connection between the half-space and the exterior DtN maps that we analyze further the half-space DtN map, and use insights obtained from this analysis to treat the exterior DtN map.

\subsection{Solution to the half-space problem in a uniform medium using the Green's function}\label{sec:hsG}

In section \ref{sec:Gfsol} we used the free-space Green's function $G$ to obtain the solution $u$ anywhere from integrating $G$ with the right-hand side $f$. We can do the same for the half-space problem. First, we define the Green's function of the half-space problem:
\begin{definition}\label{def:Greenhalf}
The half-space \emph{Green's function} for the Helmholtz equation \ref{eq:hsHE} is the unique function $G(\x,\y)$ which solves the half-space problem with zero boundary data, that is $g=0$ in \eqref{eq:bchalf}, and with right-hand side the delta function $f(\x)=\delta(|\x-\y|)$, for every fixed $\y$.
\end{definition}
This half-space Green's function, which we shall call $G_\text{half}$, is once again well-known, and can be obtained from the free-space Green's function by the reflection principle (p. 110 of \cite{FollandIntroPDE}), with $ \x=(x_1,x_2)$, $\y=(y_1,y_2)$ and $\x'=(-x_1,x_2)$:
\begin{equation}\label{eq:Greenhalf}
G_\text{half}(\x,\y)=G(\x,\y)- G(\x',\y).
\end{equation}
Then, the solution $u$ to the half-space problem with $g=0$ is as expected:
\begin{equation}\label{eq:Gsolhalf}
u(\x)=\int_S G_\text{half}(|\x-\y|)f(\y) \ dy,  \ \x \in \mathbb{R}^2 \setminus S.
\end{equation}
where $S=\left\{(x_1,x_2):x_2=0\right\}$ is the $x_1$ axis. This half-space Green's function will be useful to us, and in particular, we are interested in the following result of \cite{Hsweep}, slightly reformulated for our purposes:
\begin{lemma}\label{Glowrank}
\emph{Theorem 2.3 of \cite{Hsweep}.} Let $G_\text{half}$ be the half-space Green's function as defined above. Let $n \in \mathbb{N}$ be some discretization parameter, $n$ even, and let $h=1/n$. Let $Y=\left\{\y_j=(jh,-h), j=1, \ldots, n/2 \right\}$ and $X=\left\{\x_j=(jh,-h), j=n/2+1, \ldots, n \right\}$. Then $\left(G_\text{half}(\x,\y)\right)_{\x \in X, \ \y \in Y}$ is numerically low-rank. More precisely, for any $\eps >0$, there exist a constant $J=O(\log k |\log\eps |^2)$, functions $\left\{\alpha_p(\x)\right\}_{1\leq p \leq J}$ for $\x \in X$ and functions $\left\{\beta_p(\y)\right\}_{1\leq p \leq J}$ for $\y \in Y$ such that
\[ \left| G(\x,\y)-\sum_{p=1}^J \alpha_p(\x)\beta_p(\y) \right| \leq \eps \ \text{for} \ \x \in X, \y \in Y.\]
\end{lemma}
We provide this lemma here because the concept of an operator being \emph{numerically low-rank} away from its diagonal will be useful later.

For now, we define the half-space DtN map. We first want to find a kernel that will give us the solution on the bottom half-space when multiplied with the boundary data $g$ and not with the right-hand side $f$, which we will soon assume to be 0.
Such a kernel exists, and is related to $G_\text{half}$. From equation (2.37) of \cite{FollandIntroPDE}, we may write
$$ u(\x)=\int_{S} \partial_{\nu_\y}G_\text{half}(\x,\y) \ u(\y) \ dS_\y, \ \x \in \mathbb{R}^2 \setminus S $$
where $\nu$ is the outward- (thus downward-) pointing normal to $S$. Hence we obtain
\begin{equation}\label{eq:halfsol}
u(\x)=-\int \left. \partial_{y_2}G_\text{half}\left(\x,(y_1,y_2)\right)\right|_{y_2=0} \ u(y_1,0) \ dy_1,  \ \x \in \mathbb{R}^2 \setminus S.
\end{equation}
We have found a different way of expressing the solution $u$ to the half-space problem using the half-space Green's function $G_\text{half}$.

\subsection{The analytical half-space DtN map}

Since we wish to consider the exterior Dirichlet-to-Neumann map for $u$ on $x_2=0$, $Du(x_1,0)=\left. \partial_{\nu_\x}u(x_1,x_2) \right|_{x_2=0}=-\left. \partial_{x_2} u(x_1,x_2) \right|_{x_2=0}$, we take the normal derivative in \eqref{eq:halfsol} and evaluate at $x_2=0$ (using the fact that $u(y_1,0)=g(y_1,0)$ is the boundary data): 
$$Du(x_1,0)= \int \left. \partial_{x_2} \partial_{y_2} G_\text{half}\left((x_1,x_2),(y_1,y_2)\right) \right|_{x_2=0,y_2=0} \ g(y_1,0) \  dy_1, \ x_1 \in \mathbb{R}.$$
Hence we have found the kernel of the half-space DtN map to be, essentially, two derivatives of the half-space Green's function. Since we know that Green's function, we can use it to obtain an analytical expression for the half-space DtN map. Notice that $\frac{\partial}{\partial z} H_0^{(1)}(z)=-H_1^{(1)}(z)$, %
We have
$$\partial_{x_2} \partial_{y_2} G(\x,\y) = \frac{i}{4} \partial_{x_2} \partial_{y_2} H_0^{(1)}(z)$$
with
$$z=k|\x-\y|=k\sqrt{(x_1-y_1)^2+(x_2-y_2)^2},$$
so that
$$\partial_{y_2} z = k^2\frac{y_2-x_2}{z}, \  \partial_{x_2} z = k^2\frac{x_2-y_2}{z}.$$
Hence
$$\partial_{y_2} G(\x,\y) = \frac{i}{4} \frac{\partial z}{\partial y_2} \partial_z H_0^{(1)}(z) =  \frac{i k^2 }{4} \frac{x_2-y_2}{z} H_1^{(1)}(z)$$
and 
\begin{equation}\label{eq:Gxy}
\partial_{x_2} \partial_{y_2} G(\x,\y) =  \frac{i k^2}{4} \left( \frac{z-k^2(x_2-y_2)^2/z}{z^2} H_1^{(1)}(z) + k^2\left(\frac{x_2-y_2}{z}\right)^2 \partial_z H_1^{(1)}(z) \right).
\end{equation}
Also, we have 
$$\partial_{x_2} \partial_{y_2} G(\x',\y) = \frac{i}{4} \partial_{x_2} \partial_{y_2} H_0^{(1)}(z')$$
with
$$z'=k\sqrt{(x_1-y_1)^2+(-x_2-y_2)^2},$$
so that
$$\partial_{y_2} z' = k^2\frac{y_2+x_2}{z'}, \  \partial_{x_2} z' = k^2\frac{y_2+x_2}{z'}.$$
Hence
$$\partial_{y_2} G(\x',\y) = \frac{i}{4} \frac{\partial z'}{\partial y_2} \partial_{z'} H_0^{(1)}(z') =  \frac{i k^2}{4} \frac{-x_2-y_2}{z'} H_1^{(1)}(z')$$
and 
\begin{equation}\label{eq:Gxpy}
\partial_{x_2} \partial_{y_2} G(\x',\y) = \frac{i k^2}{4} \left( \frac{-z'+k^2(x_2+y_2)^2/z'}{z'^2} H_1^{(1)}(z') -k^2\left(\frac{x_2+y_2}{z'}\right)^2 \partial_{z'} H_1^{(1)}(z') \right). 
\end{equation}
Now let $x_2=0,y_2=0$ but $0 \neq k|\x-\y|= k|x_1-y_1|=\left. z \right|_{x_2=0,y_2=0}=\left. z' \right|_{x_2=0,y_2=0}$ in \eqref{eq:Gxy} and \eqref{eq:Gxpy}, so that 
\begin{eqnarray*}
\partial_{x_2} \partial_{y_2} \left. G_\text{half}(\x,\y) \right|_{x_2=0,y_2=0} &=& \partial_{x_2} \partial_{y_2} \left. G(\x,\y) \right|_{x_2=0,y_2=0} - \partial_{x_2} \partial_{y_2} \left.  G(\x',\y) \right|_{x_2=0,y_2=0} \\
&=& \left. \frac{i k^2}{2k|\x-\y|} H_1^{(1)}(k|\x-\y|) \right|_{x_2=0,y_2=0} \\
&=&  \frac{i k^2}{2k|x_1-y_1|} H_1^{(1)}(k|x_1-y_1|)
\end{eqnarray*}
Thus we have found that the half-space DtN map kernel is:
\begin{equation}\label{eq:hsDtNkernel}
K(r)=\frac{i k^2}{2kr} H_1^{(1)}(kr)
\end{equation}
and thus the half-space DtN map is:
\begin{equation}\label{eq:hsDtNmap}
Du(x_1,0)= \int K(|x_1-y_1|) g(y_1,0) \ dy_1.
\end{equation}

As we shall prove in chapter \ref{ch:plr}, the half-space DtN map kernel \eqref{eq:hsDtNkernel} for a uniform medium is numerically low-rank away from its singularity, just as the half-space Green's function is from Lemma \ref{Glowrank}. This means that, if $x$ and $y$ are coordinates along the infinite boundary $S$, and $|x-y|\geq r_0$ for some constant quantity $r_0$, then the DtN map, a function of $x$ and $y$, can be approximated up to error $\eps$ as a sum of products of functions of $x$ only with functions of $y$ only (\emph{separability}) and that this sum is finite and in fact involves only a small numbers of summands (\emph{low-rank}). Hence, a numerical realization of that DtN map in matrix form should be compressible. In particular, blocks of that matrix which are away from the diagonal should have low column rank. We shall make all of this precise in chapter \ref{ch:plr}.

Recall that the goal of this thesis, to compress ABCs, will be achieved by approximating the DtN map in more general cases than the uniform half-space problem. Before we explain our approach for compressing an ABC in the next two chapters, we first explain the most straightforward way of obtaining a DtN map from an ABC, by eliminating the unknowns in the absorbing layer in order to obtain a reduced system on the interior nodes. This solution, however, is computationally impractical. It only serves to make explicit the relationship between ABCs and the DtN map.

\section{Eliminating the unknowns: from any ABC to the DtN map}\label{sec:strip}

We now wish to explain the fundamental relationship between discrete ABCs and the discrete DtN map: any discrete ABC can be reformulated as a discrete DtN map on $\pd \Omega$. We present two similar ways of obtaining that relationship using Schur complements: one for the half-space problem, and one for the free-space problem.

\subsection{Eliminating the unknowns in the half-space problem}

We consider the half-space problem, in which we care about the solution in the top half-plane. We assume from the SRC that the solution far away from $x_2=0$ in the bottom half is small. We want to eliminate unknowns from far away on the bottom side, where the solution is so small we ignore it as zero, towards the positive $x_2$ direction in order to obtain an outgoing Dirichlet-to-Neumann map in the $-x_2$ direction. We assume $f=0$ everywhere in the bottom half-plane. Let $u_1$ denote the first line of unknowns at the far bottom, $u_2$ the next line, and so on. We first use a Schur complement to eliminate $u_1$ from the discretized (using the standard five-point stencil) Helmholtz system which is as follows:
\[ \frac{1}{h^2}
\begin{pmatrix} \\
 & S_1 & P_1 & \\ 
 & & & \\
& P_1^T& C_1 & \\
& & & \\
\end{pmatrix}  \quad 
\begin{pmatrix} \\ u_{1} \\  \vdots \\ \\ \end{pmatrix}
= \begin{pmatrix} \\ 0 \\  \vdots  \\ \\ \end{pmatrix}.
\]
We then define similarly the matrices $S_k$, $P_k$ and $C_k$ corresponding to having eliminated lines $u_1$ through $u_{k-1}$ from the system. Only the upper block of $C_1$, or $C_k$ when subsequently eliminating line $u_k$, will be modified from the Schur complement. Indeed, since we eliminate from bottom to top and we have a five-point stencil, the matrices $P_k$ will be
\[ P_k=
\begin{pmatrix} & I & 0 & \cdots \\
\end{pmatrix} 
\]
and so we may write the following recursion for the Schur complements:
\begin{equation}\label{Srel}
 S_{k+1}=M-S_k^{-1},
\end{equation}
where $M$ is the main block of the 2D Helmholtz operator multiplied by $h^2$, that is,
\[ M=
\begin{pmatrix} 
& -4+h^2k^2 & 1              &        & \ & \\ 
& 1             &  -4+h^2k^2 & \ddots & \ & \\
&                & \ddots          & \ddots & \ & \\
& & & & \ & \\
\end{pmatrix}  = -2I + Lh^2.\]
Here $L$ is the 2nd order discretization of the 1D Helmholtz operator in the $x_1$ direction, $\pd_{x_1}\pd_{x_1}+k^2$. Now, at each step we denote
\begin{equation}\label{StoD}
 S_k=hD_k-I.
\end{equation}
Indeed, looking at the first block row of the reduced system at step $k-1$ we have
\[S_ku_k +Iu_{k+1}=0\]
or
\[(hD_k-I)u_k+Iu_{k+1}=0.\]
From this we obtain the DtN map $D_k$ from a backward difference:
\[ \frac{u_k-u_{k+1}}{h}=D_ku_k. \]
Now we use (\ref{StoD}) inside (\ref{Srel}) to obtain
\[ hD_{k+1}-I=M+(I-hD_k)^{-1}=M+I+hD_k+h^2D_k^2+O(h^3)\]
 or 
\[hD_{k+1}-hD_k=Lh^2+h^2D_k^2+ O(h^3), \]
which we may rewrite to obtain a discretized Riccati equation (something similar was done in \cite{keller}):
\begin{equation}\label{DR}
 \frac{D_{k+1}-D_k}{h}=L+D_k^2+ O(h),
\end{equation}
of which we may take the limit as $h \rightarrow 0$ to obtain the Riccati equation for the DtN map $D$ in the $-x_2$ direction:
\begin{equation}\label{R}
 D_{x_2}=(\pd_{x_1}\pd_{x_1}+k^2)+D^2.
\end{equation}
This equation is to be evolved in the $+x_2$ direction, starting far away in the bottom half-space. Looking at the steady state, $D_{x_2}=0$, we get back $D^2=-\pd_{x_1}\pd_{x_1}-k^2$, which is the Helmholtz equation with 0 right-hand side $f$ (which we have assumed to hold in the bottom half-space). Hence we conclude that the Riccati equation for the DtN map could be used to obtain the DtN map in the half-space case, and maybe even more complicated problems. We leave this to future work, and turn to a very similar way of eliminating unknowns, but for the exterior DtN map with $\Omega=[0,1]^2$ this time. This technique will not give rise to a Riccati equation, but will help us understand how the DtN map can be used numerically to solve the free-space Helmholtz equation reformulated on $\Omega$.

\subsection{Eliminating the unknowns in the exterior problem}

In this subsection, we assume we need an absorbing layer of large width, $w \geq N$ in number of points. We write the system for the discrete Helmholtz equation as
\begin{equation}\label{HEsys}
\begin{pmatrix} \\
 & A & P & \\ 
 & & & \\
& P^T& C & \\
& & &
\end{pmatrix}  \quad 
\begin{pmatrix} \\ u_{out} \\ \\ u_{\Omega} \\ \\ \end{pmatrix}
= \begin{pmatrix} \\ 0 \\ \\ f_\Omega \\ \\ \end{pmatrix},
\end{equation}
with $A=\Delta_{layer} + k^2 I$ and $C=\Delta + k^2 I$, with $\Delta$ overloaded to denote discretization of the Laplacian operator, and $\Delta_{layer}$ the discretization of the Laplacian operator inside the absorbing layer. We wish to eliminate the exterior unknowns $u_{out}$ from this system in order to have a new system which only depends on the interior unknowns $u_{\Omega}$. The most obvious way of eliminating those unknowns is to form the Schur complement $S=C-P^TA^{-1}P$ of $A$ by any kind of Gaussian elimination. For instance, in the standard raster scan ordering of the unknowns, the computational cost of this method\footnote{The cost of the Schur complement procedure is dominated by that of Gaussian elimination to apply $A^{-1}$ to $P$. Gaussian elimination on a sparse banded matrix of size $s$ and band $b$ is $O(sb^2)$, as can easily be infered from Algorithm 20.1 of \cite{tref}.} is $O(w^4)$ --- owing from the fact that $A$ is a sparse banded matrix of size $4Nw+4w^2$ which is $O(w^2)$, and band $N+2w$. Alternatively, elimination of the unknowns can be performed by layer-stripping, starting with the outermost unknowns from $u_{out}$, until we eliminate the layer of points that is just outside of $\pd \Omega$. The computational cost will be $O(w^4)$ in this case as well. To see this, let $u_{w}$ be the points on the outermost layer, $u_{w-1}$ the points in the layer just inside of $u_{w}$, etc. Then we have the following system:
\[
\begin{pmatrix} \\
 & A_w & P_w & \\ 
 & & & \\
& P_w^T& C_w & \\
& & & \\
\end{pmatrix}  \quad 
\begin{pmatrix} \\ u_{w} \\ \\ \vdots \\ \\ \end{pmatrix}
= \begin{pmatrix} \\ 0 \\ \\ \vdots  \\ \\ \end{pmatrix}
\]

Note that, because of the five-point stencil, $P_w$ has non-zeros exactly on the columns corresponding to $u_{w-1}$. Hence the matrix $P_w^TA_w^{-1}P_w$ in the first Schur complement $S_w=C_w-P_w^TA_w^{-1}P_w$ is non-zero exactly at the entries corresponding to $u_{w-1}$. It is then clear that, in the next Schur complement, to eliminate the next layer of points, the matrix $A_{w-1}$ (the block of $S_w$ corresponding to the points $u_{w-1}$) to be inverted will be full. For the same reason, every matrix $A_j$ to be inverted thereafter, for every subsequent layer to be eliminated, will be a full matrix. Hence at every step the cost of forming the corresponding Schur complement is at least on the order of $m^3$, where $m$ is the number of points in that layer. Hence the total cost of eliminating the exterior unknowns by layer stripping is approximately
\[ \sum_{j=1}^{w} (4(N+2j))^3 = O(w^4). \]

Similar arguments can be used for the Helmholtz equation in 3 dimensions. In this case, the computational complexity of the Schur complement or layer-stripping methods would be $O(w^3 (w^2)^2)=O(w^7)$ or $~\sum_{j=1}^{w} (6(N+2j)^2)^3=O(w^7)$. Therefore, direct elimination of the exterior unknowns is quite costly. Some new insight will be needed to construct the DtN map more efficiently.

We now remark that, whether we eliminate exterior unknowns in one pass or by layer-stripping, we obtain a reduced system. It looks just like the original Helmholtz system on the interior unknowns $u_{\Omega}$, except for the top left block, corresponding to $u_{0}$ the unknowns on $\pd \Omega$, which has been modified by the elimination procedure. Hence with the help of some dense matrix $D$ we may write the reduced, $N^2$ by $N^2$ system as
\begin{equation}\label{HEred}
Lu=
\begin{pmatrix} 
(hD-I)/h^2  & I/h^2 & 0 & & \cdots & \\ 
 &  &  & \\
 I/h^2 &  & & \\
       & & & \\
0 & & [ \; \Delta + k^2 I \; ] & \\
& & & \\
\vdots & & & \\
\\
\end{pmatrix} \quad 
\begin{pmatrix} u_{0} \\  \\   u_{-1}  \\ \\   u_{-2}  \\ \\ \vdots \\ \\ \end{pmatrix} 
= \begin{pmatrix}  0  \\  \\ f_{-1} \\ \\ f_{-2} \\ \\ \vdots \\ \\ \end{pmatrix} 
\end{equation}
and we have thus obtained an absorbing boundary condition which we may use on the boundary of $\Omega$, independent of the right-hand side $f$. Indeed, if we call $u_{-1}$ the first layer of points inside $\Omega$, we have $ (I-hD)u_{0} = u_{-1} $, or
\[
 \frac{u_{0} - u_{-1}}{h}=Du_{0} ,
\]
a numerical realization of the DtN map in \eqref{eq:dtn-abc}, using the ABC of choice. Indeed, elimination can be used to reformulate any computationally intensive ABC, not just absorbing layers, into a realization of \eqref{eq:dtn-abc}. Any ABC is equivalent to a set of equations relating unknowns on the surface to unknowns close to the surface, and possibly auxiliary variables. Again, elimination can reduce those equations to relations involving only unknowns on the boundary and on the first layer inside the boundary, to obtain a numerical DtN map $D$. A drawback is that forming this matrix $D$ by elimination is prohibitive, as we have just seen.%

As for the complexity of solving the Helmholtz equation, reducing the ABC confers the advantage of making the number of nonzeros in the matrix $L$ (of Section \ref{sec:strip}) independent of the width of the absorbing layer or complexity of the ABC. After elimination of the layer, it is easy to see that $L$ has about $20N^2$ nonzero entries, instead of the $5N^2$ one would expect from a five-point stencil discretization of the Helmholtz equation, because the matrix $D$ (part of a small block of $L$) is full. Although obtaining a fast matrix-vector product for our approximation of $D$ could reduce the application cost of $L$ from $20N^2$ to something closer to $5N^2$, it should be noted that the asymptotic complexity does not change -- only the constant does.%

This thesis addresses those two problems, obtaining the DtN map and applying it fast. The next chapter, chapter \ref{ch:probing}, suggests adapting the framework of matrix probing in order to obtain $D$ in reasonable complexity. Subsequently, chapter \ref{ch:plr} presents a compression method which leads to a fast application of the DtN map.

%% file: probing2.tex
\chapter{Matrix probing for expanding the Dirichlet-to-Neumann map}\label{ch:probing}

Recall that the goal of this thesis is to introduce a new compression scheme for ABCs. This scheme consists of two steps:
\begin{enumerate}
\item a precomputation sets up an expansion of the Dirichlet-to-Neumann map, then 
\item a fast algorithm is used to apply the DtN map in a Helmholtz solver.
\end{enumerate}
This chapter is concerned with the first step of this procedure, namely, setting up an expansion of the exterior DtN map in a precomputation. This will pave the way for compression in step two, presented in the next chapter.

The main strategy we use in this chapter is matrix probing, we introduce it in section \ref{sec:introprobe}. For matrix probing to be an efficient expansion scheme, we need to carefully choose the basis for this expansion. We present our choices and their rationales in section \ref{sec:basis}. In particular, inverse powers multiplied by a complex exponential work quite well as kernels for the basis. We then present a detailed study of using matrix probing to expand the DtN map in various different media, use that expansion to solve the Helmholtz equation, and document the complexity of the method, all in section \ref{sec:numexp}. Then, we prove in section \ref{sec:BasisPf} a result on approximating the half-space DtN map with the particular set of functions mentioned before, inverse powers multiplied by a complex exponential. We also present a numerical confirmation of that result in section \ref{sec:NumPf}. %

\section{Introduction to matrix probing}\label{sec:introprobe}

The idea of matrix probing is that a matrix $D$ with adequate structure can sometimes be recovered from the knowledge of a fixed, small number of matrix-vector products $Dg_j$, where $g_j$ are typically random vectors. In the case where $D$ is the numerical DtN map (with a slight abuse of notation), each $g_j$ consists of Dirichlet data on $\pd \Omega$, and each application $Dg_j$ requires solving an exterior Helmholtz problem to compute the derivative of the solution normal to $\pd \Omega$. We first explain how to obtain the matrix-vector multiplication of the DtN map with any vector, without having to use the costly procedure of layer-stripping. We then introduce matrix probing.

\subsection{Setup for the exterior problem}\label{sec:ext}

Recall the exterior problem of section \ref{sec:extprob}: solving the heterogeneous-medium Helmholtz equation at frequency $\omega$, outside $\Omega=[0,1]^2$, with Dirichlet boundary condition $u=g$ on $\pd \Omega$. This problem is solved numerically with the five-point stencil of finite differences (FD), using $h$ to denote the grid spacing and $N$ the number of points across one dimension of $\Omega$. We use a Perfectly Matched Layer (PML) or pPML, introduced in section \ref{sec:layers}, as our ABC. The layer starts at a fixed, small distance away from $\Omega$, so that we keep a small strip around $\Omega$ where the equations are unchanged. Recall that the width of the layer is in general as large as $O(\omega^{1.5})$ grid points. We number the edges of $\pd \Omega$ counter-clockwise starting from $(0,0)$, hence side 1 is the bottom edge $(x,0)$, $0\leq x \leq 1$, side 2 is the right edge, etc. The exterior DtN map for this problem is defined from $\pd \Omega$ to itself. Thus its numerical realization, which we also call $D$ by a slight abuse of notation, has a $4\times 4$ block structure. Hence the numerical DtN map $D$ has 16 sub-blocks, and is $n \times n$ where $n=4N$. As an integral kernel, $D$ would have singularities at the junctions between these blocks (due to the singularities in $\pd \Omega)$, so we shall respect this feature by probing $D$ sub-block by sub-block. We shall denote a generic such sub-block by $M$, or as the $(i_M,j_M)$ sub-block of $D$, referring to its indices in the $4 \times 4$ sub-block structure.

The method by which the system for the exterior problem is solved is immaterial in the scope of this paper, though for reference, the experiments in this paper use UMFPACK's sparse direct solver \cite{UMFPACK}. For treating large problems, a better solver should be used, such as the sweeping preconditioner of Engquist and Ying \cite{Hsweep,Msweep}, the shifted Laplacian preconditioner of Erlangga \cite{erlangga}, the domain decomposition method of Stolk \cite{stolk}, or the direct solver with spectral collocation of Martinsson, Gillman and Barnett \cite{dirfirst,dirstab}. This in itself is a subject of ongoing research which we shall not discuss further. 

For a given boundary condition $g$, we solve the system and obtain a solution $u$ in the exterior computational domain. In particular we consider $u_{1}$, the solution in the layer just outside of $\pd \Omega$. We are using the same notation as in section \ref{sec:strip}, where as we recall $u_0$ was the solution on the boundary, hence here $u_0=g$. We know from Section \ref{sec:strip} that $u_1$ and $g$ are related by
\begin{equation}\label{eq:D}
 \frac{u_{1} - g}{h}=Dg
\end{equation}
The matrix $D$ that this relation defines needs not be interpreted as a first-order approximation of the continuous DtN map: it is the algebraic object of interest that will be ``probed" from repeated applications to different vectors $g$.

Similarly, for probing the $(i_M,j_M)$ block $M$ of $D$, one needs matrix-vector products of $D$ with vectors $g$ of the form $[z, 0, 0, 0]^T$, $[0, z, 0, 0]^T$, etc., to indicate that the Dirichlet boundary condition is $z$ on the side indexed by $j_M$, and zero on the other sides. The application $Dg$ is then restricted to side $i_M$.

\subsection{Matrix probing}\label{sec:probe}

The dimensionality of $D$ needs to be limited for recovery from a few $Dg_j$ to be possible, but matrix probing is \emph{not} an all-purpose low-rank approximation technique. Instead, it is the property that $D$ has an efficient representation in some adequate pre-set basis that makes recovery from probing possible. As opposed to the randomized SVD method which requires the number of matrix-vector applications to be greater than the rank \cite{Halko-randomSVD}, matrix probing can recover interesting structured operators from a single matrix-vector application \cite{Chiu-probing, Demanet-probing}.

We now describe a model for $M$, any $N \times N$ block of $D$, that will sufficiently lower its dimensionality to make probing possible. Assume we can write $M$ as
\begin{equation}\label{eq:Dexp}
M \approx \sum_{j=1}^p c_j B_j
\end{equation}
where the $B_j$'s are fixed, known basis matrices, that need to be chosen carefully in order to give an accurate approximation of $M$. In the case when the medium $c$ is uniform, we typically let $B_j$ be a discretization of the integral kernel
\begin{equation}\label{eq:Bj}
B_j(x,y)= \frac{e^{ik|x-y|}}{(h+|x-y|)^{j/2}},
\end{equation}
where again $h=1/N$ is the discretization parameter. We usually add another index to the $B_j$, and a corresponding multiplicative factor, to allow for a smooth dependence on $x+y$ as well. We shall further detail our choices and discuss their rationales in Section \ref{sec:basis}. For now, we note that the advantage of the specific choice of basis matrix \eqref{eq:Bj}, and its generalizations explained in Section \ref{sec:basis}, is that it results in accurate expansions with a number of parameters $p$ which is ``essentially independent" of $N$, namely that grows either logarithmically in $N$, or at most like a very sublinear fractional power law (such as $N^{0.12}$, see section \ref{sec:pwN}). This is in sharp contrast to the scaling for the layer width, $w = O(N)$ grid points, discussed earlier. The form of $B_j$ suggested in equation \eqref{eq:Bj} is motivated by the fact that they provide a good expansion basis for the uniform-medium half-space DtN map in $\R^2$. This will be proved in section \ref{sec:BasisPf}.

Given a random vector $z^{(1)} \sim N(0,I_N)$ (other choices are possible), the product $w^{(1)}=Mz^{(1)}$ and the expansion \eqref{eq:Dexp}, we can now write
\begin{equation}\label{mp}
w^{(1)}=Mz^{(1)} \approx \sum_{j=1}^p c_j B_j z^{(1)} = \Psi_{z^{(1)}} \, \bc.
\end{equation}
Multiplying this equation on the left by the pseudo-inverse of the $N$ by $p$ matrix $\Psi_{z^{(1)}}$ will give an approximation to $\bc$, the coefficient vector for the expansion \eqref{eq:Dexp} of $M$. More generally, if several applications $w^{(j)} = M z^{(j)}$, $j = 1,\ldots, q$ are available, a larger system is formed by concatenating the $\Psi_{z^{(j)}}$ into a tall-and-thin $Nq$ by $p$ matrix ${\bm \Psi}$. The computational work is dominated, here and in other cases \cite{Chiu-probing, Demanet-probing}, by the matrix-vector products $Dg^{(1)}$, or $Mz^{(j)}$. Note that both $\Psi_{z^{(j)}}$ and the resulting coefficient vector $\bc$ depend on the vectors $z^{(j)}$. In the sequel we let $z^{(j)}$ be gaussian iid random.

In a nutshell, recovery of $\bc$ works under mild assumptions on $B_j$, and when $p$ is a small fraction of $Nq$ up to log factors. In order to improve the conditioning in taking the pseudo-inverse of the matrix $\Psi_z$ and reduce the error in the coefficient vector $\bc$, one may use $q > 1$ random realizations of $M$. %
There is a limit to the range of $p$ for which this system is well-posed: past work by Chiu and Demanet \cite{Chiu-probing} covers the precise conditions on $p$, $N$, and the following two parameters, called \emph{weak condition numbers}, for which recoverability of $\bc$ is accurate with high probability.

\begin{definition}
\emph{Weak condition number $\lambda$.} 
\[ \lambda = \max_j \frac{\| B_j \|_2 \sqrt{N}}{\| B_j \|_F} \]
\end{definition}
\begin{definition}\label{kap}
\emph{Weak condition number $\kappa$.} 
\[ \kappa = \mbox{cond}( {\bf B}), \ {\bf B}_{j \ell} = \mbox{Tr} \, (B_j^T B_\ell)\]
\end{definition}

It is desirable to have a small $\lambda$, which translates into a high rank condition on the basis matrices, and a small $\kappa$, which translates into a Riesz basis condition on the basis matrices. Having small weak condition numbers will guarantee a small failure probability of matrix probing and a bound on the condition number of ${\bf \Psi}$, i.e. guaranteed accuracy in solving for $\bc$. Also, using $q > 1$ allows to use a larger $p$, to achieve greater accuracy. These results are contained in the following theorem.
\begin{theorem} (Chiu-Demanet, \cite{Chiu-probing}) Let $z$ be a Gaussian i.i.d. random vector of length $qN$, and ${\bf \Psi}$ as above. Then $\mbox{cond}({\bf \Psi}) \leq 2\kappa + 1$ with high probability provided that $p$ is not too large, namely 
\[
q N \geq C \, p \, (\kappa \lambda \log N)^2,
\]
for some number $C > 0$.
\end{theorem}

As noted previously, the work necessary for probing the matrix $M$ is on the order of $q$ solves of the original problem. Indeed, computing $Mz^{(1)}, \ldots , Mz^{(q)}$ means solving $q$ times the exterior problem with the AL. This is roughly equivalent to solving the original Helmholtz problem with the AL $q$ times, assuming the AL width $w$ is at least as large as $N$. Then, computing the $qp$ products of the $p$ basis matrices with the $q$ random vectors amounts to a total of at most $qpN^2$ work, or less if the basis matrices have a fast matrix-vector product. And finally, computing the pseudo-inverse of ${\bf \Psi}$ has cost $Nqp^2$. Hence, as long as $p,q \ll N$, the dominant cost of matrix probing\footnote{We will see later that we also need to perform a QR factorization on the basis matrices, and this has cost $N^2p^2$. This precomputation has a cost similar or smaller to the cost of an exterior solve using current Helmholtz solvers. It might also be possible to not need a QR factorization if basis matrices closer to orthonormal are used.} comes from solving $q$ times the exterior problem with a random Dirichlet boundary condition. In our experiments, $q=O(1)$ and $p$ can be as large as a few hundreds for high accuracy.

Finally, we note that the information from the $q$ solves can be re-used for any other block which is in the same block column as $M$. However, if it is needed to probe blocks of $D$ which are not all in the same block column, then another $q$ solves need to be performed, with a Dirichlet boundary condition on the appropriate side of $\partial \Omega$. This of course increases the total number of solves. Another option would be to probe all of $D$ at once, using a combination of basis matrices that have the same size as $D$, but that are 0 except on the support of each distinct block in turn. In this case, $\kappa$ remains the same because we still orthogonalize our basis matrices, but $\lambda$ doubles ($\| B_j \|_2 $ and $\| B_j \|_F$ do not change but $N \rightarrow 4N$) and this makes the conditioning worse, in particular a higher value of $q$ is needed for the same accuracy, given by $p$. Hence we have decided not to investigate further this approach, which might become more advantageous in the case of a more complicated polygonal domain.

\subsection{Solving the Helmholtz equation with a compressed ABC}

Once we have obtained approximations $\tilde{M}$ of each block $M$ in compressed form through the coefficients $\bc$ using matrix probing, we construct block by block the approximation $\tilde{D}$ of $D$ and use it in a solver for the Helmholtz equation on the domain $\Omega=[0,1]^2$, with the boundary condition
$$\frac{\pd u}{\pd \nu}  = \tilde{D}u , \qquad x \in \pd \Omega.$$

\section{Choice of basis matrices for matrix probing}\label{sec:basis}

The essential information of the DtN map needs to be summarized in broad strokes in the basis matrices $B_j$, with the details of the numerical fit left to the probing procedure. In the case of $D$, most of its physics is contained in its \emph{diagonal singularity} and \emph{oscillations}, as predicted by geometrical optics.

A heuristic argument to obtain the form of $D$ starts from the Green's formula \eqref{eq:GRF}, that we differentiate one more time in the normal direction. After accounting for the correct jump condition, we get an alternative Steklov-Poincare identity, namely
\[
D = (T^* + \frac{1}{2} I)^{-1} H,
\]
where $H$ is the hypersingular integral operator with kernel $\frac{\pd^2 G}{\pd \nu_{\x} \pd \nu_{\y}}$, again $G(\x,\y)$ is the free-space Green's function and $\nu_{\x}$, $\nu_{\y}$ are the normals to $\pd \Omega$ in $\x$ and $\y$ respectively. The presence of $(T^* + \frac{1}{2} I)^{-1}$ is somewhat inconsequential to the form of $D$, as it involves solving a well-posed second-kind integral equation. As a result, the properties of $D$ are qualitatively similar to those of $H$. (The exact construction of $D$ from $G$ is of course already known in a few special cases, such as the uniform medium half-space problem considered earlier.)

\subsection{Oscillations and traveltimes for the DtN map}

Geometrical optics will reveal the form of $G$. In a context where there is no multi-pathing, that is, where there is a single traveltime $\tau(\x,\y)$ between any two points $\x,\y \in \Omega$, one may write a high-$\omega$ asymptotic series for $G$ as
\begin{equation}\label{eq:geoopts}
 G(\x,\y) \sim e^{i\omega \tau(\x,\y)} \sum_{j\geq 0} A_j(\x,\y) \omega^{-j},
\end{equation}
$\tau(\x,\y)$ is the traveltime between points $\x$ and $\y$, found by solving the Eikonal equation
\begin{equation} \label{eq:tau}
 \| \nabla_{\x} \tau(\x,\y) \| = \frac{1}{c(\x)},
\end{equation}
and the amplitudes $A_j$ satisfy transport equations. In the case of multi-pathing (possible multiple traveltimes between any two points), the representation \eqref{eq:geoopts} of $G$ becomes instead
 \[ 
 G(\x,\y) \sim \sum_j e^{ i \omega \tau_j(\x,\y)} \sum_{k \geq 0} A_{jk}(\x,\y) \omega^{-k}, 
 \]
where the $\tau_j$'s are the traveltimes, each obeying \eqref{eq:tau} away from caustic curves. The amplitudes are singular at caustic curves in addition to the diagonal $\x=\y$, and contain the information of the Maslov indices. Note that traveltimes are symmetric: $\tau_j(\x,\y)=\tau_j(\y,\x)$, and so is the kernel of $D$.

The singularity of the amplitude factor in \eqref{eq:geoopts}, at $\x = \y$, is $O \left( \log | \x - \y| \right)$ in 2D and $O \left( | \x - \y |^{-1} \right)$ in 3D. After differentiating twice to obtain $H$, the homogeneity on the diagonal becomes $O \left( | \x - \y|^{-2} \right)$ in 2D and $O \left( | \x - \y |^{-3} \right)$ in 3D. For the decay at infinity, the scalings are different and can be obtained from Fourier analysis of square root singularities; the kernel of $H$ decays like $O \left(| \x - \y|^{-3/2} \right)$ in 2D, and $O \left(| \x - \y|^{-5/2} \right)$ in 3D. In between, the amplitude is smooth as long as the traveltime is single-valued.

As mentioned before, much more is known about DtN maps, such as boundedness and coercivity theorems. Again, we did not attempt to leverage these properties of $D$ in the scheme presented here.

For all these reasons, we define the basis matrices $B_j$ as follows. Assume $\tau$ is single-valued. In 1D, denote the tangential component of $\x$ by $x$, and similarly that of $\y$ by $y$, in coordinates local to each edge with $0 \leq x,y \leq 1$. Each block $M$ of $D$ relates to a couple of edges of the square domain. Let $j = (j_1, j_2)$ with $j_1, j_2$ nonnegative integers. The general forms that we consider are
\[
\beta_j(x,y) = e^{i \omega \tau(x,y)} (h + |x-y|)^{-\frac{j_1}{\alpha}} (h + \theta(x,y))^{-\frac{j_2}{\alpha}}
\]
and
\[
\beta_j(x,y) = e^{i \omega \tau(x,y)} (h + |x-y|)^{-\frac{j_1}{\alpha}} (h + \theta(x,y))^{j_2},
\]
where again $h$ is the grid spacing of the FD scheme, and $\theta(x,y)$ is an adequate function of $x$ and $y$ that depends on the particular block of interest. The more favorable choices for $\theta$ are those that respect the singularities created at the vertices of the square; we typically let $\theta(x,y) = \min(x+y, 2-x-y)$. The parameter $\alpha$ can be taken to be equal to 2, a good choice in view of the numerics and in the light of the asymptotic behaviors on the diagonal and at infinity discussed earlier. 

If several traveltimes are needed for geometrical reasons, then different sets of $\beta_j$ are defined for each traveltime. (More about this in the next subsection.) The $B_j$ are then obtained from the $\beta_j$ by QR factorization within each block\footnote{Whenever a block of $D$ has symmetries, we enforce those in the QR factorization by using appropriate weights on a subset of the entries of that block. This also reduces the complexity of the QR factorization.}, where orthogonality is defined in the sense of the Frobenius inner product $\< A, B \> = \tr(A B^T)$. This automatically sets the $\kappa$ number of probing to 1.

In many of our test cases it appears that the ``triangular" condition $j_1 + 2 j_2 < $ \emph{constant} works well. The number of couples $(j_1,j_2)$ satisfying this relation will be $p/T$, where $p$ is the number of basis matrices in the matrix probing algorithm and $T$ is the number of distinct traveltimes. The eventual ordering of the basis matrices $B_j$ respects the increase of $j_1 + 2 j_2$.

\subsection{More on traveltimes}\label{sec:tt}

Determining the traveltime(s) $\tau(\x,\y)$ is the more ``supervised" part of this method, but is needed to keep the number $p$ of parameters small in the probing expansion. A few different scenarios can arise.

\begin{itemize}
\item In the case when $\nabla c(\x)$ is perpendicular to a straight segment of the boundary, locally, then this segment is itself a ray and the waves can be labeled as interfacial, or ``creeping". The direct traveltime between any two points $\x$ and $\y$ on this segment is then simply given by the line integral of $1/c(\x)$. An infinite sequence of additional interfacial waves result from successive reflections at the endpoints of the segment, with traveltimes predicted as follows.

We still consider the exterior problem for $[0,1]^2$. We are interested in the traveltimes between points $\x, \y$ on the same side of $\pd \Omega$ -- for illustration, let $\x=(x,0)$ and $\y=(y,0)$ on the bottom side of $\Omega=[0,1]^2$, with $x \leq y$ (this is sufficient since traveltimes are symmetric). Assume that all the waves are interfacial. The first traveltime $\tau_1$ corresponds to the direct path from $\x$ to $\y$. The second arrival time $\tau_2$ will be the minimum traveltime corresponding to: either starting at $\x$, going left, reflecting off of the $(0,0)$ corner, and coming back along the bottom side of $\pd \Omega$, past $\x$ to finally reach $\y$; or starting at $\x$, going past $\y$, reflecting off of the $(1,0)$ and coming straight back to $\y$. The third arrival time $\tau_3$ is the maximum of those two choices. The fourth arrival time then corresponds to starting at $\x$, going left, reflecting off of the $(0,0)$ corner, travelling all the way to the $(1,0)$ corner, and then back to $\y$. The fifth arrival time corresponds to leaving $\x$, going to the $(1,0)$ corner this time, then back to the $(0,0)$ corner, then on to $\y$. And so on. To recap, we have the following formulas:
\begin{eqnarray*}
\tau_1(\x,\y)&=& \int_x^y \frac{1}{c(t,0)} \ dt, \\
\tau_2(\x,\y)&=& \tau_1(\x,\y) + 2\min \left( \int_0^x \frac{1}{c(t,0)} \ dt, \int_y^1 \frac{1}{c(t,0)} \ dt \right), \\
\tau_3(\x,\y)&=& \tau_1(\x,\y) + 2\max \left( \int_0^x \frac{1}{c(t,0)} \ dt, \int_y^1 \frac{1}{c(t,0)} \ dt \right) = 2\int_0^1 \frac{1}{c(t,0)} \ dt - \tau_2(\x,\y), \\
\tau_4(\x,\y)&=& 2\int_0^1 \frac{1}{c(t,0)} \ dt - \tau_1(\x,\y), \\
\tau_5(\x,\y)&=& 2\int_0^1 \frac{1}{c(t,0)} \ dt + \tau_1(\x,\y), \qquad \mbox{etc.} \\
\end{eqnarray*}
All first five traveltimes can be expressed as a sum of $\pm \tau_1$, $\pm \tau_2$ and the constant phase $2\int_0^1 \frac{1}{c(t,0)} \ dt$, which does not depend on $\x$ or $\y$. In fact, one can see that any subsequent traveltime corresponding to traveling solely along the bottom boundary of $\pd \Omega$ should be again a combination of those quantities. This means that if we use $\pm \tau_1$ and $\pm \tau_2$ in our basis matrices, we are capturing all the traveltimes relative to a single side, which helps to obtain higher accuracy for probing the diagonal blocks of $D$.

This simple analysis can be adapted to deal with creeping waves that start on one side of the square and terminate on another side, which is important for the nondiagonal blocks of $D$.

\item In the case when $c(\x)$ increases outward in a smooth fashion, we are also often in presence of body waves, going off into the exterior and coming back to $\pd \Omega$. The traveltime for these waves needs to be solved either by a Lagrangian method (solving the ODE for the rays), or by an Eulerian method (solving the Eikonal PDE shown earlier). In this paper we used the fast marching method of Sethian \cite{sethart} to deal with these waves in the case that we label ``slow disk" in the next section. 
\item In the case when $c(\x)$ has singularities in the exterior domain, each additional reflection creates a traveltime that should (ideally) be predicted. Such is the case of the ``diagonal fault" example introduced in the next section, where a straight jump discontinuity of $c(\x)$ intersects $\pd \Omega$ at a non-normal angle:  we can construct by hand the traveltime corresponding to a path leaving the boundary at $\x$, reflecting off of the discontinuity and coming back to the boundary at $\y$.  More precisely, we consider again $\x=(x,0)$, $\y=(y,0)$ and $x \leq y$, with $x$ larger than or equal to the $x$ coordinate of the point where the reflector intersects the bottom side of $\pd \Omega$. We then  reflect the point $\y$ across the discontinuity into the new point $\y'$, and calculate the Euclidean distance between $\x$ and $\y'$. To obtain the traveltime, we then divide this distance by the value $c(\x)=c(\y)$ of $c$ on the right side of the discontinuity, assuming that value is constant. This body traveltime is used in the case of the ``diagonal fault", replacing the quantity $\tau_2$ that was described above. This increased accuracy by an order of magnitude, as mentioned in the numerical results of the next section.
\end{itemize}

\section{Numerical experiments}\label{sec:numexp}

\begin{figure}[H]
\begin{minipage}[t]{0.30\linewidth}
\includegraphics[scale=.33,trim = 7mm 2mm 15mm 0mm, clip]{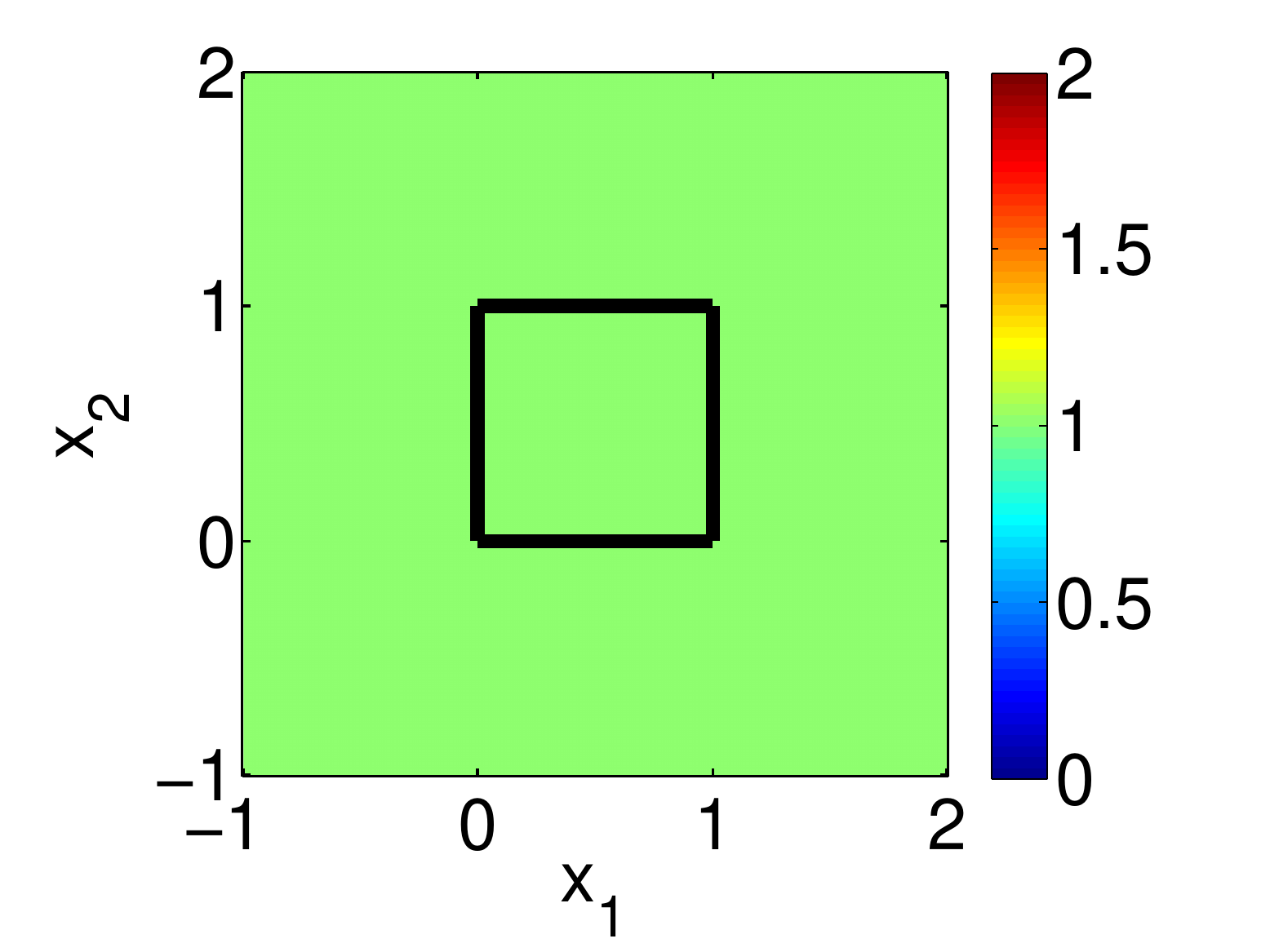}
\caption{Color plot of $c(\x)$ for the uniform medium.}\label{c1}
\includegraphics[scale=.33,trim = 7mm 2mm 15mm 0mm, clip]{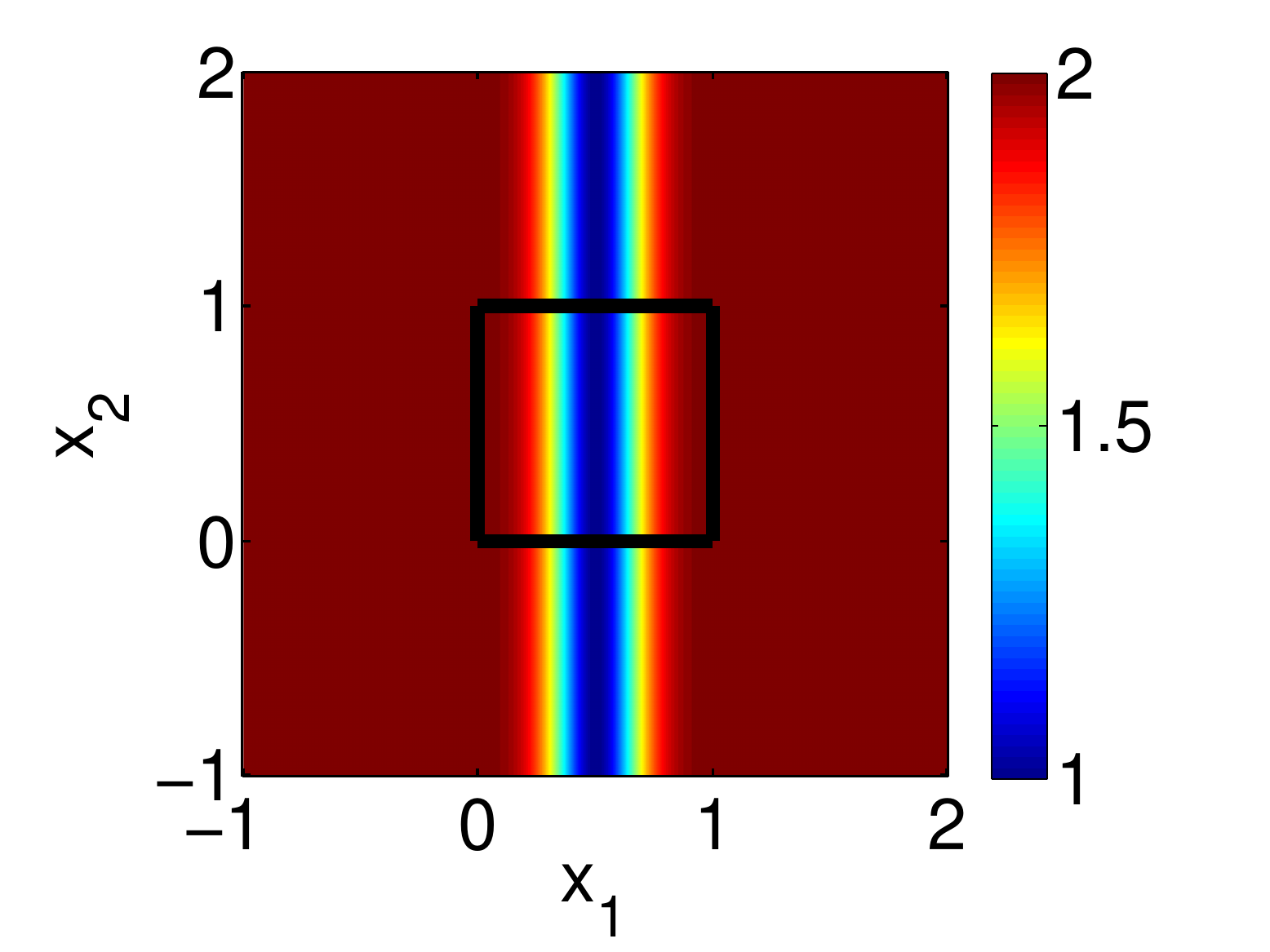}
\caption{Color plot of $c(\x)$ for the Gaussian waveguide.}\label{wg}
\end{minipage}
\begin{minipage}[t]{0.05\linewidth}
\end{minipage}
\begin{minipage}[t]{0.30\linewidth}
\includegraphics[scale=.33,trim = 7mm 2mm 15mm 0mm, clip]{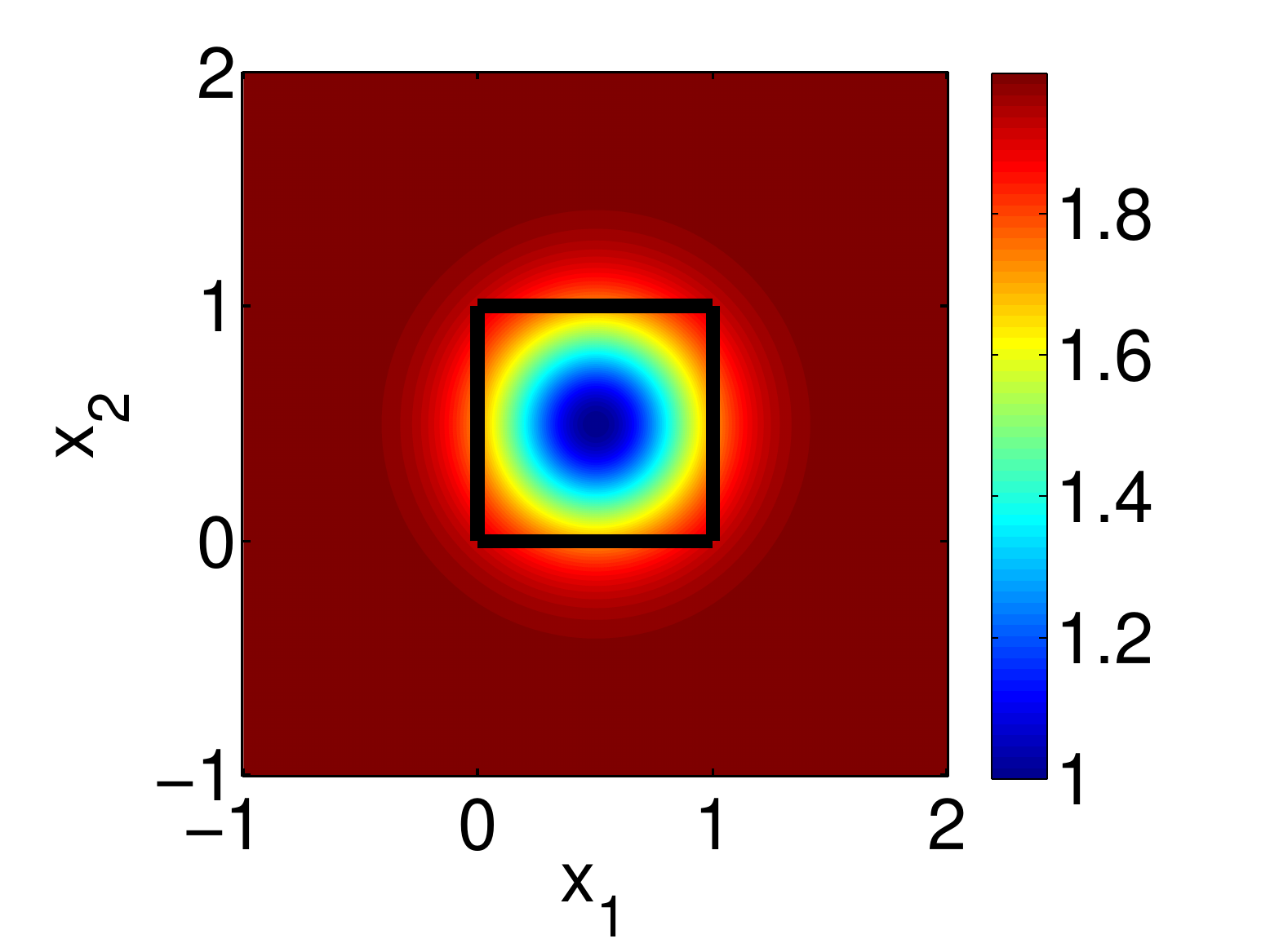}
\caption{Color plot of $c(\x)$ for the Gaussian slow disk.}\label{slow}
\includegraphics[scale=.33,trim = 7mm 2mm 15mm 0mm, clip]{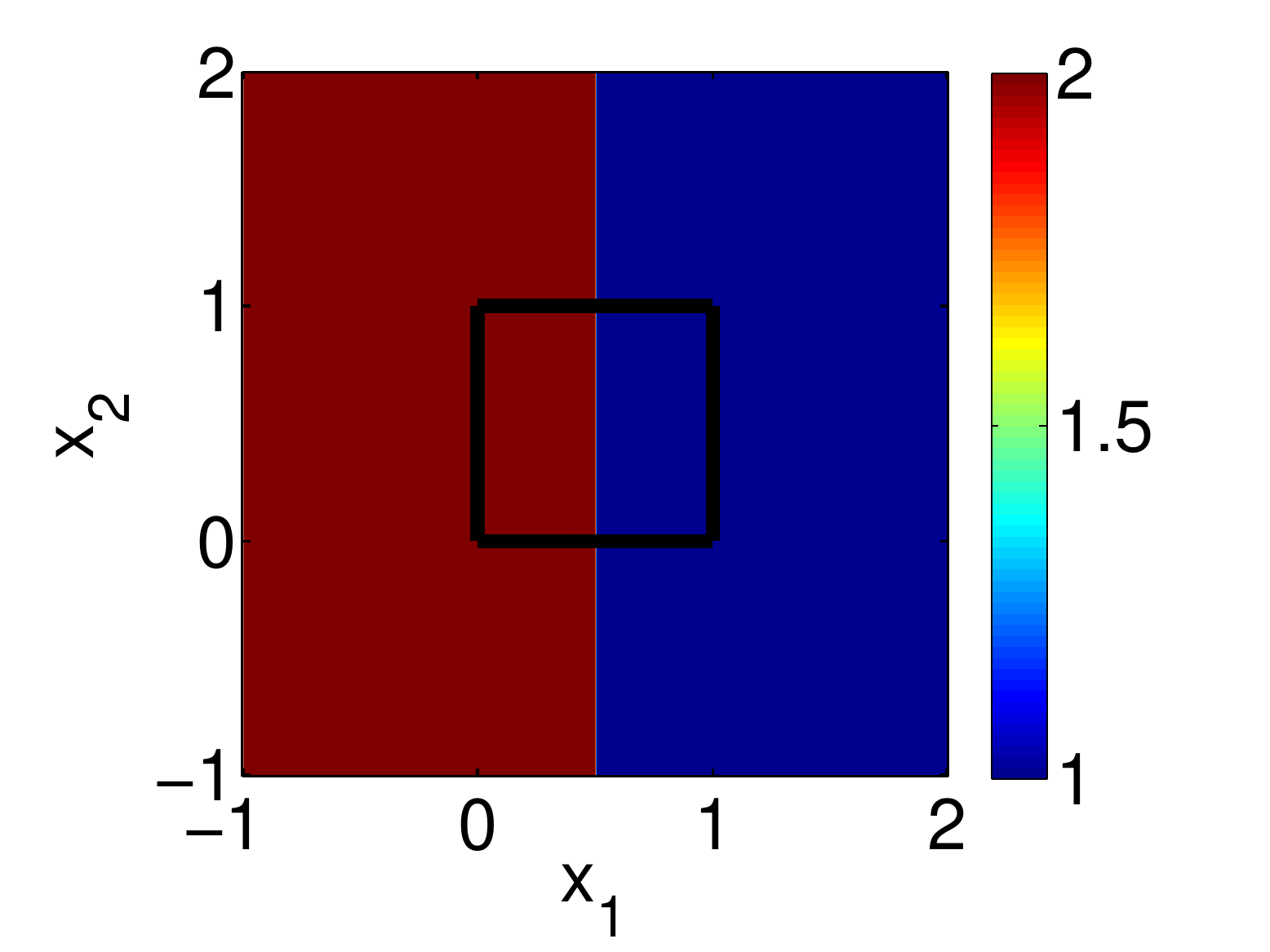}
\caption{Color plot of $c(\x)$ for the vertical fault.}\label{fault}
\end{minipage}
\begin{minipage}[t]{0.05\linewidth}
\end{minipage}
\begin{minipage}[t]{0.30\linewidth}
\includegraphics[scale=.33,trim = 7mm 2mm 15mm 0mm, clip]{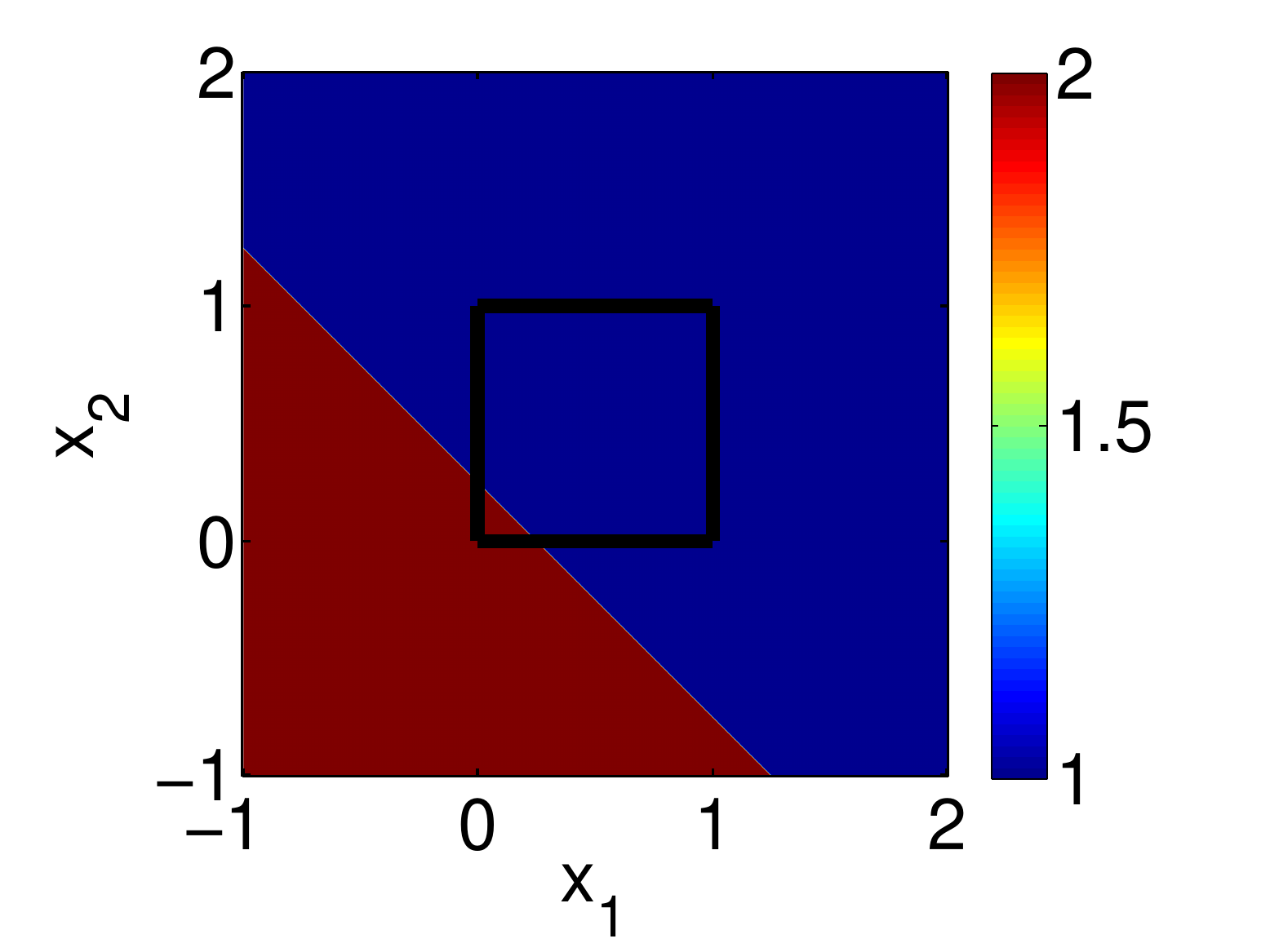}
\caption{Color plot of $c(\x)$ for the diagonal fault.}\label{diagfault}
\includegraphics[scale=.33,trim = 7mm 2mm 15mm 0mm, clip]{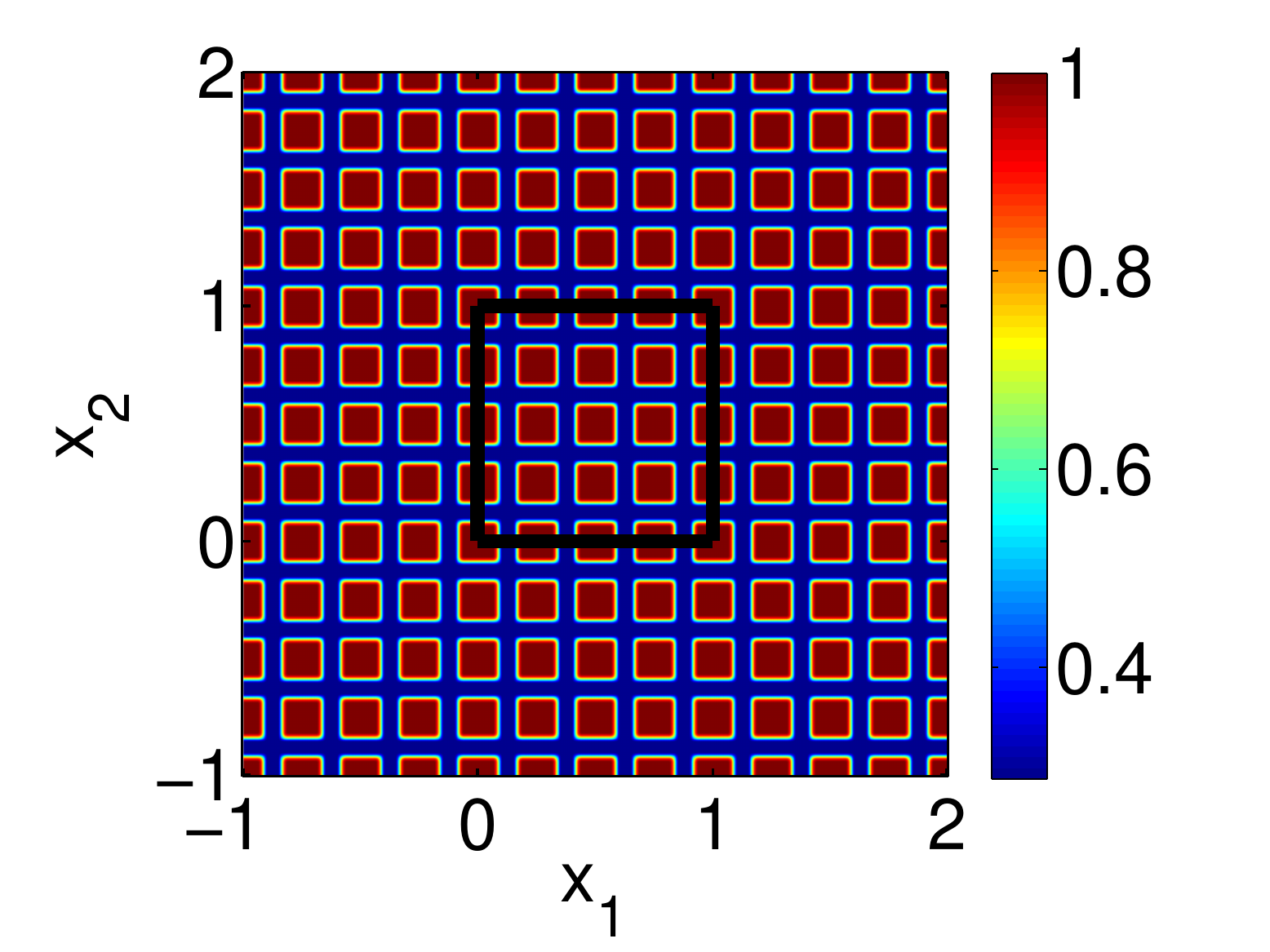}
\caption{Color plot of $c(\x)$ for the periodic medium.}\label{period}
\end{minipage}
\end{figure}

Our benchmark media $c(\x)$ are as follows:

\begin{enumerate}
\item a uniform wave speed of 1, $c \equiv 1$ (Figure \ref{c1}),
\item a ``Gaussian waveguide" (Figure \ref{wg}),
\item a ``Gaussian slow disk" (Figure \ref{slow}) large enough to encompass $\Omega$ -- this will cause some waves going out of $\Omega$ to come back in, 
\item a ``vertical fault" (Figure \ref{fault}), 
\item a ``diagonal fault" (Figure \ref{diagfault}), 
\item and a discontinuous periodic medium (Figure \ref{period}). The periodic medium consists of square holes of velocity 1 in a background of velocity $1/\sqrt{12}$. 
\end{enumerate}

All media used are continued in the obvious way (i.e., they are \emph{not} put to a uniform constant) outside of the domain in which they are shown in the figures if needed. The outline of the $[0,1]^2$ box is shown in black.

We can use a standard Helmholtz equation solver to estimate the relative error in the Helmholtz equation solution caused by the Finite Difference discretization (the \emph{FD error}\footnote{To find this FD error, we use a large pseudo-PML, and compare the solution $u$ for different values of $N$. What we call the FD error is the relative $\ell_2$ error in $u$ inside $\Omega$.}), and also the error caused by using the specified pPML width\footnote{To obtain the error caused by the absorbing layer, we fix $N$ and compare the solution $u$ for different layer widths $w$, and calculate the relative $\ell_2$ error in $u$ inside $\Omega$.}. Those errors are presented in Table \ref{FDPMLerr}, along with the main parameters used in the remaining of this section, including the position of the point source or right-hand side $f$. We note that, whenever possible, we try to use an AL with error smaller than the precision we seek with matrix probing, so with a width $w$ greater than that showed in Table \ref{FDPMLerr}. This makes probing easier, i.e. $p$ and $q$ can be smaller.

\begin{table}
\begin{center} \footnotesize
\begin{tabular}{|l|l|l|l|l|l|l|} \hline  
Medium			&$N$	&$\omega/2\pi$	&FD error	&$w$		&$P$	&Source position	\\ \hline  
$c \equiv 1$		&1023	&51.2		&{$2.5e-01$}	&{$4$}		&8	&$(0.5,0.25)$		\\ \hline  
waveguide		&1023	&51.2		&{$2.0e-01$}	&{$4$}		&56	&$(0.5,0.5)$		\\ \hline  
slow disk		&1023	&51.2		&{$1.8e-01$}	&{$4$}		&43	&$(0.5,0.25)$		\\ \hline
fault, left source	&1023	&51.2		&{$1.1e-01$}	&{$4$}		&48	&$(0.25,0.5)$		\\ \hline  
fault, right source	&1023	&51.2		&{$2.2e-01$}	&{$4$}		&48	&$(0.75,0.5)$		\\ \hline
diagonal fault		&1023	&51.2		&{$2.6e-01$}	&{$256$}	&101	&$(0.5,0.5)$		\\ \hline
periodic medium		&319	&6		&{$1.0e-01$}	&{$1280$}	&792	&$(0.5,0.5)$		\\ \hline
\end{tabular}
\end{center} 
\caption{For each medium considered, we show the parameters $N$ and $\omega/2\pi$, along with the resulting discretization error caused by the Finite Difference (FD error) formulation. We also show the width $w$ of the pPML needed, in number of points, to obtain an error caused by the pPML of less than $1e-1$. Furthermore, we show the total number $P$ of basis matrices needed to probe the entire DtN map with an accuracy of about $1e-1$ as found in Section \protect\ref{sec:tests}. Finally, we show the position of the point source used in calculating the solution $u$.} 
\label{FDPMLerr} 
\end{table}

Consider now a block $M$ of $D$, corresponding to the restriction of $D$ to two sides of $\partial \Omega$. We note that some blocks in $D$ are the same up to transpositions or flips (inverting the order of columns or rows) if the medium $c$ has symmetries. 

\begin{definition}
\emph{Multiplicity of a block of $D$.} Let $M$ be a block of $D$, corresponding to the restriction of $D$ to two sides of $\partial \Omega$. The \emph{multiplicity} $m(M)$ of $M$ is the number of copies of $M$ appearing in $D$, up to transpositions or flips. 
\end{definition}

Only the distinct blocks of $D$ need to be probed. Once we have chosen a block $M$, we may calculate the \emph{true probing coefficients}.

\begin{definition}
\emph{True probing coefficients of block $M$.} Let $M$ be a block of $D$, corresponding to the restriction of $D$ to two sides of $\partial \Omega$. Assume orthonormal probing basis matrices $\left\{B_j \right\}$. The true coefficients $c^t_j$ in the probing expansion of $M$ are the inner products $c^t_j = \< B_j, M\>$.
\end{definition}

We may now define the \emph{$p$-term approximation error} for the block $M$.%
\begin{definition}
\emph{The $p$-term approximation error of block $M$.} Let $M$ be a block of $D$, corresponding to the restriction of $D$ to two sides of $\partial \Omega$. For orthonormal probing basis matrices $\left\{B_j \right\}$, we have the true coefficients $c^t_j$ in the probing expansion of $M$. Let $M_p =\sum_{j=1}^p c^t_j B_j$ be the probing $p$-term approximation to $M$. The $p$-term approximation error for $M$ is
\begin{equation}\label{apperr}
\sqrt{m(M)} \frac{\|M-M_p\|_F}{\|D\|_F},
\end{equation}
using the matrix Frobenius norm.
\end{definition}

Because the blocks on the diagonal of $D$ have a singularity, their Frobenius norm can be a few orders of magnitude greater than that of other blocks, and so it is more important to approximate those well. This is why we consider the error relative to $D$, not to the block $M$, in the $p$-term approximation error. Also, we multiply by the square root of the multiplicity of the block to give us a better idea of how big the total error on $D$ will be. For brevity, we shall refer to (\ref{apperr}) simply as the approximation error when it is clear from the context what $M$, $p$ $\left\{B_j\right\}$, $D$ are.

Then, using matrix probing, we will recover a coefficient vector $\bc$ close to $\bc^t$, which gives an approximation $\tilde{M}=\sum_{j=1}^p c_j B_j$ to $M$. %
We now define the \emph{probing error} (which depends on $q$ and the random vectors used), for the block $M$.%

\begin{definition}
\emph{Probing error of block $M$.} Let $\bc$ be the probing coefficients for $M$ obtained with $q$ random realizations $z^{(1)}$ through $z^{(q)}$. Let $\tilde{M}=\sum_{j=1}^p c_j B_j$ be the probing approximation to $M$. The probing error of $M$ is
\begin{equation}\label{acterr}
\sqrt{m(M)}\frac{\|M-\tilde{M}\|_F}{\|D\|_F}. 
\end{equation}
\end{definition}

Again, for brevity, we refer to (\ref{acterr}) as the probing error when other parameters are clear from the context. Once all distinct blocks of $D$ have been probed, we can consider the \emph{total probing error}.

\begin{definition}
\emph{Total probing error.} The total probing error is defined as the total error made on $D$ by concatenating all probed blocks $\tilde{M}$ to produce an approximate $\tilde{D}$, and is equal to
\begin{equation}\label{eq:Derr}
\frac{\|D-\tilde{D}\|_F}{\|D\|_F}. 
\end{equation}
\end{definition}

In order to get a point of reference for the accuracy benchmarks, for small problems only, the actual matrix $D$ is computed explicitly by solving the exterior problem $4N$ times using the standard basis as Dirichlet boundary conditions, and from this we can calculate \eqref{eq:Derr} exactly. For larger problems, we only have access to a black-box that outputs the product of $D$ with some input vector by solving the exterior problem. We can then estimate \eqref{eq:Derr} by comparing the products of $D$ and $\tilde{D}$ with a few random vectors different from those used in matrix probing.

We shall present results on the approximation and probing errors for various media, along with related condition numbers, and then we shall verify that using an approximate $\tilde{D}$ (constructed from approximate $\tilde{M}$'s for each block $M$ in $D$) does not affect the accuracy of the new solution to the Helmholtz equation, using the \emph{solution error from probing}.
\begin{definition}
\emph{Solution error from probing.} Once we have obtained an approximation $\tilde{D}$ to $D$ from probing the distinct blocks of $D$, we may use this $\tilde{D}$ in a Helmholtz solver to obtain an approximate solution $\tilde{u}$, and compare that to the true solution $u$ using $D$ in the solver. The solution error from probing is the $\ell_2$ error on $u$ inside $\Omega$:
\begin{equation}\label{eq:solerr}
\frac{\|u-\tilde{u}\|_2}{\|u\|_2} \text{ in } \Omega. 
\end{equation}
\end{definition}

\subsection{Probing tests}\label{sec:tests}

As we saw in Section \ref{sec:probe}, randomness plays a role in the value of $\mbox{cond}({\bf \Psi})$ and of the probing error. Hence, whenever we show plots for those quantities in this section, we have done 10 trials for each value of $q$ used. The error bars show the minimum and maximum of the quantity over the 10 trials, and the line is plotted through the average value over the 10 trials. As expected, we will see in all experiments that increasing $q$ gives a better conditioning, and consequently a better accuracy and smaller failure probability. The following probing results will then be used in Section \ref{sec:insolver} to solve the Helmholtz equation.

\subsubsection{Uniform medium}

\begin{figure}[ht]
\begin{minipage}[t]{0.48\linewidth}
\includegraphics[scale=.5]{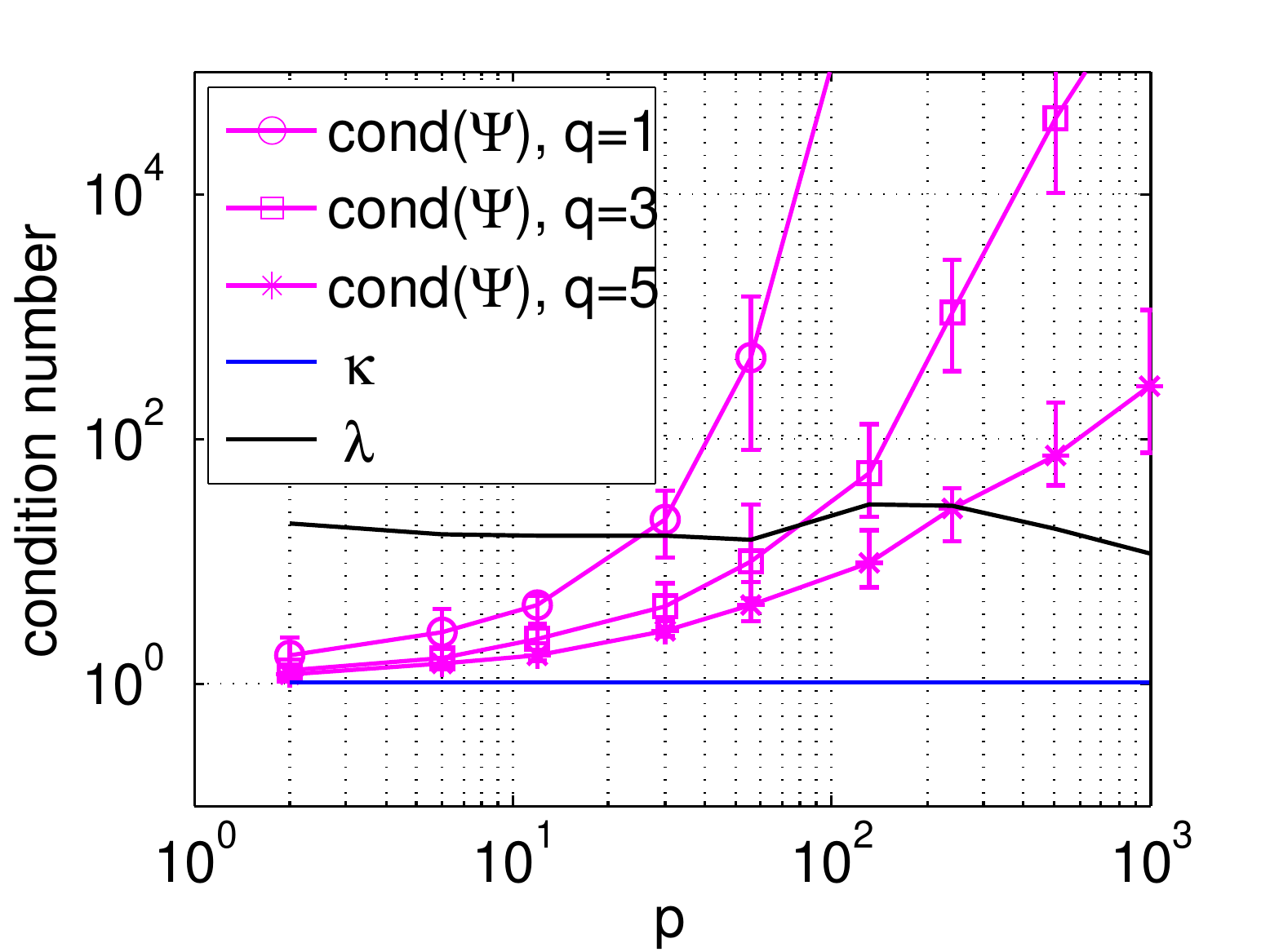}
\caption{Condition numbers for the $(1,1)$ block, $c\equiv 1$.}
\label{cond11_1024_c1}
\end{minipage}
\hspace{0.1cm}
\begin{minipage}[t]{0.48\linewidth}
\includegraphics[scale=.5]{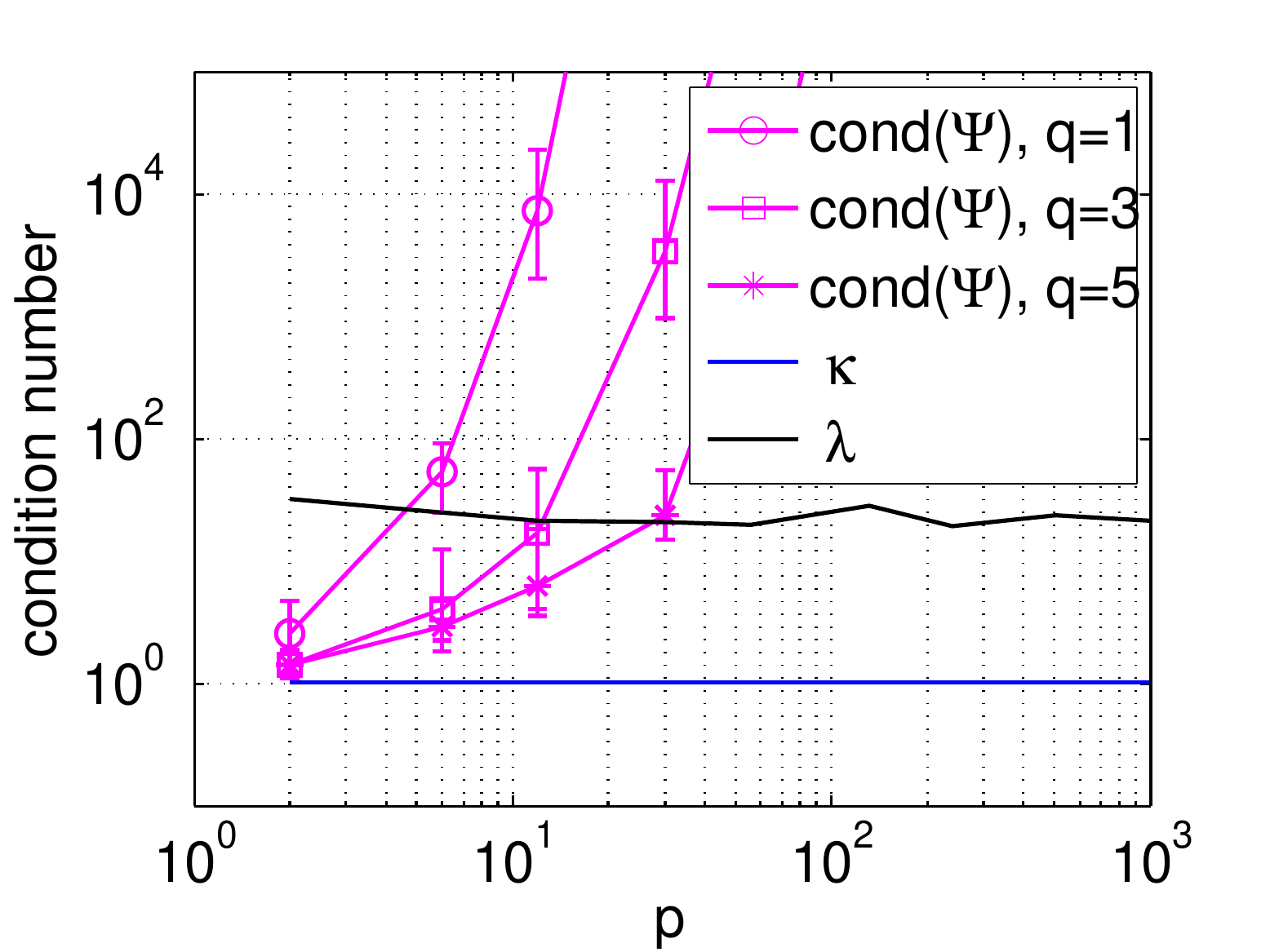}
\caption{Condition numbers for the $(2,1)$ block, $c\equiv 1$.}
\label{cond21_1024_c1}
\end{minipage}
\end{figure}

For a uniform medium, $c \equiv 1$, we have three blocks with the following multiplicities: $m((1,1))=4$ (same edge), $m((2,1))=8$ (neighboring edges), and $m((3,1))=4$ (opposite edges). Note that we do not present results for the $(3,1)$ block: this block has negligible Frobenius norm\footnote{We can use probing with $q=1$ and a single basis matrix (a constant multiplied by the correct oscillations) and have a probing error of less than $10^{-6}$ for that block.} compared to $D$.  First, let us look at the conditioning for blocks $(1,1)$ and $(2,1)$. Figures \ref{cond11_1024_c1} and \ref{cond21_1024_c1} show the three relevant conditioning quantities: $\kappa$, $\lambda$ and $\mbox{cond}({\bf \Psi})$ for each block. As expected, $\kappa=1$ because we orthogonalize the basis functions. Also, we see that $\lambda$ does not grow very much as $p$ increases, it remains on the order of 10. As for $\mbox{cond}({\bf \Psi})$, it increases as $p$ increases for a fixed $q$ and $N$, as expected. This will affect probing in terms of the failure probability (the odds that the matrix ${\bf \Psi}$ is far from the expected value) and accuracy (taking the pseudo-inverse will introduce larger errors in $\bc$). We notice these two phenomena in Figure \ref{erb1023_c1}, where we show the approximation and probing errors in probing the $(1,1)$ block for various $p$, using different $q$ and making 10 tests for each $q$ value as explained previously. As expected, as $p$ increases, the variations between trials get larger. Also, the probing error, always larger than the approximation error, becomes farther and farther away from the approximation error. Comparing Figure \ref{erb1023_c1} with Table \ref{c1solve} of the next section, we see that in Table \ref{c1solve} we are able to achieve higher accuracies. This is because we use the first two traveltimes (so four different types of oscillations, as explained in Section \ref{sec:basis}) to obtain those higher accuracies. But we do not use four types of oscillations for lower accuracies because this demands a larger number of basis matrices $p$ and of solves $q$ for the same error level.

\subsubsection{The waveguide}

\begin{figure}[ht]
\begin{minipage}[t]{0.48\linewidth}
\includegraphics[scale=.5]{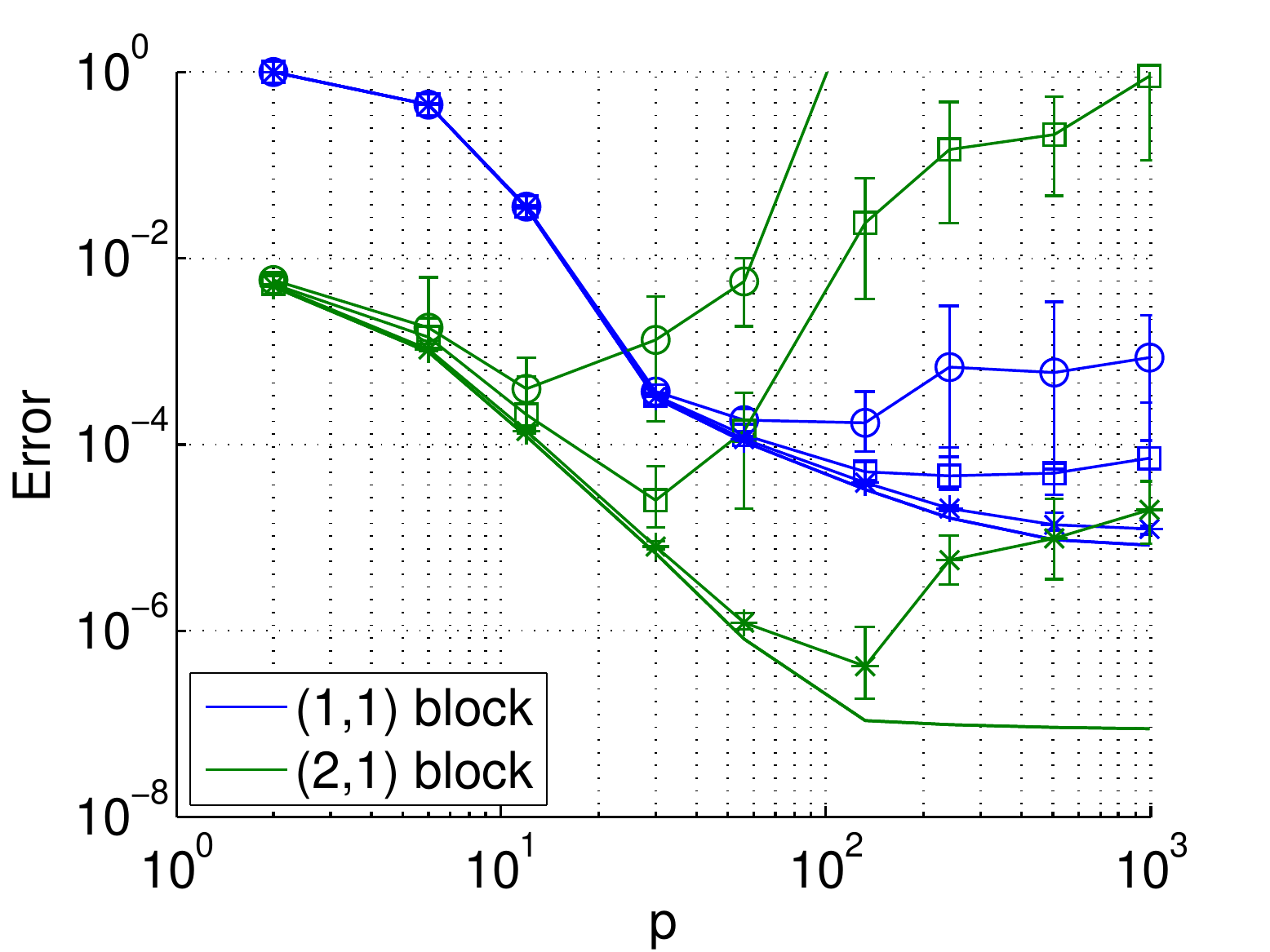}
\caption{Approximation error (line) and probing error (with markers) for the blocks of $D$, $c\equiv 1$. Circles are for $q=3$, squares for $q=5$, stars for $q=10$. }
\label{erb1023_c1}
\end{minipage}
\hspace{0.1cm}
\begin{minipage}[t]{0.48\linewidth}
\includegraphics[scale=.5]{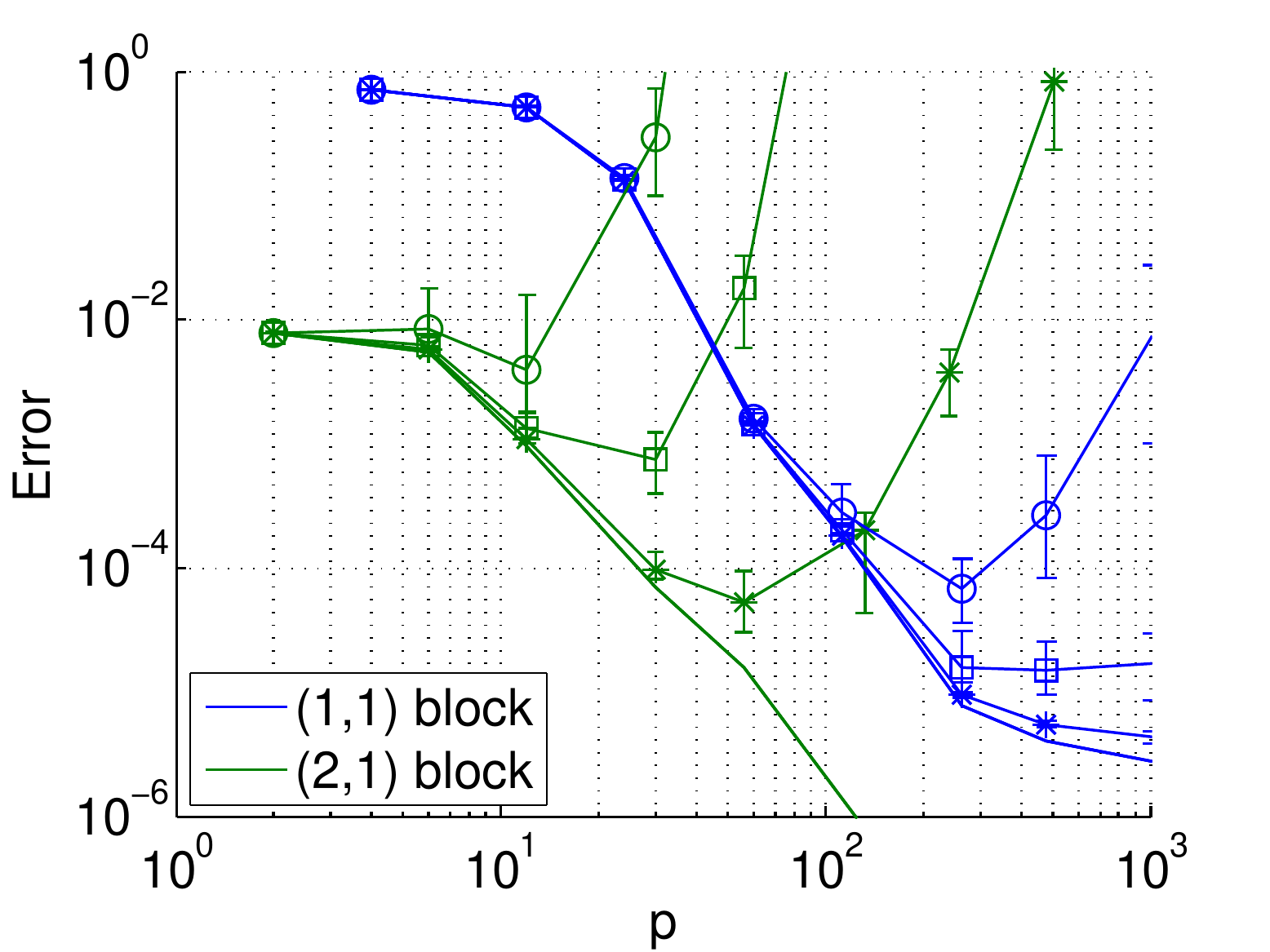}
\caption{Approximation error (line) and probing error (with markers) for the blocks of $D$, $c$ is the waveguide. Circles are for $q=3$, squares for $q=5$, stars for $q=10$.}
\label{erb1023_c2}
\end{minipage}
\end{figure}

For a waveguide as a velocity field, we have more blocks compared to the uniform medium case, with different multiplicities: $m((1,1))=2$, $m((2,2))=2$, $m((2,1))=8$, $m((3,1))=2$, $m((4,2))=2$. Note that block $(2,2)$ will be easier to probe than block $(1,1)$ since the medium is smoother on that interface. Also, we can probe blocks $(3,1)$ and $(4,2)$ with $q=1$, $p=2$ and have a probing error less than $10^{-7}$. Hence we only show results for the probing and approximation errors of blocks $(1,1)$ and $(2,1)$, in Figure \ref{erb1023_c2}. Results for using probing in a solver can be found in Section \ref{sec:insolver}.

\subsubsection{The slow disk}

\begin{figure}[ht]
\begin{minipage}[t]{0.48\linewidth}
\includegraphics[scale=.5]{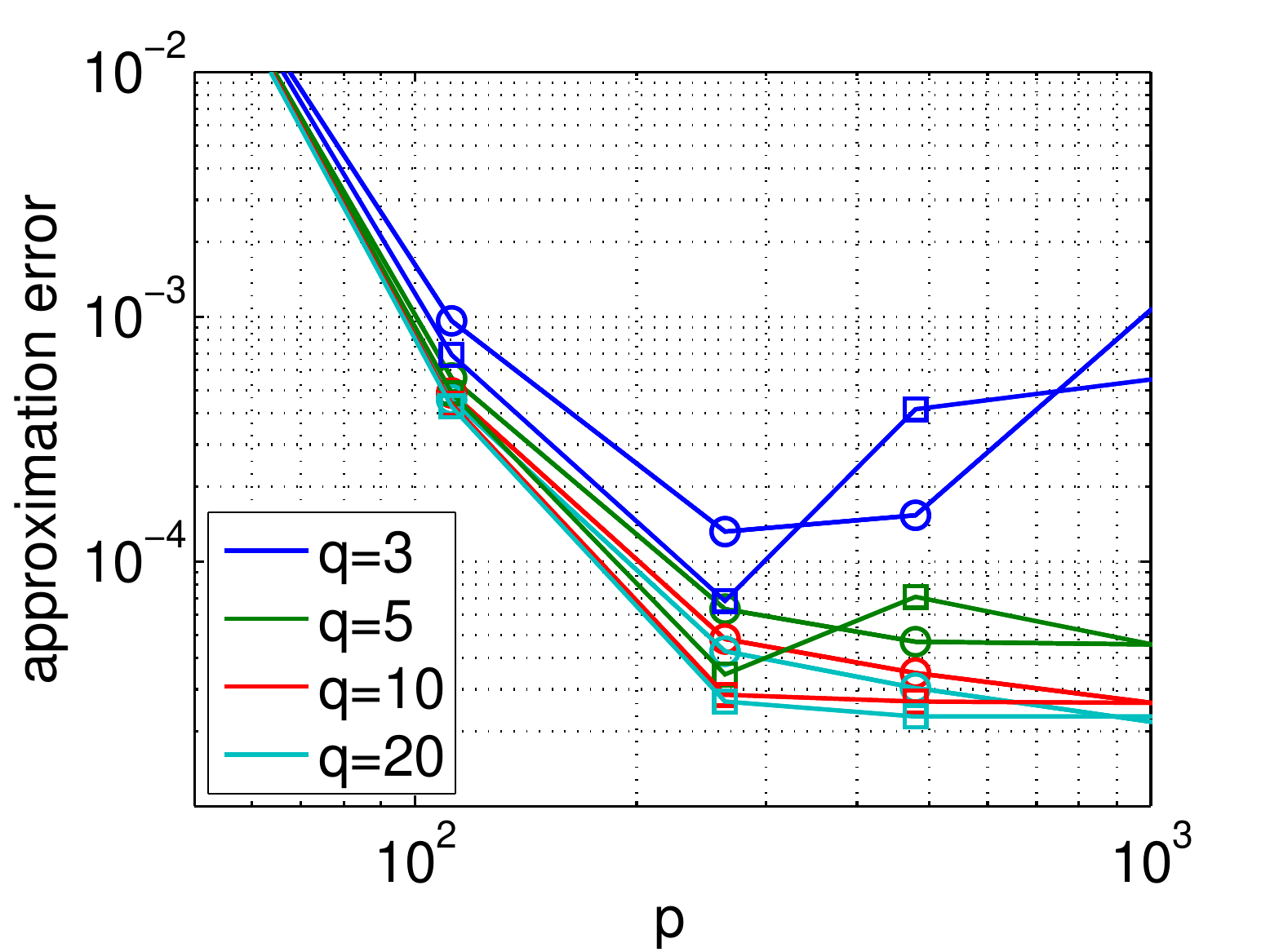}
\caption{Approximation error for the $(1,1)$ blocks of $D$, $c$ is the slowness disk, comparing the use of the normal traveltime (circles) to the fast marching traveltime (squares).}
\label{c5fastvsnorm1}
\end{minipage}
\hspace{0.1cm}
\begin{minipage}[t]{0.48\linewidth}
\includegraphics[scale=.5]{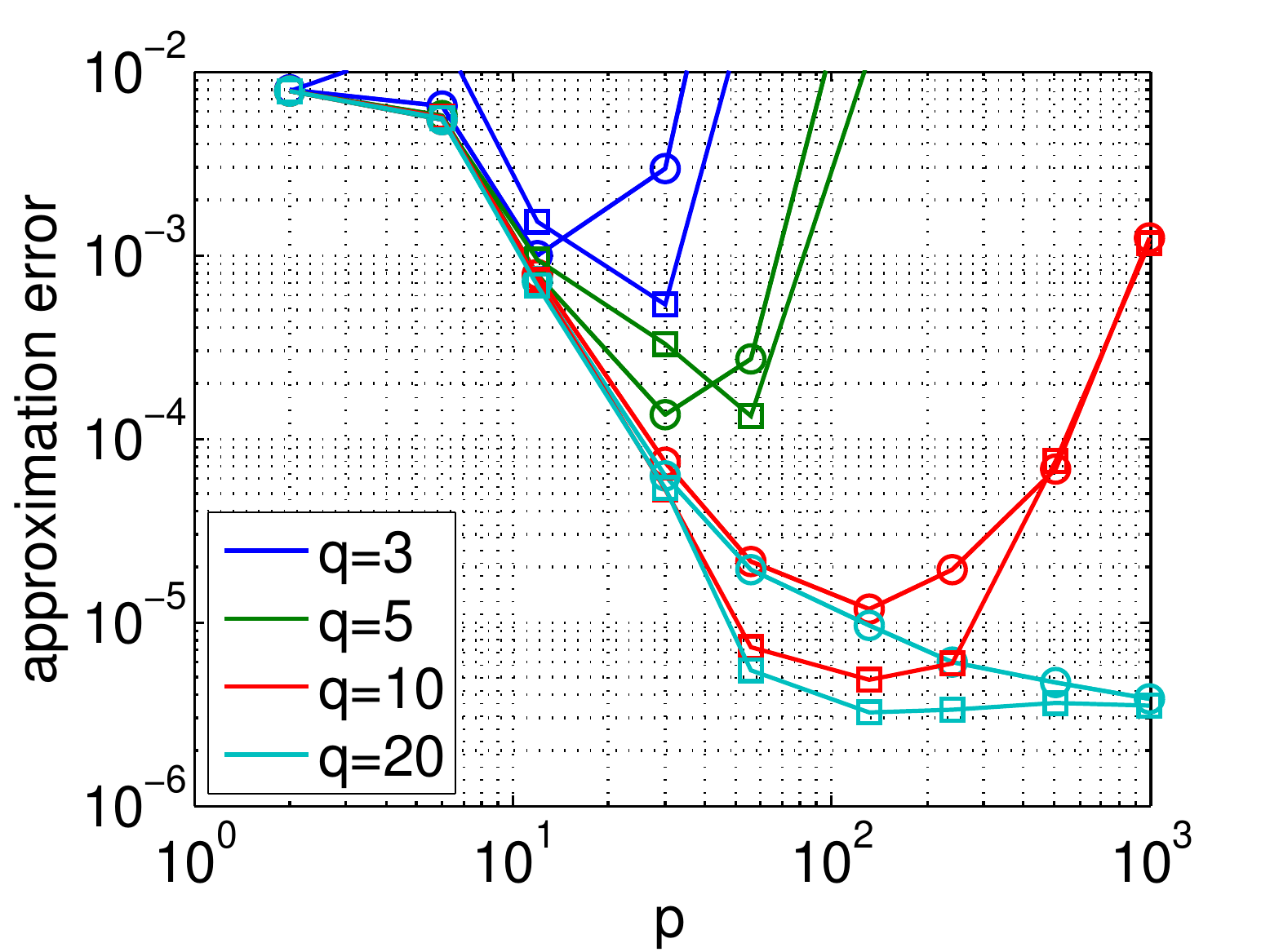}
\caption{Approximation error for the $(2,1)$ blocks of $D$, $c$ is the slowness disk, comparing the use of the normal traveltime (circles) to the fast marching traveltime (squares).}
\label{c5fastvsnorm2}
\end{minipage}
\end{figure}

Next, we consider the slow disk. Here, we have a choice to make for the traveltime upon which the oscillations depend. We may consider interfacial waves, traveling in straight line segments along $\partial \Omega$, with traveltime $\tau$. There is also the first arrival time of body waves, $\tau_f$, which for some points on $\partial \Omega$ involve taking a path that goes away from $\partial \Omega$, into the exterior where $c$ is higher, and back towards $\partial \Omega$. We have approximated this $\tau_f$ using the fast marching method of Sethian \cite{sethart}. For this example, it turns out that using either $\tau$ or $\tau_f$ to obtain oscillations in our basis matrices does not significantly alter the probing accuracy or conditioning, although it does seem that, for higher accuracies at least, the fast marching traveltime makes convergence slightly faster. Figures \ref{c5fastvsnorm1} and \ref{c5fastvsnorm2} demonstrate this for blocks $(1,1)$ and $(2,1)$ respectively. We omit plots of the probing and approximation errors, and refer the reader to Section \ref{sec:insolver} for final probing results and using those in a solver.

\subsubsection{The vertical fault}

\begin{figure}[ht]
\begin{minipage}[t]{0.48\linewidth}
\includegraphics[scale=.5]{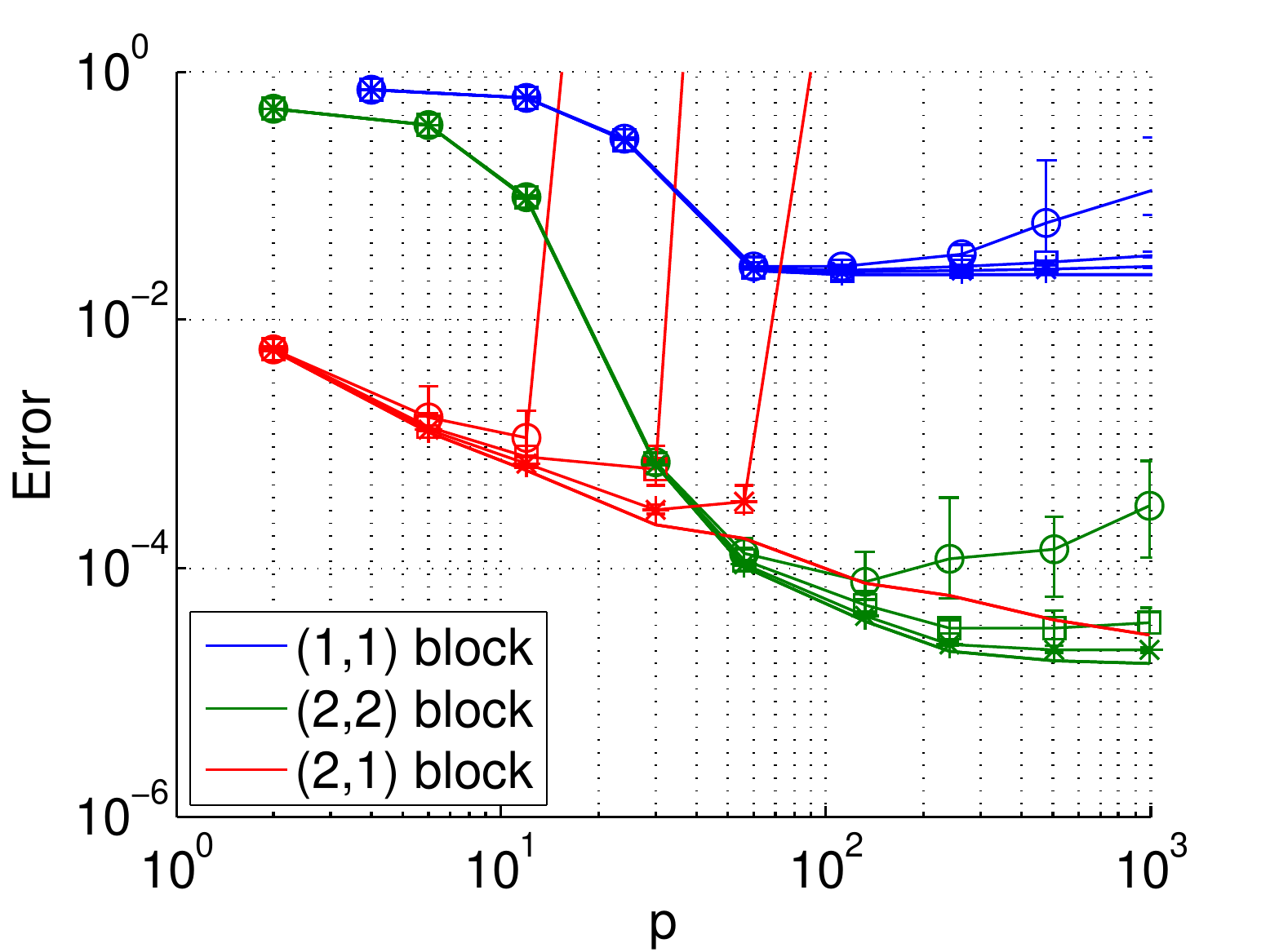}
\caption{Approximation error (line) and probing error (with markers) for the blocks of $D$, $c$ is the fault. Circles are for $q=3$, squares for $q=5$, stars for $q=10$.}
\label{erb1023_c16}
\end{minipage}
\hspace{0.1cm}
\begin{minipage}[t]{0.48\linewidth}
\includegraphics[scale=.5]{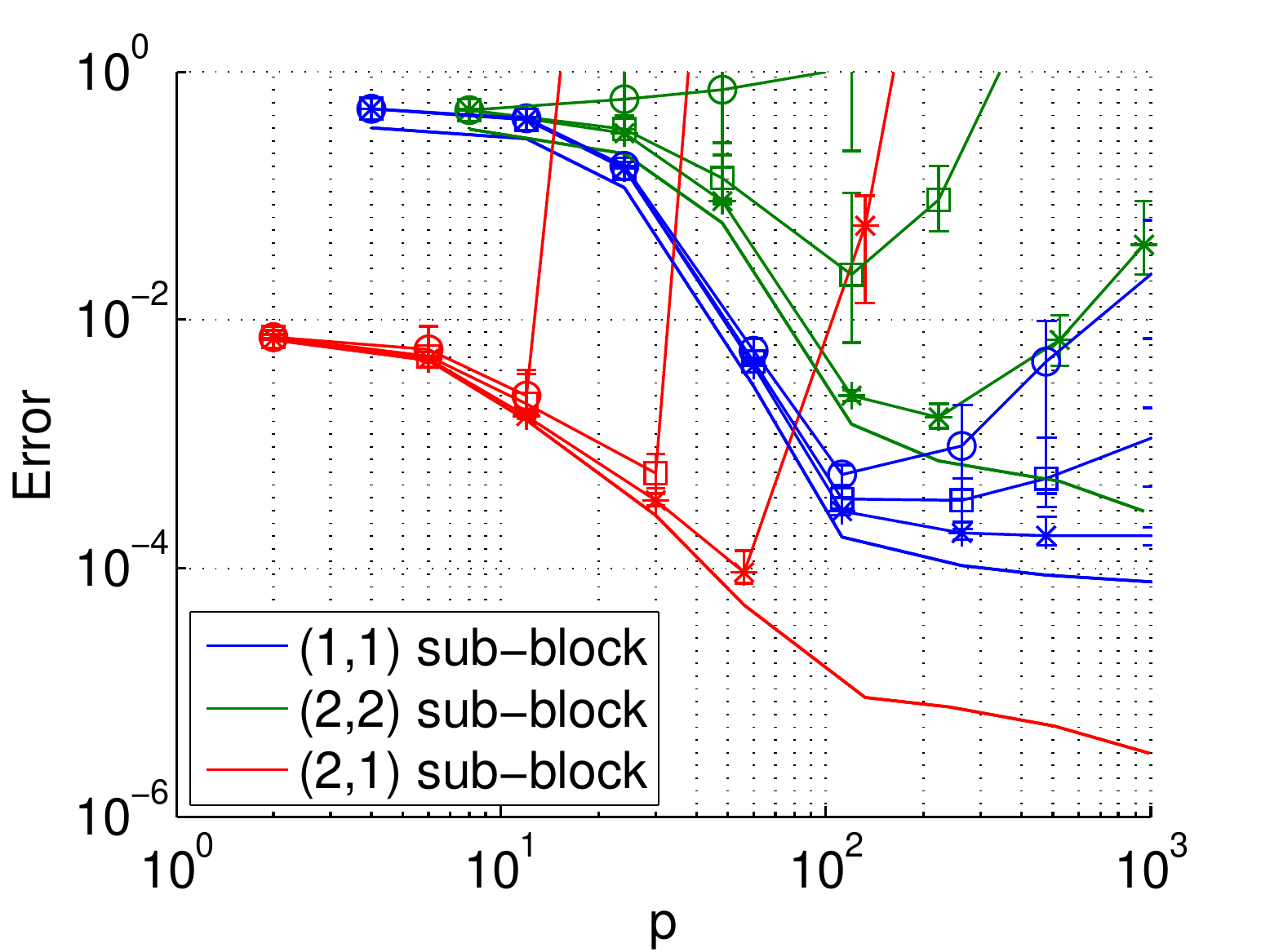}
\caption{Approximation error (line) and probing error (with markers) for the sub-blocks of the $(1,1)$ block of $D$, $c$ is the fault. Circles are for $q=3$, squares for $q=5$, stars for $q=10$.}
\label{erb1023_c16_sub}
\end{minipage}
\end{figure}

Next, we look at the case of the medium $c$ which has a vertical fault. We note that this case is harder because some of the blocks will have themselves a 2 by 2 or 1 by 2 structure caused by the discontinuity in the medium. Ideally, as we shall see, each sub-block should be probed separately. There are 7 distinct blocks, with different multiplicities: $m((1,1))=2$, $m((2,2))=1$, $m((4,4))=1$, $m((2,1))=4$, $m((4,1))=4$, $m((3,1))=2$, $m((4,2))=2$. Blocks $(2,2)$ and $(4,4)$ are easier to probe than block $(1,1)$ because they do not exhibit a sub-structure. Also, since the velocity is smaller on the right side of the fault, the frequency there is higher, which means that blocks involving side 2 are slightly harder to probe than those involving side 4. Hence we first present results for the blocks $(1,1)$, $(2,2)$ and $(2,1)$ of $D$. In Figure \ref{erb1023_c16} we see the approximation and probing errors for those blocks. Then, in Figure \ref{erb1023_c16_sub}, we present results for the errors related to probing the 3 distinct sub-blocks of the $(1,1)$ block of $D$. We can see that probing the $(1,1)$ block by sub-blocks helps achieve greater accuracy. We could have split other blocks too to improve the accuracy of their probing (for example, block $(2,1)$ has a 1 by 2 structure because side 1 has a discontinuity in $c$) but the accuracy of the overall DtN map was still limited by the accuracy of probing the $(1,1)$ block, so we do not show results for other splittings.

\subsubsection{The diagonal fault}

Now, we look at the case of the medium $c$ which has a diagonal fault. Again, some of the blocks will have themselves a 2 by 2 or 1 by 2 structure. There are 6 distinct blocks, with different multiplicities: $m((1,1))=2$, $m((2,2))=2$, $m((2,1))=4$, $m((4,1))=2$, $m((3,2))=2$, $m((3,1))=4$. Again, we split up block $(1,1)$ in 4 sub-blocks and probe each of those sub-blocks separately for greater accuracy, but do not split other blocks. We then use two traveltimes for the $(2,2)$ sub-block of block $(1,1)$. Using as the second arrival time the geometrical traveltime consisting of leaving the boundary and bouncing off the fault, as mentioned in Section \ref{sec:tt}, allowed us to increase accuracy by an order of magnitude compared to using only the first arrival traveltime, or compared to using as a second arrival time the usual bounce off the corner (or here, bounce off the fault where it meets $\delta \Omega$). We omit plots of the probing and approximation errors, and refer the reader to Section \ref{sec:insolver} for final probing results and using those in a solver.

\subsubsection{The periodic medium}

Finally, we look at the case of the periodic medium presented earlier. There are 3 distinct blocks, with different multiplicities: $m((1,1))=4$, $m((2,1))=8$, $m((3,1))=4$. We expect the corresponding DtN map to be harder to probe because its structure will reflect that of the medium, i.e. it will exhibit sharp transitions at points corresponding to sharp transitions in $c$ (similarly as with the faults). First, we notice that, in all the previous mediums we tried, plotting the norm of the anti-diagonal entries of diagonal blocks (or sub-blocks for the faults) shows a rather smooth decay away from the diagonal. However, that is not the case for the periodic medium: it looks like there is decay away from the diagonal, but variations from that decay can be of relative order 1. This prevents our usual strategy, using basis matrices containing terms that decay away from the diagonal such as $(h+|x-y|)^{-j_1/\alpha}$, from working adequately. Instead, we use polynomials along anti-diagonals, as well as polynomials along diagonals as we previously did.

It is known that solutions to the Helmholtz equation in a periodic medium are Bloch waves with a particular structure \cite{JohnsonPhot}. However, using that structure in the basis matrices is not robust. Indeed, using a Bloch wave structure did not succeed very well, probably because our discretization was not accurate enough and so $D$ exhibited that structure only to a very rough degree. Hence we did not use Bloch waves for probing the periodic medium. Others have successfully used the known structure of the solution in this setting to approximate the DtN map. In \cite{Fliss}, the authors solve local cell problems and Ricatti equations to obtain discrete DtN operators for media which are a perturbation of a periodic structure. In \cite{antoine}, the authors develop a DtN map eigenvalue formulation for wave propagation in periodic media. We did not attempt to use those formulations here.

For this reason, we tried basis matrices with no oscillations, but with polynomials in both directions as explained previously, and obtained the results of Section \ref{sec:insolver}.

Now that we have probed the DtN map and obtained compressed blocks to form an approximation $\tilde{D}$ of $D$, we may use this $\tilde{D}$ in a Helmholtz solver as an absorbing boundary condition.

\subsection{Using the probed DtN map in a Helmholtz solver}\label{sec:insolver}

\begin{figure}[ht]
\begin{minipage}[t]{0.45\linewidth}
\includegraphics[scale=.45]{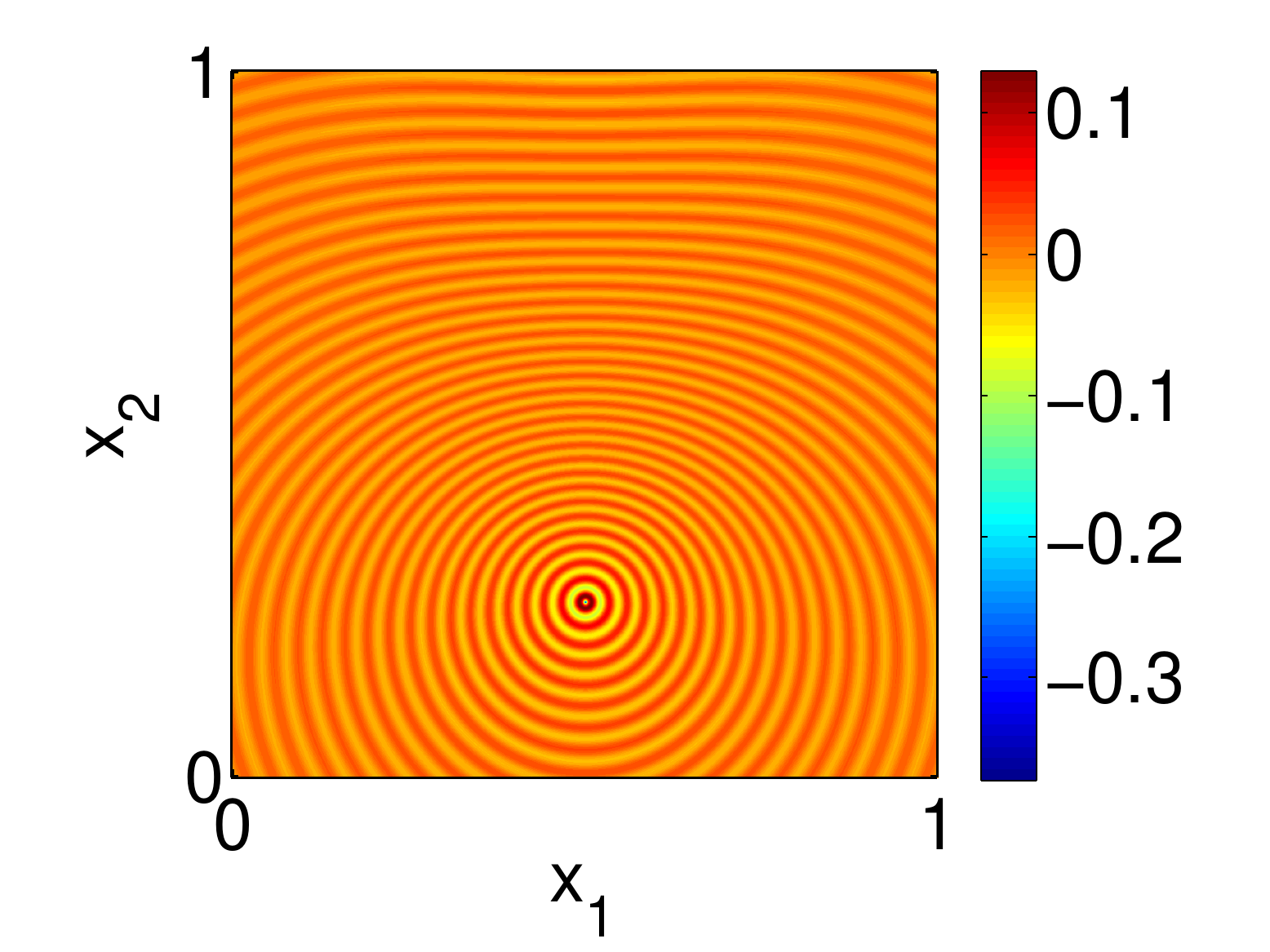}
\caption{Real part of the solution, $c$ is the slow disk.}
\label{solc5}
\end{minipage}
\hspace{1cm}
\begin{minipage}[t]{0.45\linewidth}
\includegraphics[scale=.45]{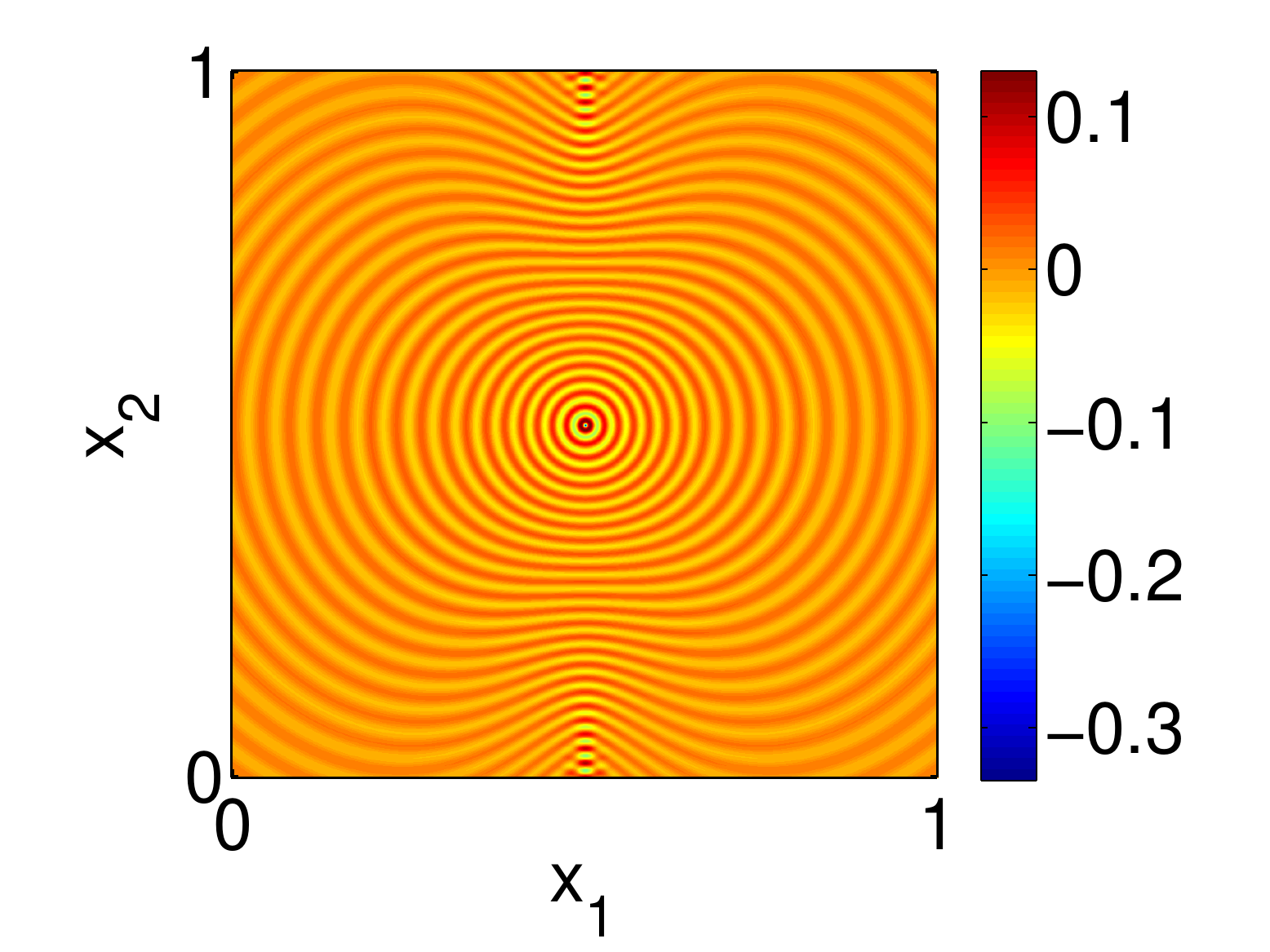}
\caption{Real part of the solution, $c$ is the waveguide.}
\label{solc3}
\end{minipage}
\end{figure}

\begin{figure}[ht]
\begin{minipage}[t]{0.45\linewidth}
\includegraphics[scale=.45]{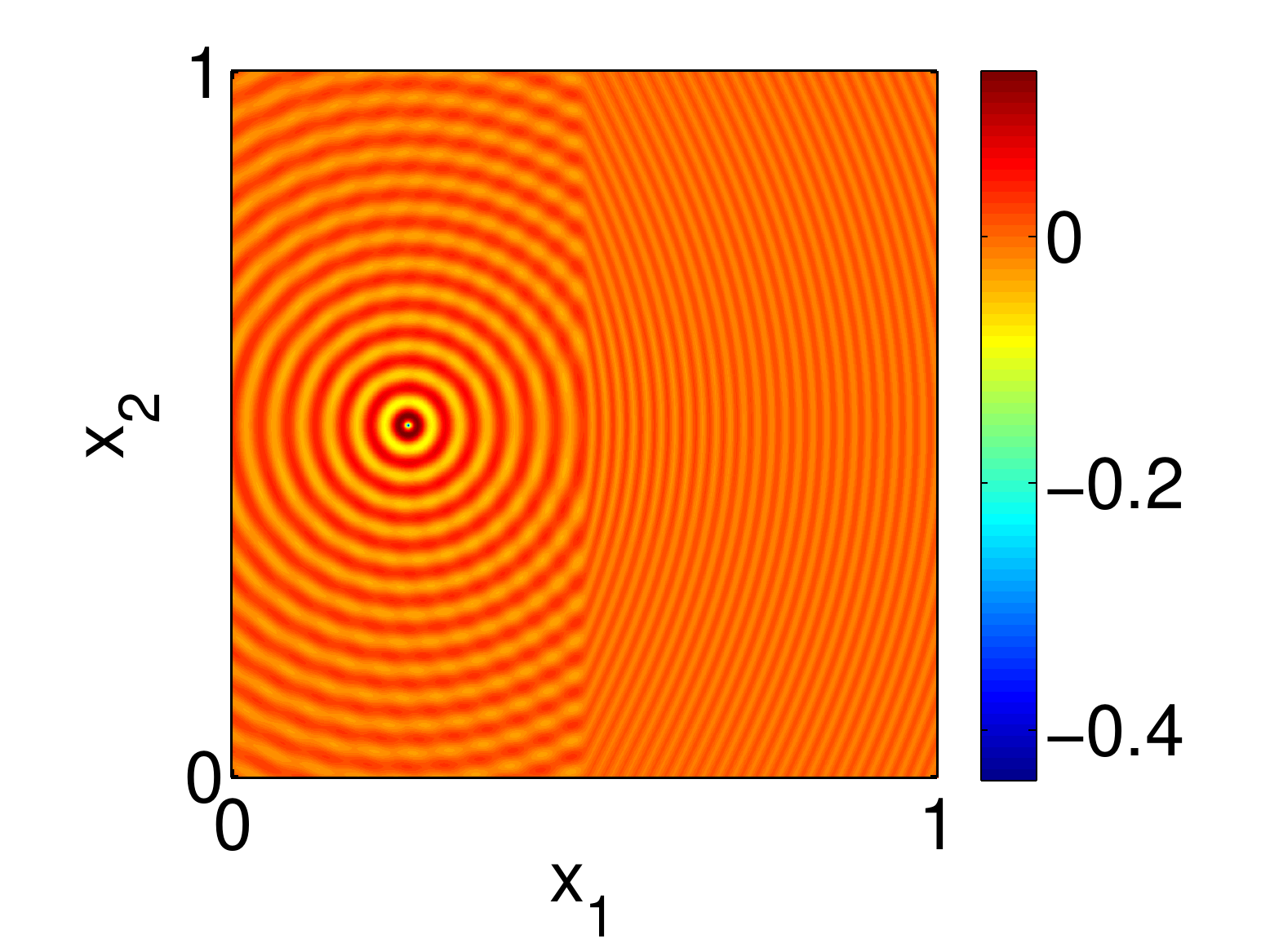}
\caption{Real part of the solution, $c$ is the vertical fault with source on the left.}
\label{solc16l}
\end{minipage}
\hspace{1cm}
\begin{minipage}[t]{0.45\linewidth}
\includegraphics[scale=.45]{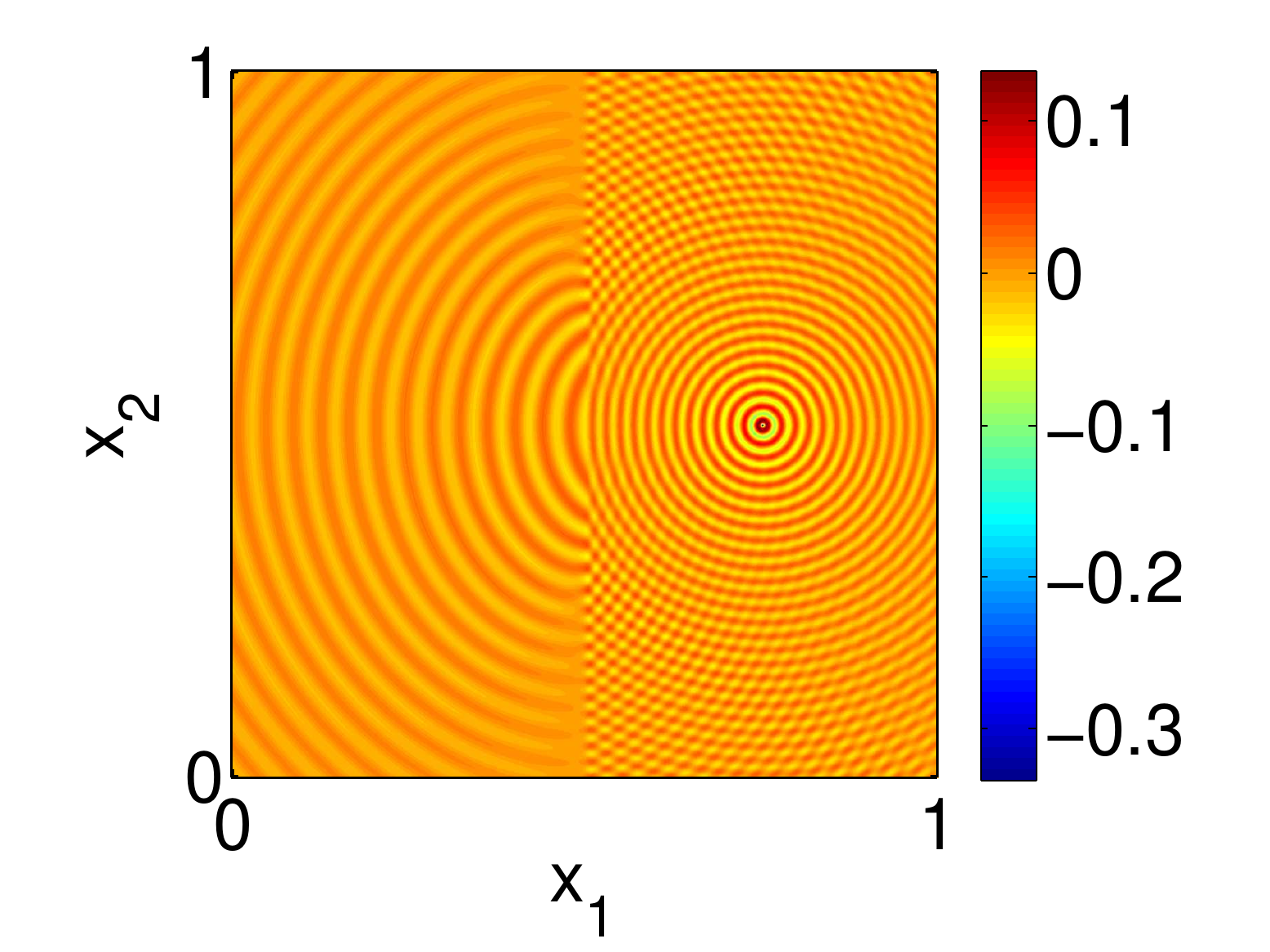}
\caption{Real part of the solution, $c$ is the vertical fault with source on the right.}
\label{solc16r}
\end{minipage}
\end{figure}

\begin{figure}[ht]
\begin{minipage}[t]{0.45\linewidth}
\includegraphics[scale=.45]{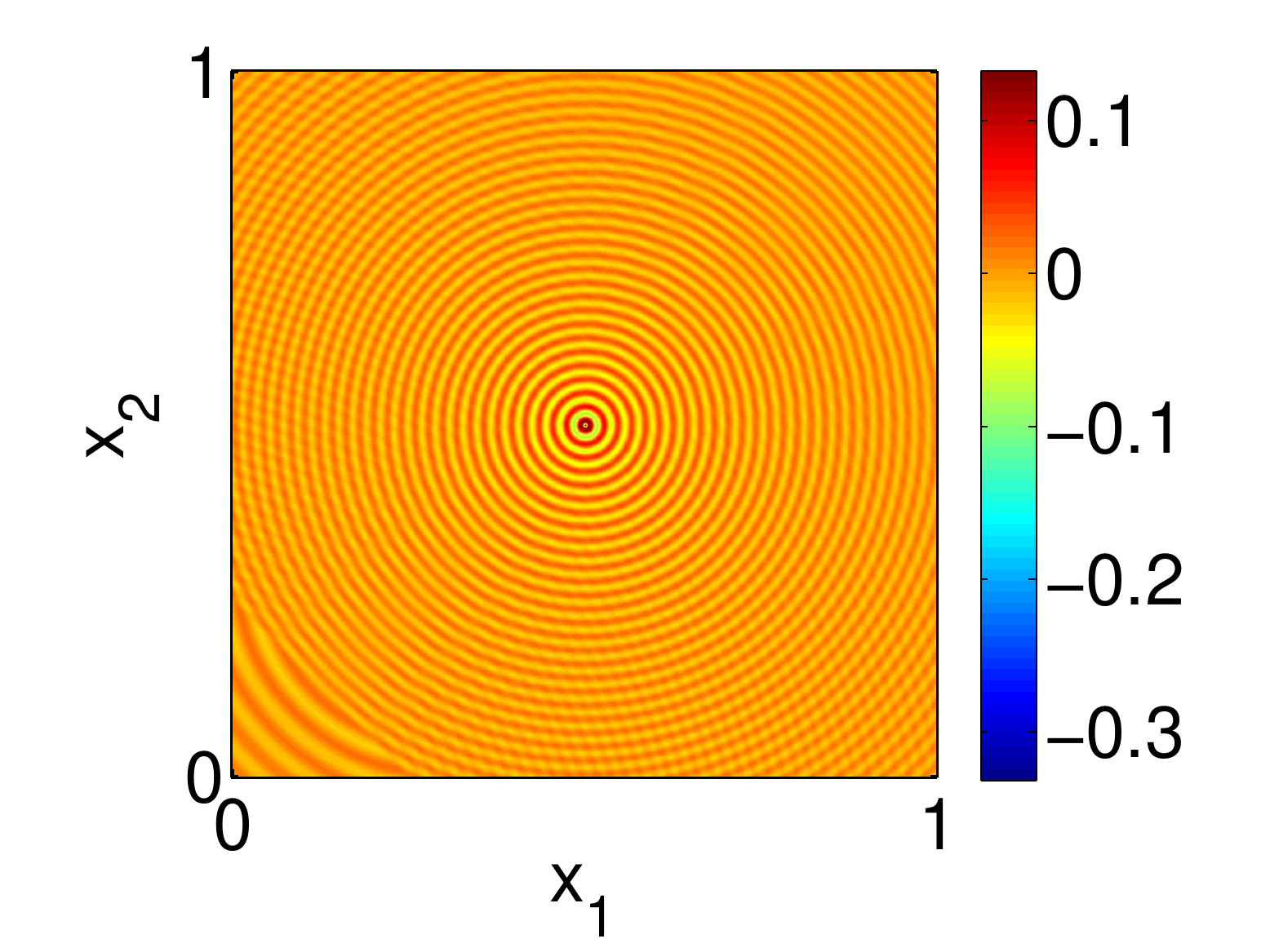}
\caption{Real part of the solution, $c$ is the diagonal fault.}
\label{solc18}
\end{minipage}
\hspace{1cm}
\begin{minipage}[t]{0.45\linewidth}
\includegraphics[scale=.45]{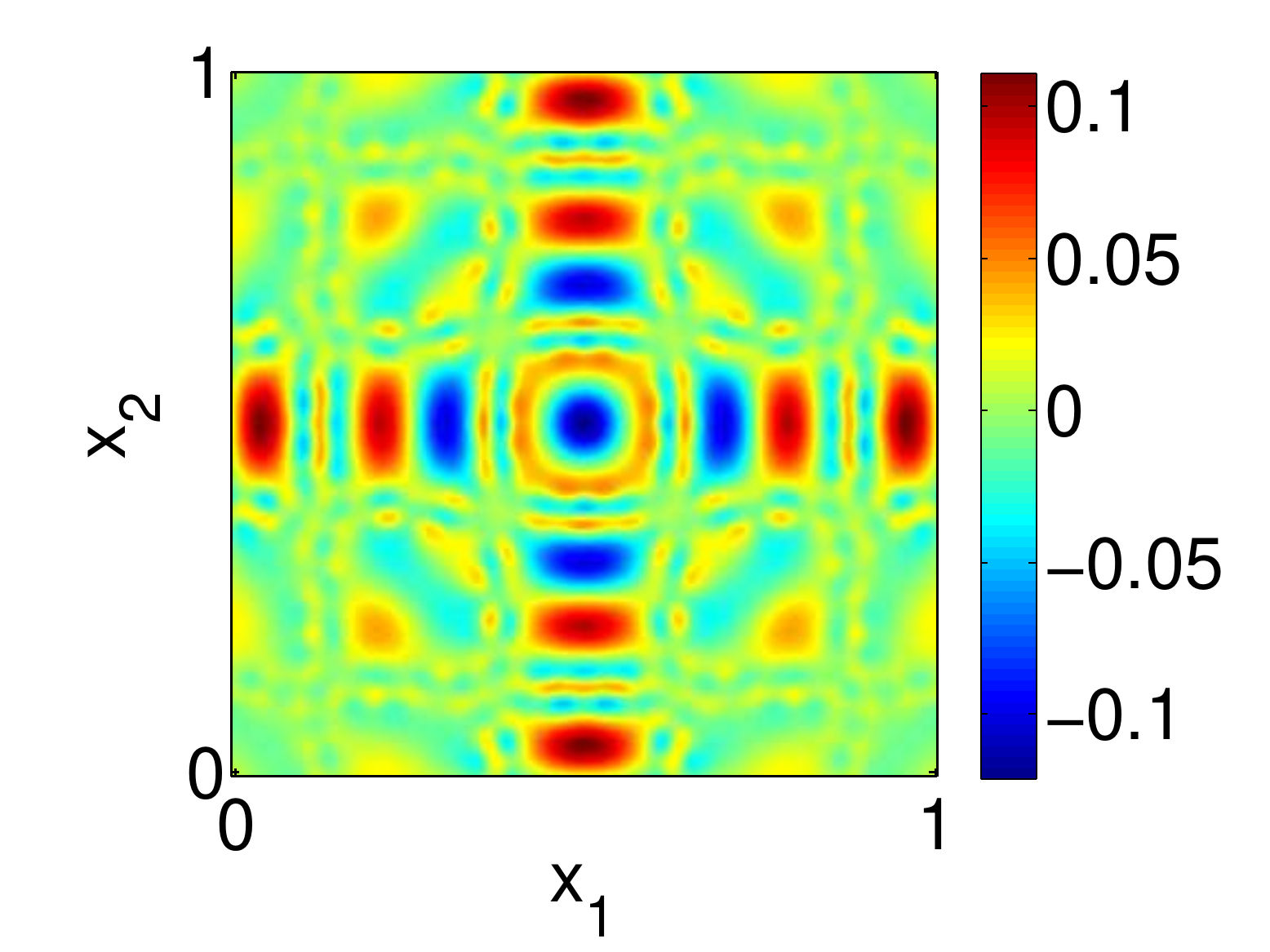}
\caption{Imaginary part of the solution, $c$ is the periodic medium.}
\label{solc33}
\end{minipage}
\end{figure}

In Figures \ref{solc5}, \ref{solc3}, \ref{solc16l}, \ref{solc16r}, \ref{solc18} and \ref{solc33} we can see the standard solutions to the Helmholtz equation on $[0,1]^2$ using a large PML or pPML for the various media we consider, except for the uniform medium, where the solution is well-known. We use those as our reference solutions.

We have tested the solver with the probed $\tilde{D}$ as an absorbing boundary condition with success. See Tables \ref{c1solve}, \ref{c5solve}, \ref{c3solve} and \ref{c16solve} for results corresponding to each medium. We show the number $p$ of basis matrices required for some blocks for that tolerance, the number of solves $q$ of the exterior problem for those blocks, the total number of solves $Q$, the total probing error \eqref{eq:Derr} in $D$ and the solution error from probing \eqref{eq:solerr}. %
As we can see from the tables, the solution error from probing \eqref{eq:solerr} in the solution $u$ is no more than an order of magnitude greater than the total probing error \eqref{eq:Derr} in the DtN map $D$, for a source position as described in Table \ref{FDPMLerr}. Grazing waves, which can arise when the source is close to the boundary of the computational domain, will be discussed in the next subsection, \ref{sec:graz}. We note again that, for the uniform medium, using the second arrival traveltime as well as the first for the $(1,1)$ block allowed us to achieve accuracies of 5 and 6 digits in the DtN map, which was not possible otherwise. Using a second arrival time for the cases of the faults was also useful. Those results show that probing works best when the medium $c$ is rather smooth. For non-smooth media such as a fault, it becomes harder to probe the DtN map to a good accuracy, so that the solution to the Helmholtz equation also contains more error.

\begin{table}
\caption{$c\equiv 1$} 
\begin{center} \footnotesize
\begin{tabular}{|l|l|l|l|l|l|} \hline  
$p$ for $(1,1)$		& $p$ for $(2,1)$	& $q=Q$		& $\frac{\|D-\tilde{D}\|_F}{\|D\|_F}$ 	& $\frac{\|u-\tilde{u}\|_2}{\|u\|_2}$ \\ \hline  
{$6$}			& {$1$}			& {$1$}		& {$2.0130e-01$}			& {$3.3191e-01$} \\ \hline  
{$12$}			& {$2$}			& {$1$}		& {$9.9407e-03$}			& {$1.9767e-02$} \\ \hline  
{$20$}			& {$12$}		& {$3$}		& {$6.6869e-04$}			& {$1.5236e-03$} \\ \hline  
{$72$}			& {$20$}		& {$5$}		& {$1.0460e-04$}			& {$5.3040e-04$} \\ \hline  
{$224$}			& {$30$}		& {$10$}	& {$8.2892e-06$}			& {$9.6205e-06$} \\ \hline  
{$360$}			& {$90$}		& {$10$}	& {$7.1586e-07$}			& {$1.3044e-06$} \\ \hline  
\end{tabular}
\end{center} 
\label{c1solve} 
\end{table}

\begin{table}
\caption{$c$ is the waveguide} 
\begin{center} \footnotesize
\begin{tabular}{|l|l|l|l|l|l|l|l|l|} \hline  
$p$ for $(1,1)$	&$p$ for $(2,1)$&$q$	&$p$ for $(2,2)$&$q$  &$Q$	&$\frac{\|D-\tilde{D}\|_F}{\|D\|_F}$	&$\frac{\|u-\tilde{u}\|_2}{\|u\|_2}$ \\ \hline  
$40$		&$2$		&$1$	&$12$		&$1$  &$2$	&$9.1087e-02$				&$1.2215e-01$ \\ \hline  
$40$		&$2$		&$3$	&$20$		&$1$  &$4$	&$1.8685e-02$				&$7.6840e-02$       \\ \hline  
$60$		&$20$		&$5$	&$20$		&$3$  &$8$	&$2.0404e-03$				&$1.3322e-02$       \\ \hline  
$112$		&$30$		&$10$	&$30$		&$3$  &$13$	&$2.3622e-04$				&$1.3980e-03$      \\ \hline  
$264$		&$72$		&$20$	&$168$		&$10$ &$30$	&$1.6156e-05$				&$8.9911e-05$      \\ \hline  
$1012$		&$240$		&$20$	&$360$		&$10$ &$30$	&$3.3473e-06$				&$1.7897e-05$      \\ \hline  
\end{tabular}
\end{center} 
\label{c3solve} 
\end{table}

\begin{table}
\caption{$c$ is the slow disk} 
\begin{center} \footnotesize
\begin{tabular}{|l|l|l|l|l|l|} \hline  
$p$ for $(1,1)$	& $p$ for $(2,1)$	&$q=Q$	&$\frac{\|D-\tilde{D}\|_F}{\|D\|_F}$ 	& $\frac{\|u-\tilde{u}\|_2}{\|u\|_2}$ \\ \hline  
{$40$}		& {$2$}			&{$3$}	& {$1.0730e-01$}			& {$5.9283e-01$} \\ \hline  
{$84$}		& {$2$}			&{$3$}	& {$8.0607e-03$}			& {$4.5735e-02$} \\ \hline  
{$180$}		& {$12$}		&{$3$}	& {$1.2215e-03$}			& {$1.3204e-02$} \\ \hline  
{$264$}		& {$30$}		&{$5$}	& {$1.5073e-04$}			& {$7.5582e-04$} \\ \hline  
{$1012$}	& {$132$}		&{$20$}	& {$2.3635e-05$}			& {$1.5490e-04$} \\ \hline  
\end{tabular}
\end{center} 
\label{c5solve} 
\end{table}

\begin{table}
\caption{$c$ is the fault} 
\begin{center} \footnotesize
\begin{tabular}{|l|l|l|l|l|l|l|l|l|} \hline  
Q 	&$\frac{\|D-\tilde{D}\|_F}{\|D\|_F}$	&$\frac{\|u-\tilde{u}\|_2}{\|u\|_2}$, left source 	&$\frac{\|u-\tilde{u}\|_2}{\|u\|_2}$, right source\\ \hline  
{$5$}	&{$2.8376e-01$}				&{$6.6053e-01$} 					&{$5.5522e-01$} \\ \hline  
{$5$}	&{$8.2377e-03$}				&{$3.8294e-02$} 					&{$2.4558e-02$} \\ \hline  
{$30$}	&{$1.1793e-03$}				&{$4.0372e-03$} 					&{$2.9632e-03$} \\ \hline  
\end{tabular}
\end{center} 
\label{c16solve} 
\end{table}

\begin{table}
\caption{$c$ is the diagonal fault} 
\begin{center} \footnotesize
\begin{tabular}{|l|l|l|l|l|l|l|l|l|} \hline  
Q 	&$\frac{\|D-\tilde{D}\|_F}{\|D\|_F}$	&$\frac{\|u-\tilde{u}\|_2}{\|u\|_2}$ \\ \hline  
{$4$}	&{$1.6030e-01$}				&{$4.3117e-01$}   \\ \hline  
{$6$}	&{$1.7845e-02$}				&{$7.1500e-02$}   \\ \hline  
{$23$}	&{$4.2766e-03$}				&{$1.2429e-02$}   \\ \hline  
\end{tabular}
\end{center} 
\label{c18solve} 
\end{table}

\begin{table}
\caption{$c$ is the periodic medium} 
\begin{center} \footnotesize
\begin{tabular}{|l|l|l|l|l|l|l|l|l|} \hline  
Q 	&$\frac{\|D-\tilde{D}\|_F}{\|D\|_F}$	&$\frac{\|u-\tilde{u}\|_2}{\|u\|_2}$ \\ \hline  
{$50$}	&{$1.8087e-01$}				&{$1.7337e-01$}   \\ \hline  
{$50$}	&{$3.5714e-02$}				&{$7.1720e-02$}   \\ \hline  
{$50$}	&{$9.0505e-03$}				&{$2.0105e-02$}   \\ \hline  
\end{tabular}
\end{center} 
\label{c33solve}
\end{table}

\subsection{Grazing waves}\label{sec:graz}

\begin{figure}[ht]
\begin{center}
\includegraphics[scale=.5]{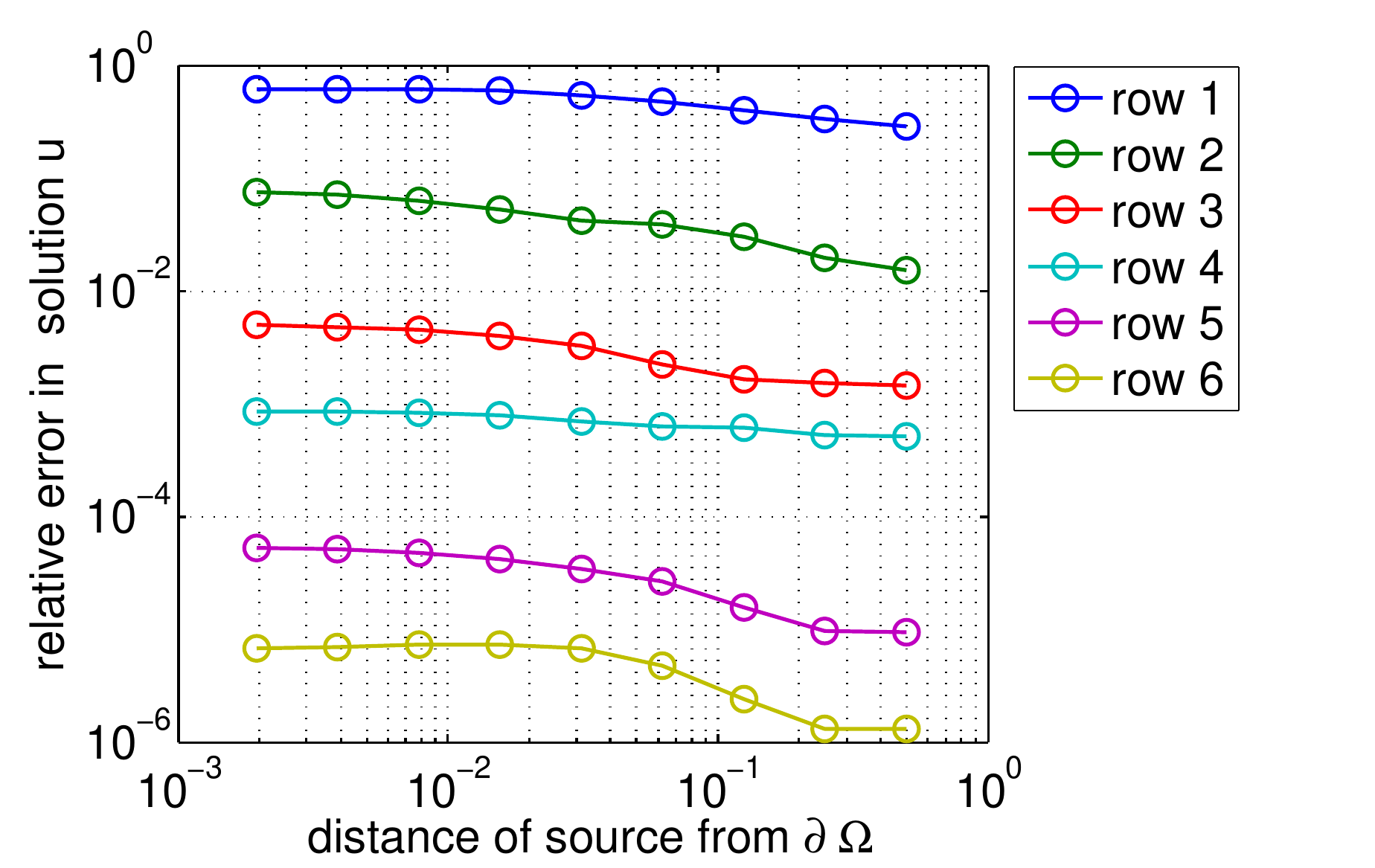}
\caption{Error in solution $u$, $c\equiv 1$, moving point source. Each line in the plot corresponds to using $\tilde{D}$ from a different row of Table \ref{c1solve}.}
\label{c1graz}
\end{center}
\end{figure}

It is well-known that ABCs often have difficulties when a source is close to a boundary of the domain, or in general when waves incident to the boundary are almost parallel to it. We wish to verify that the solution $\tilde{u}$ using the result $\tilde{D}$ of probing $D$ does not degrade as the source becomes closer and closer to some side of $\partial \Omega$. For this, we use a right-hand side $f$ to the Helmholtz equation which is a point source, located at the point $(x_0,y_0)$, where $x_0=0.5$ is fixed and $y_0>0$ becomes smaller and smaller, until it is a distance $2h$ away from the boundary (the point source's stencil has width $h$, so a source at a distance $h$ from the boundary does not make sense). We see in Figure \ref{c1graz} that, for $c \equiv 1$, the solution remains quite good until the source is a distance $2h$ away from the boundary. In this figure, we have used the probed maps we obtained in each row of Table \ref{c1solve}. %
We obtain very similar results for the waveguide, slow disk and faults (for the vertical fault we locate the source at $(x_0,y_0)$, where $y_0=0.5$ is fixed and $x_0$ goes to $0$ or $1$). This shows that the probing process itself does not significantly affect how well grazing waves are absorbed.

\subsection{Variations of $p$ with $N$}\label{sec:pwN}

\begin{figure}[ht]
\begin{minipage}[t]{0.45\linewidth}
\includegraphics[scale=.45]{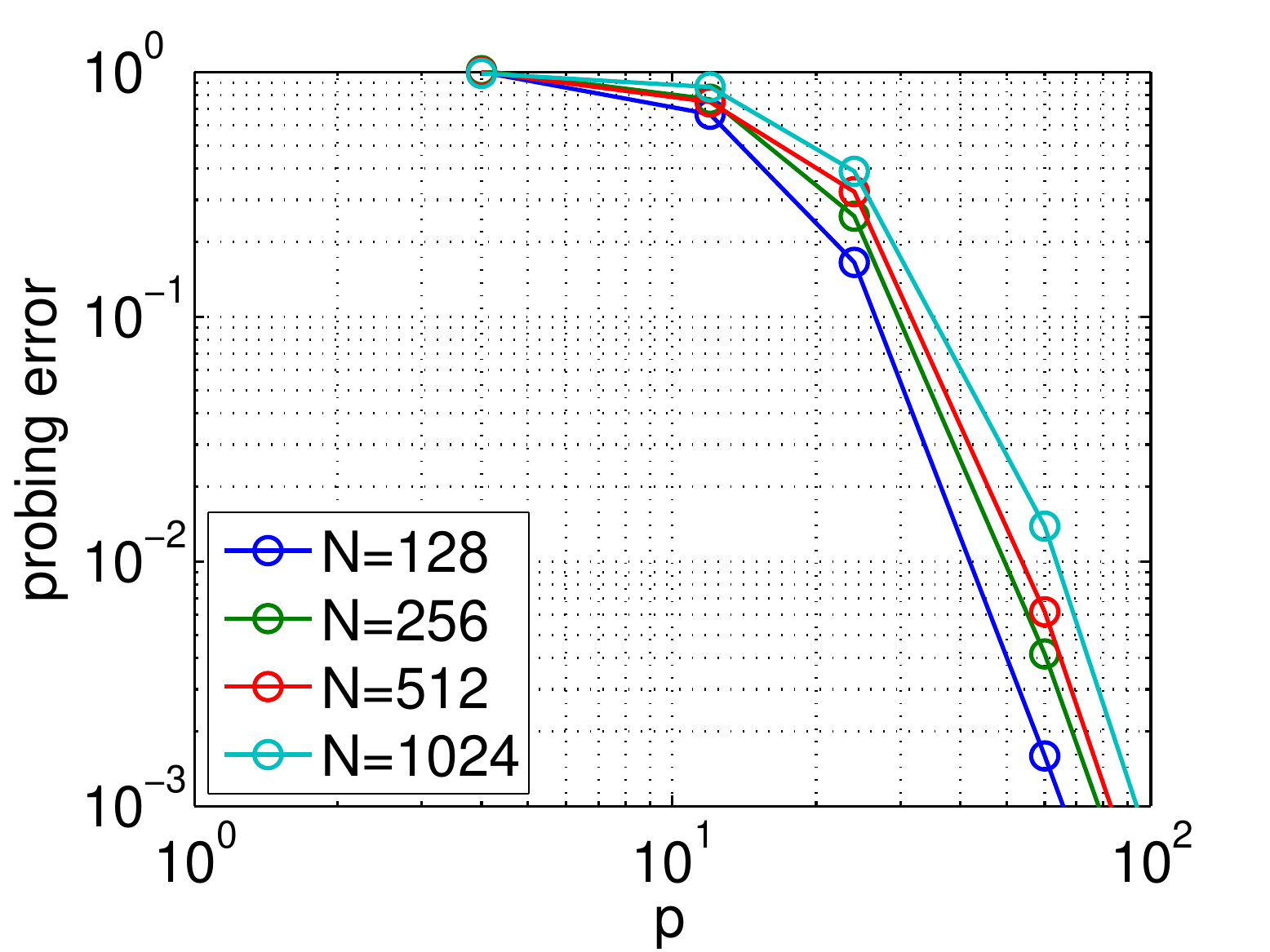}
\caption{Probing error of the $(1,1)$ block of the DtN map for the slow disk, fixed FD error level of $10^{-1}$, increasing $N$. This is the worst case, where $p$ follows a weak power law with $N$.}
\label{fig:c5pvsn}
\end{minipage}
\hspace{1cm}
\begin{minipage}[t]{0.45\linewidth}
\includegraphics[scale=.45]{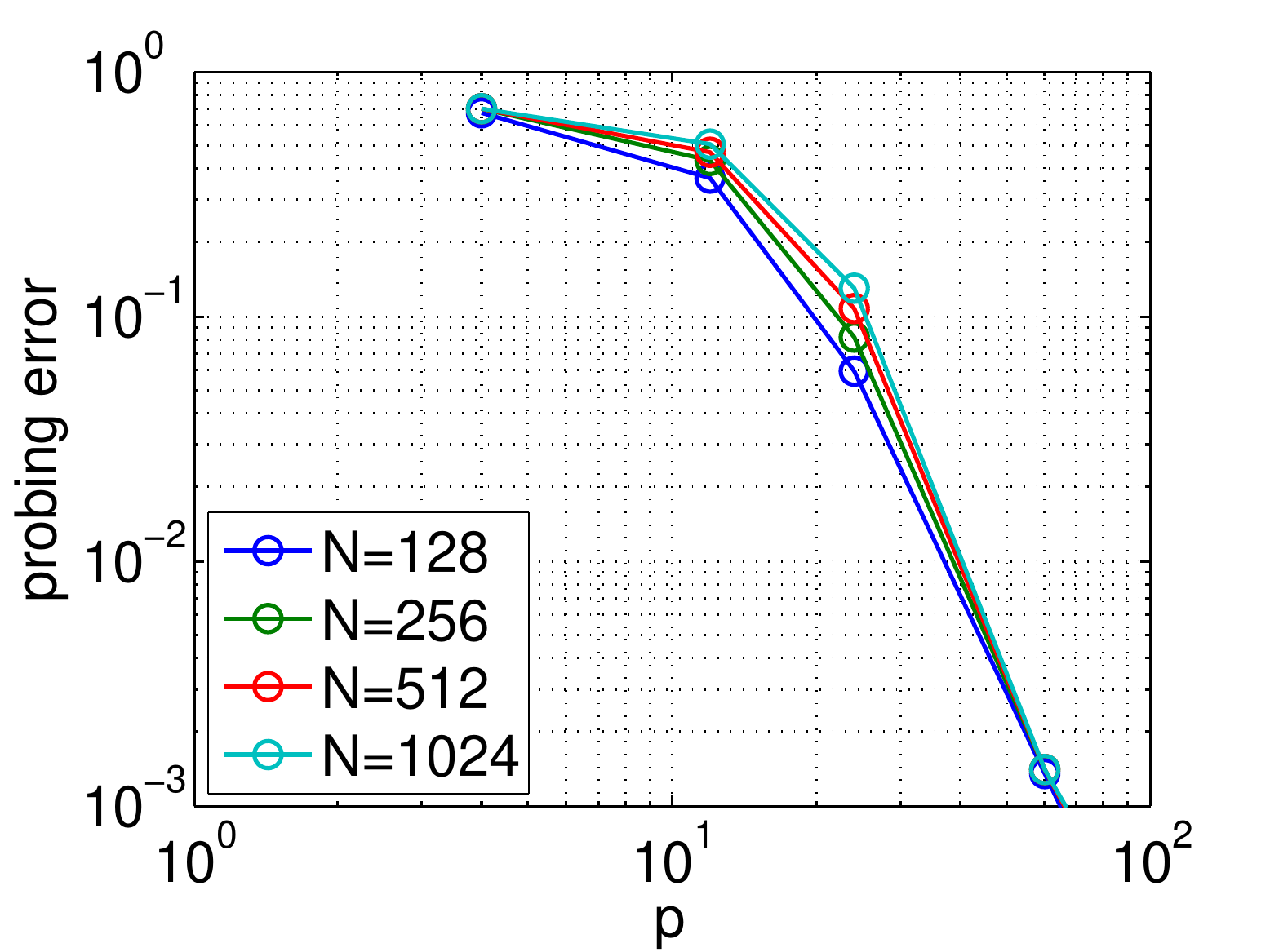}
\caption{Probing error of the $(1,1)$ block of the DtN map for the waveguide, fixed FD error level of $10^{-1}$, increasing $N$. Here $p$ follows a logarithmic law.}
\label{fig:c3pvsn}
\end{minipage}
\end{figure}

We now discuss how the number of basis matrices $p$ needed to achieve a desired accuracy depends on $N$ or $\omega$. To do this, we pick 4 consecutive powers of 2 as values for $N$, and find the appropriate $\omega$ such that the finite discretization error remains constant at $10^{-1}$, so that in fact $N \sim \omega^{1.5}$ as we have previously mentioned. We then probe the $(1,1)$ block of the corresponding DtN map, using the same parameters for all $N$, and observe the required $p$ to obtain a fixed probing error. The worst case we have seen in our experiments came from the slow disk. As we can see in Figure \ref{fig:c5pvsn}, $p$ seems to follow a very weak power law with $N$, close to $p \sim 15N^{.12}$ for a probing error of $10^{-1}$ or $p \sim 15N^{.2}$ for an probing error of $10^{-2}$. In all other cases, $p$ is approximately constant with increasing $N$, or seems to follow a logarithmic law with $N$ as for the waveguide (see Figure \ref{fig:c3pvsn}).

\section{Convergence of probing for the half-space DtN map: theorem}\label{sec:BasisPf}

In this section, we consider the half-space DtN map kernel in uniform medium
$K(r) = \frac{ik}{2r} H_1^{(1)}(kr)$
that we found in section \ref{sec:hsG}. We wish to approximate this kernel for values of $r$ that are relevant to our numerical scheme. Because we take $\Omega$ in our numerical experiments to be the $[0,1]^2$ box, $r=|x-y|$ will be between 0 and 1, in increments of $h$, as coordinates $x$ and $y$ along edges of $\pd \Omega$ vary between 0 and 1 in increments of $h$. However, as we know, $K(r)$ is singular at $r=0$, and since discretization effects dominate near the diagonal in the matrix representation of the DtN map, we shall consider only values of $r$ in the range $r_0 \leq r \leq 1$, with $0 < r_0 \leq 1/k$ (hence $r_0$ can be on the order of $h$). Since we know to expect oscillations $e^{ikr}$ in this kernel, we can remove those from $K$ to obtain%
\beq\label{eq:Hofr}
H(r) = \frac{ik}{2r} H_1^{(1)}(kr) e^{-ikr}
\eeq
(not to be confused with the hypersingular kernel $H$ of section \ref{sec:basis}), a smoother function which will be easier to approximate. Equivalently, we can add those oscillations to the terms in an approximation of $H$, to obtain an approximation of $K$.

For this section only, we denote by $\tilde{D}$ the corresponding operator with integral kernel $H$, while we use $D$ for the half-space Dirichlet-to-Neumann map, that is, the operator with kernel $K$.

\begin{theorem}\label{teo:main}
Let $\alpha > \frac{2}{3}$, $0 < r_0 < 1/k$, and let $K_p(r)$ be the best uniform approximation of $K(r)$ in 
$$\mbox{span} \{ \frac{e^{ikr}}{r^{j/\alpha}} : j = 1, \ldots, p, \mbox{ and } r_0 \leq r \leq 1 \}.$$
Denote by $D_p$ the operator defined with $K_p$ in place of $K$. Then, in the operator norm,
$$\| D - D_p \| \leq C_\alpha \, p^{1 - \lfloor 3\alpha/2 \rfloor} \, \| \tilde{K} \|_\infty,$$
for some $C_\alpha > 0$ depending on $\alpha$ and $r_0 \leq r \leq 1$.
\end{theorem}

The important point of the theorem is that the quality of approximation is otherwise independent of $k$, i.e., the number $p$ of basis functions does not need to grow like $k$ for the error to be small. In other words, it is unnecessary to ``mesh at the wavelength level" to spell out the degrees of freedom that go in the representation of the DtN map's kernel.

\begin{remark}
Growing $\alpha$ does not automatically result in a better approximation error, because a careful analysis of the proof shows that $C_\alpha$ grows factorially with $\alpha$. This behavior translates into a slower onset of convergence in $p$ when $\alpha$ is taken large, as the numerics show in the next section. This can in turn be interpreted as the result of ``overcrowding" of the basis by very look-alike functions.
\end{remark}

\begin{remark}
It is easy to see that the operator norm of $D$ grows like $k$, for instance by applying $D$ to the function $e^{-ik\x}$. %
The uniform norms of $K$ and $H$ once we cut out the diagonal, however, grow like $k^{1/2}/r_0^{3/2}$, so the result above shows that we incur an additional factor $k^{-1/2}r_0^{-3/2}$ in the error (somewhat akin to numerical pollution) in addition to the factor $k$ that we would have gotten from $\| D \|$.
\end{remark}

The result in Theorem \ref{teo:main} points the way for the design of basis matrices to be used in matrix probing, for the more general case of the exterior DtN map in heterogeneous media. We prove \ref{teo:main} in the next subsections, and present a numerical verification in the next section.

\subsection{Chebyshev expansion}

We mentioned the domain of interest for the $r$ variable is $[r_0,1]$. Again, expanding $K(r)$ in the system of Theorem \ref{teo:main} is equivalent to expanding $H(r)$ in polynomials of $r^{-1/\alpha}$ over $[r_0,1]$. It will be useful to perform the affine rescaling 
\[
\xi(r) = \frac{2}{r_0^{-1/\alpha} - 1} (r^{-1/\alpha} - 1 ) - 1 \qquad \Leftrightarrow \qquad r(\xi) = \left( \frac{\xi+1}{2} (r_0^{-1/\alpha} - 1) + 1 \right)^{-\alpha}
\]
so that the bounds $r \in [r_0,1]$ turn into $\xi \in [-1,1]$. We further write $\xi = \cos \theta$ with $\theta \in [0,\pi]$. Our strategy is to expand $H$ in Chebyshev polynomials $T_n(\xi)$. By definition, the best $p$-term approximation of $H(r)$ in polynomials of $r^{-1/\alpha}$ (best in a uniform sense over $[r_0,1]$) will result in a lower uniform approximation error than that associated with the $p$-term approximation of $H(r(\xi))$ in the $T_n(\xi)$ system. Hence in the sequel we overload notations and let $H_p$ for the $p$-term approximant of $H$ in our Chebyshev system.

We write out the Chebyshev series for $H(r(\xi))$ as 
$$H(r(\xi)) = \sum^{\infty}_{j=0} c_j T_j(\xi), \qquad c_j = \frac{2}{\pi} \int_{-1}^1 \frac{H(r(\xi)) T_j(\xi)}{(1-\xi^2)^{1/2}} \ d\xi, $$
with $T_j(\xi)=\cos{(j(\cos^{-1}\xi))}$, and $c_j$ alternatively written as
$$ c_j = \frac{2}{\pi} \int_0^\pi H(r(\cos{\theta})) \cos{j\theta} \ d \theta = \frac{1}{\pi} \int_0^{2\pi} H(r(\cos{\theta})) \cos{j\theta} \ d \theta. $$
The expansion will converge fast because we can integrate by parts in $\theta$ and afford to take a few derivatives of $H$, say $M$ of them, as done in \cite{tadmor}. After noting that the boundary terms cancel out because of periodicity in $\theta$, we express the coefficients $c_j$ for $j > 0$, up to a sign, as
$$ c_j  = \pm \frac{1}{\pi j^M} \int_0^{2\pi} \sin{j\theta} \frac{d^M}{d\theta^M} H(r(\cos{\theta})) \ d\theta, \qquad \ M \ \text{odd,} $$
$$ c_j  = \pm \frac{1}{\pi j^M} \int_0^{2\pi} \cos{j\theta} \frac{d^M}{d\theta^M} H(r(\cos{\theta})) \ d\theta, \qquad \ M \ \text{even.} $$
It follows that, for $j > 0$, and for all $M > 0$,
\[
\left| c_j \right| \leq \frac{2}{j^M} \max_\theta \left| \frac{d^M}{d\theta^M} H(r(\cos{\theta})) \right|.
\]
Let $B_M$ be a bound on this $M$-th order derivative. The uniform error we make by truncating the Chebyshev series to $H_p=\sum^{p}_{j=0} c_j T_j$ is then bounded by %
\begin{equation}\label{eq:bndderiv}
\left\| H-H_p \right\|_{L^\infty[r_0,1]} \leq \sum_{j=p+1}^\infty \left|c_j \right| \leq 2 B_M \sum_{j=p+1}^\infty \frac{1}{j^{M}} \leq \frac{2B_M}{(M-1) p^{M-1}}, \qquad \ p > 1.
\end{equation}
The final step is a simple integral comparison test.

\subsection{Bound on the derivatives of the DtN map kernel with oscillations removed}

The question is now to find a favorable estimate for $B_M$, from studying successive $\theta$ derivatives of $H(r)$ in \eqref{eq:Hofr}.  From the bound for the derivatives of Hankel functions in Lemma 1 of \cite{flatland}: given any $C > 0$, we have
\begin{equation}\label{eq:derivH}
\left| \frac{d^m}{dr^m} \left( H_1^{(1)}(kr) e^{-ikr} \right) \right| \leq C_m (kr)^{-1/2} r^{-m} \qquad \ \text{for} \ kr \geq C.
\end{equation}
The change of variables from $r$ to $\theta$ results in
\begin{eqnarray*}
\frac{dr}{d\theta}&=& \frac{d\xi}{d\theta} \frac{dr}{d\xi} = \left(-\sin \theta \right) \left( -\alpha  \left( \frac{\xi+1}{2} (r_0^{-1/\alpha} - 1) + 1 \right)^{-\alpha-1} \frac{\left(r_0^{-1/\alpha}-1\right)}{2} \right) \\
&=&\left(-\sin \theta \right) \left( -\alpha \ r^{1+1/\alpha} \ \frac{r_0^{-1/\alpha}(1-r_0^{1/\alpha})}{2} \right).
\end{eqnarray*}
Hence
\begin{equation}\label{eq:drdtheta}
\frac{dr}{d\theta} = r (r/r_0)^{1/\alpha} \ \frac{\alpha \sin \theta (1-r_0^{1/\alpha})}{2}.
\end{equation}
Derivatives of higher powers of $r$ are handled by the chain rule, resulting in
\begin{equation}\label{eq:derivr}
\frac{d}{d \theta} (r^p) = p r^p (r/r_0)^{1/\alpha} \ \frac{\alpha \sin \theta (1-r_0^{1/\alpha})}{2}.
\end{equation}
We see that the action of a $\theta$ derivative is essentially equivalent to multiplication by $(r/r_0)^{1/\alpha}$. As for higher derivatives of powers of $r$, it is easy to see by induction that the product rule has them either hit a power of $r$, or a trigonometric polynomial of $\theta$, resulting in a growth of at most $(r/r_0)^{1/\alpha}$ for each derivative:
\[
| \frac{d^m}{d\theta^m} r^p | \leq C_{m,p,\alpha} \, r^p (r/r_0)^{m/\alpha}.
\]
These estimates can now be combined to bound $\frac{d^m}{d \theta^m} \left( H_1^{(1)}(kr) e^{-ikr} \right)$. One of two scenarios occur when applying the product rule:
\bit
\item either $\frac{d}{d\theta}$ hits $\frac{d^{m_2}}{d\theta^{m_2}} \left( H_1^{(1)}(kr) e^{-ikr} \right)$ for some $m_2 < m$. In this case, one negative power of $r$ results from $\frac{d}{dr}$ as we saw in \eqref{eq:derivH}, and a factor $r (r/r_0)^{1/\alpha}$ results from $\frac{dr}{d\theta}$ as we saw in \eqref{eq:drdtheta};
\item or $\frac{d}{d\theta}$ hits some power of $r$, possibly multiplied by some trigonometric polynomial in $\theta$, resulting in a growth of an additional factor $(r/r_0)^{1/\alpha}$ as we saw in \eqref{eq:derivr}.
\eit
Thus, we get at most a $(r/r_0)^{1/\alpha}$ growth factor per derivative in every case. The situation is completely analogous when dealing with the slightly more complex expression $\frac{d^m}{d \theta^m} \left( \frac{1}{r} H_1^{(1)}(kr) e^{-ikr} \right)$. The number of terms is itself at most factorial in $m$, hence we get
\begin{equation}\label{eq:derivbnd}
| \frac{d^m}{d \theta^m} \frac{k}{r} \left( H_1^{(1)}(kr) e^{-ikr} \right) | \leq C_{m, \alpha} \ \frac{k}{r} \left( \frac{r}{r_0} \right)^{\frac{m}{\alpha} - \frac{1}{2}} \leq C_{m,\alpha} \ \frac{k}{r_0} \left( \frac{r}{r_0} \right)^{\frac{m}{\alpha} - \frac{3}{2}}.
\end{equation}
We now pick $m \leq M = \lfloor 3 \alpha /2 \rfloor$, so that the max over $\theta$ is realized when $r = r_0$, and $B_M$ is on the order of $k/r_0$. It follows from \eqref{eq:bndderiv} and \eqref{eq:derivbnd} that
$$ \left\| H-H_p \right\|_{L^\infty[r_0,1]} \leq  C_\alpha \ \frac{k}{r_0} \ \frac{1}{p^{\lfloor 3\alpha/2 \rfloor - 1}}, \qquad \ p > 1, \ \alpha > 2/3.$$
The kernel of interest, $K(r) = H(r) e^{ikr}$ obeys the same estimate if we let $K_p$ be the $p$-term approximation of $K$ in the Chebyshev system modulated by $e^{ikr}$.

\subsection{Bound on the error of approximation}

For ease of writing, we now let $D^0$, $D^0_p$ be the operators with respective kernels $K^0(r) = K(r) \chi_{[r_0,1]}(r)$ and $K^0_p(r) = K_p(r) \chi_{[r_0,1]}(r)$. We now turn to the operator norm of $D^0 - D^0_p$ with kernel $K^0 -K^0_p$:
$$(D^0-D^0_p)g(x)=\int_0^1 (K^0-K^0_p)(|x-y|)g(y) \ dy, \qquad x \in [0,1]. $$
We use the Cauchy-Schwarz inequality to bound
\begin{eqnarray*}
\|(D^0-D^0_p)g\|_2 &=& \left(\int_{0\leq x \leq 1}  \left|\int_{0\leq y \leq 1, \ |x-y|\geq r_0} (K^0-K^0_p)(|x-y|)g(y) \ dy \right|^2 dx \right)^{1/2} \\
& \leq & \left(\int_{0\leq x \leq 1} \int_{0\leq y \leq 1, \ |x-y|\geq r_0}\left|(K^0-K^0_p)(|x-y|)\right|^2 \ dy dx \right)^{1/2} \|g\|_2 \\
& \leq & \left( \int_{0\leq x \leq 1} \int_{0\leq y \leq 1, \ |x-y|\geq r_0} 1 \ dy \ dx \right)^{1/2} \|g\|_2 \ \max_{0\leq x,y \leq 1, \ |x-y| \geq r_0} |(K^0-K^0_p)(|x-y|)|  \\
& \leq &  \|g\|_2 \ \| K^0-K^0_p \|_{L^{\infty}[r_0,1]}.
\end{eqnarray*}
Assembling the bounds, we have
\[
\| D^0-D^0_p \|_2 \leq \| K^0-K^0_p \|_{L^{\infty}[r_0,1]} \leq C_{\alpha} \, p^{1 - \lfloor 3 \alpha / 2 \rfloor} \, \frac{k}{r_0}.
\]
It suffices therefore to show that $\| K^0 \|_\infty = \| K \|_{L^\infty[r_0,1]}$ is bigger than $k/r_0$ to complete the proof. Letting $z = kr$, we see that
$$\max_{r_0 \leq r \leq 1} |K(r)|=\frac{k}{2r_0} \max_{kr_0 \leq z \leq k} \left|H_1^{(1)}(z) \right| \geq C \frac{k^{1/2}}{r_0^{3/2}}.$$
The last inequality follows from the fact that there exist a positive constant $c_1$ such that $c_1 z^{-1/2} \leq \left|H_1^{(1)}(z) \right| $, from Lemma 3 of \cite{flatland}. But $k^{1/2}/r_0^{3/2} \geq k/r_0$ precisely when $r_0 \leq 1/k$. Hence we have proved the statement of Theorem \ref{teo:main}.

In the next section, we proceed to a numerical confirmation of Theorem \ref{teo:main}.

\section{Convergence of probing for the half-space DtN map: numerical confirmation}\label{sec:NumPf}

\begin{figure}[ht]
\begin{minipage}[t]{0.48\linewidth}
\includegraphics[scale=.5]{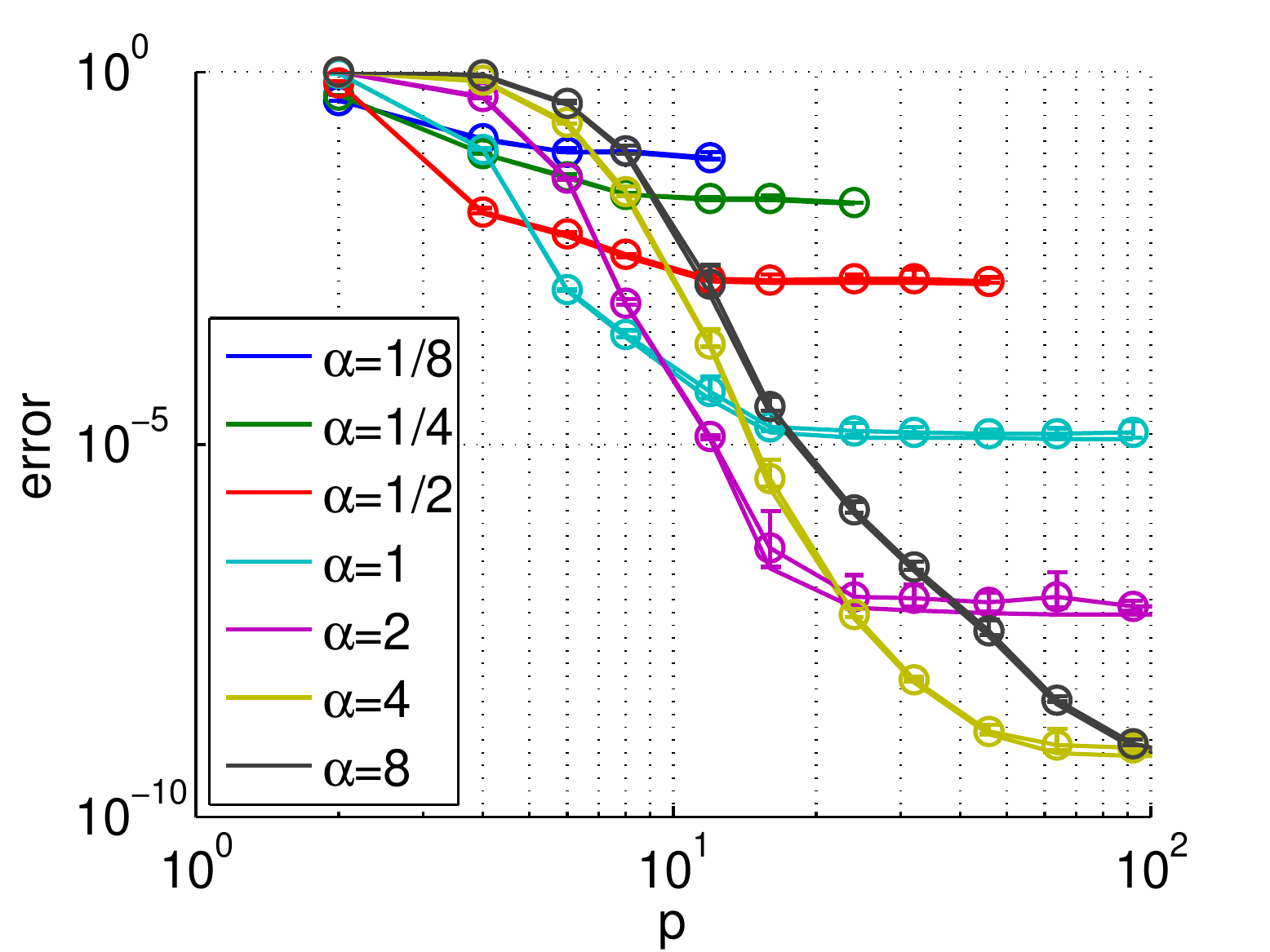}
\caption{Probing error of the half-space DtN map($q=1$, 10 trials, circle markers and error bars) compared to the approximation error (line), $c\equiv 1$, $L=1/4$, $\alpha=2$, $n=1024$, $\omega=51.2$.}
\label{q1}
\end{minipage}
\hspace{0.1cm}
\begin{minipage}[t]{0.48\linewidth}
\includegraphics[scale=.5]{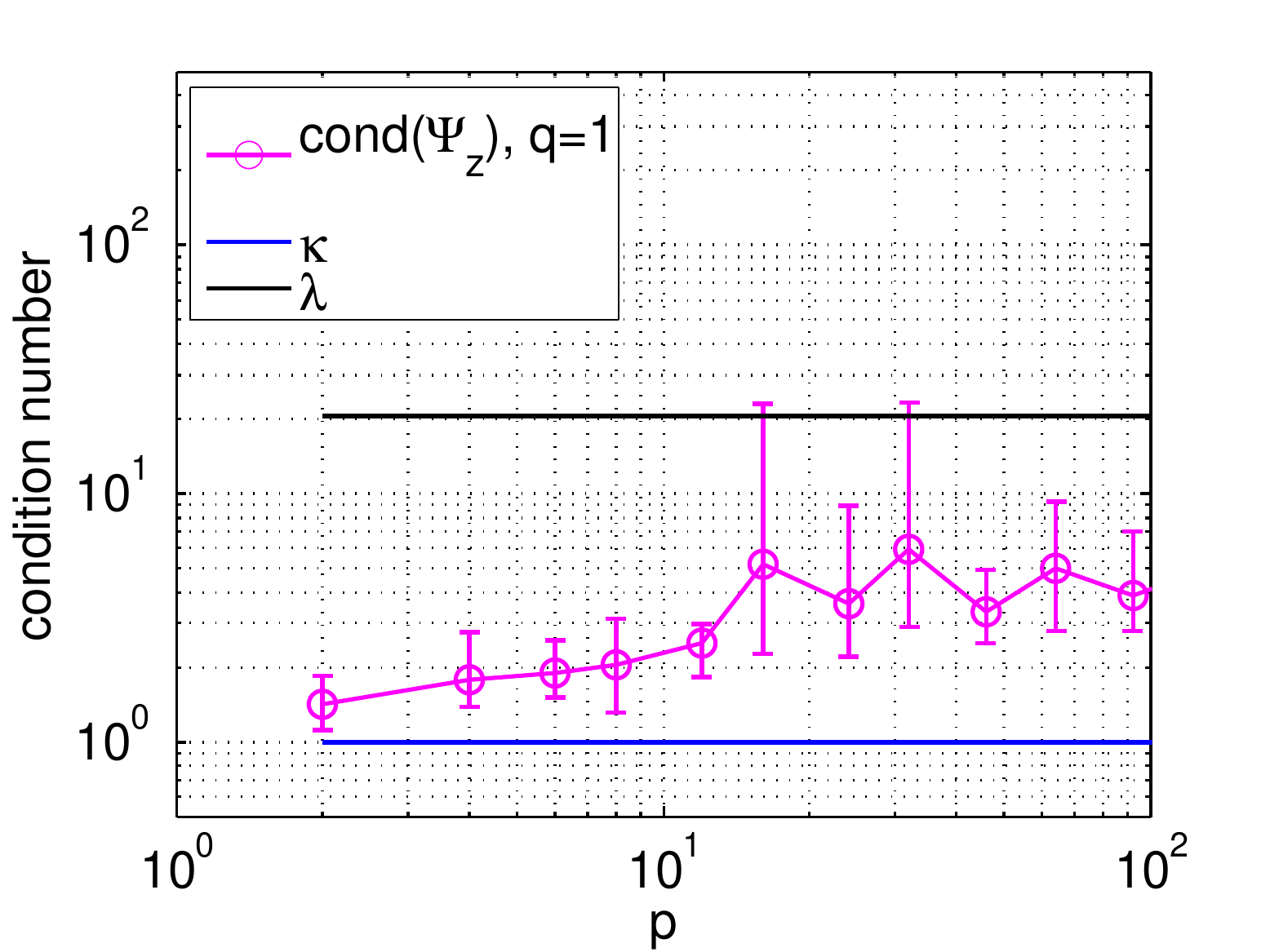}
\caption{Condition numbers for probing the half-space DtN map, $c\equiv 1$, $L=1/4$, $\alpha=2$, $n=1024$, $\omega=51.2$, $q=1$, 10 trials.}
\label{p2}
\end{minipage}
\end{figure}

In order to use Theorem \ref{teo:main} to obtain convergent basis matrices, we start with the set $\left\{(r)^{-j/\alpha}  \right\}_{j=0}^{p-1}$. We have proved the theorem for the interval $r_0 \leq r \leq 1$, but here we consider $h \leq r \leq 1$, which is a larger interval than in the theorem. We then put in oscillations, orthonormalize, and use this new set as a basis for probing the DtN map. Thus we have pre-basis matrices ($0 \leq j \leq p$)
$$(\beta_j)_{\ell m}=\frac{e^{ikh|\ell-m|}}{|\ell-m|^{j/\alpha}} \ \text{for} \ \ell \neq m,$$
with $(\beta_j)_{\ell \ell}=0$. We add to this set the identity matrix in order to capture the diagonal of $D$, and orthonormalize the resulting collection to get the $B_j$. Alternatively, we have noticed that orthonormalizing the set of $\beta_j$'s with
\begin{equation}\label{halfbasis}
 (\beta_j)_{\ell m}=\frac{e^{ikh|\ell-m|}}{(h+h|\ell-m|)^{j/\alpha}}
\end{equation}
works just as well, and is simpler because there is no need to treat the diagonal separately. We use this same technique for the probing basis matrices of the exterior problem.

The convergent basis matrices in \eqref{halfbasis} have been used to obtain a numerical confirmation of Theorem \ref{teo:main}, again for the half-space DtN map. To obtain the DtN map in this setup, instead of solving the exterior problem with a PML or pPML on all sides, we solve a problem on a thin strip, with a random Dirichlet boundary condition (for probing) one of the long edges, and a PML or pPML on the other three sides. This method for a numerical approximation of the solution to the half-space problem was introduced in section \ref{sec:half} and Figure \ref{fig:half}.

\subsection{Uniform medium}

In Figure \ref{q1}, we show the approximation error, which we expect will behave as in Theorem \ref{teo:main}. We also plot error bars for the probing error, corresponding to ten trials of probing, with $q=1$. The probing results are about as good as the approximation error, because the relevant condition numbers are all well-behaved as we see in Figure \ref{p2} for the value of choice of $\alpha=2$. Back to the approximation error, we notice in Figure \ref{q1} that increasing $\alpha$ delays the onset of convergence as expected, because of the factor which is factorial in $\alpha$ in the error of Theorem \ref{teo:main}. And we can see that, for small $\alpha$, we are taking very high inverse powers of $r$, an ill-conditioned operation. Hence the appearance of a convergence plateau for smaller $\alpha$ is explained by ill-conditioning of the basis matrices, and the absence of data points is because of computational overflow. 

Finally, increasing $\alpha$ from $1/8$ to $2$ gives a higher rate of convergence, as it should because in the error we have the factor $p^{-3\alpha/2}$, which gives a rate of convergence of $3\alpha/2$. This is roughly what we obtain numerically. As discussed, further increasing $\alpha$ is not necessarily advantageous since the constant $C_\alpha$ in \ref{teo:main} grows fast in $\alpha$.

%% file: plr2.tex
\chapter{Partitioned low-rank matrices for compressing the Dirichlet-to-Neumann map}\label{ch:plr}

In the previous chapter, we explained in detail the first step of two in our numerical scheme for compressing ABC's. This first step consisted in approximating the DtN map $D$ by $\tilde{D}$ using matrix probing. To do this, we considered each block $M$ of $D$ separately, corresponding to an edge-to-edge restriction of $D$. We approximated each $M$ by a matrix probing expansion. %
We saw how we could obtain an accurate approximation of $M$ using appropriate basis matrices $B_j$. We were left with the task of producing a fast algorithm for applying the resulting $\tilde{M}$ to vectors, since this is an operation that a Helmholtz solver needs. This is what we do in the current chapter, by compressing $\tilde{M}$ into a new $\overline{M}$ which can then be applied fast.

We also explained in the previous chapter why we needed probing: to obtain an explicit approximation of $D$, to be compressed in order to obtain a fast matrix-vector product. Indeed, we do not have direct access to the entries of $D$, but rather we need to solve the costly exterior problem every time we need a multiplication of $D$ with a vector. We have already mentioned how the approach of Lin et al. \cite{Hmatvec}, for example, would require $O(\log N)$ such solves, with a large constant.

We alluded that we might be able to compress each block $M$ (or $\tilde{M}$) of $D$ (or $\tilde{D}$) when we presented background material in chapter \ref{ch:back}. Indeed, we discussed the fact that the half-space Green's function $G_\text{half}$ in constant medium is separable and low-rank away from its singularity. Because the half-space DtN map kernel $K$ is simply two derivatives of $G_\text{ext}$, we expect $K$ to also be separable and low-rank, and we prove this at the end of the present chapter, in section \ref{sec:septheo}. See also the numerical verification of that theorem in section \ref{sec:sepnum}. Because the half-space DtN map is strongly related to the exterior DtN map as we mentioned in chapter \ref{ch:back}, we expect the exterior DtN map kernel to also be separable and low rank, at least in favorable conditions such as a constant medium. But first, as in the previous chapter, we begin by explaining the technique we use, partitioned low-rank (PLR) matrices, in section \ref{sec:introplr}. Compression of an $N$ by $N$ matrix into the PLR framework is nearly linear in $N$, and so is matrix-vector multiplication. We then show the results of using this PLR technique on test cases in section \ref{sec:usingplr}. %

\section{Partitioned low-rank matrices}\label{sec:introplr}

As we have discussed in chapter \ref{ch:back}, when an operator is separable and low-rank, we expect its numerical realization to have low-rank blocks under certain conditions. In our case, the DtN map $K(x-y)$ is separable and low-rank away from the singularity $x=y$ and so we expect its numerical realization to have low-rank blocks away from its diagonal. This calls for a compression scheme such as the hierarchical matrices of Hackbush et al. \cite{hmat1}, \cite{hmat2}, \cite{hmatlect}, to compress off-diagonal blocks. However, because we expect higher ranks away from the singularity in variable media, and because different blocks of the DtN map will show a singularity elsewhere than on the diagonal, we decide to use a more flexible scheme called \emph{partitioned low rank} matrices, or PLR matrices, from \cite{Jones}.

\subsection{Construction of a PLR matrix}

PLR matrices are constructed recursively, using a given tolerance $\epsilon$ and a given maximal rank $R_{\text{max}}$. We start at the top level, level 0, with the matrix $M$ which is $N$ by $N$ and $N$ is a power of two\footnote{Having a square matrix with dimensions that are powers of two is not necessary, but makes the discussion easier.}. We wish to compress $M$ (in the next sections we will use this compression scheme on probed blocks $\tilde{M}$ of the DtN map, but we use $M$ here for notational simplicity). We first ask for the numerical rank $R$ of $M$. The numerical rank is defined by the Singular Value Decomposition and the tolerance $\epsilon$, as the number $R$ of singular values that are larger than or equal to the tolerance. If $R>R_{\text{max}}$, we split the matrix in four blocks and recurse to the next level, level 1 where blocks are $N/2^1$ by $N/2^1$. If the numerical rank of $M$ is less than or equal to $R_{\text{max}}$, $R \leq R_{\text{max}}$, we express $M$ in its low-rank form by truncating the SVD of $M$ after $R$ terms. That is, the SVD of $M$ is $M=U\Sigma V^*=\sum_{j=1}^{N} U_j \sigma_j V_j^*$ where $U$ and $V$ are orthonormal matrices with columns $\left\{U_j\right\}_{j=i}^{N}$ and $\left\{V_j\right\}_{j=i}^{N}$, and $\Sigma$ is the diagonal matrix of decreasing singular values: $\Sigma=\text{diag}(\sigma_1, \sigma_2, \ldots, \sigma_N)$. Then, if $R \leq R_{\text{max}}$, we compress $M$ to $\overline{M}=\sum_{j=1}^{R} U_j \sigma_j V_j^*$ by truncating the SVD of $M$ after $R$ terms.

If we need to split $M$ and recurse down to the next level, we do the following. First, we split $M$ in four square blocks of the same size: take the first $N/2$ rows and columns to make the first block, then taking the first $N/2$ rows and last $N/2$ columns to make the second block, etc. And now we apply the step described in the previous paragraph to each block of $M$, checking the block's numerical rank and compressing it or splitting it depending on that numerical rank. Whenever we split up a block, we label it as ``hierarchical'', and call its four sub-blocks its \emph{children}. Whenever a block was not divided, and hence compressed instead, we label it as ``compressed'', and we may call it a ``leaf'' as well.

If a block has dimension $R_{\text{max}}$ by $R_{\text{max}}$, then its numerical rank will be $R_{\text{max}}$, and so once blocks have dimensions smaller than or equal to the maximal desired rank $R_{\text{max}}$, we can stop recursing and store the blocks directly. However, especially if $R_{\text{max}}$ is large, we might still be interested in compressing those blocks using the SVD. This is what we do in our code, and we label such blocks as ``compressed'' as well. When we wish to refer to how blocks of a certain matrix $M$ have been divided when $M$ was compressed in the PLR framework, or in particular to the set of all leaf blocks and their positions in $M$, we refer to the ``structure'' of $M$. We see then that the structure of a PLR matrix will have at most $L$ levels, where $N/R_\text{max}=2^L$ so $L=\log_2{N/R_\text{max}}$.

\subsubsection{Implementation details}

Algorithm \ref{alg:PLR_matrix} presents pseudocode for the construction of a PLR matrix from a dense matrix. In practice, when we compute the SVD of a block, we use the randomized SVD\footnote{Theoretically, this randomized SVD has a failure probability, but we can choose a parameter to make this probability on the order of $10^{-16}$, and so we ignore the fact that the randomized SVD could fail.} \cite{randomSVD}. This allows us to use only a few matrix-vector multiplies between the block (or its transpose) and random vectors to form an approximate reduced SVD. This is a faster way of producing the SVD, and thus also of finding out whether the numerical rank of the block is larger than $R$. The randomized SVD requires about 10 more random matrix-vector multiplies than the desired maximal rank $R$. This is why, in algorithm \ref{alg:PLR_matrix}, the call to \emph{svd} has two arguments: the block we want to find the SVD of, and the maximal desired rank $R_{\text{max}}$. The randomized SVD algorithm then uses 10 more random vectors than the quantity $R_{\text{max}}$ and returns an SVD of rank $R_{\text{max}}+1$. We use the $R_{\text{max}}+1^{\text{st}}$ singular value in $\Sigma$ to test whether we need to split the block and recurse or not.

\begin{algorithm} Compression of matrix $M$ into Partitioned Low Rank form, with maximal rank $R_{\text{max}}$ and tolerance $\epsilon$ \label{alg:PLR_matrix}
\begin{algorithmic}[1]
\Function{H = PLR}{$M$, $R_{\text{max}}$, $\epsilon$}
\State  $[U, \Sigma ,V] = \texttt{svd}(M, R_{\text{max}})$ \Comment{Randomized SVD}
\If{ $\exists R \in \left\{1,2,\ldots,R_{\text{max}}\right\} : \Sigma(R+1,R+1) < \epsilon$}
\State Let $R$ be the smallest such integer.
\State \texttt{H.data = \{$U(:,1:R)\cdot \Sigma(1:R,1:R)$, $V(:,1:R)^{*}$\} }
\State \texttt{H.id = 'c' }    \Comment{This block is ``compressed''}
\Else \Comment{The $M_{ij}$'s are defined in the text}
\For{i = 1:2}
\For{j = 1:2}
\State \texttt{H.data\{i,j\} = } PLR($M_{ij}$, $R_{\text{max}}$, $\epsilon$ ) \Comment Recursive call
\EndFor
\EndFor
\State \texttt{H.id = 'h' }    \Comment{This block is ``hierarchical''}
\EndIf
\EndFunction
\end{algorithmic}
\end{algorithm}

\subsubsection{Complexity}
The complexity of this algorithm depends on the complexity of the SVD algorithm we use. The randomized SVD has complexity $O(N_B \Rm^2)$ where $N_B$ is the dimension of block $B$ whose SVD we are calculating. This can be much better, especially for larger blocks, than standard SVD algorithms which have complexity $O(N_B^3)$. The total complexity of the compression algorithm will then depend on how many blocks of which size and rank we find the randomized SVD of. We shall discuss this in more detail when we discuss also the complexity of a matrix-vector product for special structures.

\subsubsection{Error analysis}
To understand the error we make by compressing blocks in the PLR framework, we first note that compressing a block $M$ to its $R$-rank approximation $\overline{M}=\sum_{j=1}^{R} U_j \sigma_j V_j^*$ by truncating its SVD, we make an error in the $L_2$ norm of $\sigma_{R+1}$. That is,
\[ \|M-\overline{B} \|_2 = \| \sum_{j=R+1}^{N} U_j \sigma_j V_j^* \|_2 = \sigma_{R+1} . \]
Hence, by compressing a block, we make an $L_2$ error for that block of at most $\eps$ because we make sure that $\sigma_{R+1} \leq \eps$. Of course, errors from various blocks will compound to affect the total error we make on matrix $M$. We shall mention this in more detail when we discuss particular structures. 

The relative Frobenius error we make between $M$ and $\overline{M}$, the compressed approximation of $M$, will usually be larger than $\epsilon$ because of two factors. First of all, as we just saw, the PLR compression algorithm uses the $L_2$ norm. To use the Frobenius norm when deciding whether to compress or divide a blcok, we would need access to all the dingular values of each block. This would be possible using typical SVD algorithm, but quite costly. Hence we use the randomized SVD, which is faster but with which we are forced to use the $L_2$ norm. Another factor we have yet to mention is that errors from different blocks will compound to make the total error between $M$ and $\overline{M}$ larger than the error between any individual blocks of $M$ and $\overline{M}$. This of course depends on how many blocks there are in the structure of nay particular PLR matrix. We will talk more about this in subsection \ref{sec:structncomp}, where we explore the complexity of the compression and matrix-vector algorithms of the PLR framework. But first, we introduce matrix-vector products.

\subsection{Multiplication of a PLR matrix with a vector}

To multiply a PLR matrix $M$  with a vector $v$, we again use recursion. Starting at the highest level block, the whole matrix $M$ itself, we ask whether this block has been divided into sub-blocks. If not, we multiply the block directly with the vector $v$. If the block has been subdivided, then we ask for each of its children whether those have been subdivided. If not, we multiply the sub-block with the appropriate restriction of $v$, and add the result to the correct restriction of the output vector. If so, we recurse again. Algorithm \ref{alg:PLR_matvec} presents the pseudocode for multiplying a vector by a PLR matrix. The algorithm to left-multiply a vector by a matrix is similar, we do not show it here.

\begin{algorithm} Multiplication of a PLR matrix $H$ with column vectors $x$ \label{alg:PLR_matvec}
\begin{algorithmic}[1]
\Function{y = matvec}{H,x}
\If{ \texttt{H.id == 'c' }  }
\State \texttt{y =  H.data\{1\}$\cdot$(H.data\{2\}$\cdot$x) }
\Else 
\State $\texttt{y}_{1}\texttt{ = }\text{matvec}\texttt{(H.data\{1,1\},}\texttt{x}(\texttt{1:end/2,:}))$
\State $\qquad +\text{matvec}\texttt{(H.data\{1,2\},}\texttt{x}(\texttt{end/2:end,:}))$
\State $\texttt{y}_{2}\texttt{ = }\text{matvec}\texttt{(H.data\{2,1\},}\texttt{x}(\texttt{1:end/2,:}))$
\State $\qquad +\text{matvec}\texttt{(H.data\{2,2\},}\texttt{x}(\texttt{end/2:end,:}))$
\State \texttt{y = $ \left[ \begin{array}{c} 
										\texttt{y}_1\\
										\texttt{y}_2
										\end{array}
							  \right ]$ }
\EndIf
\EndFunction
\end{algorithmic}
\end{algorithm}

\subsubsection{Complexity}
The complexity of this algorithm is easily understood.  Recall that we store blocks that are not divided not as a full matrix, but as column vectors corresponding to the orthonormal matrix $U$ of the SVD, and to the product $\Sigma V^*$ of the SVD. Every time we multiply such an $N_B$ by $N_B$ block $B$ that has not been subdivided with the corresponding restriction of $v$, we first multiply the restriction $\tilde{v}$ of $v$ with $\Sigma V^*$, and then multiply the result with $U$. Let $R_B \leq R_\text{max}$ be the numerical rank of that block. Then, we first make $N_BR_B$ multiplication operations and $(N_B-1)R_B$ addition operations for the product $\Sigma V^* \tilde{v}$, and then $R_BN_B$ multiplication operations and $(R_B-1)N_B$ addition operations for the product of that with $U$. Hence we make approximately $4N_BR_B$ operations per block, where again $N_B$ is the dimension of the block $B$ and $R_B$ is its numerical rank.

The total number of operations for multiplying a vector $v$ with a PLR matrix $M$, then, is about 
\begin{equation}\label{eq:comp}
\sum_{B \text{ is compressed}} 4 N_B R_B ,
\end{equation}
where we sum over all ``compressed'' blocks.

Evidently there is a trade-off here. Asking for a small maximal rank $R_\text{max}$ may force blocks to be subdivided a lot. We then have small $R_B$'s, and small $N$'s, but a lot of blocks. On the other hand, having a larger $R_\text{max}$ means a lot of blocks can be remain large. We then have large $R_B$'s and $N_B$'s, but very few blocks. We shall use the complexity count \eqref{eq:comp} later on to decide on which $R_\text{max}$ to choose for each block of the DtN map. To have an idea of whether such a matrix-vector multiplication might result in a fast algorithm, we need to introduce new terminology regarding the structure of PLR matrices, and we do so in the next subsection.

\subsection{Structure and complexity}\label{sec:structncomp}

As we did with matrix probing, we again split up the probed exterior DtN map $\tilde{D}$ in submatrices corresponding to the different sides of $\pd \Omega$. We called those \emph{blocks} of $\tilde{D}$ in the previous chapter, but we now call them \emph{submatrices}, to differentiate them from the blocks in the structure we obtain from compressing to a PLR matrix. So $\tilde{D}$ is split up in submatrices $\tilde{M}$ that represent edge-to-edge restrictions of $\tilde{D}$. We then only have to compress once each unique submatrix $\tilde{M}$ to obtain an approximation $\overline{M}$ of $\tilde{M}$. We can use those compressed $\overline{M}$'s to define $\overline{D}$, an approximation of $\tilde{D}$. If we need to multiply our final approximation $\overline{D}$ by a vector, we may then split that vector in blocks corresponding to the sides of $\pd \Omega$ and use the compressed submatrices and the required PLR matrix algebra to obtain the result with low complexity.

What is particular about approximating the DtN map on a square boundary $\pd \Omega$ is that distinct submatrices are very different. Those that correspond to the restriction of $D$ from one edge to that same edge, and as such are on the diagonal of $D$, are harder to probe as we saw in the previous chapter because of the diagonal singularity. And, because of the diagonal singularity, they might be well-suited for compression by hierarchical matrices \cite{bebendorf}.

However, submatrices of $D$ corresponding to two edges that are side by side (for example, the bottom and right edges of the boundary of $[0,1]^2$) see the effects of the diagonal of $D$ in their upper-right or lower-left corners, and entries of such submatrices decay in norm away from that corner. Thus a hierarchical matrix would be ill-suited to compress such a submatrix. This is why the PLR framework is so useful to us: it automatically adapts to the submatrix at hand, and to whether there is a singularity in the submatrix, and where that singularity might be. 

Similarly, when dealing with a submatrix of $D$ corresponding to opposite edges of $\pd \Omega$, we see that entries with higher norm are in the upper-right and bottom-left corners, so again PLR matrices are more appropriate than hierarchical ones. However, note that because such submatrices have very small relative norm compared to $D$, and were probed with only one or two basis matrices in the previous chapter, their PLR structure is often trivial.

In order to help us understand the complexity of PLR compression and matrix-vector products, we first study typical structures of hierarchical \cite{hmatlect}, \cite{Hweak} and PLR matrices.

\subsubsection{Weak hierarchical matrices}
\begin{definition}
A matrix is said to have \emph{weak hierarchical structure} when a block is compressed if and only if its row and column indices do not overlap.
\end{definition}

\begin{figure}[H]
\begin{minipage}[t]{0.32\linewidth}
\includegraphics[scale=.2]{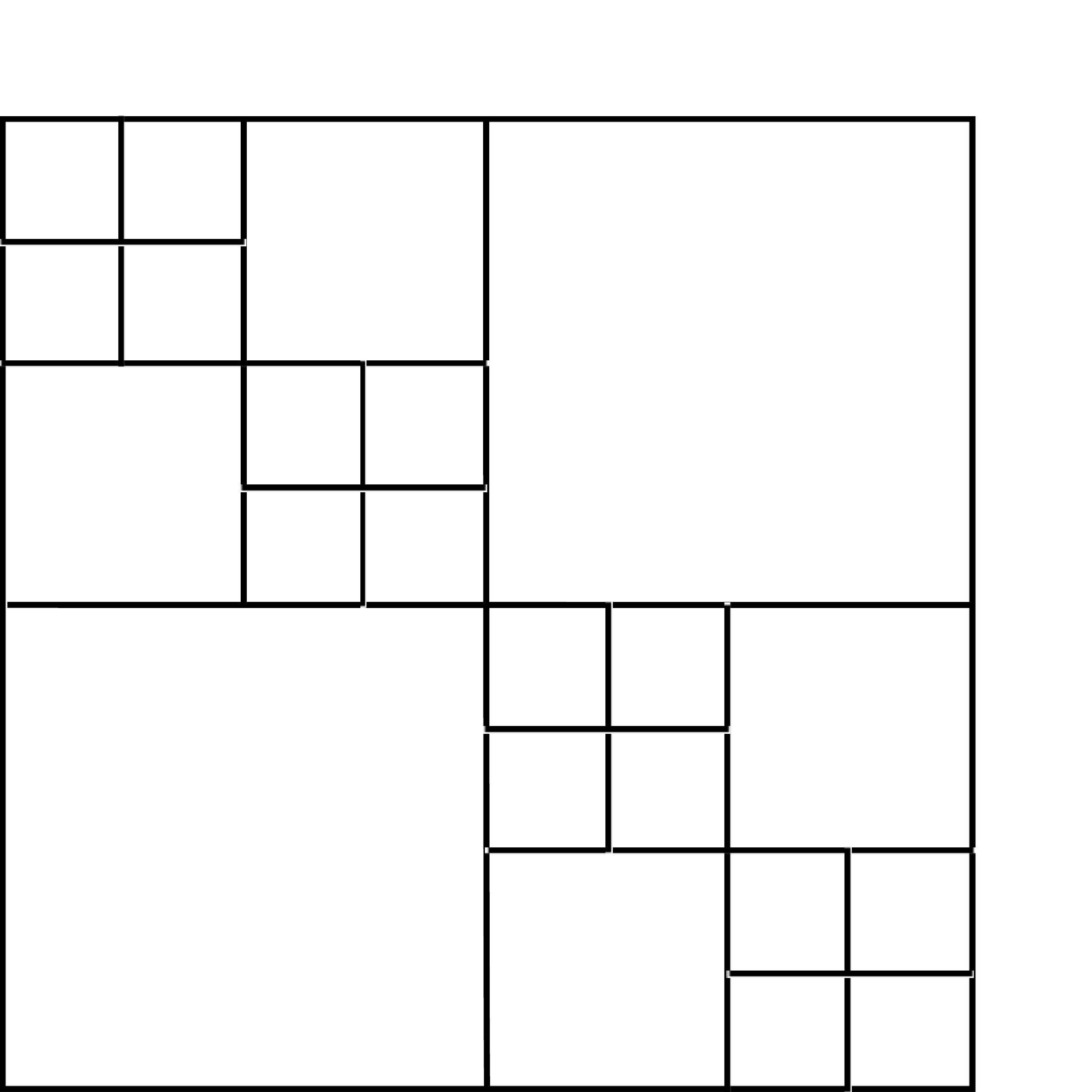}
\caption{Weak hierarchical structure, $\frac{N}{\Rm}=8$.}\label{fig:weak}
\end{minipage}
\begin{minipage}[t]{0.32\linewidth}
\includegraphics[scale=.2]{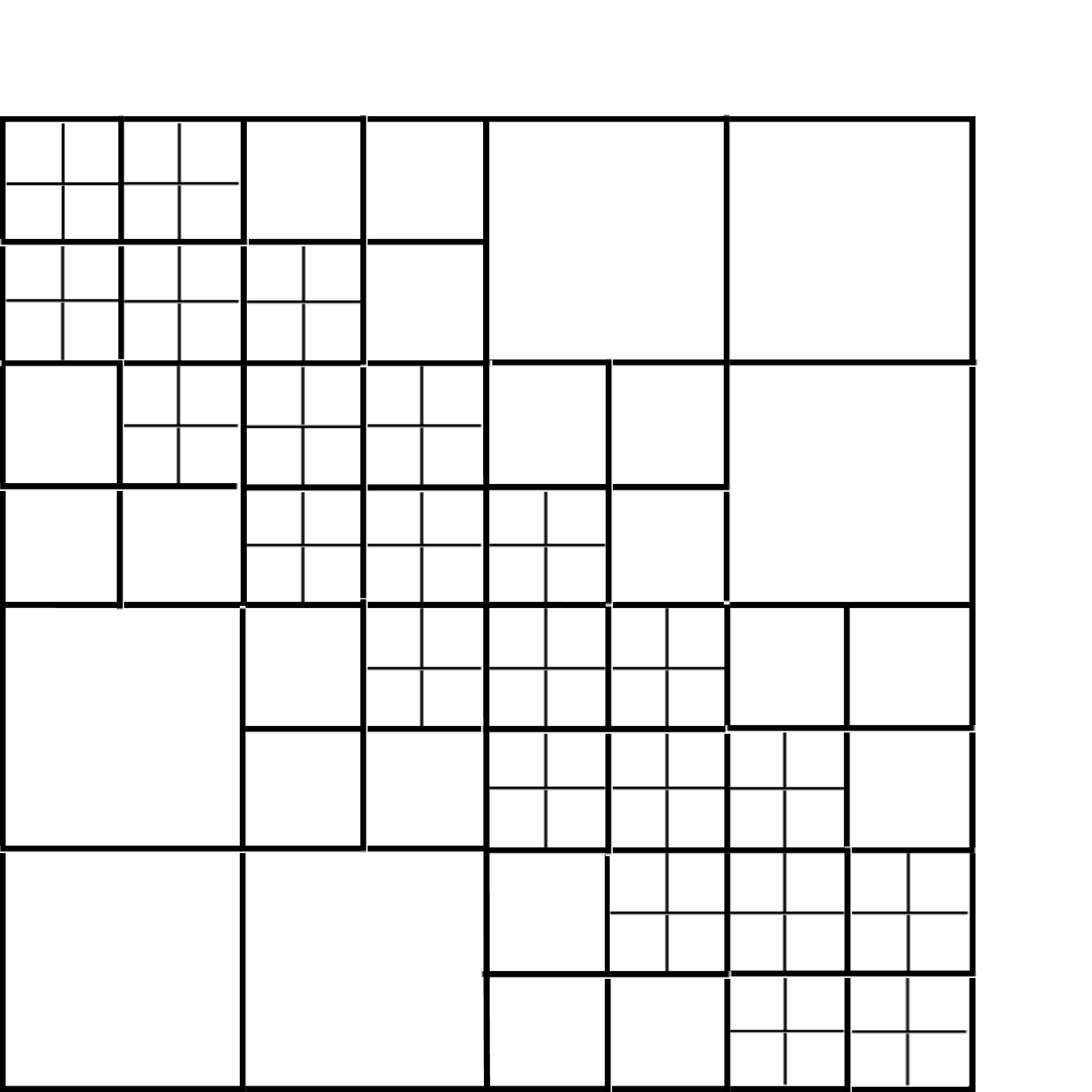}
\caption{Strong hierarchical structure, $\frac{N}{\Rm}=16$.}\label{fig:strong}
\end{minipage}
\begin{minipage}[t]{0.32\linewidth}
\includegraphics[scale=.2]{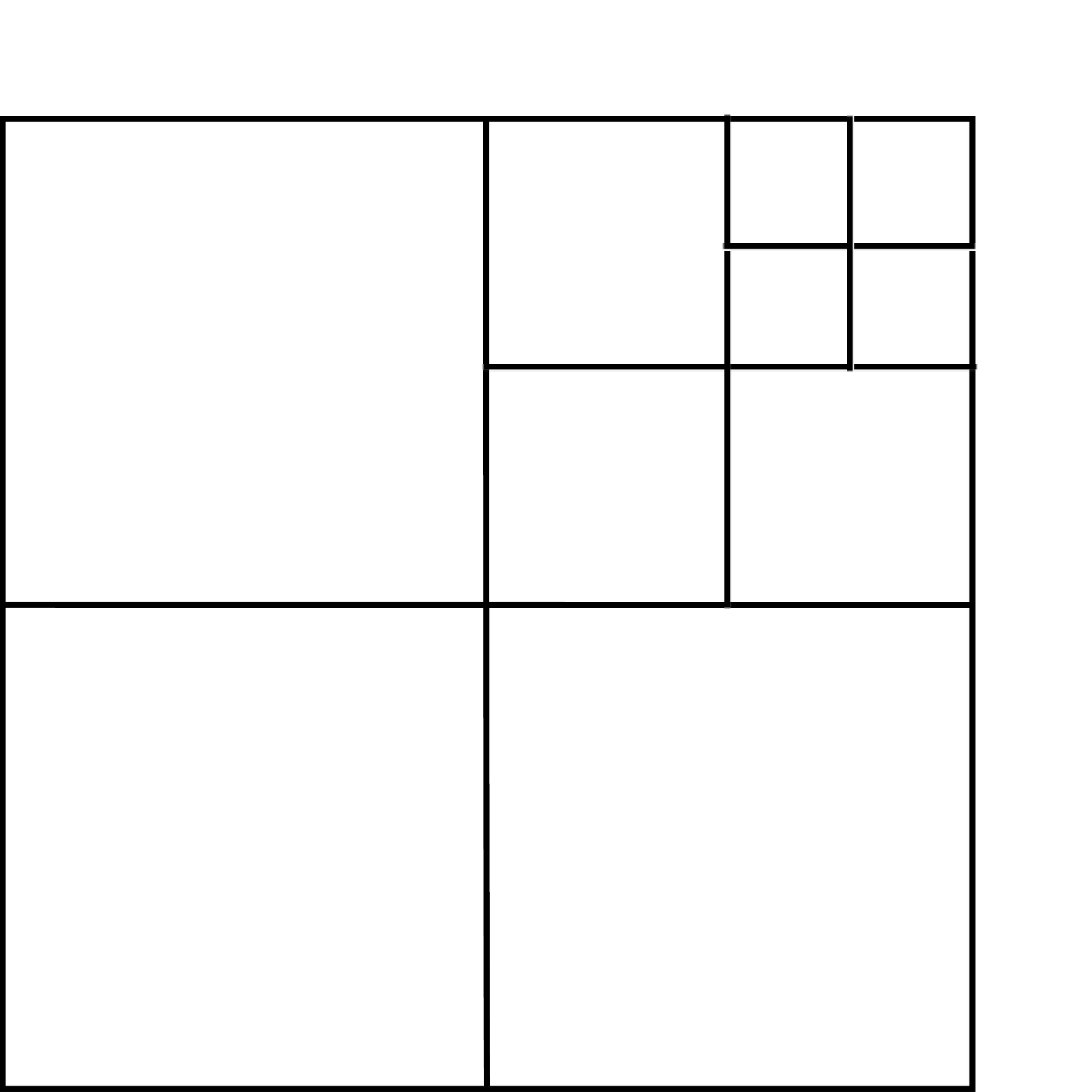}
\caption{Corner PLR structure, $\frac{N}{\Rm}=8$.}\label{fig:corner}
\end{minipage}
\end{figure}

The weak hierarchical structure of a matrix is shown in Figure \ref{fig:weak}. For example, let the matrix $M$ be $8 \times 8$. Then the block at level 0 is $M$ itself. The row indices of $M$ are $\left\{1, 2, \ldots, 8\right\}$, and so are its column indices. Since those overlap, we divide the matrix in four. We are now at level 1, with four blocks. The $(1,1)$ block has row indices $\left\{1, 2, 3, 4\right\}$, and its column indices are the same. This block will have to be divided. The same holds for block $(2,2)$. However, block $(1,2)$ has row indices $\left\{1, 2, 3, 4\right\}$ and column indices $\left\{5, 6, 7, 8\right\}$. Those two sets do not overlap, hence this block is compressed. The same will be true of block $(2,1)$.

We note that, if the matrix $M$ has a weak hierarchical structure, we have a fast matrix-vector product. We may use the heuristic in \ref{eq:comp} to obtain the complexity of that product, assuming for simplicity that all $R_B$'s are $R_\text{max}$. Hence we account for all compressed blocks, starting from the 2 larger blocks on level 1, of size $N/2$ by $N/2$ (one on each side of the diagonal): they correspond to a maximum of $4 N/2 \times R_\text{max}$ operations (multiplications and divisions) each, and there is two of them, so they correspond to a total of $4 N R_\text{max}$ operations. Then, the next larger blocks are of size $N/4$ by $N/4$, and there is 4 of them (two on each side of the diagonal). Hence they correspond to a total of $4 \times 4 \times N/4 \times R_\text{max}=4 N R_\text{max}$ operations. Since we have $L=\log_2{N/R_\text{max}}$ levels, or different block sizes, and as we can see each of those block sizes will contribute at most $4 N R_\text{max}$ operations, we have about $4 N \log{N/R_\text{max}} R_\text{max}$ operations from off-diagonal blocks. We are left with the diagonal blocks. Those have size $R_\text{max}$ and there are $N/R_\text{max}$ of them, so the complexity of multiplying them by a vector is at most $ 4 N/R_\text{max} R^2_\text{max} = 4N R_\text{max}$ operations. Hence the total complexity of a matrix-vector multiplication with a weak hierarchical matrix is
\begin{equation}\label{eq:comweak}
4N R_\text{max} \log{\frac{2N}{R_\text{max}}} .
\end{equation}
This is clearly faster, asymptotically, than the typical complexity of a dense matrix-vector product which is of $2N^2$ operations.

As to the complexity of the compression algorithm, we do a similar calculation, but here the cost per block is $O(N_B \Rm^2)$, so all we need is to sum the dimensions of all blocks that we used the SVD on. We start by taking the SVD of the matrix itself, then of all the blocks on level 1, then half of the blocks on level 2, then a quarter of the blocks on level 3, etc. Hence the complexity of compression is
\[ \Rm^2 \left( N + \sum_{l=1}^{L-1} \frac{N}{2^l} \frac{4^l}{2^{l-1}} \right) =  \Rm^2 (N+2N(L-1)) \leq 2N \Rm^2 \log{\frac{N}{\Rm}}\]
which is nearly linear in $N$.

Finally, we address briefly the error in the matrix $M$ that is made when it is compressed. We especially care about the error in a matrix-vector multiplication $w=Mv$. We can see in this case that, for any entry $j$ in $w$, there will be error coming from all multiplications of the appropriate restriction of $v$ with the corresponding block intersecting row $j$ of $M$. Since there are about $\log\frac{N}{\Rm}$ such blocks in row $j$, we can estimate that by giving a tolerance $\eps$ to the PLR compression algorithm, we will obtain an error in matrix-vector multiplications of about $\eps \log N$. As we will see in section \ref{sec:usingplr}, dividing the ``desired'' error by a factor of 1 to 25 to obtain the necessary $\eps$ will work quite well for our purposes.

\subsubsection{Strong hierarchical matrices}
Next, we define a matrix with a \emph{strong hierarchical structure}. This will be useful for matrices with a singularity on the diagonal.
\begin{definition}\label{def:strong}
A matrix is said to have \emph{strong hierarchical structure} when a block is compressed if and only if its row and column indices are separated by at least the width of the block.
\end{definition}

The strong hierarchical structure of a matrix is shown in Figure \ref{fig:strong}. We can see that, the condition for compression being stronger than in the weak case, more blocks will have to be divided. For example, let the matrix $M$ be $8 \times 8$ again. Then the block at level 0 is $M$ itself, and again its row and column indices overlap, so we divide the matrix in four. We are now at level 1, with four blocks. The $(1,1)$ block will still have to be divided, its row and column indices being equal. The same holds for block $(2,2)$. Now, block $(1,2)$ has row indices $\left\{1, 2, 3, 4\right\}$ and column indices $\left\{5, 6, 7, 8\right\}$. Those two sets do not overlap, but the distance between them, defined as the minimum of $|i-j|$ over all row indices $i$ and column indices $j$ for that block, is 1. Since the width of the block is 4, which is greater than 1, we have to divide the block following Definition \ref{def:strong}. However, at level 2 which has 16 blocks of width 2, we can see that multiple blocks will be compressed: $(1,3), (1,4), (2,4), (3,1), (4,1), (4,2)$.

The matrix-vector multiplication complexity of matrices with a strong hierarchical structure is
\begin{equation}\label{eq:compstrong}
12N R_\text{max} \log{\frac{N}{2 R_\text{max}}}.
\end{equation}
Again this is faster, asymptotically, than the typical $2N^2$ operations of a dense matrix-vector product. We can obtain this number once again by accounting for all the blocks and using \ref{eq:comp}. More precisely, we have $3\sum_{j=1}^{l-1} 2^l=6\frac{1-2^{l-1}}{1-2}=6(2^{l-1}-1)$ compressed blocks at level $l$, hence those blocks have size $N/2^l$. This is true for $l=2, \ldots L-1$, where again $L=\log_2{N/R_\text{max}}$ is the number of levels. Notice that, as expected, we do not have any compressed blocks, or leaves, at level 1. The contribution of those blocks to the matrix-vector complexity will be
\begin{eqnarray*}
4R_\text{max} \sum_{l=2}^{L-1} \left(6(2^{l-1}-1) \ \frac{N}{2^l} \right) &=&12 N R_\text{max} \sum_{l=2}^{L-1} \frac{2^l-2}{2^l}\\
&=& 12 N R_\text{max} (L-2 - \sum_{l=0}^{L-3} \frac{1}{2} \frac{1}{2^l} ) \\
&=& 12 N R_\text{max} (L-2 - (1-1/2^{L-2}) ) \\
&=& 12 N \Rm (L-3+1/2^{L-2}).
\end{eqnarray*}
We need to add to this quantity the complexity coming from the smallest blocks, of size $N/2^L$. There are
\[ 6(2^{L-1}-1) + 2^L + 2(2^L-1)=6\cdot 2^L-8\]
such blocks, and so the corresponding complexity is
\[ 4R_\text{max} (6 \cdot 2^L-8)(N/2^L)=4 \Rm N (6-8/2^L).\]
Adding this to our previous result, we obtain the final complexity of a matrix-vector multiplication:
\begin{eqnarray*}
& & 12 N \Rm (L-3+1/2^{L-2})+ 4 \Rm N (6-8/2^L)\\
&=& 12 N \Rm (L-1+\frac{1}{2^{L-2}}-\frac{2}{3\cdot2^{L-2}}) \\
&\leq& 12 N \Rm L,
\end{eqnarray*}
as stated previously.

For the complexity of the compression algorithm, again we sum the dimensions of all blocks whose SVD we calculated: the matrix itself, the 4 blocks of level 1, the 16 blocks of level 2, 40 blocks in level 3, etc. Hence the complexity of compression is
\begin{eqnarray*}
\Rm^2 \left( N + \frac{N}{2} 4 + \sum_{l=2}^L \frac{N}{2^l} (6\cdot 2^l-8) \right) &=& \Rm^2 \left( N + 2N + N \sum_{l=2}^{L-1} \left(6-\frac{8}{2^l} \right)  \right) \\
&=&  \Rm^2 N \left(3+6(L-2)-8 \ \frac{1}{4}\ \frac{1-1/2^{L-2}}{1-1/2} \right) \\
&=&  \Rm^2 N \left(6L-9-4(1-1/2^{L-2}) \right) \\
&=&  \Rm^2 N \left(6L-13+1/2^{L-2} \right) \\
&\leq & \Rm^2 N \left(6L-12 \right)
\end{eqnarray*}
or
\[6N \Rm^2 \log{\frac{N}{4\Rm}}. \]
This again is nearly linear. Using similar arguments as in the weak case, we can estimate that by giving a tolerance $\eps$ to the PLR compression algorithm, we will obtain an error in matrix-vector multiplications of about $\eps \log N$ again.

\subsubsection{Corner PLR matrices}
One final structure we wish to define, now useful for matrices with a singularity in a corner, is the following:
\begin{definition}\label{def:corner}
A matrix is said to have \emph{corner PLR structure}, with reference to a specific corner of the matrix, when a block is divided if and only if both its row and column indices contain the row and column indices of the entry corresponding to that specific corner.
\end{definition}

Figure \ref{fig:corner} shows a top-right corner PLR structure. Again, we take an $8 \times 8$ matrix $M$ as an example. The top-right entry has row index 1 and column index 8. We see that the level 0 block, $M$ itself, certainly contains the indices $(1,8)$, so we divide it. On the next level, we have four blocks. Block $(1,2)$ is the only one that has both row indices that contain the index 1, and column indices that contain the index 8, so this is the only one that is divided. Again, on level 2, we have 16 blocks of size 2, and block $(1,4)$ is the only one divided.

As for the corner PLR matrices, their matrix-vector multiplication complexity is:
\begin{equation}\label{eq:compcorn}
8N R_\text{max}.
\end{equation}
Indeed, we see we have 3 blocks of size $N/2^L$ for $l=1,2, \ldots, L-1$. This is a constant number of blocks per level, which means that matrix-vector multiplication will be even faster. We also have 4 blocks at the lowest level, of size $N/2^L$. The complexity is then
\begin{eqnarray*}
4 \Rm 3 \sum_{l=1}^{L-1} N/2^l \ +4\Rm 4 N/2^L &=& 4 \Rm N (3(1-1/2^{L-1}) +2/2^{L-1})\\
&=& 8 N \Rm.
\end{eqnarray*}

For the complexity of the compression algorithm, we sum the dimensions of all blocks whose SVD we calculated: the matrix itself, the 4 blocks of level 1, 4 blocks in level 2, 4 blocks in level 3, etc. Hence the complexity of compression is
\begin{eqnarray*}
\Rm^2 \left( N + 4 \sum_{l=1}^L \frac{N}{2^l} \right) &=& \Rm^2 N \left( 1+4 \ \frac{1}{2} \ \frac{1-1/2^L}{1-1/2} \right) \\
&=&  \Rm^2 N \left(1+4-4/2^L \right) \\
&\leq & \Rm^2 N \left(4 \right)
\end{eqnarray*}
or
\[4N \Rm^2 . \]
Hence the complexity of compression for corner PLR matrices is linear. And again, we estimate that by giving a tolerance $\eps$ to the PLR compression algorithm, we will obtain an error in matrix-vector multiplications of about $\eps \log N$.
\newline
\newline
Now that we have explained these three special structures, and how they provide a fast matrix-vector product, we are ready to discuss using PLR matrices specifically for the exterior DtN map.

\section{Using PLR matrices for the DtN map's submatrices}\label{sec:usingplr}

As we recall, obtaining the full DtN map from solving the exterior problem $4N$ times is too costly, and so we use matrix probing to approximate the DtN map $D$ by $\tilde{D}$ using only a few exterior solves. If we were to try to use PLR matrices directly on $D$, we would have to find the SVD of many blocks. Since we do not have access to the blocks themselves, we would need to use the randomized SVD, and hence to solve the exterior problem, on random vectors restricted to the block at hand. As we mentioned before, Lin et al. have done something similar in \cite{Hmatvec}, which required $O(\log{N}$ matrix-vector multiplies with a large constant, or in our case exterior solves. This is too costly, and this is why we use matrix probing first to obtain an approximate $\tilde{D}$ with a small, nearly constant number of exterior solves.

Now that we have access to $\tilde{D}$ from matrix probing, we can approximate it using PLR matrices. Indeed, we have access to the full matrix $\tilde{D}$, and so finding the SVD of a block is not a problem. In fact, we use the randomized SVD for speed, not because we only have access to matrix-vector multiplies.

Compressing one of those edge-to-edge submatrices under the PLR matrix framework requires that we pick both a tolerance $\epsilon$ and a maximal desired rank $R_\text{max}$. We explain in the next subsections how to choose appropriate values for those parameters.

\subsection{Choosing the tolerance}

Because our submatrices come from probing, they already have some error attached to them, that is, the relative probing error as defined in equation \ref{acterr} of chapter \ref{ch:probing}. Therefore, it would be wasteful to ask for the PLR approximation to do any better than that probing error.

Also, when we compress blocks in the PLR compression algorithm, we make an absolute error in the $L_2$ norm. However, because of the high norm of the DtN map, it makes more sense to consider the relative error. We can thus multiply the relative probing error we made on each submatrix $\tilde{M}$ by the norm of the DtN map $D$ to know the absolute error we need to ask of the PLR compression algorithm. And since the $L_2$ norm is smaller than the Frobenius norm, and errors from each block compound, we have found empirically that asking for a tolerance $\epsilon$ which is a factor of $1$ to $1/100$ of the absolute probing error of a submatrix works for obtaining a similar Frobenius error from the PLR approximation. As a rule of thumb, this factor needs to be smaller for diagonal submatrices $M$ of $D$, but can be equal to 1 for submatrices corresponding to opposite edges of $\pd \Omega$. 

Of course, we do not want to use an $\eps$ which is too smaller either. That might force the PLR compression algorithm to divide blocks more than needed, and make the matrix-vector multiplications slower than needed.

\subsection{Minimizing the matrix-vector application time}

Our main objective in this chapter is to obtain a fast algorithm. To this end, we try to compress probed submatrices of the DtN map using various values of the parameter $R_\text{max}$, and choose the value that will give us the fastest matrix-vector multiplies. We use the known complexity of a matrix-vector multiplication \eqref{eq:comp} to find the rank $R_\text{max}$ that minimizes the complexity, from doing a few tests, and we use the compressed submatrix corresponding to that particular maximal rank in our Helmholtz solver. A different choice of complexity might be used depending on the operating system and coding language used, since slow downs might occur because of cache size, operations related to memory, matrix and vector operations, etc. 

However, we note that we may compare the actual complexity from the particular structure obtained by PLR compression to the ``ideal'' complexities coming from the special structures we have mentioned before. Indeed, for a submatrix on the diagonal, we can compare its matrix-vector complexity to that of weak and strong hierarchical matrices. That will give us an idea of whether we have a fast algorithm. One thing we notice in most of our experiments is that, for diagonal blocks, the actual complexity usually becomes smaller as $\Rm$ increases, until we arrive at a minimum with the $\Rm$ that gives us the better compromise between too many and blocks of too high a rank. Then, the complexity increases again. However, the complexity increases slower than that of both weak and strong matrices. 

\begin{figure}[ht]
\begin{minipage}[t]{0.45\linewidth}
\includegraphics[scale=.45]{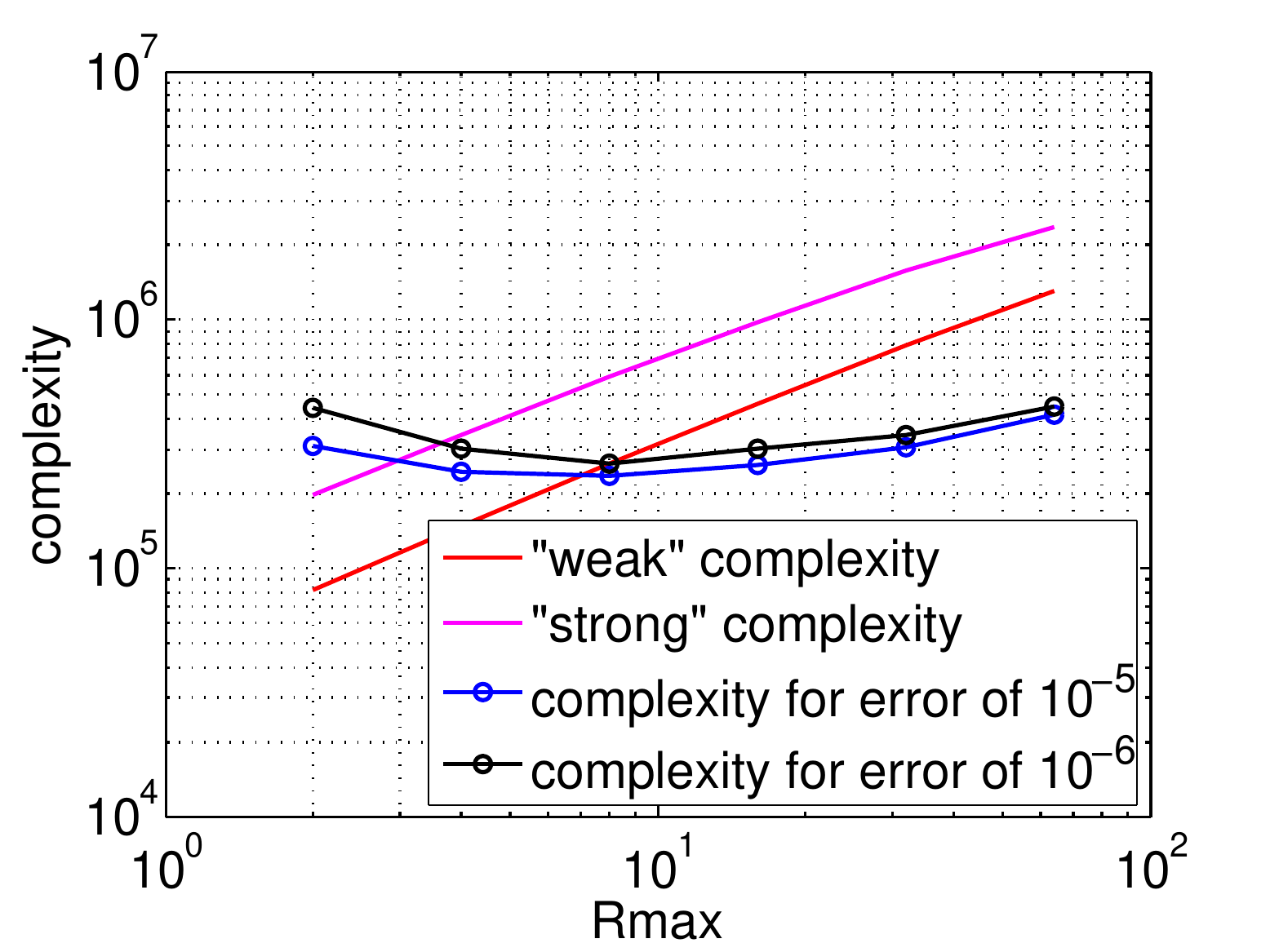}
\caption{Matrix-vector complexity for submatrix $(1,1)$ for $c \equiv 1$, various $\Rm$. Probing errors of $10^{-5}$, $10^{-6}$.}
\label{fig:compscomp1}
\end{minipage}
\hspace{1cm}
\begin{minipage}[t]{0.45\linewidth}
\includegraphics[scale=.45]{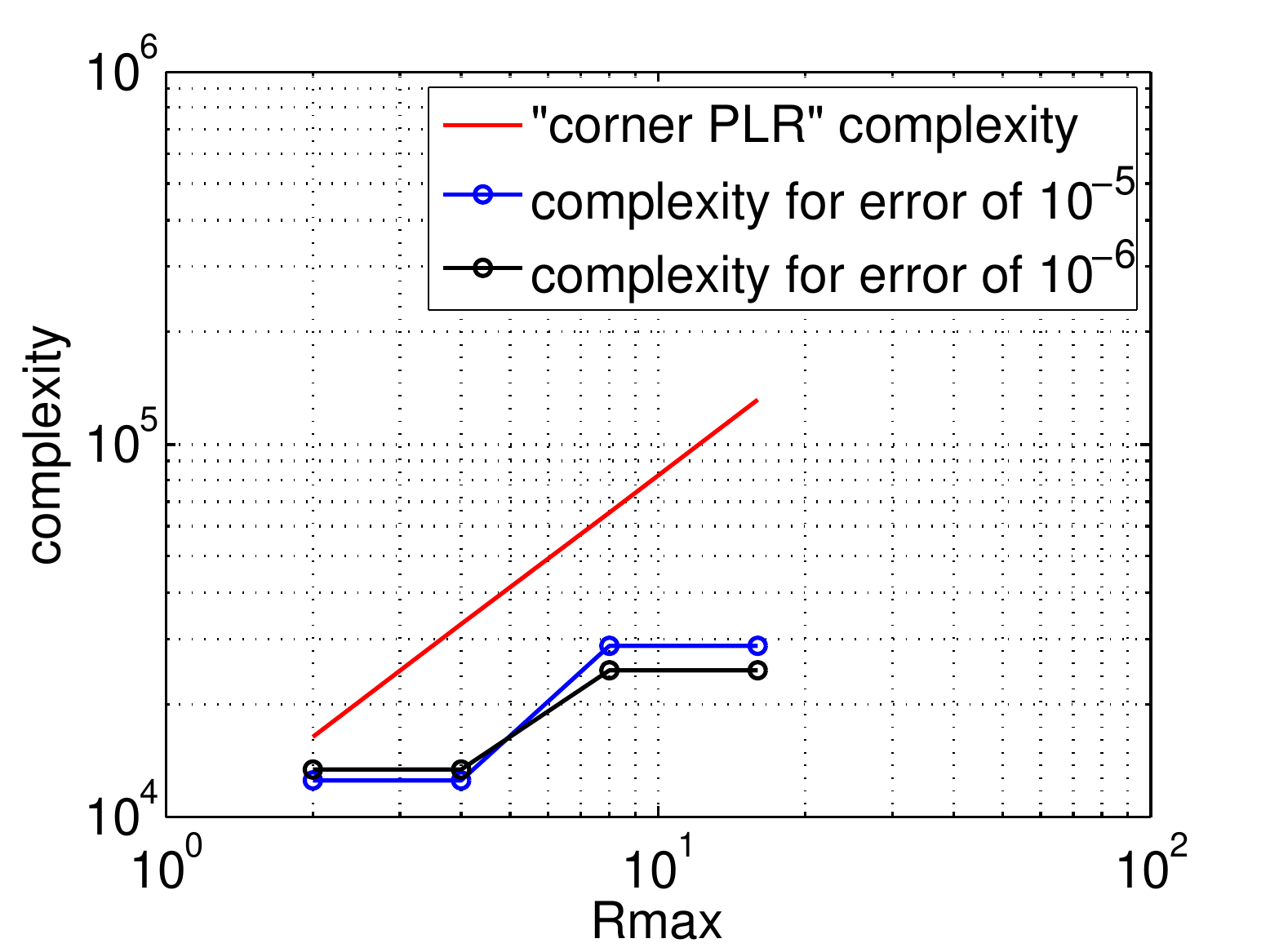}
\caption{Matrix-vector complexity for submatrix $(2,1)$ for $c \equiv 1$, various $\Rm$. Probing errors of $10^{-5}$, $10^{-6}$.}
\label{fig:compscomp2}
\end{minipage}
\end{figure}

Figure \ref{fig:compscomp1} confirms this phenomenon for the $(1,1)$ block of the constant medium. From this figure, we would then pick $\Rm=8$ since this is the value of $\Rm$ that corresponds to the smallest actual complexity of a matrix-vector product, both for a relative probing error of $10^{-5}$ and $10^{-6}$. Figure \ref{fig:compscomp2} confirms this phenomenon as well for the $(2,1)$ block of the constant medium. From this figure, we would then pick $\Rm=4$ again both for a relative probing error of $10^{-5}$ and $10^{-6}$, since this is the value of $\Rm$ that corresponds to the smallest actual complexity of a matrix-vector product (it is hard to tell from the figure, but the complexity for $\Rm=2$ is just larger than that for $\Rm=4$ in both cases).

\begin{figure}[ht]
\begin{minipage}[t]{0.45\linewidth}
\includegraphics[scale=.45]{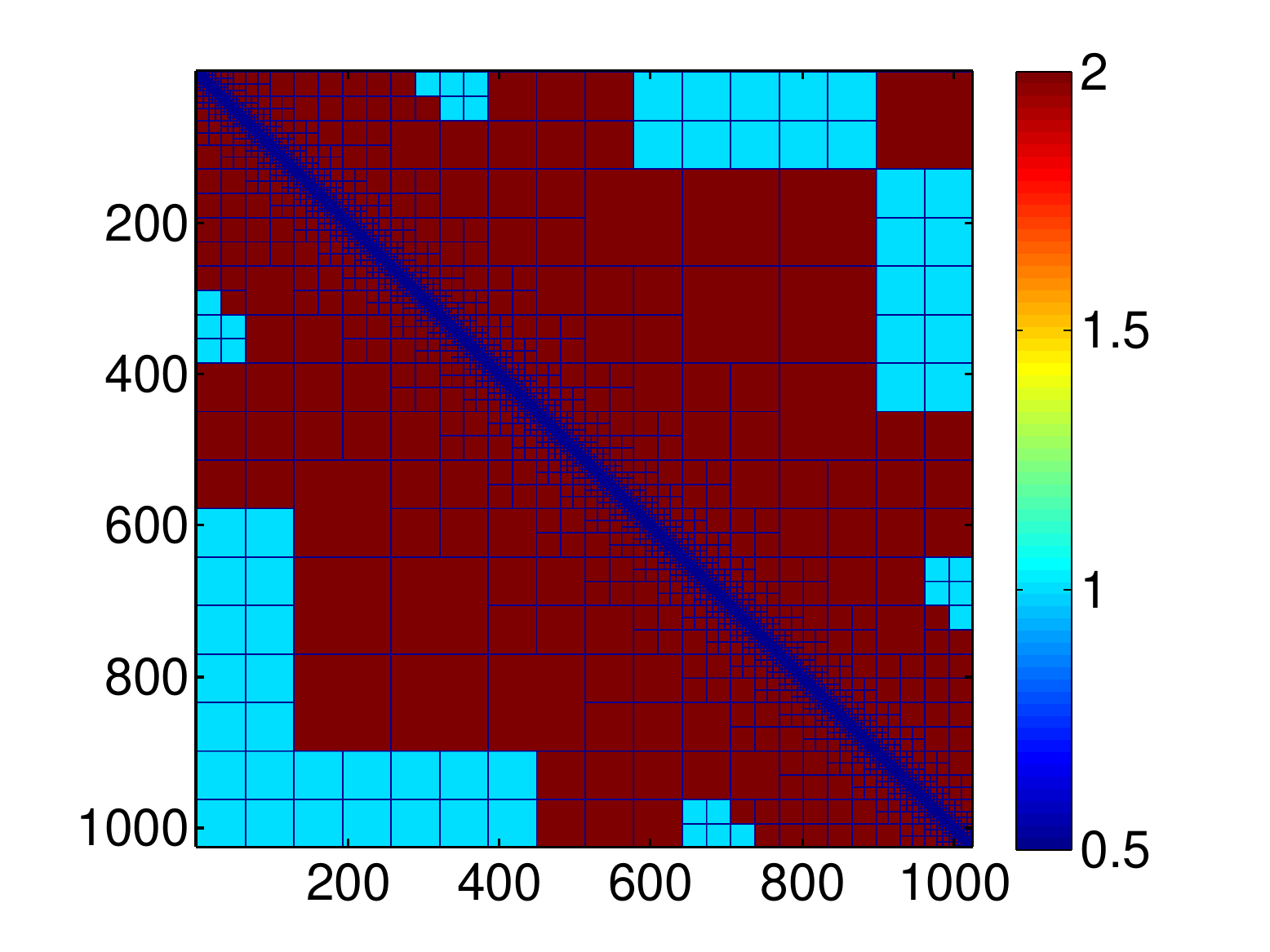}
\caption{PLR structure of the probed submatrix $(1,1)$ for $c \equiv 1$, $R_\text{max}=2$. Each block is colored by its numerical rank.}
\label{fig:c1d1b1}
\end{minipage}
\hspace{1cm}
\begin{minipage}[t]{0.45\linewidth}
\includegraphics[scale=.45]{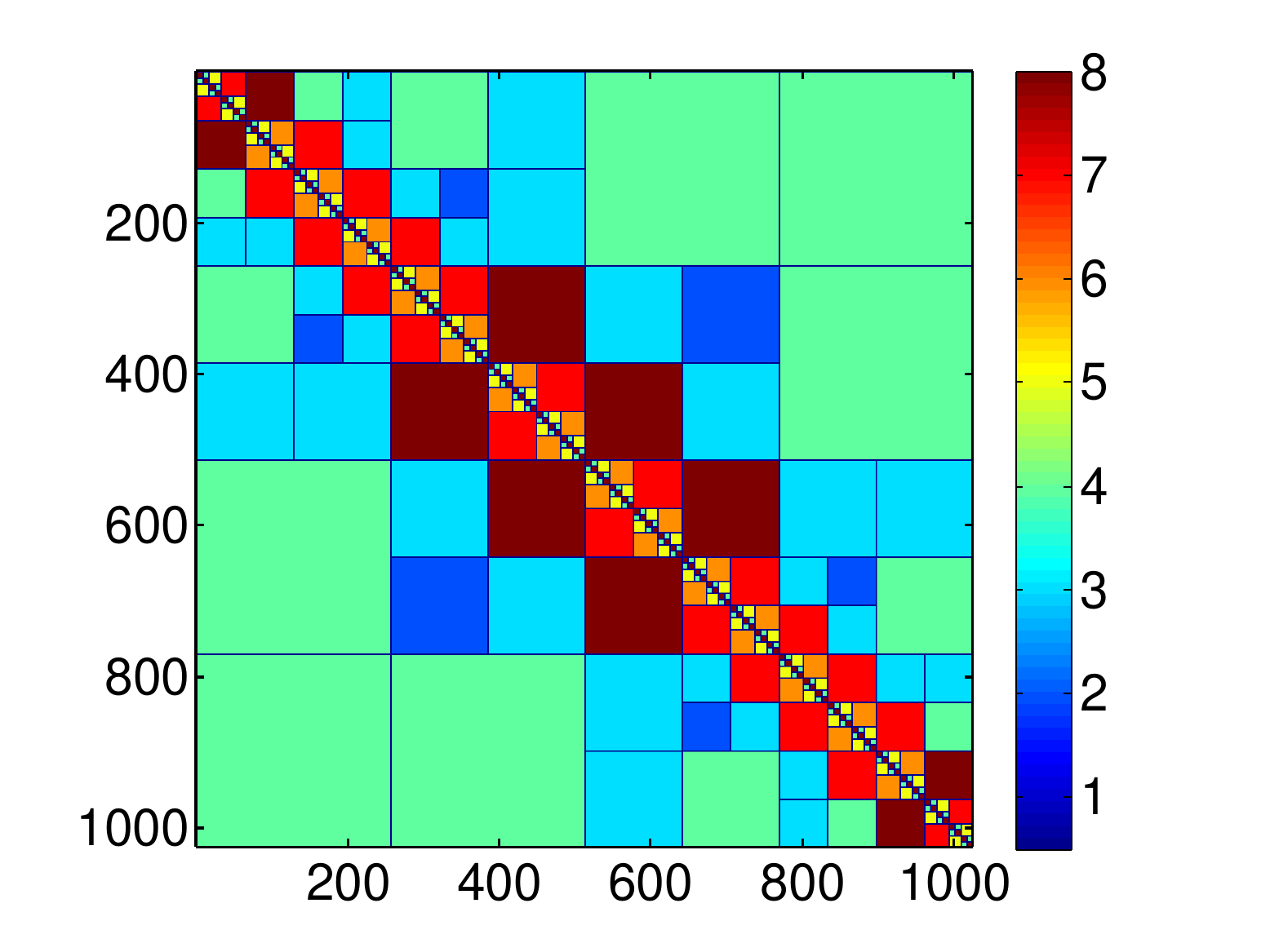}
\caption{PLR structure of the probed submatrix $(1,1)$ for $c \equiv 1$, $R_\text{max}=8$. Each block is colored by its numerical rank.}
\label{fig:c1d3b1}
\end{minipage}
\end{figure}
\begin{figure}[ht]
\begin{minipage}[t]{0.45\linewidth}
\includegraphics[scale=.45]{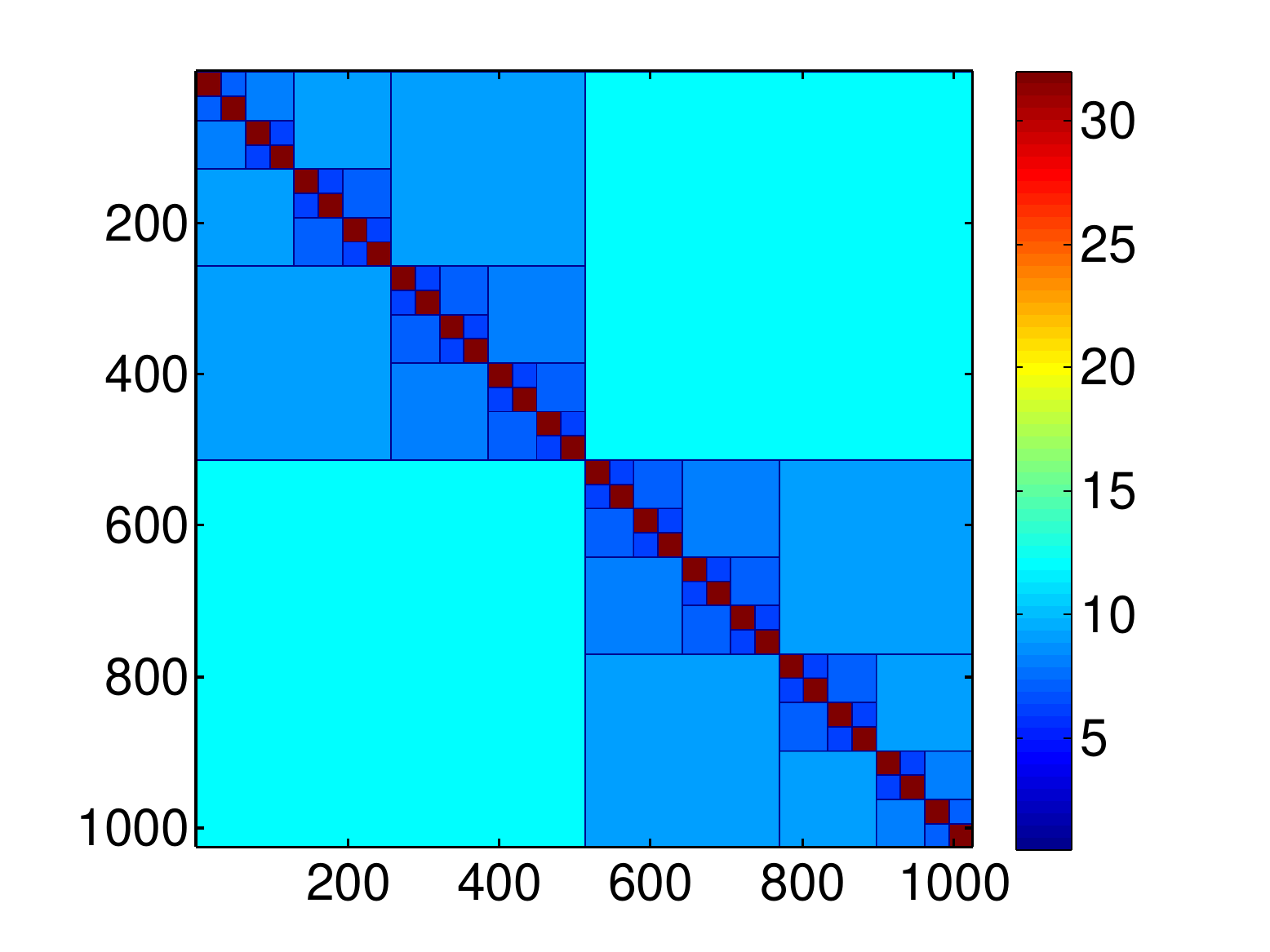}
\caption{PLR structure of the probed submatrix $(1,1)$ for $c \equiv 1$, $R_\text{max}=32$. Each block is colored by its numerical rank.}
\label{fig:c1d5b1}
\end{minipage}
\hspace{1cm}
\begin{minipage}[t]{0.45\linewidth}
\includegraphics[scale=.45]{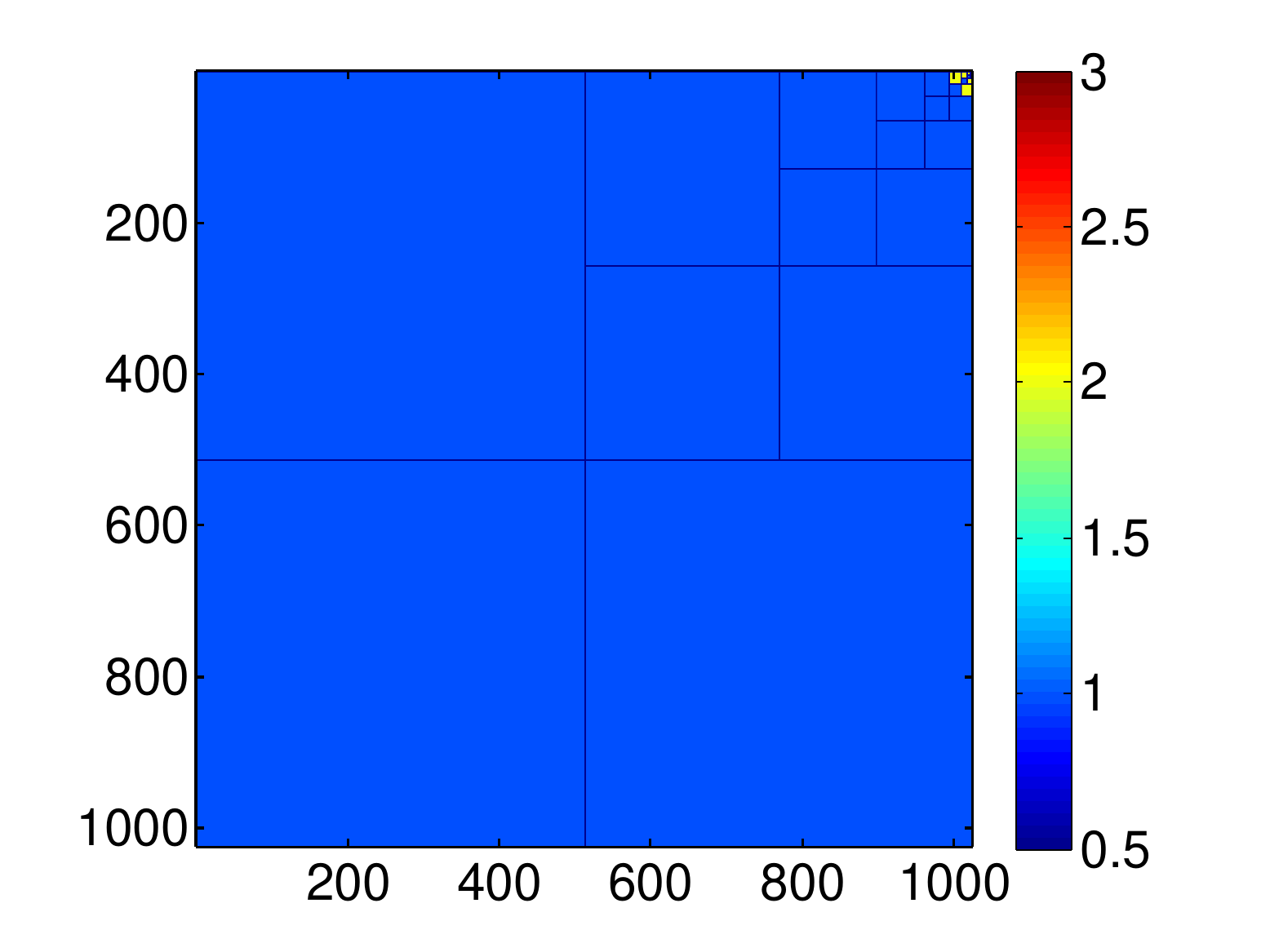}
\caption{PLR structure of the probed submatrix $(2,1)$ for $c \equiv 1$, $R_\text{max}=4$. Each block is colored by its numerical rank.}
\label{fig:c1d2b2}
\end{minipage}
\end{figure}

It is also informative to have a look at the structures we obtain, and the ranks of each block. Figures \ref{fig:c1d1b1}, \ref{fig:c1d3b1} and \ref{fig:c1d5b1} refer to block $(1,1)$ again of $c \equiv 1$, for different values of $\Rm$. As we expect, having $\Rm=2$ in Figure \ref{fig:c1d1b1} in this case forces blocks to be very small, which is wasteful. On the other hand, a larger $\Rm=32$ in Figure \ref{fig:c1d5b1} is not better because then we have fewer blocks, but they have large rank. Still, because the structure is that of a weak hierarchical matrix, with blocks that have in fact rank smaller than $\Rm$, we obtain a fast matrix-vector product. However, the ideal is $\Rm=8$ in Figure \ref{fig:c1d3b1}, which minimizes the complexity of a matrix-vector product by finding the correct balance between fewer blocks but small rank. We see it almost has a strong hierarchical structure, but in fact with more large blocks and fewer small blocks. As for the $(2,1)$ block, we see its PLR structure in Figure \ref{fig:c1d2b2}: it is actually a corner PLR structure but the numerical ranks of blocks are always lower than $\Rm=4$ so the matrix-vector multiplication of that submatrix will be even faster than for a corner PLR structure, as we knew from Figure \ref{fig:compscomp2}.

\subsection{Numerical results}

\begin{table}[ht]
\caption{PLR compression results, $c\equiv 1$} 
\begin{center} \footnotesize
\begin{tabular}{|l|l|l|l|l|} \hline  
$R_\text{max}$ for $(1,1)$	& $R_\text{max}$ for $(2,1)$	& $\|D-\overline{D}\|_F/\|D\|_F$ & $\|u-\overline{u}\|_F/\|u\|_F$ & Speed-up\\ \hline 
{$2$} & {$2$} & {$4.2126e-01$} &{$6.5938e-01$} & $115$ \\ \hline 
{$2$} & {$2$} & {$4.2004e-02$} &{$7.3655e-02$} & $93$ \\ \hline 
{$2$} & {$2$} & {$1.2517e-03$} &{$2.4232e-03$}  & $55$ \\ \hline 
{$4$} & {$2$} & {$1.1210e-04$} &{$4.0003e-04$}  & $42$ \\ \hline 
{$8$} & {$4$} & {$1.0794e-05$} &{$1.4305e-05$}  & $32$ \\ \hline 
{$8$} & {$4$} & {$6.5496e-07$} &{$2.1741e-06$}  & $29$ \\ \hline 
\end{tabular}
\end{center} 
\label{c1solveplr} 
\end{table}

\begin{table}[ht]
\caption{PLR compression results, $c$ is the Gaussian waveguide.} 
\begin{center} \footnotesize
\begin{tabular}{|l|l|l|l|l|l|} \hline  
$R_\text{max}$ for $(1,1)$	& $R_\text{max}$ for $(2,2)$	& $R_\text{max}$ for $(2,1)$	& $\|D-\overline{D}\|_F/\|D\|_F$ & $\|u-\overline{u}\|_F/\|u\|_F$  & Speed-up \\ \hline 
{$2$} & {$2$} & {$2$} & {$6.6034e-02$} &{$1.4449e-01$} & $105$ \\ \hline 
{$2$} & {$2$} & {$2$} & {$1.8292e-02$} &{$7.4342e-02$} & $74$ \\ \hline 
{$2$} & {$2$} & {$2$} & {$2.0948e-03$} &{$1.1014e-02$} & $59$ \\ \hline 
{$4$} & {$4$} & {$2$} & {$2.3740e-04$} &{$1.6023e-03$} & $47$ \\ \hline 
{$8$} & {$4$} & {$4$} & {$1.5369e-05$} &{$8.4841e-05$} & $36$ \\ \hline 
{$8$} & {$8$} & {$4$} & {$3.4148e-06$} &{$1.7788e-05$} & $30$ \\ \hline 
\end{tabular}
\end{center} 
\label{c3solveplr} 
\end{table}

\begin{table}[ht]
\caption{PLR compression results, $c$ is the Gaussian slow disk.} 
\begin{center} \footnotesize
\begin{tabular}{|l|l|l|l|l|} \hline  
$R_\text{max}$ for $(1,1)$	& $R_\text{max}$ for $(2,1)$	& $\|D-\overline{D}\|_F/\|D\|_F$ & $\|u-\overline{u}\|_F/\|u\|_F$ & Speed-up\\ \hline 
{$2$} & {$2$} & {$9.2307e-02$} &{$1.2296e+00$} & $97$ \\ \hline 
{$2$} & {$2$} & {$8.1442e-03$} &{$4.7922e-02$} & $69$ \\ \hline 
{$4$} & {$2$} & {$1.2981e-03$} &{$3.3540e-02$} & $44$ \\ \hline 
{$4$} & {$2$} & {$1.1680e-04$} &{$1.0879e-03$} & $39$ \\ \hline 
{$4$} & {$2$} & {$2.5651e-05$} &{$1.4303e-04$} & $37$ \\ \hline 
\end{tabular}
\end{center} 
\label{c5solveplr} 
\end{table}

We have compressed probed DtN maps and used them in a Helmholtz solver with success. We have used the same probed matrices as in the previous chapter, and so we refer the reader to Tables \ref{FDPMLerr}, \ref{c1solve}, \ref{c3solve}, \ref{c5solve}, \ref{c16solve}, \ref{c18solve}, \ref{c33solve} for all the parameters we used then. 

\begin{table}[ht]
\caption{PLR compression results, $c$ is the vertical fault, sources on the left and on the right.} 
\begin{center} \footnotesize
\begin{tabular}{|l|l|l|l|l|l|} \hline  
$R_\text{max}$ for $(1,1)$	& $R_\text{max}$ for $(2,2)$	& $\frac{\|D-\overline{D}\|_F}{\|D\|_F}$ & $\frac{\|u-\overline{u}\|_F}{\|u\|_F}$, left  & $\frac{\|u-\overline{u}\|_F}{\|u\|_F}$, right  & Speed-up\\ \hline
{$2$} & {$2$} & {$2.6972e-01$} &{$5.8907e-01$} &{$4.6217e-01$} & $105$\\ \hline 
{$2$} & {$2$} & {$9.0861e-03$} &{$3.9888e-02$} &{$2.5051e-02$} & $67$ \\ \hline 
{$1$} & {$4$} & {$8.7171e-04$} &{$3.4377e-03$} &{$2.4279e-03$} & $53$\\ \hline 
\end{tabular}
\end{center} 
\label{c16solveplrl} 
\end{table}
\begin{table}[ht]
\caption{PLR compression results, $c$ is the diagonal fault.} 
\begin{center} \footnotesize
\begin{tabular}{|l|l|l|l|l|} \hline  
$R_\text{max}$ for $(1,1)$	& $R_\text{max}$ for $(2,2)$	& $\|D-\overline{D}\|_F/\|D\|_F$ & $\|u-\overline{u}\|_F/\|u\|_F$ & Speed-up\\ \hline 
{$2$} & {$2$} & {$1.4281e-01$} &{$5.3553e-01$} & $98$ \\ \hline 
{$2$} & {$2$} & {$1.9108e-02$} &{$7.8969e-02$} & $76$ \\ \hline 
{$2$} & {$4$} & {$2.5602e-03$} &{$8.7235e-03$} & $49$ \\ \hline 
\end{tabular}
\end{center} 
\label{c18solveplr} 
\end{table}
\begin{table}[ht]
\caption{PLR compression results, $c$ is the periodic medium.} 
\begin{center} \footnotesize
\begin{tabular}{|l|l|l|l|l|} \hline  
$R_\text{max}$ for $(1,1)$	& $R_\text{max}$ for $(2,1)$	& $\|D-\overline{D}\|_F/\|D\|_F$ & $\|u-\overline{u}\|_F/\|u\|_F$ & Speed-up\\ \hline 
{$2$} & {$2$} & {$1.2967e-01$} &{$2.1162e-01$} & $32$ \\ \hline 
{$2$} & {$2$} & {$3.0606e-02$} &{$5.9562e-02$} & $22$ \\ \hline 
{$8$} & {$2$} & {$9.0682e-03$} &{$2.6485e-02$} & $11$ \\ \hline 
\end{tabular}
\end{center} 
\label{c33solveplr} 
\end{table}

We now present results for PLR compression in a Helmholtz solver in Tables \ref{c1solveplr}, \ref{c3solveplr}, \ref{c5solveplr}, \ref{c16solveplrl}, \ref{c18solveplr}, \ref{c33solveplr}. For each medium, we show the chosen $\Rm$ of the most important (in norm) submatrices. For all other submatrices, $\Rm \leq 2$. We then show the relative norm of the error between the PLR compression $\overline{D}$ and the actual DtN map $D$. We also show the relative error between the solution $\overline{u}$ computed using $\overline{D}$ in the Helmholtz solver and the actual solution $u$ using $D$ as the DtN map. Finally, we show the ``speed-up'' obtained from taking the ratio of the complexity of using a dense matrix-vector product for $\tilde{D}$, which would be of about $2\times 16N^2$, to the total complexity of a matrix-vector product of $\bar{D}$. This ratio tells us how much faster than a dense product the PLR compression is. We see that this ratio ranges from about 50 to 100 for all media expect the periodic media, with a smaller ratio associated with asking for a higher accuracy, as expected. The speed-up ratio is between 10 and 30 for the periodic media, but as we recall the value of $N$ here is smaller: $N=320$. Larger values of $N$ should lead to a better speed-up.

\section{The half-space DtN map is separable and low rank: theorem}\label{sec:septheo}

As we have mentioned before, the Green's function for the half-space Helmholtz equation is separable and low rank \cite{Hsweep}. We investigate here the half-space DtN map kernel, which is related to the Green's function through two derivatives, as we saw in section \ref{sec:hsG}, and we obtain a similar result to that of \cite{Hsweep}. We state the result here, and prove it in the next subsections. We then end this section with a discussion on generalizing our theorem for heterogeneous media.

Let $\x=(x,0)$ and $\y=(y,0)$ be points along the half-space boundary, $x\neq y$. Recall the Dirichlet-to-Neumann map kernel for the half-space Helmholtz equation \eqref{eq:hsHE} with homogeneous medium $c \equiv 1$ and $\omega=k/c$ is:
\begin{equation}\label{eq:hsDtN}
K(|x-y|)= \frac{ik^2}{2} \frac{H_1^{(1)}(k|x-y|)}{k|x-y|}.
\end{equation}

\begin{theorem}\label{theo:sep}
Let $0 <\epsilon \leq 1/2$, and $0<r_0<1$, $r_0=\Theta(1/k)$. There exists an integer $J$, functions $\left\{\Phi_j,\chi_j\right\}_{j=1}^J$ and a number $C$ such that we can approximate $K(|x-y|)$ for $r_0\leq |x-y| \leq 1$ with a short sum of smooth separated functions:
\begin{equation}\label{eq:sep}
K(|x-y|)=\sum_{j=1}^J \Phi_j(x)\chi_j(y) + E(x,y)
\end{equation}
where $|E(x,y)| \leq \epsilon$, and $J \leq C \left(\log k \max(|\log\epsilon|,\log k\right))^2$ with $C$ which does not depend on $k$, or $\epsilon$. $C$ does depend weakly on $r_0$ through the constant quantity $r_0k$, for more details see Remark \ref{rem:C_0}.
\end{theorem}

To prove this, we shall first consider the Hankel function $H_1^{(1)}$ in \eqref{eq:hsDtN}, and see that it is separable and low rank. Then, we shall look at the factor of $1/k|x-y|$ in \eqref{eq:hsDtN}, and see the need for using quadratures on a dyadic partition of the interval $\left[r_0,1\right]$ in order to prove that this factor is also separable and low rank. Finally, we prove Theorem \ref{theo:sep} and make a few remarks.

\subsection{Treating the Hankel function}

From Lemmas 7 and 8 of \cite{mr} we know that $H_1^{(1)}(k|x-y|)$ is separable and low-rank as a function of $x$ and $y$. Looking in particular at Lemma 7 from \cite{mr}, we make a slight modification of the proof to obtain the following lemma.
\begin{lemma}\label{lem:7n8}
Slight modification of Lemmas 7 and 8, \cite{mr}. Let $0<\epsilon \leq 1/2$, $k>0$ and $r_0>0$, $r_0=\Theta(1/k)$. Let $|x-y|>r_0$. Then there exists an integer $J_1$, a number $C_1$, and functions $\left\{\Phi^{(1)}_j,\chi^{(1)}_j\right\}_{j=1}^{J_1}$ such that
\begin{equation}\label{eq:h1}
H_1^{(1)}(k|x-y|)=\sum_{j=1}^{J_1} \Phi^{(1)}_j(x)\chi^{(1)}_j(y) + E^{(1)}_{J_1}(x,y),
\end{equation}
where $|E^{(1)}_{J_1}(x,y)| \leq \epsilon$, and $J_1 \leq C_1 \log k |\log\epsilon|$ with $C_1$ which does not depend on $k$, or $\epsilon$. $C_1$ does depend weakly on $r_0$ through the quantity $r_0k$. Again, see Remark \ref{rem:C_0}.
\end{lemma}

\subsection{Treating the $1/kr$ factor}

We now show that the $1/kr$ factor is also separable and low rank\footnote{A different technique than the one presented here would be to expand $1/|x-y|$ using the Taylor expansion for $x>y>0$: $$\frac{1}{x} \frac{1}{1-y/x} = \frac{1}{x} \left(1+\frac{y}{x}+\left(\frac{y}{x}\right)^2 + \ldots \right).$$ However, making an error of $\eps$ requires $k\log\eps$ terms in the expansion because the error is large when $y/x \approx 1$ or $y\approx x$ or $r\approx r_0$.}. Notice:
\begin{equation}\label{eq:1okr}
\int_0^\infty e^{-krt} dt = \left. \frac{e^{-krt}}{-kr} \right|_0^\infty = 0 + \frac{e^{-kr0}}{kr} = \frac{1}{kr}
\end{equation}
and 
\begin{equation}\label{eq:split}
\int_0^\infty e^{-krt} dt = \int_0^{T} e^{-krt} dt + \int_{T}^\infty e^{-krt} dt,
\end{equation}
where 
\[ \int_{T}^\infty e^{-krt} dt = \frac{e^{-krT}}{kr}. \]
Equation \eqref{eq:1okr} means we can write $1/kr$ as an integral, and equation \eqref{eq:split} means we can split this integral in two, one part being a definite integral, the other indefinite. But we can choose $T$ so that the indefinite part is smaller than our error tolerance:
\[  \left| \int_{T}^\infty e^{-krt} dt \right| \leq \epsilon . \]
For this, we consider that $\frac{e^{-krT}}{kr} \leq \frac{e^{-krT}}{C_0} \leq \epsilon $ and so we need $krT \geq |\log C_0| + |\log \epsilon |$ or $T \geq (|\log C_0| + |\log \epsilon|)/C_0$ or
\begin{equation}\label{eq:T}
T = O(| \log \epsilon |).
\end{equation}
If we assume \eqref{eq:T} holds, then we can use a Gaussian quadrature to obtain a low-rank, separable expansion of $1/kr$:
\[ \frac{1}{kr} \approx \int_0^{T} e^{-krt} dt \approx \sum_{p=1}^n w_p e^{-krt_p} \]
where the $w_p$ are the Gaussian quadrature weights and the $t_p$ are the quadrature points. To determine $n$, the number of quadrature weights and points we need for an accuracy of order $\epsilon$, we can use the following Gaussian quadrature error estimate \cite{na} on the interval $[a,b]$:
\begin{equation}\label{eq:quaderr}
\frac{(b-a)^{2n+1} (n!)^4 }{(2n+1)[(2n)!]^3} f^{(2n)}(\xi)
\end{equation}
where $f$ is the integrand, and $\xi$ is within the bounds of integration: $f(\xi)=e^{-kr\xi}$ and $a \leq \xi \leq b$ where here $a=0$ and $b=T$. Clearly
\[ f^{(2n)}(\xi) = (-kr)^{2n}e^{-kr\xi}. \]
The worst case will be when $\xi=0$ and $r=1$ so 
\[\max_{0\leq \xi \leq T}\left| f^{(2n)}(\xi) \right|= (k)^{2n}.\]
We can put this back in the error estimate \eqref{eq:quaderr}, using Stirling's approximation \cite{ans}
$$ \sqrt{2\pi n}\ n^n e^{-n} \leq n! \leq e \sqrt{n} \ n^n e^{-n}, $$
for the factorials, to get:
\begin{eqnarray*}
\left| \frac{(b-a)^{2n+1} (n!)^4 }{(2n+1)[(2n)!]^3} f^{(2n)}(\xi) \right| &\leq& \frac{T^{2n+1} (n!)^4}{(2n+1)[(2n)!]^3} (k)^{2n} \\
&\leq& \frac{T^{2n+1} e^4(n)^2 (n)^{4n} e^{6n}}   {(2n+1)e^{4n}(2\pi 2n)^{3/2} (2n)^{6n} } (k)^{2n} \\
&\leq& \frac{T^{2n+1} e^4 (n)^{1/2} (ke)^{2n}} {(2n+1) \pi^{3/2} (2)^{6n+3} n^{2n}}  \\
&\leq& \frac{Te^4}{16\sqrt{n} \pi^{3/2}} \left( \frac{Tke}{8n}\right)^{2n} .
\end{eqnarray*}
This is problematic because in order for the quadrature scheme to converge, we are forced to have $n > Tke/8 \approx k \log \epsilon$, which is prohibitively large. This can be understood as the difficulty of integrating accurately a function which has large higher derivatives such as this sharp exponential over a large domain such as the interval $[0,\log{\epsilon}]$. To solve this problem, we make a dyadic partition of the $[0,T]$ interval in $O(\log{k})$ subintervals, each of which will require $O(\log{\epsilon})$ quadrature points.
Before we get to the details, let us redo the above error analysis for a dyadic interval $[a,2a]$. The maximum of $\left| f^{(2n)}(\xi) \right|= (kr)^{2n}e^{-kr\xi}$ as a function of $kr$ occurs when $kr=2n/\xi$, and the maximum of that as a function of $\xi$ is when $\xi=a$, so the maximum is $\left| f^{(2n)}(a) \right|= (2n/a)^{2n}e^{-(2n/a)*a}$. We can put this back in the error estimate \eqref{eq:quaderr} to get%
\begin{eqnarray*}
\left| \frac{(2a-a)^{2n+1} (n!)^4 }{(2n+1)[(2n)!]^3} f^{(2n)}(\xi) \right| &\leq& \frac{a^{2n+1} (n!)^4}{(2n+1)[(2n)!]^3} (2n/a)^{2n} e^{-2n} \\
&\leq& \frac{a e^4(n)^2 (n)^{4n} e^{6n}}   {(2n+1)e^{4n}(2\pi 2n)^{3/2} (2n)^{6n} } (2n)^{2n} e^{-2n} \\
&\leq& \frac{ 1.22 a (n)^{1/2} } {(2n+1) (2)^{4n} }  \\
&\leq& \frac{a } {\sqrt{n} (2)^{4n} } .
\end{eqnarray*}
To have this error less than $\epsilon$, we thus need
\begin{equation}\label{eq:n}
4n \log 2 \geq |\log\epsilon | + \log{a/\sqrt{n}},
\end{equation}
with $a \leq T \approx |\log\epsilon |$, and we see that in fact
\begin{equation}\label{eq:nval}
n=|\log \epsilon|
\end{equation}
will work.
\begin{remark}
We found the maximum of $\left| f^{(2n)}(\xi) \right|= (kr)^{2n}e^{-kr\xi}$ to be when both $kr=2n/\xi$ and $\xi=a$. However, we neek $kr\leq k$, so that we need $a\geq 2n/k=2|\log\eps|/k$. %
In the next subsection, we make sure that $a \geq 2|\log\eps|/k$ by having the lower of the interval $I_1$ equal to $2|\log\eps|/k$.%
\end{remark}

\subsection{Dyadic interval for the Gaussian quadrature}

Now we are ready get into the details of how we partition the interval. The subintervals are:
\begin{eqnarray*}
I_0&=&\left[0,\frac{2|\log\epsilon|}{k}\right] \\
I_j&=&\left[\frac{2^{j}|\log\epsilon|}{k},\frac{2^{j+1}|\log\epsilon|}{k}\right], \qquad j=1, \dots , M-1
\end{eqnarray*}
where $T=\frac{2^M|\log\epsilon|}{k}=O(|\log{\epsilon}|)$ which implies that
\begin{equation}\label{eq:M}
M=O(\log{k}).
\end{equation}
Then, for each interval $I_j$ with $j \geq 1$, we apply a Gaussian quadrature as explained above, and we need $n=|\log \epsilon|$ quadrature points to satisfy the error tolerance of $\epsilon$.
As for interval $I_0$, we return to the Gaussian quadrature error analysis, where this time again $k^{2n}$ is the maximum of $\left| f^{(2n)}(\xi) \right|$ for $\xi \in I_0$. Thus we have that the quadrature error is:
\begin{eqnarray*}
\left| \frac{(2|\log\epsilon|/k-0)^{2n+1} (n!)^4 }{(2n+1)[(2n)!]^3} f^{(2n)}(\xi) \right| &\leq& \frac{(2|\log\epsilon|/k)^{2n+1} (n!)^4}{(2n+1)[(2n)!]^3} k^{2n} \\
&\leq& \frac{(2|\log\epsilon|/k)^{2n+1} e^4 (n)^2 (n)^{4n} e^{6n}}   {(2n+1)e^{4n}(2\pi 2n)^{3/2} (2n)^{6n} } k^{2n} \\
&\leq& \frac{2^{2n+1}|\log\epsilon|^{2n+1} e^{2n} e^4 n^{1/2} } {k(2n+1) (2)^{6n+3} n^{2n}}  \\
&\leq& \frac{2|\log\epsilon| e^4 \sqrt{n}}{8k(2n+1)} \left( \frac{2|\log\epsilon| e}{8n} \right)^{2n} 
\end{eqnarray*}
and $n=O(|\log\epsilon|)$ will satisfy the error tolerance.

To recap, we have approximated the function $1/kr$, as a function of $r$, by a low-rank separable expansion with error $\epsilon$:
\[ \frac{1}{k|x-y|} = \sum_{j=1}^{J_2} w_j e^{-k|x-y|t_j} + E^{(2)}_{J_2}(x,y),\]
where $J_2=O(\log k |\log\epsilon|)$ (again, from using $O(\log k)$ intervals with $O(|\log\epsilon|)$ quadrature points on each interval), $C_0 \leq k|x-y| \leq k$, and $|E^{(2)}_{J_2}(x,y)|<\epsilon$.
Clearly this expansion is separable: depending on the sign of $(x-y)$, we have $e^{-k|x-y|t_j}=e^{-kxt_j}e^{kyt_j}$ or $e^{-k|x-y|t_j}=e^{kxt_j}e^{-kyt_j}$. Either way, the exponential has been expressed as a product of a function of $x$ only and a function of $y$ only. Thus we have the following lemma.

\begin{lemma}\label{lem:1overkr}
Let $0<\epsilon $, $k>0$ and $r_0>0$, $r_0=\Theta(1/k)$. Let $|x-y|>r_0$. Then there exists an integer $J_2$, a number $C_2$, and functions $\left\{\Phi^{(2)}_j,\chi^{(2)}_j\right\}_{j=1}^{J_2}$ such that
\begin{equation}\label{eq:kr}
\frac{1}{k|x-y|}=\sum_{j=1}^{J_2} \Phi^{(2)}_j(x)\chi^{(2)}_j(y) + E^{(2)}_{J_2}(x,y),
\end{equation}
where $|E^{(2)}_{J_2}(x,y)| \leq \epsilon$, and $J_2 \leq C_2 \log k |\log\epsilon|$ with $C_2$ which does not depend on $k$, or $\epsilon$. $C_2$ does depend weakly on $r_0$ through the constant quantity $r_0k$. Again, see Remark \ref{rem:C_0}.
\end{lemma}

\subsection{Finalizing the proof}

We now come back to the DtN map kernel $K$ in \eqref{eq:hsDtN}. %
Using Lemmas \ref{lem:7n8} and \ref{lem:1overkr}, we can write each factor of $K$ in its separable expansion:

\begin{eqnarray*}
\frac{K(|x-y|)}{ik^2/2}&=&\ H_1^{(1)}(k|x-y|) \ \frac{1}{k|x-y|} \\
&=&\left( \sum_{j=1}^{J_1} \Phi^{(1)}_j(x)\chi^{(1)}_j(y) + E^{(1)}_{J_1}(x,y) \right) \left( \sum_{j=1}^{J_2} \Phi^{(2)}_j(x) \chi^{(2)}_j(y) + E^{(2)}_{J_2}(x,y) \right)\\
&=& \frac{K_{(J_1,J_2)}(|x-y|)}{ik^2/2} + E^{(2)}_{J_2}(x,y)\sum_{j=1}^{J_1} \Phi^{(1)}_j(x)\chi^{(1)}_j(y) + E^{(1)}_{J_1}(x,y) \sum_{j=1}^{J_2} \Phi^{(2)}_j(x) \chi^{(2)}_j(y) \\
&+& E^{(1)}_{J_1}(x,y) E^{(2)}_{J_2}(x,y)
\end{eqnarray*}
where
\begin{eqnarray}\label{eq:hssepDtN}
K_{(J_1,J_2)}(|x-y|)&=&\frac{ik^2}{2}\left( \sum_{j=1}^{J_1} \Phi^{(1)}_j(x)\chi^{(1)}_j(y) \right) \left( \sum_{j=1}^{J_2} \Phi^{(2)}_j(x) \chi^{(2)}_j(y) \right)\\
&=&  \frac{ik^2}{2}\sum_{j=1}^{J_1 J_2} \Phi_j(x)\chi_j(y) .
\end{eqnarray}
It follows that
\begin{equation}\label{eq:errDs}
\left| \frac{K-K_{(J_1,J_2)}}{ik^2/2} \right| \leq \left|E^{(2)}_{J_2}\right| \left| \sum_{j=1}^{J_1} \Phi^{(1)}_j(x)\chi^{(1)}_j(y)\right| + \left|E^{(1)}_{J_1}\right| \left|\sum_{j=1}^{J_2} \Phi^{(2)}_j(x) \chi^{(2)}_j(y)\right|+ \left|E^{(1)}_{J_1}\right| \left|E^{(2)}_{J_2}\right|
\end{equation}

Now clearly
\[ \max_{C_0 \leq kr \leq k} \frac{1}{kr} = \frac{1}{C_0}. \]
We also have from Lemma 3 of \cite{flatland} that
\[ \left| e^{-ikr}H_1^{(1)}(kr) \right| \leq C(kr)^{-1/2}, \qquad kr \geq C_0 \]
for some constant $C$ which does not depend on $k$. Then, we have that
\[\max_{kr \geq C_0} \left| H_1^{(1)}(kr) \right| \leq \frac{C}{C_0^{1/2}}. \]
What we have shown is that the quantities $1/kr$ and $H_1^{(1)}(kr)$ are bounded by some constant, call it $\tilde{C}$, for the range of $kr$ we are interested, that is, $C_0 \leq kr \leq k$. We can now go back to \eqref{eq:errDs}, using our approximations from Lemmas \ref{lem:7n8} and \ref{lem:1overkr} which make an absolute error of no more than $\epsilon$, and see that
\begin{equation*}
\left| \frac{K-K_{(J_1,J_2)}}{ik^2/2} \right| \leq 2\epsilon (\tilde{C}+\epsilon) + \epsilon^2
\end{equation*}
or
\begin{equation}\label{eq:errDsrel}
\left| \frac{K-K_{(J_1,J_2)}}{ik^2/2} \right|=O(\epsilon).
\end{equation}
Note that the expansion of $K_{(J_1,J_2)}$ in \eqref{eq:hssepDtN} now contains $J_1 J_2 = O\left((\log k |\log\epsilon| )^2\right)$ terms. In order to obtain an absolute error bound, we need to multiply through with $ik^2/2$ in \eqref{eq:errDsrel}. With $\epsilon \rightarrow \epsilon k^2$, we have now finally shown that
\begin{equation*}
\left| K-K_{(J_1,J_2)}\right| \leq \epsilon
\end{equation*}
with the expansion of $K_{(J_1,J_2)}$ in \eqref{eq:hssepDtN} containing
\[J=J_1 J_2 =O\left( (\log k (|\log\epsilon|+2\log k) )^2 \right)\]
terms. We thus conclude that the DtN map is low-rank and separable away from the diagonal, with the prescriptions of Theorem \ref{theo:sep}:
\begin{equation}%
K(|x-y|)= \frac{ik^2}{2} \frac{H_1^{(1)}(k|x-y|)}{k|x-y|}= \sum_{j=1}^J \Phi_j(x)\chi_j(y) + E(x,y)
\end{equation}
where $|E(r)| \leq \epsilon$ for $C_0 \leq kr \leq k$ and there is a number $C$ which does not depend on $k$ or $\epsilon$ such that $J \leq C (\log k \max(|\log\epsilon|,\log k))^2$.

\begin{remark}\label{rem:C_0}
We can understand the number $J$ of ranks as made of two terms, one which is $(C_1\log k) (C_2|\log\epsilon|)$, the other $(C_1\log k)(\log C_3 k)$. The numbers $C_1$and $C_3$, but not $C_2$, also weakly depend on the separation $C_0$, in the sense that the larger the separation is, the smaller those numbers are. First, we note from the discussion before equation $\eqref{eq:T}$ that $T$ is smaller when $C_0$ is bigger. Then, a smaller $C_1$ comes from the discussion before equation $\eqref{eq:M}$. The fact that $C_2$ does \emph{not} depend on $C_0$ comes from the discussion after equation \eqref{eq:n}. As for $C_3$, we can understand its \emph{very weak} dependence on $C_0$ by looking at equation \eqref{eq:n} and plugging in $a=T$, remembering how $T$ depends on $C_0$. Thus both terms in $J$ should behave somewhat similarly as $C_0$ changes. Physically, we know a greater separation means we are farther away from the singularity of the DtN map, and so we expect the map to be smoother there, and hence have lower rank. 
\end{remark}
\begin{remark}\label{rem:highpow}
In our numerical verifications, we have not noticed the square power in $J$. Rather, we observe in general that $J \sim \log k |\log\epsilon|$. The only exceptions to this behavior that we observed were for larger $\epsilon$, such as $\epsilon=1/10$ or sometimes $1/100$, where $J \sim \log k \log k$. From Remark \ref{rem:C_0}, we know $J$ is made up of two terms, and it makes sense that the term $\sim \log k \log k$ might become larger than the term $\sim \log k |\log \eps|$ when $\eps$ is large. %
\end{remark}
\begin{remark}\label{rem:h}
We also note that in our numerical verifications, we use $r_0$ as small as $h$, which is smaller than the $r_0\sim 1/k \sim h^{2/3}$ we prove the theorem with. If we used $r_0\sim h$ in the theorem, this would mean $C_0\sim n^{-1/3}<1$. By Remark \ref{rem:C_0}, this would affect $J$: the power of the $|\log\epsilon |$ factor would go from 2 to 3. Again, we do not notice such a higher power in the numerical simulations.
\end{remark}

\subsection{The numerical low-rank property of the DtN map kernel for heterogeneous media}

We would like to know as well if we can expect the half-space DtN map kernel to be numerically low-rank in heterogeneous media.  We saw in section \ref{sec:BasisPf} how the half-space DtN map, in constant medium, consists of the amplitude $H(r)$, singular at $r=0$, multiplied by the complex exponential $e^{ikr}$. Because of the geometrical optics expansion \eqref{eq:geoopts} of the Green's function $G$ for the Helmholtz equation free-space problem in heterogeneous media, we expect $G$ to have an amplitude $A(\x,\y)$, which is singular at $\x=\y$, multiplied by a complex exponential with a phase corresponding to the traveltime between points $\x$ and $\y$: $e^{i\omega \tau(\x,\y)}$. We can expect to be able to treat the amplitude in the same way as we did before, and approximate it away from the singularity with a low-rank separable expansion. However, the complex exponential is harder to analyze because of the phase $\tau(\x,\y)$, which is not so simple as $|\x-\y|$.

However, a result of \cite{butterflyFIO} still allows us to find a separable low-rank approximation of a function such as $e^{i\omega \tau(\x,\y)}$. We refer the reader to Theorem 3.1 of \cite{butterflyFIO} for the details of the proof, and simply note here the main result. Let $X$, $Y$ be bounded subsets of $\mathbf{R}^2$ such that we only consider $\x \in X$ and $\y \in Y$. The width of $X$ (or $Y$) is defined as the maximal Euclidean distance between any two points in that set. Then Theorem 3.1 of \cite{butterflyFIO} cites that $e^{i\omega \tau(\x,\y)}$, $\x \in X$ and $\y \in Y$, is numerically low-rank with rank $O(|\log^4 \eps|)$, given that the product of the widths of $X$ and $Y$ is less than $1/k$. 

This translates into a restriction on how large can the off-diagonal blocks of the matrix $D$ be while still being low-rank. Since we use square blocks in the PLR compression algorithm, we expect blocks to have to remain smaller than $1/\sqrt{k}$, equivalent to $N/\sqrt{k}$ points, in variable media. If $N^{2/3} \sim k$ as we have in this thesis because of the pollution effect, this translates into a maximum expected number of blocks of $N^{1/3}$. If we kept instead $N \sim k$, then the maximal number of blocks would be $\sqrt{N}$.

This is why, again, using PLR matrices for compressing the DtN map makes more sense than using hierarchical matrices: the added flexibility means blocks will be divided only where needed, in other words only where the traveltime $\tau$ requires blocks to be smaller in order to have low ranks. And as we saw in section \ref{sec:usingplr}, where we presented our numerical results, ranks indeed remain very small, between 2 and 8, for off-diagonal blocks of the submatrices of the exterior DtN map, even in heterogeneous media.

\section{The half-space DtN map is separable and low rank: numerical verification}\label{sec:sepnum}

We first compute the half-space DtN map for various $k \sim N^{2/3}$, which ensures a constant error from the finite difference discretization (FD error) as we saw in section \ref{sec:compABC}. We also choose a pPML width consistent with the FD error level. Then, we compute the maximum off-diagonal ranks for various fixed separations from the diagonal, that is, various $r_0$ such that $r \geq r_0$. To compute the ranks of a block, we fix a tolerance $\epsilon$, find the Singular Value Decomposition of that block, and discard all singular values smaller than that tolerance. The number of remaining singular values is the numerical rank of that block with tolerance $\epsilon$ (the error we make in Frobenius norm is not larger than $\epsilon$). Then, the maximum off-diagonal rank for a given separation $r_0$ is the maximum rank of any block whose entries correspond to $r\geq r_0$. %
Hence we consider all blocks that have $|i-j| \geq r_0/h$, or $i-j \geq r_0/h$ with $i>j$ since the DtN map is symmetric (and so is its numerical realization, up to machine precision).

\subsection{Slow disk and vertical fault}

\begin{figure}[h]
\begin{minipage}[t]{0.48\linewidth}
\includegraphics[scale=.5]{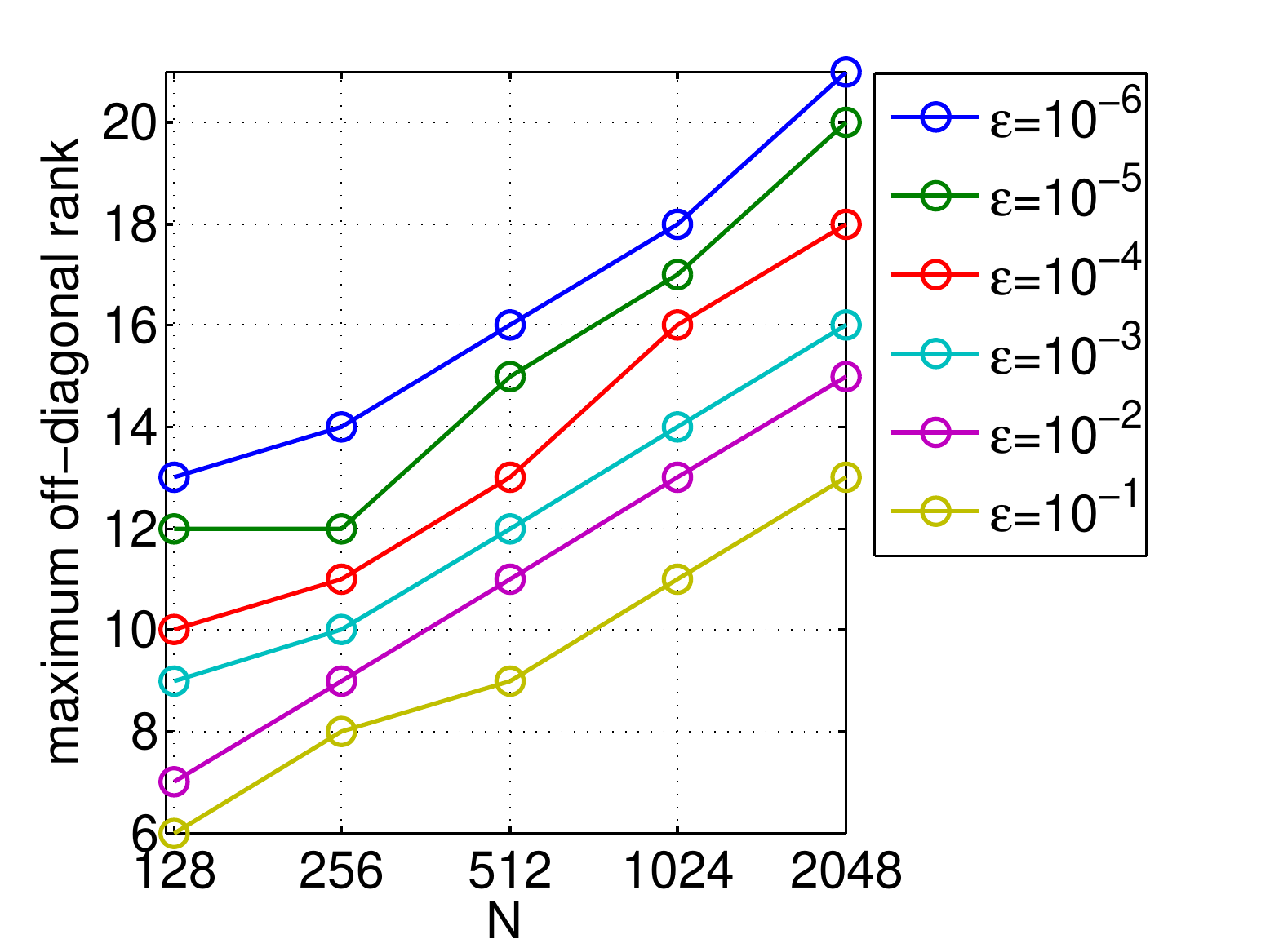}
\caption{Maximum off-diagonal ranks with $N$ for the slow disk, various $\epsilon$. FD error of $10^{-3}$, $r_0=h$.}
\label{c5-vsN-r1}
\end{minipage}
\hspace{0.1cm}
\begin{minipage}[t]{0.48\linewidth}
\includegraphics[scale=.5]{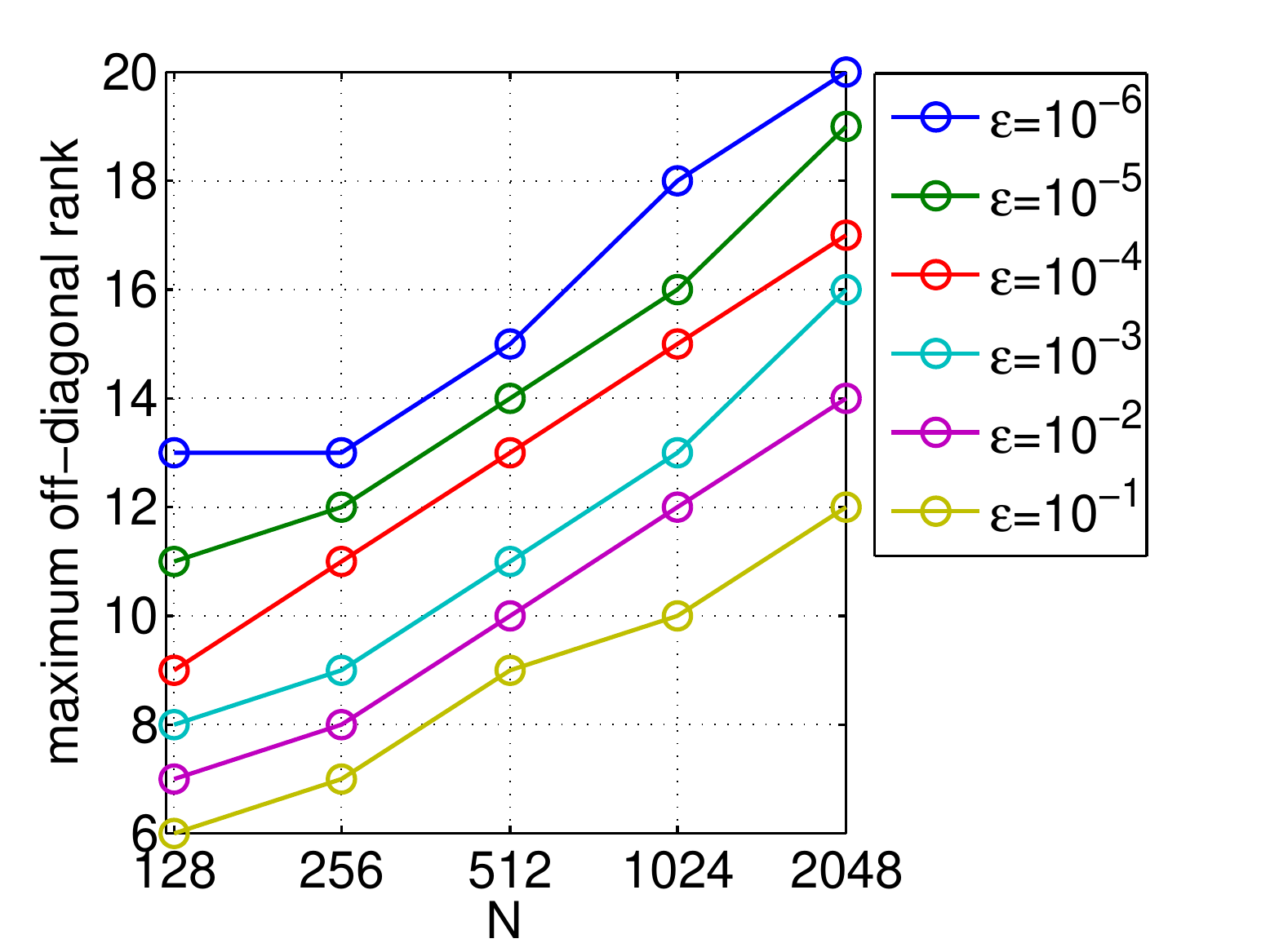}
\caption{Maximum off-diagonal ranks with $N$ for the slow disk, various $\epsilon$. FD error of $10^{-3}$, $r_0=4h$.}
\label{c5-vsN-r4}
\end{minipage}
\end{figure}

\begin{figure}[h]
\begin{minipage}[t]{0.48\linewidth}
\includegraphics[scale=.5]{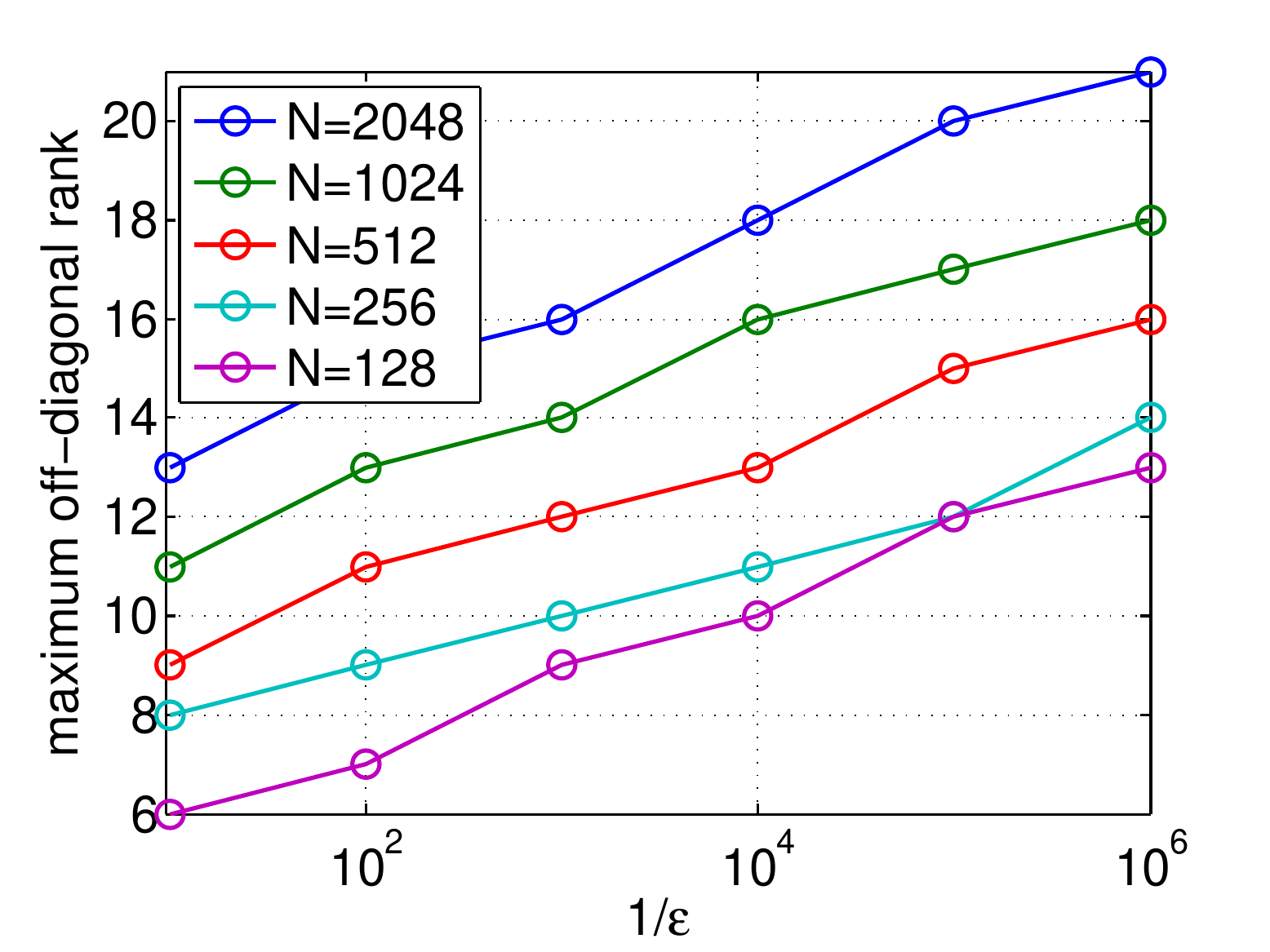}
\caption{Maximum off-diagonal ranks with $\eps$ for the slow disk, various $N$. FD error of $10^{-3}$, $r_0=h$.}
\label{c5-vsEps-r1}
\end{minipage}
\hspace{0.1cm}
\begin{minipage}[t]{0.48\linewidth}
\includegraphics[scale=.5]{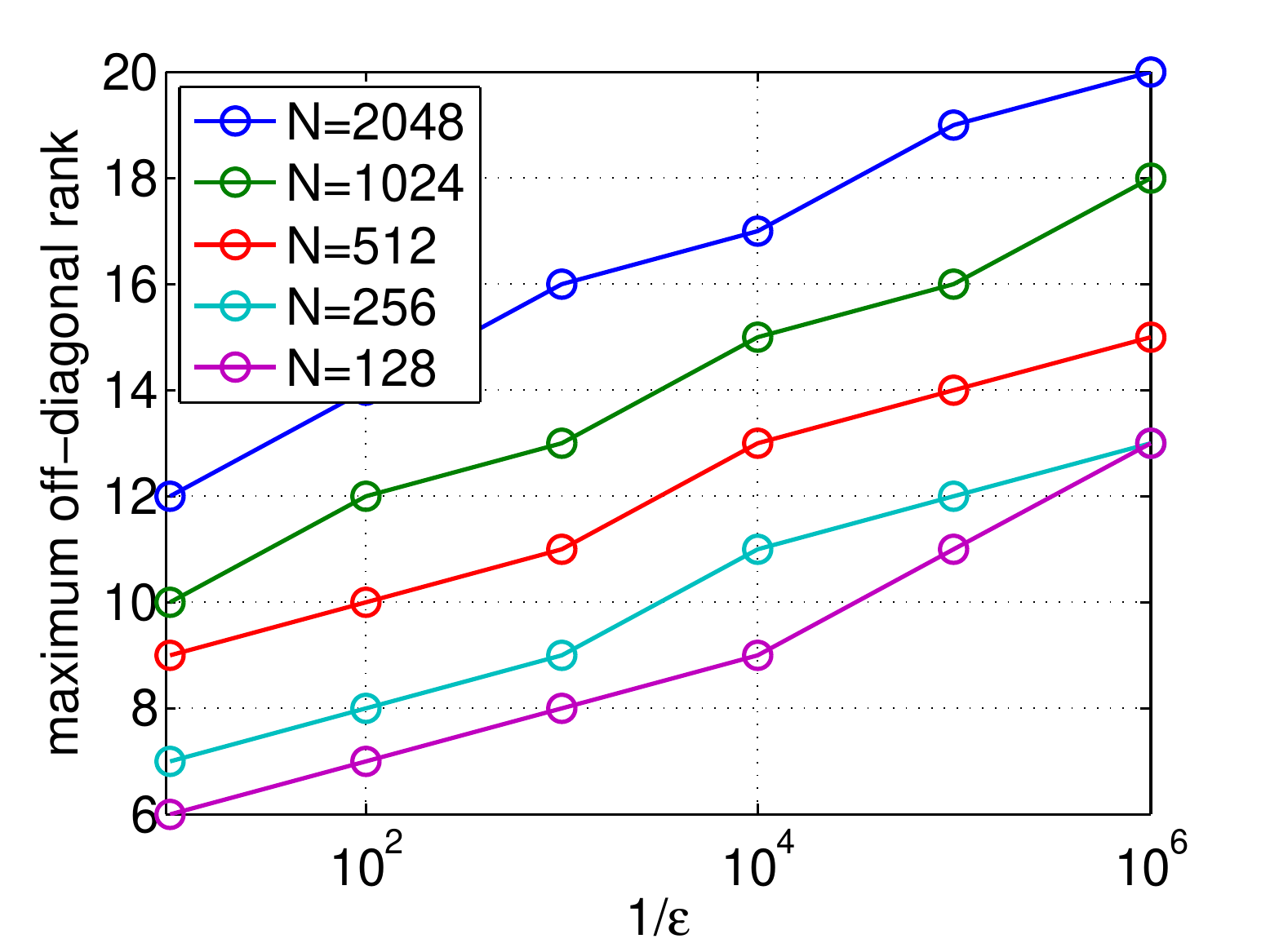}
\caption{Maximum off-diagonal ranks with $\eps$ for the slow disk, various $N$. FD error of $10^{-3}$, $r_0=4h$.}
\label{c5-vsEps-r4}
\end{minipage}
\end{figure}

The Figures \ref{c5-vsN-r1}, \ref{c5-vsN-r4}, \ref{c5-vsEps-r1}, \ref{c5-vsEps-r4} show the relationship between the ranks and $N$ or $\eps$ for the slow disk, FD error of $10^{-3}$ and separations of $r_0=1$ and $r_0=4$.

We expect ranks to be slightly smaller for a larger separation $r_0$ (hence larger $C_0$) because of Remark \ref{rem:C_0}. This is indeed the case in our numerical simulations, as we can see by comparing Figures \ref{c5-vsN-r1} ($r_0=h$) and \ref{c5-vsN-r4} ($r_0=4h$), or Figures \ref{c5-vsEps-r1} ($r_0=h$) and \ref{c5-vsEps-r4} ($r_0=4h$). We can clearly see also how the maximum ranks behave as in the previous theorem, except for the missing square power in $J$, as alluded to in Remark \ref{rem:highpow}: they vary logarithmically with $k$ (or $N$) when the tolerance $\epsilon$ is fixed. We expect the slope in a graph of the ranks versus $\log N$ to increase slowly as $\epsilon$ becomes smaller, and we see that in fact the slope barely does (from slightly smaller than $2$ to slightly larger than $2$) as $\epsilon$ goes from $10^{-1}$ to $10^{-6}$. Similarly, when we fix $N$, we expect the ranks to grow logarithmically with $1/\epsilon$, and this is the case. Once again, the slope of the graph with a logarithmic scale for $1/\epsilon$ grows, but only from $1$ to $2$ or so, as $N$ goes from $128$ to $2048$.

The off-diagonal ranks of the DtN map for the slow disk behave very similarly as the above for an FD error of $10^{-2}$, and also for various other separations $r_0$. The same is true for the vertical fault, and so we do not show those results.

\subsection{Constant medium, waveguide, diagonal fault}

\begin{figure}[ht]
\begin{minipage}[t]{0.48\linewidth}
\includegraphics[scale=.5]{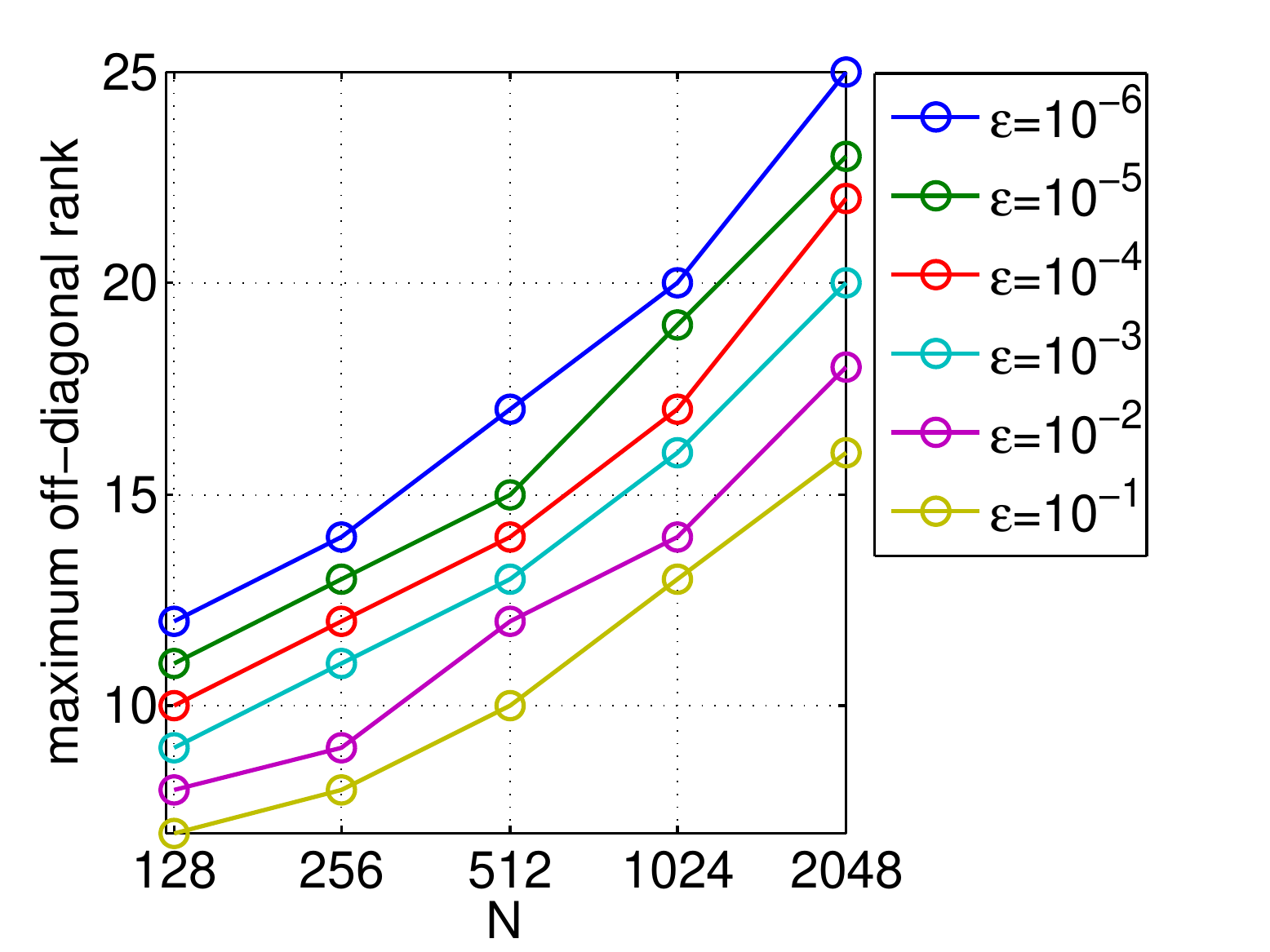}
\caption{Maximum off-diagonal ranks with $N$ for the diagonal fault, various $\epsilon$. FD error of $10^{-2}$, $r_0=h$. }
\label{c18-vsN-r1}
\end{minipage}
\hspace{0.1cm}
\begin{minipage}[t]{0.48\linewidth}
\includegraphics[scale=.5]{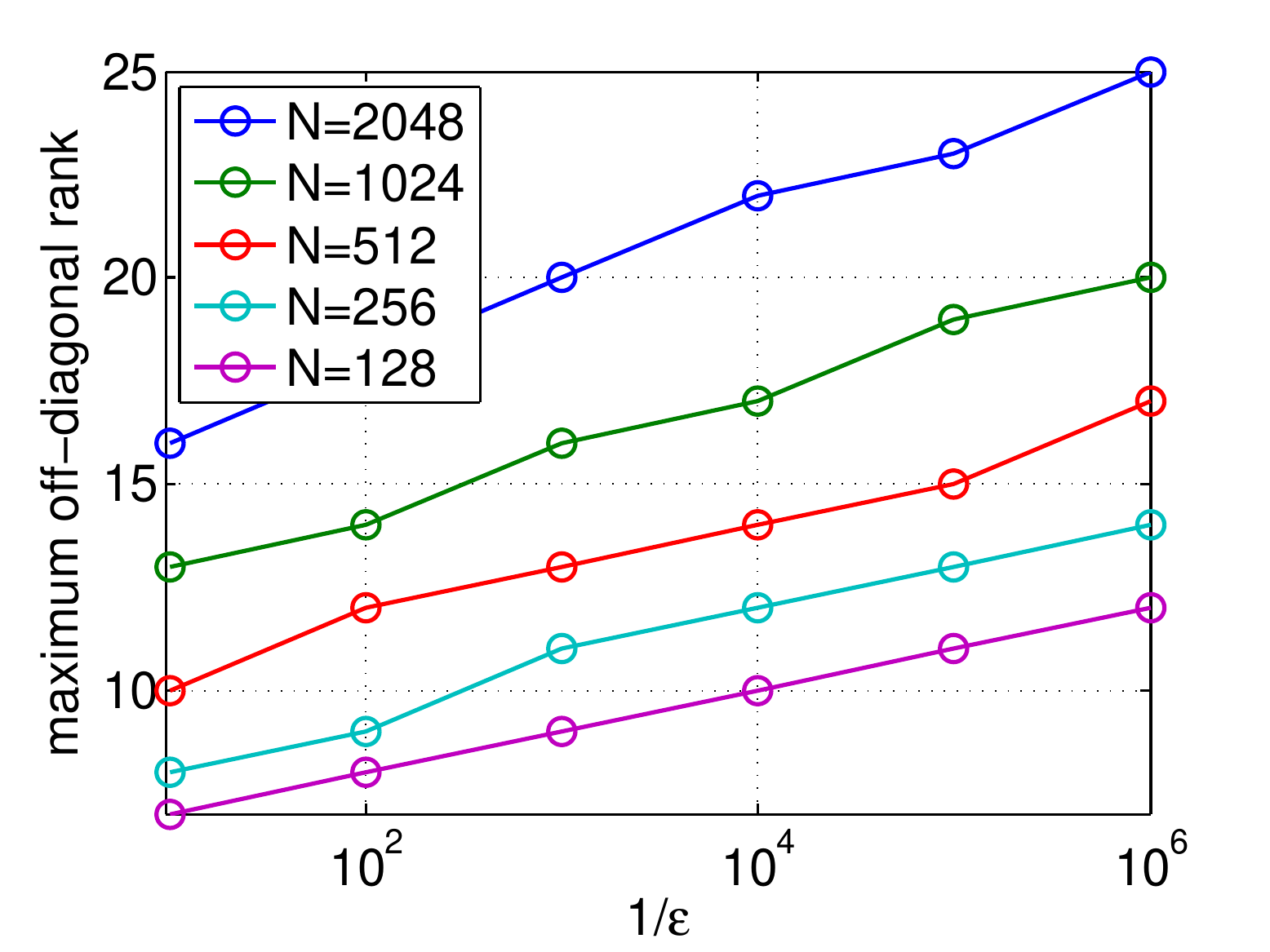}
\caption{Maximum off-diagonal ranks with $\eps$ for the diagonal fault, various $N$. FD error of $10^{-2}$, $r_0=h$.}
\label{c18-vsEps-r1}
\end{minipage}
\end{figure}

\begin{figure}[ht]
\begin{minipage}[t]{0.48\linewidth}
\includegraphics[scale=.5]{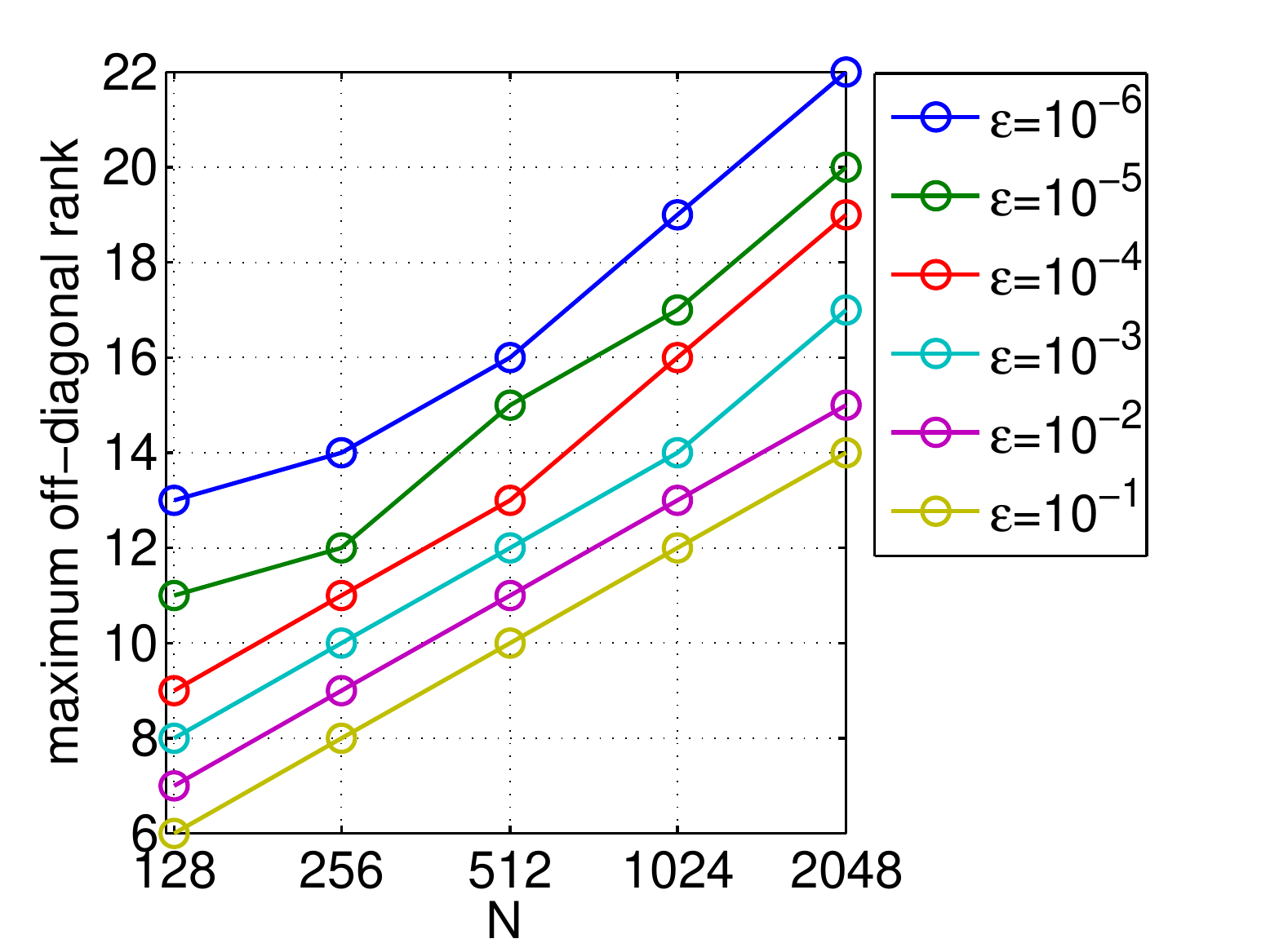}
\caption{Maximum off-diagonal ranks for the diagonal fault as a function of $N$, various tolerances $\epsilon$. The separation is $r_0=h$. FD error of $10^{-3}$.}
\label{c18-vsN-r1-FDm3}
\end{minipage}
\hspace{0.1cm}
\begin{minipage}[t]{0.48\linewidth}
\includegraphics[scale=.5]{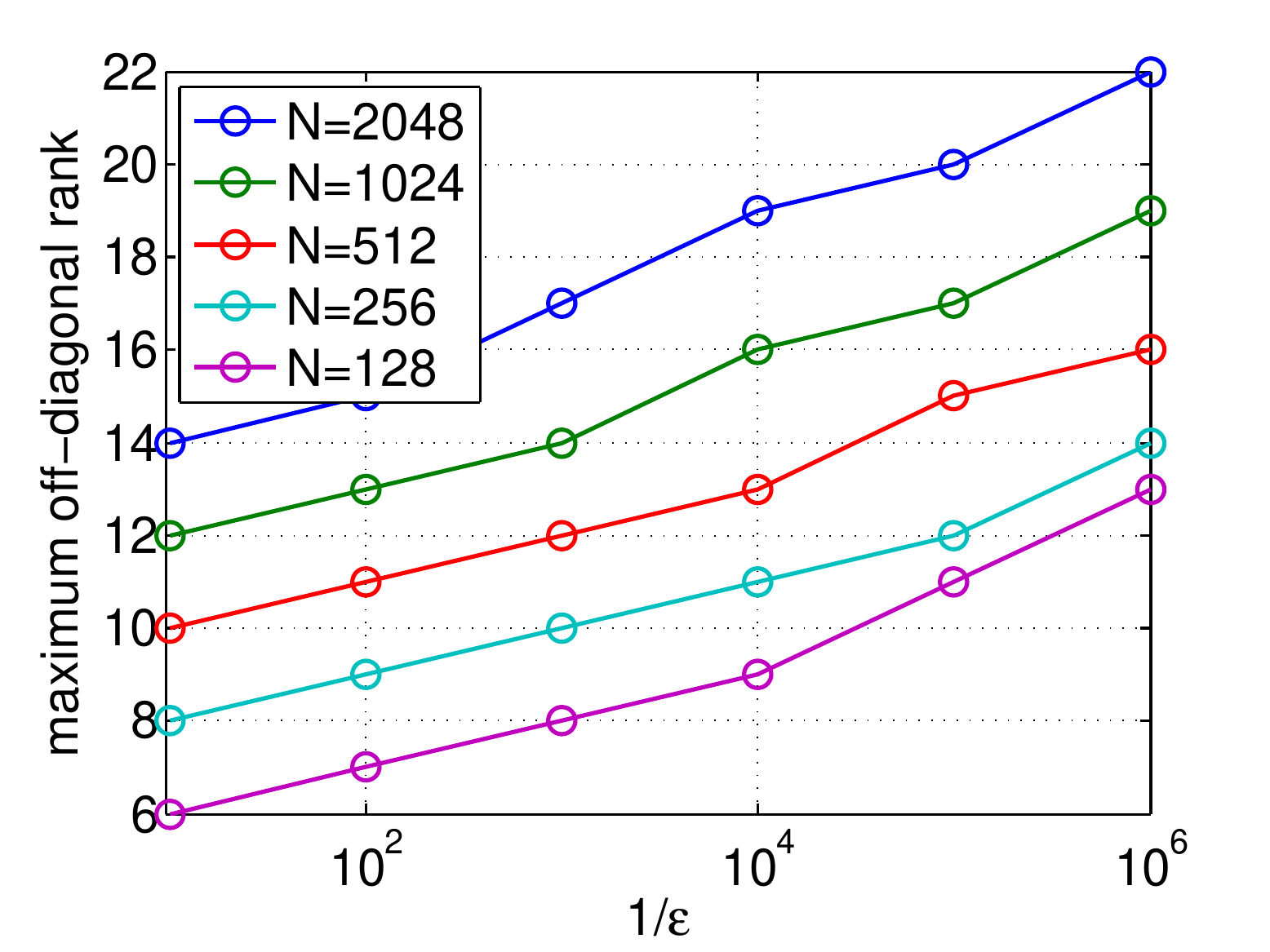}
\caption{Maximum off-diagonal ranks for the diagonal fault as a function of the tolerance $\epsilon$, various $N$. The separation is $r_0=h$. FD error of $10^{-3}$.}
\label{c18-vsEps-r1-FDm3}
\end{minipage}
\end{figure}

As for the constant medium, waveguide and diagonal fault, it appears that the term $O(\log^2 k)$ we expect in the size of the ranks is larger than the term $O(\log k |\log\epsilon|)$, especially when the FD error is $10^{-2}$. This was mentioned in Remark \ref{rem:highpow}. As we can see in Figure \ref{c18-vsN-r1} for the diagonal fault, the dependance of the ranks with $\log N$ seems almost quadratic, not linear. This can also be seen in Figure \ref{c18-vsEps-r1}: here we still see a linear dependence of the ranks with $\log\epsilon$, but we can see that the ranks jump up more and more between different $N$, as $N$ grows, than they do for the slow disk for example (compare to Figure \ref{c5-vsEps-r1}). %
This phenomenon disappears for a smaller FD error (Figures \ref{c18-vsN-r1-FDm3}, \ref{c18-vsEps-r1-FDm3}).%

\begin{figure}[ht]
\begin{minipage}[t]{0.48\linewidth}
\includegraphics[scale=.5]{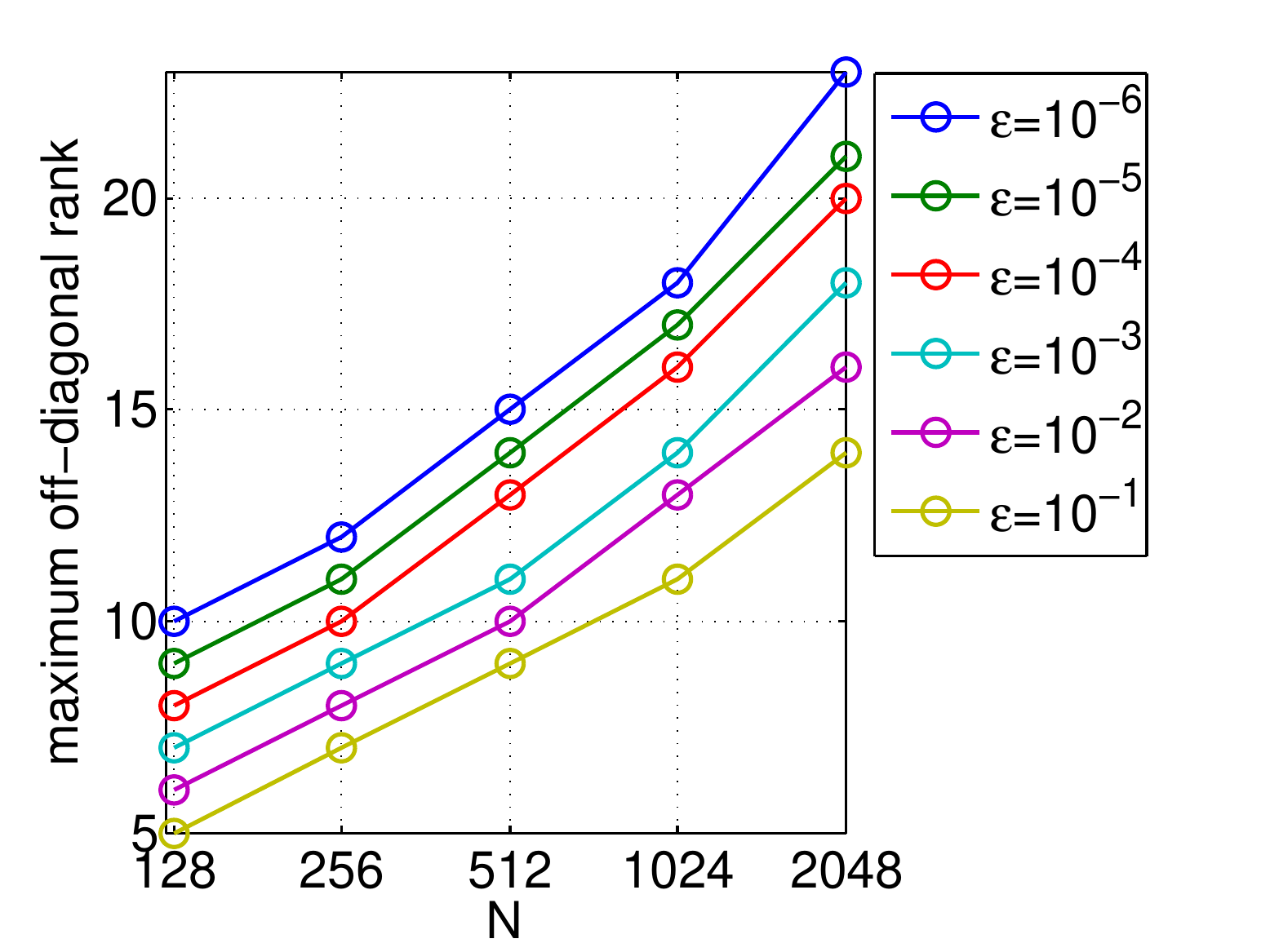}
\caption{Maximum off-diagonal ranks for the diagonal fault as a function of $N$, various tolerances $\epsilon$. The separation is $r_0=8h$. FD error of $10^{-2}$.}
\label{c18-vsN-r8}
\end{minipage}
\hspace{0.1cm}
\begin{minipage}[t]{0.48\linewidth}
\includegraphics[scale=.5]{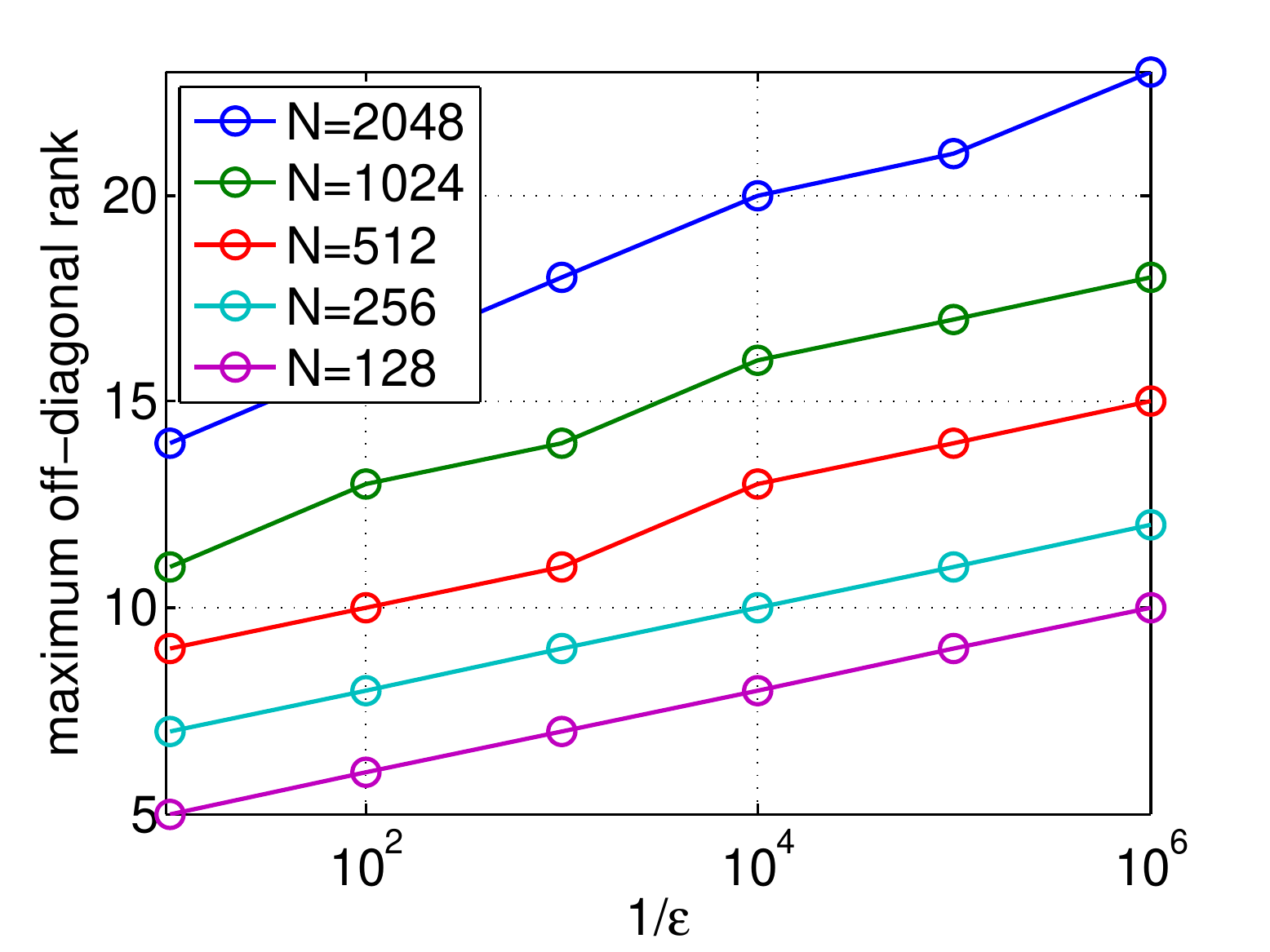}
\caption{Maximum off-diagonal ranks for the diagonal fault as a function of the tolerance $\epsilon$, various $N$. The separation is $r_0=8h$. FD error of $10^{-2}$.}
\label{c18-vsEps-r8}
\end{minipage}
\end{figure}

Finally, we also notice that the term $O(\log^2 k)$ seems to remain important compared to the term $O(\log k |\log\epsilon|)$ as the separation $r_0$ (or $C_0$) grows, as is somewhat expected from Remark \ref{rem:C_0}. This can be seen by comparing Figures \ref{c18-vsN-r8} and \ref{c18-vsN-r1}, or Figures \ref{c18-vsEps-r8} and \ref{c18-vsEps-r1}.

\subsection{Focusing and defocusing media}

We have also tested a smooth defocusing medium, that is, a medium where the value of $c$ decreases away from the half-space boundary, tending to a small value far away from the boundary. The equation we have used is $c(x,y)=1+\frac{1}{\pi}\arctan(4(y-1/2))$. This means that as $y \rightarrow \infty$, $c \rightarrow 3/2$ and that as $y \rightarrow -\infty$, $c \rightarrow 1/2$. Choosing the half-space to be $y<0$, we see that $c$ decreases away from $y=0$ into the negative $y$'s: this is a defocusing medium. We expect the waves in this case to never come back, and so we expect off-diagonal ranks of the DtN map to be very nice, just as in the constant medium case, and this is indeed what happens. We could not see any significant difference between the defocusing medium and the constant medium in terms of off-diagonal ranks.

We have also looked at a focusing medium, that is, one in which $c$ increases away from the interface. This forces waves to come back toward the interface. With the same medium $c(x,y)=1+\frac{1}{\pi}\arctan(4(y-1/2))$ as above, but choosing now $y>1$ as our half-space, we see that $c$ increases away from $y=1$ into the large positive $y$'s. This is a focusing medium. We have noticed that the off-diagonal ranks of the DtN map for this medium are the same or barely larger than for the constant medium.

\begin{figure}[h]
\begin{minipage}[t]{0.48\linewidth}
\includegraphics[scale=.5]{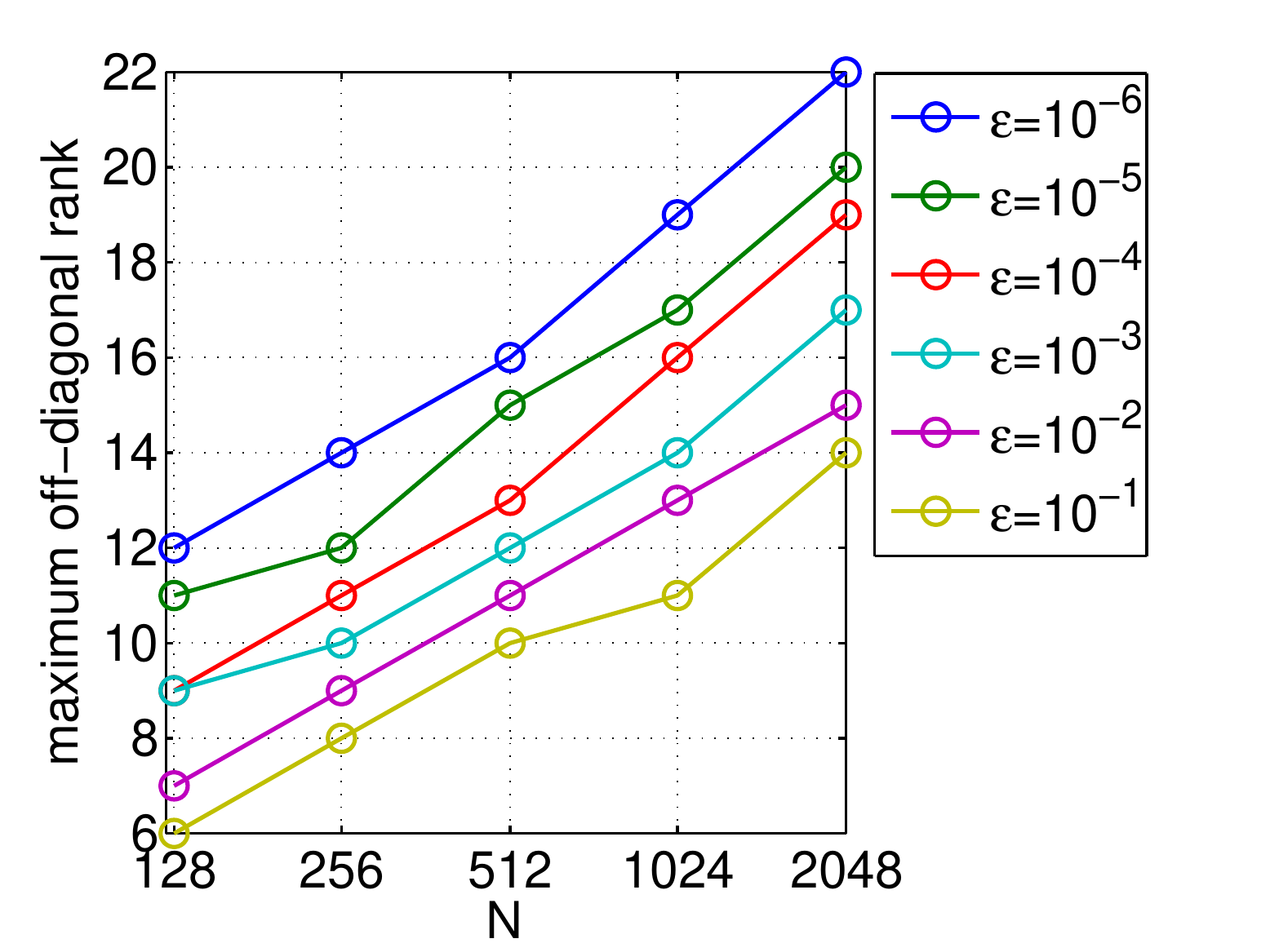}
\caption{Maximum off-diagonal ranks for the focusing medium \eqref{eq:focus} as a function of $N$, various tolerances $\epsilon$. Separation is $r_0=h$, FD error of $10^{-2}$.}
\label{c8-vsN-r1-FD2}
\end{minipage}
\hspace{0.1cm}
\begin{minipage}[t]{0.48\linewidth}
\includegraphics[scale=.5]{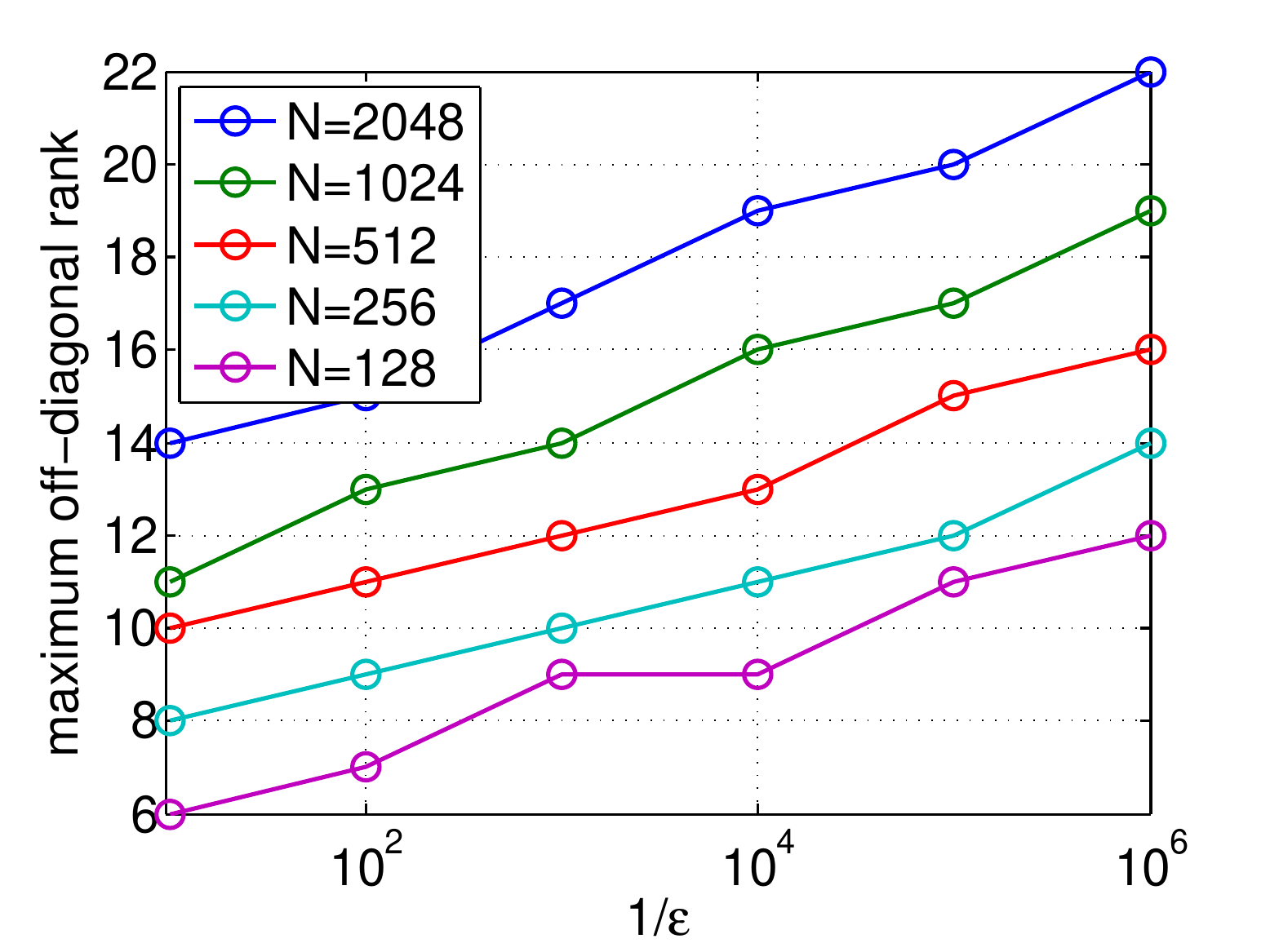}
\caption{Maximum off-diagonal ranks for the focusing medium \eqref{eq:focus} as a function of the tolerance $\epsilon$, various $N$. Separation is $r_0=h$, FD error of $10^{-2}$.}
\label{c8-vsEps-r1-FD2}
\end{minipage}
\end{figure}

\begin{figure}[h]
\begin{minipage}[t]{0.48\linewidth}
\includegraphics[scale=.5]{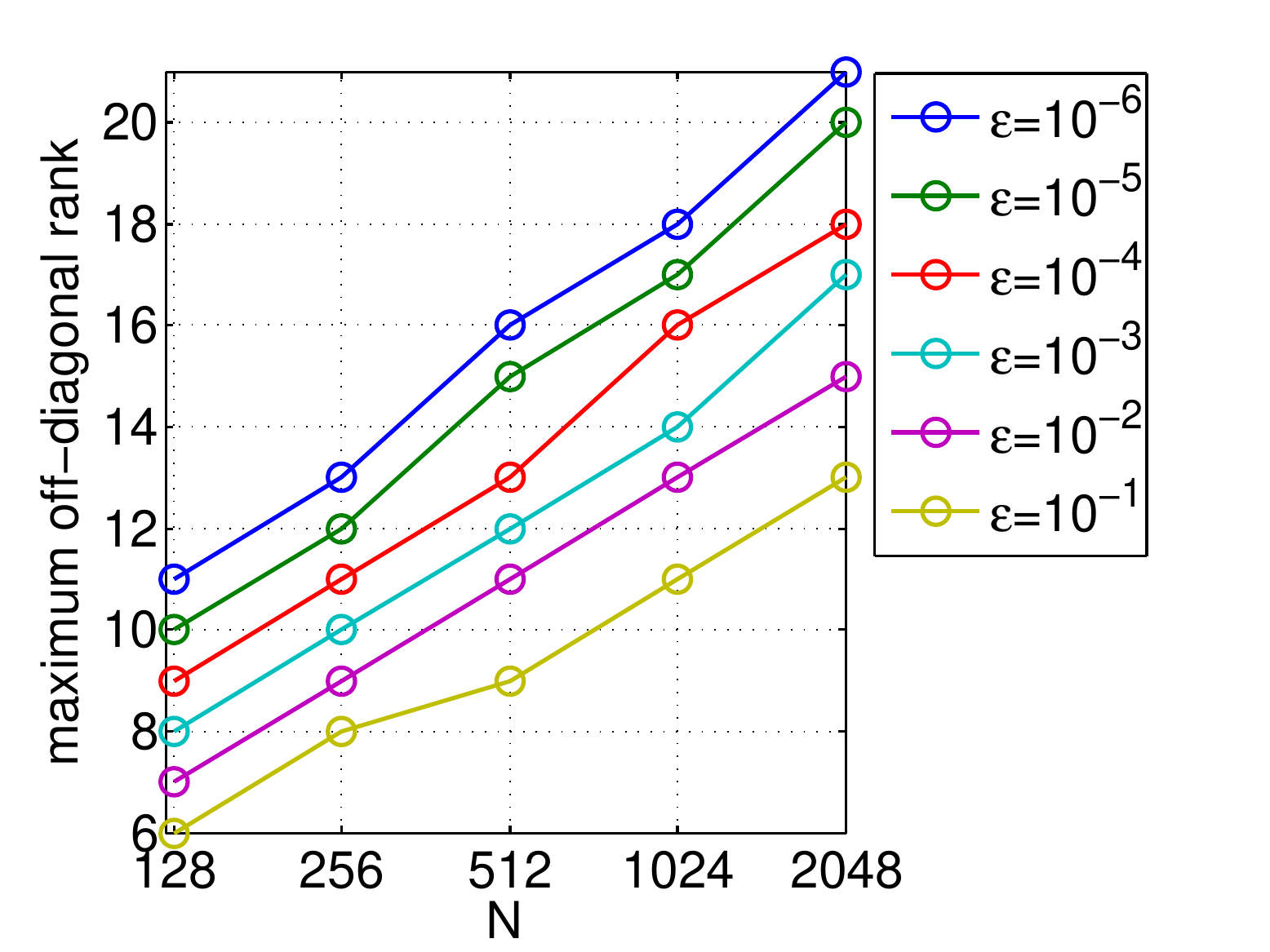}
\caption{Maximum off-diagonal ranks for the focusing medium \eqref{eq:focus} as a function of $N$, various tolerances $\epsilon$. Separation is $r_0=h$, FD error of $10^{-3}$.}
\label{c8-vsN-r1-FD3}
\end{minipage}
\hspace{0.1cm}
\begin{minipage}[t]{0.48\linewidth}
\includegraphics[scale=.5]{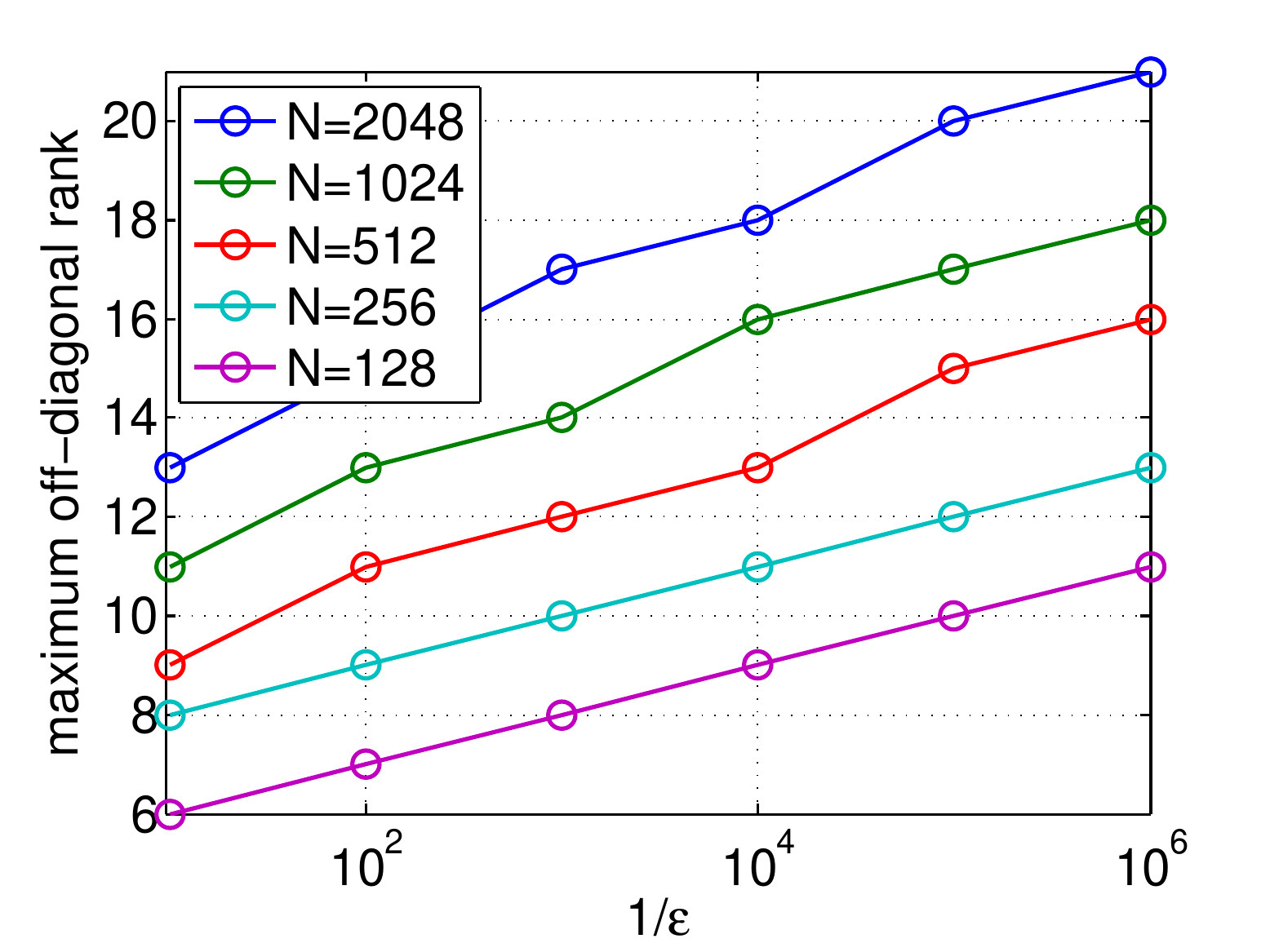}
\caption{Maximum off-diagonal ranks for the focusing medium \eqref{eq:focus} as a function of the tolerance $\epsilon$, various $N$. Separation is $r_0=h$, FD error of $10^{-3}$.}
\label{c8-vsEps-r1-FD3}
\end{minipage}
\end{figure}

This might only mean that the medium we chose did not have many returning waves. A more interesting medium is the following:
\begin{equation}\label{eq:focus}
c(x,y)=1/2+|y-1/2|.
\end{equation}
This linear $c$ has a first derivative bounded away from 0. Of course, this means that solving the Helmholtz equation in this case is much harder, and in particular, the pPML layer needs to be made thicker than for other media. Still, we notice that the ranks are very similar to the previous cases, as we can see in Figures \ref{c8-vsN-r1-FD2}, \ref{c8-vsEps-r1-FD2}, \ref{c8-vsN-r1-FD3}, \ref{c8-vsEps-r1-FD3}.

This could be explained by the fact that we do not go to high enough accuracies to really see those returning waves. Saying this another way, the relative amplitude of the returning waves might be too small to notice at lower accuracies. We were not able to construct an example with returning waves which affected the off-diagonal ranks of the DtN map. Thus we conclude that Theorem \ref{theo:sep} holds in a broad range of contexts, at least when $\epsilon$ is not very small.

%% file: conclusion2.tex
\chapter{Conclusion}\label{ch:conclusion}

In this thesis, we have compressed the Dirichlet-to-Neumann map for the Helmholtz equation in two steps, using matrix probing followed by the partitioned low-rank matrix framework. This approach is useful for applications in heterogeneous media because absorbing boundary conditions can be very costly. Especially when one needs to solve the Helmholtz equation with many different sources, it makes sense to perform the precomputation required of our two-step scheme, to speed up subsequent solves.

Probing the DtN map $D$ ultimately makes sense in conjunction with a fast algorithm for its application. In full matrix form, $D$ costs $N^2$ operations to apply. With the help of a compressed representation, this count becomes $p$ times the application complexity of any atomic basis function $B_j$, which may or may not be advantageous depending on the particular expansion scheme. The better solution for a fast algorithm, however, is to post-process the compressed expansion from probing into a slightly less compressed, but more algorithmically favorable one, such as hierarchical or partitioned low-rank matrices. These types of matrix structures are not parametrized enough to lend themselves to efficient probing -- see for instance \cite{Hmatvec} for an illustration of the large number of probing vectors required -- but give rise to faster algorithms for matrix-vector multiplication. Hence the feasibility of probing, and the availability of a fast algorithm for matrix-vector multiplication, are two different goals that require different expansion schemes.

We found that in order to obtain an efficient compression algorithm, we need to perform some precomputations. The leading cost of those is equivalent to a small number of solves of the original problem, which we can afford if we plan to make many solves in total. Then, a matrix-vector application of the DtN map in a Helmholtz solver is nearly linear in the dimension of the DtN map. The worst-case complexity of a matrix-vector application is in fact super-linear. General oscillatory integrals can often be handled in optimal complexity with the butterfly algorithm \cite{butterflyFIO}. We did not consider this algorithm because currently published methods are not adaptive, nor can they handle diagonal singularities in a kernel-dependent fashion.
\newline
\newline
Let us summarize the different complexity scalings of our method, recalling that without our compressed ABC, a thousand solves of the free-space problem would cost a thousand times the complexity of one solve, and that complexity depends on the solver used and on the fact that the ABC used might be very costly. By contrast, the cost of precomputation for matrix probing is of a few (1 to 50) solves of the exterior problem -- roughly equivalent to a solve of the free-space problem with the costly ABC. The precomputation cost for PLR compression is $O(N R^2 |\mathcal{B}|)$ where $|\mathcal{B}|$ is the number of leaves in the compressed DtN map -- empirically $|\mathcal{B}|=O(\log N)$ with a worst case of $|\mathcal{B}|=O(\sqrt{N})$. Finally, the application cost of the compressed ABC is $O(N R |\mathcal{B}|)$, which typically leads to a speed-up of a factor of 30 to 100 for $N \approx 1000$. Using the compressed ABC in a solver, the computational cost of making a thousand solves is then reduced to that of a thousand solves where the ABC is not so costly anymore, but only amounts to a matrix-vector multiplication of nearly linear cost.

Of course, the method presented here has some limitations. The main ones come from the first step of the scheme, matrix probing, which requires a careful design of basis matrices. This is harder to do when the wavelength is comparable to the size of features in the medium, or when the medium has discontinuities.
\newline
\newline
In addition to the direct applications the proposed two-step scheme has in making a typical (iterative) Helmholtz solver faster when in presence of a heterogeneous media, it could also improve solvers based on Domain Decomposition. Some interesting recent work on fast solvers for the Helmholtz equation has been along the direction of the Domain Decomposition Method (DDM). This method splits up the domain $\Omega$ in multiple subdomains. Depending on the method, the subdomains might or might not overlap. In non-overlapping DDM, transmission conditions are used to transfer information about the solution from one domain to the next. Transmission conditions are crucial for the convergence of the algorithm. After all, if there is a source in one subdomain, this will create waves which should travel to other parts of $\Omega$, and this information needs to be transmitted to other subdomains as precisely as possible.%

In particular, Stolk in \cite{stolk} has developed a nearly linear complexity solver for the Helmholtz equation, with transmission conditions based on the PML. Of course, this relies on the PML being an accurate absorbing boundary condition. One may assume that, for faster convergence in heterogeneous media, thicker PMLs might be necessary. In this case, a precomputation to compress the PML might prove useful. Other solvers which rely on transmission conditions are those of Engquist and Ying \cite{Hsweep,Msweep}, Zepeda-N\'u\~nez et al. \cite{zepeda}, Vion and Geuzaine \cite{Vion}.

Another way in which the two-step numerical scheme presented in this thesis could be useful is for compressing the Helmholtz solution operator itself. Indeed, we have compressed ABCs to the DtN map, by using knowledge of the DtN map kernel. If one knows another kernel, one could expand that kernel using matrix probing then compress it with PLR matrices, and obtain a fast application of that solution operator.

%% file: steps2.tex
\chapter{Summary of steps}\label{ch:steps}

In this appendix, we summarize the various operations needed in each step of the numerical scheme presented in this thesis.

\subsubsection{Preliminary remarks}
We need two solvers for the numerical scheme, one which does exterior problem solves with boundary condition $g$ on $\pd \Omega$, and one which solves the reformulated problem inside $\Omega$ using $D$ or $\overline{D}$ as a boundary condition. Both solvers should be built with the other in mind, so their discretization points agree on $\pd \Omega$.

For the exterior solves, we shall impose a boundary condition $g$ on $\pd \Omega$ (called $u_0$ in our discussion of layer-stripping), and find the solution on the layer of points just outside of $\Omega$ (called $u_1$ in our discussion of layer-stripping). Note that $u_1$ has eight more points than $u_0$ does. However, the four corner points of $u_1$ are not needed in the normal derivative of $u$ on $\pd \Omega$ (again, because we use the five-point stencil). Also, the four corner points of $u_0$ need to be used twice each. For example, the solution $u_0$ at point $(0,0)$ is needed for the normal derivative in the negative $x_1$ direction (going left) and the normal derivative in the negative $x_2$ direction (going down). Hence we obtain a DtN operator which takes the $4N$ solution points $u_0$ (with corners counted twice) to the $4N$ normal derivatives $(u_1-u_0)/h$ (with corners omitted in $u_1$). In this way, one can impose a (random) boundary condition on the $N$ points of any one side of $\pd \Omega$ and obtain the resulting Neumann data on the $N$ points of any side of $\pd \Omega$. 

Once we have probed and compressed the DtN map $D$ into $\overline{D}$, we will need to use this $\overline{D}$ in a Helmholtz solver. We do this using ghost points $u_1$ just outside of $\pd \Omega$ (but not at the corners), which we can eliminate using $\overline{D}$. For best results, it is important to use $\overline{D}$ for the same solution points it was obtained. In other words, here we defined $\overline{D}$ as
$$\overline{D}u_0=\frac{u_1-u_0}{h}.$$
If instead we use the DtN map as 
$$\overline{D}u_1=\frac{u_1-u_0}{h} \qquad \text{ or } \qquad \overline{D}u_0=\frac{u_0-u_{-1}}{h},$$
we lose accuracy. In a nutshell, one needs to be careful about designing the exterior solver and the Helmholtz solver so they agree with each other. We are now ready to consider the steps of the numerical scheme.

\section{Matrix probing: $D \rightarrow \tilde{D}$}

\begin{enumerate}
\item \label{steporg} Organize information about the various submatrices of $D$, their multiplicity and the location and ``orientation'' of each of their copies in $D$. For example, %
in the waveguide medium case, the $(1,1)$ submatrix has multiplicity 2, and appears as itself (because of symmetries in $c$) also in position $(3,3)$. However, the $(2,1)$ submatrix has multiplicity 8 but appears as itself in $D$ only in position $(4,3)$ as well. It appears as its transpose (\emph{not} conjugate transpose) in positions $(1,2)$ and $(3,4)$. It appears with column order flipped in positions $(4,1)$ and $(2,3)$, and finally as its transpose with column order flipped in positions $(1,4)$ and $(3,2)$.

\item \label{steprep}Pick a representative for each distinct submatrix of $D$. To do this, think of which block columns of $D$ will be used for probing. Minimizing the distinct block columns used will minimize the cost of matrix probing by minimizing $Q$, the sum of all solves needed. See step \ref{stepq} as well.

\item If the medium $c$ has discontinuities, it might be necessary to split up submatrices further, and to keep track of the ``sub-submatrices'' and their positions, multiplicities and orientations inside the representative submatrix.

\item \label{stepq} Pick a $q$ for each block column, keeping in mind that diagonal submatrices are usually the hardest to probe (hence need a higher $p$ and $q$), and that submatrices the farthest from the diagonal are typically very easy to probe. It might be wise to pick representatives, in step \ref{steprep}, knowing that some will require a higher $q$ than others.

\item \label{stepext} Solve the exterior problem $q$ times on each block column, saving the restriction of the result to the required block rows depending on the representative submatrices you chose. Also, save the random vectors used.

\item \label{steperror} For error checking in step \ref{steperrcheck}, also solve the exterior problem a fixed number of times, say 15, with different random vectors. Again, save both the random vectors and the exterior solves. Use those results to approximate the norm of $D$ as well.
\item For each representative submatrix $M$ (and representative sub-submatrix, if needed), do the following:
\begin{enumerate}
\item Pick appropriate basis matrices $B_j$.
\item Orthonormalize the basis matrices if necessary. If this step is needed, it is useful to use symmetries in the basis matrices to both reduce the complexity of the orthonormalization and enforce those symmetries in the orthonormalized basis matrices.
\item Multiply each basis matrix by the random vectors used in step \ref{stepext} in solving the exterior problem corresponding to the correct block column of $M$. Organize results in the matrix ${\bf \Psi}$.
\item Take the pseudo-inverse of ${\bf \Psi}$ on the results of the exterior solves from step \ref{stepext}, corresponding to the correct block row of $M$, to obtain the probing coefficients $\bc$ and $\tilde{M}=\sum c_j B_j$.
\item \label{steperrcheck} To check the probing error, multiply $\tilde{M}$ with the random vectors used in the exterior solves for error checking purposes, in step \ref{steperror}. Compare to the results of the corresponding exterior solves. Multiply that error by the square root of the multiplicity, and divide by the estimated norm of $D$.
\end{enumerate}
\item If satisfied with the probing errors of each submatrix, move to next step: PLR compression.
\end{enumerate}

\section{PLR compression: $\tilde{D} \rightarrow \overline{D}$}

\begin{enumerate}
\item For each probed representative submatrix $\tilde{M}$ (and representative sub-submatrix, if needed), do the following:
\begin{enumerate}
\item Pick a tolerance $\eps$ which is smaller than the probing error for that $\tilde{M}$. A factor of 25 smaller works well usually. Also, pick a maximal desired rank $\Rm$. Usually $\Rm \leq 8$ for a diagonal submatrix, $\Rm \leq 4$ for a submatrix just off of the diagonal, and $\Rm =2$ for a submatrix furthest away from the diagonal work well.
\item Compress $\tilde{M}$ using the PLR compression algorithm. Keep track of the dimensions and ranks of each block, to compare the matrix-vector complexity with that of a dense product.
\item Check the error made by the PLR compression by comparing $\tilde{M}$ and $\overline{M}$, again multiply that error by the square root of the multiplicity, and divide by the estimated norm of $D$.
\end{enumerate}
\item  If satisfied with the PLR errors of each submatrix, move to next step: using the PLR compression of probed submatrices into a Helmholtz solver.
\end{enumerate}

\section{Using $\overline{D}$ in a Helmholtz solver}

Using the Helmholtz solver described in the preliminary remarks of this appendix, obtain the approximate solution $\overline{u}$ from solving with the appropriate boundary condition using $\overline{D}$. Every time a product of a vector $v$ with $\overline{D}$ is needed, build the result from all submatrices of $\overline{D}$. For each submatrix, multiply the correct restriction of that vector $v$ by the correct probed and compressed representative submatrix $\overline{M}$, taking into account the orientation of the submatrix as discussed in step \ref{steporg} of the matrix probing part of the numerical scheme.

%
%

%% file: biblio.tex
\begin{singlespace}
\bibliography{main2}
\bibliographystyle{plain}
\end{singlespace}